\documentclass{svmono}

\usepackage{makeidx}     % allows index generation
\usepackage{graphicx}    % standard LaTeX graphics tool
\usepackage{multicol}    % used for the two-column index
\makeindex             % used for the subject index
\usepackage{amscd}
\usepackage{amsmath}
\usepackage{amssymb}
\usepackage{latexsym}
\usepackage{theorem}

\usepackage[all]{xy}

\newenvironment{pf}
{\noindent {\it Proof.}}
{\hfill $\Box$}

\newcommand{\N}{{\mathbb N}}
\newcommand{\Z}{{\mathbb Z}}

\newcommand{\R}{{\mathbb R}}
\newcommand{\C}{{\mathbb C}}

\newcommand{\CP}{{\mathbb C}{\mathbb P}}
\newcommand{\K}{{\mathbb K}}
\newcommand{\vbbK}{\check{{\mathbb K}}}

\newcommand{\bbS}{{\mathbb S}}

\newcommand{\ConfA}{{\sf ConfA}}

\newcommand{\locLie}{{\sf locLie}}
\newcommand{\Mod}{{\sf Mod}}

\newcommand{\OPEA}{{\sf OPEA}}

\newcommand{\Set}{{\sf Set}}

\newcommand{\Vect}{\text{{\sf Vect}}}
\newcommand{\VertA}{{\rm {\sf VertA}}}
\newcommand{\Vir}{{\sf Vir}}
\newcommand{\Witt}{{\sf Witt}}

\newcommand{\cB}{{\cal B}}
\newcommand{\cC}{{\cal C}}
\newcommand{\cD}{{\cal D}}

\newcommand{\cF}{{\cal F}}

\newcommand{\cH}{{\cal H}}

\newcommand{\cL}{{\cal L}}

\newcommand{\cP}{{\cal P}}

\newcommand{\cR}{{\cal R}}
\newcommand{\cS}{{\cal S}}

\newcommand{\fg}{{\mathfrak g}}
\newcommand{\fh}{{\mathfrak h}}

\newcommand{\fsl}{{\mathfrak s}{\mathfrak l}}

\newcommand{\tfg}{\tilde{\fg}}
\newcommand{\hfg}{\hat{\fg}}

\newcommand{\End}{{\rm End}}

\newcommand{\GL}{{\rm GL}}

\newcommand{\id}{{\rm id}}

\newcommand{\PSL}{{\rm PSL}}

\newcommand{\QEnd}{{\rm QEnd}}

\newcommand{\res}{{\rm res}}

\newcommand{\rspan}{{\rm span}}
\newcommand{\SL}{{\rm SL}}

\newcommand{\iti}{{\rm (\hskip -0.7pt{\it i}\hskip -0.2pt)}}
\newcommand{\itii}{{\rm (\hskip -0.7pt{\it ii}\hskip -0.2pt)}}
\newcommand{\itiii}{{\rm (\hskip -0.7pt{\it iii}\hskip -0.2pt)}}
\newcommand{\itiv}{{\rm (\hskip -0.7pt{\it iv}\hskip -0.2pt)}}
\newcommand{\itv}{{\rm (\hskip -0.7pt{\it v}\hskip -0.2pt)}}
\newcommand{\itvi}{{\rm (\hskip -0.7pt{\it vi}\hskip -0.2pt)}}
\newcommand{\itvii}{{\rm (\hskip -0.7pt{\it vii}\hskip -0.2pt)}}
\newcommand{\itviii}{{\rm (\hskip -0.7pt{\it viii}\hskip -0.2pt)}}
\newcommand{\itip}{{\rm (\hskip -0.7pt{\it i}'\hskip -0.2pt)}}
\newcommand{\itiip}{{\rm (\hskip -0.7pt{\it ii}'\hskip -0.2pt)}}
\newcommand{\itiiip}{{\rm (\hskip -0.7pt{\it iii}'\hskip -0.2pt)}}
\newcommand{\itipp}{{\rm (\hskip -0.7pt{\it i}''\hskip -0.2pt)}}
\newcommand{\itiipp}{{\rm (\hskip -0.7pt{\it ii}''\hskip -0.2pt)}}

\newcommand{\ita}{{\rm ({\it a})}}
\newcommand{\itb}{{\rm ({\it b})}}
\newcommand{\itc}{{\rm ({\it c})}}
\newcommand{\itd}{{\rm ({\it d})}}
\newcommand{\ite}{{\rm ({\it e})}}

\newcommand{\al}{\alpha}
\newcommand{\be}{\beta}
\newcommand{\ga}{\gamma}
\newcommand{\de}{\delta}

\newcommand{\varep}{\varepsilon}

\newcommand{\io}{\iota}
\newcommand{\ka}{\kappa}

\newcommand{\la}{\lambda}
\newcommand{\si}{\sigma}

\newcommand{\bF}{\bar{F}}
\newcommand{\bH}{\bar{H}}

\newcommand{\bL}{\bar{L}}
\newcommand{\bM}{\bar{M}}

\newcommand{\bS}{\bar{S}}
\newcommand{\bT}{\bar{T}}

\newcommand{\bV}{\bar{V}}

\newcommand{\ba}{\bar{a}}
\newcommand{\bb}{\bar{b}}
\newcommand{\bc}{\bar{c}}
\newcommand{\bh}{\bar{h}}
\newcommand{\bi}{\bar{i}}

\newcommand{\bk}{\bar{k}}
\newcommand{\bm}{\bar{m}}
\newcommand{\bn}{\bar{n}}

\newcommand{\br}{\bar{r}}
\newcommand{\bs}{\bar{s}}
\newcommand{\bt}{\bar{t}}
\newcommand{\bx}{\bar{x}}
\newcommand{\bw}{\bar{w}}

\newcommand{\bz}{\bar{z}}

\newcommand{\va}{\check{a}}
\newcommand{\vc}{\check{c}}
\newcommand{\vh}{\check{h}}
\newcommand{\vi}{\check{i}}
\newcommand{\vj}{\check{j}}
\newcommand{\vk}{\check{k}}
\newcommand{\vm}{\check{m}}
\newcommand{\vn}{\check{n}}
\newcommand{\vr}{\check{r}}
\newcommand{\vs}{\check{s}}
\newcommand{\vt}{\check{t}}
\newcommand{\vx}{\check{x}}

\newcommand{\vw}{\check{w}}
\newcommand{\vz}{\check{z}}

\newcommand{\vH}{\check{H}}
\newcommand{\vL}{\check{L}}
\newcommand{\vN}{\check{N}}

\newcommand{\vT}{\check{T}}

\newcommand{\hc}{\hat{c}}
\newcommand{\hk}{\hat{k}}

\newcommand{\hV}{\hat{V}}

\newcommand{\ta}{\tilde{a}}
\newcommand{\tb}{\tilde{b}}
\newcommand{\tc}{\tilde{c}}
\newcommand{\td}{\tilde{d}}

\newcommand{\tis}{\tilde{s}}
\newcommand{\tR}{\tilde{R}}

\newcommand{\tT}{\tilde{T}}

\newcommand{\bfg}{\bar{\mathfrak{g}}}

\newcommand{\vast}{\check{\ast}}

\newcommand{\del}{\partial}

\newcommand{\uppm}{^{\pm 1}}

\newcommand{\even}{_{\bar{0}}}
\newcommand{\odd}{_{\bar{1}}}

\newcommand{\paraab}{\zeta^{\tilde{a}\tilde{b}}}

\newcommand{\lra}{\longrightarrow}
\newcommand{\tra}{\tilde{\hskip 3pt \rightarrow}}

\newcommand{\set}[1]{\{ #1\}}
\newcommand{\sqbrack}[1]{\langle #1\rangle}
\newcommand{\normord}[1]{:{\hskip -3pt #1 }{\hskip -1pt :}}
\newcommand{\pau}[1]{[\![ #1]\!]}
\newcommand{\lau}[1]{(\!( #1)\!)}

\newcommand{\SAlg}[1]{\text{$#1$-{\sf Alg}}}
\newcommand{\Alg}{\text{{\sf Alg}}}

\begin{document}

\hyphenation{ca-te-go-ry ca-te-go-ries cor-res-pon-ding cor-res-ponds di-vi-si-bi-li-ty du-a-li-ty e-qui-va-lent ge-ne-ral ge-ne-ra-lized ge-ne-ral-ly lo-ca-li-ty ma-the-ma-ti-cal-ly mo-dule mo-dules o-pe-ra-tor o-pe-ra-tors or-ga-nized pro-cee-dings pro-duct pro-ducts sa-tis-fy sa-tis-fy-ing to-po-lo-gy}

\begin{titlepage}
\topmargin 2 cm
\begin{center}  
\LARGE \bf 
OPE-Algebras\\
\end{center}
\vspace{2 cm}
\normalsize \rm \begin{center} 
{\bf Dissertation} 
\\
\vspace{.5 cm}
zur 
\\
\vspace{.5 cm}
Erlangung des Doktorgrades (Dr.~rer.nat.) 
\\
\vspace{.5 cm}
der
\\
\vspace{.5 cm}
Mathematisch-Naturwissenschaftlichen Fakult\"at 
\\
\vspace{.5 cm}
der
\\
\vspace{.5 cm}
Rheinischen Friedrich-Wilhelms-Universit\"at Bonn
\\
\vspace{2 cm}
vorgelegt von 
\\
\vspace{.5 cm}
Markus Rosellen
\\
\vspace{.5 cm}
aus Bonn
\\
\vspace{1 cm}
Bonn 2002
\end{center}
\end{titlepage}

Angefertigt mit Genehmigung der 
Mathematisch-Naturwissenschaftlichen Fakult\"at
der Rheinischen Friedrich-Wilhelms-Universit\"at Bonn

\vspace{15cm}
1.~Referent: Prof.~Dr.~Yuri Manin

2.~Referent: Prof.~Dr.~Werner Nahm

\vspace{.5cm}
Tag der Promotion:

\newpage

\author{Markus Rosellen}
\title{OPE-Algebras}
%\title{OPE-Algebras:\\
%Mathematical Foundations of\\
%An Approach to\\
%Conformal Field Theory}
%{\small SPIN Springer's internal project number, if known}
%\subtitle{Monograph -- Mathematics --}
\maketitle

\frontmatter%%%%%%%%%%%%%%%%%%%%%%%%%%%%%%%%%%%%%%%%%%%%%%%%%%%%%%

%\include{dedic}

%\contentsline {chapter}{Preface}{V}  put this into book.toc

%%%%%%%%%%%%%%%%%%%%%% pref.tex %%%%%%%%%%%%%%%%%%%%%%%%%%%%%%%%%%%%%
%
% sample preface
%
% Use this file as a template for your own input.
%
%%%%%%%%%%%%%%%%%%%%%%%% Springer-Verlag %%%%%%%%%%%%%%%%%%%%%%%%%%

\preface

This work develops some mathematical foundations 
of the operator forma\-lism of two-dimensional conformal field theory.
The approach taken is due to A.~Kapustin and D.~Orlov. 
It is an extension of the formulation of 
chiral algebras of conformal field theory
in terms of vertex algebras.

All references to the literature have been relegated to
Appendix \ref{C:biblio}.

\bigskip

{\bf Acknowledgements.}\:
The text is partly based on talks 
given during the years 2000 till 2002
at the Seminar Topological String Theory 
that was organized by 
V.~Cort\'es, D.~Huybrechts, W.~Nahm, and myself
at the Max-Planck-Institut f\"ur Mathematik in Bonn.
I am very grateful to 
the coorganizers of the seminar
and to the participants, 
T.~Bayer,
P.~Brosnan,
C.~Devchand,
O.~Grandjean,
C.~Hertling,
B.~Jurco,
M.~Lehn,
A.~Raina,
D.~Roggenkamp,
M.~Spitzweck,
J.~Stix,
S.~Stolz,
P.~Teichner,
B.~Toen,
B.-L.~Wang,
K.~Wendland, and
P.~Xu. 
Without our joint efforts this work would not exist.

My study of OPE-algebras originated from a collaboration
with D.~Huybrechts and M.~Lehn
on moduli spaces of conformal field theories.
I thank them heartily for many pleasant and fruitful
discussions.

I am grateful to D.~van Straten for inviting me to Mainz
to talk about my Ph.D.~problem
and for inviting me to give together with 
D.~Huybrechts and M.~Lehn a series of talks about
moduli spaces of toroidal CFTs at the workshop
on mirror symmetry in Kaiserslautern in 2001.

I also like to thank
M.~Bergvelt,
J.~Fuchs,
T.~Gannon,
V.~Hinich,
G.~H\"ohn,
K.~Nagatomo,
U.~Ray,
E.~Scheidegger,
C.~Schweigert,
W.~Wang, and
D.~Zagier
for helpful discussions.

I thank Prof.~Yu.~I.~Manin and Prof.~W.~Nahm
for their willingness to be my supervisors.

I am particularly grateful to D.~Orlov for valuable
email correspondence.
He showed me how to prove that the numerators $c^i(\vz,\vw)$
of the OPE are creative.
I am very thankful to him and H.~Li for writing reports 
about this text.

%% Please "sign" your preface
\vspace{1cm}
%\begin{flushright}
\noindent
Bonn, September 2002
\hfill {\it Markus Rosellen} \\

%\hfill 
%{\it First name  Surname}\\
%\end{flushright}

\newpage

{\bf Terminology.}\:
Our conventions concerning superalgebra are
in part non-standard and are
summarized in appendix \ref{SS:superalgebra}.
In the main body of this text
we work with super objects without making
this explicit in our terminology,
e.g. we will work with super vector spaces
and just call them vector spaces.
Supersigns are included.
In order to distinguish them from other signs
they are written as powers of $\zeta:=-1$.

If $\cP$ is a property defined for elements of a set $S$
then we say that a subset $T$ of $S$ satisfies $\cP$
if any element of $T$ satisfies $\cP$.
Let $\cP$ be a property 
defined for pairs of elements of a set $S$.
We say that a subset $T$ of $S$ satisfies $\cP$
if any pair of elements of $T$ satisfies $\cP$.
We say that a family $(s_i)$ in $S$ satisfies $\cP$
if $s_i$ and $s_j$ satisfy $\cP$ for any $i\ne j$.

We use the term {\bf operator} \index{operator}
synonymously with vector space endomorphism.

\bigskip

{\bf Notation.}\:
We denote by 
$\N, \Z, \R, \C$, and $\Z_<$
the monoid of non-negative integers,
the ring of integers, 
the fields of real and complex numbers, and 
the semigroup of negative integers, resp.
If $n$ is an integer then we denote by $n_+$
the maximum of $0$ and $n$.

We work over a fixed ground field $\K$ of characteristic zero.
For a finite set $S$, 
we introduce a partial order on $\K^S$ by $a\leq b$ if $b-a$ lies in $\N^S$.
We denote by $\K^{\times}$ 
the multiplicative group $\K\setminus\set{0}$.

For a non-negative integer $n$,
we denote by $\bbS_n$ the symmetric group on $n$ letters.

Let $a$ be an element of a commutative algebra and $n\in\K$. 
The binomial coefficients \index{binomial coefficient}
are defined by $\binom{a}{n}:=\big(\prod_{i=0}^{n-1}(a-i)\big)/n!$ if
$n\in\N$ and $\binom{a}{n}:=0$ otherwise.
The divided powers \index{divided power}
are defined by $a^{(n)}:=a^n/n!$ if $n\in\N$ and $a^{(n)}:=0$ otherwise.

The Kronecker symbol $\de_{a,b}$ 
is defined to be $1$ if $a=b$ and to be $0$ otherwise.
We write $a\equiv b$ if $a$ and $b$ are
by definition different notations for the same object. 

If $\cC$ is a category then we denote by $\cC(X,Y)$ 
the set of morphisms of $\cC$ from an object $X$ to an object $Y$.

\tableofcontents

\mainmatter%%%%%%%%%%%%%%%%%%%%%%%%%%%%%%%%%%%%%%%%%%%%%%%%%%%%%%%
\chapter{Introduction}
\label{C:intro}

{\bf Summary.}\:
In section \ref{S:intro overview}
we give an overview of this thesis.
In sections \ref{S:intro q fields} and \ref{S:ope alg intro}
we present the definition of an OPE-algebra 
and point out its relation to conformal field theory.
In section \ref{S:intro op alg}
we discuss the background of our first main result.
In section \ref{S:overview} 
we describe the first main result.
In section \ref{S:intro rev dual local}
we recall what is known about the axiomatic theory of vertex algebras.
In section \ref{S:intro local}
we motivate the definition of locality for OPE-algebras.
In section \ref{S:intro dual local}
we explain the second and third main result.
In section \ref{S:intro examp}
we discuss the problem of constructing further examples of OPE-algebras.
In section \ref{S:motivation}
we indicate the motivation for our work.

\section{Overview}
\label{S:intro overview}

In their paper on mirror symmetry for complex tori from 2000,
Kapustin and Orlov introduce the notion of an OPE-algebra. 
They propose that it provides an algebraic formulation
of the physical concept of a conformal field theory  
on a closed oriented surface of genus zero.
They show that toroidal conformal field theories,
i.e.~those associated to flat tori with a B-field, 
can be realized as OPE-algebras.

The holomorphic fields of a conformal field theory
form a subalgebra that is called the chiral algebra.
Its structure is formalized by the notion of a vertex algebra.
Kapustin and Orlov prove that 
the notion of an OPE-algebra generalizes that of a vertex algebra
and that 
the holomorphic fields of an OPE-algebra form a vertex subalgebra.

In this thesis we study the question which concepts and results
of the general theory of vertex algebras can be extended to OPE-algebras.

There are three main results.
First, we define products of non-holomorphic fields,
introduce the notion of multiple locality, and
prove that the space $\cF$ of fields of a multiply local OPE-algebra $V$
endowed with these products is an OPE-algebra that is isomorphic to $V$
via the state-field correspondence $Y:V\to\cF$.

Kapustin and Orlov prove the non-trivial result that locality,
which is the main axiom in the definition of an OPE-algebra,
is a generalization of locality for vertex algebras.
Haisheng Li proved in 1994 that vertex algebras can be equivalently defined 
either 
in terms of locality, or
in terms of duality and skew-symmetry, or
in terms of the Jacobi identity.
We define the notions of duality and skew-symmetry for OPE-algebras 
and prove that these two properties together imply locality.

Toroidal conformal field theories are very special.
Their OPE-algebras satisfy additive locality 
which is a special case of locality. 
We define the notions of additive duality and additive locality, 
introduce the Jacobi identity for additive OPE-algebras, 
and prove that the result of Li generalizes to additive OPE-algebras.

\section{Quantum Fields}
\label{S:intro q fields}

The crucial notion in the definition of an OPE-algebra 
is that of a {\it quantum field}.
We motivate its mathematical formulation.

\bigskip

A conformal field theory is a quantum field
theory whose group of symmetries contains the group of 
conformal spacetime transformations.
Therefore the basic data of a conformal field theory are 
that of a general quantum field theory.
They consist of 
a complex vector space $V$, called the Hilbert space or state space, 
a vector $1$ of $V$ called the vacuum, and a 
vector space $\cF$ of {\it quantum fields} on $V$.
Vectors of $V$ are called {\it states}.

A quantum field should be thought of as a distribution on
spacetime with values in operators on $V$.
At least in the case of conformal field theories
defined on a two-dimensional spacetime 
it is possible to formalize this notion in the following way.

We define a {\it distribution} with values in a complex vector space $W$
as a formal sum
$$
a(z_1, \dots, z_r)
\; =\;
\sum_{n_1, \dots, n_r\in\R}\:
a_{n_1, \dots, n_r}\, z_1^{n_1}\dots z_r^{n_r}
$$
where $z_1, \dots, z_r$ are formal variables and 
$a_{n_1, \dots, n_r}$ are vectors of $W$.
The vector space of distributions is denoted by 
$W\set{z_1, \dots, z_r}$.
It is a distinctive feature of our approach that
we only consider such formal expressions. 

Let $\C[z^{\R}]$ denote the group ring of $\R$ with basis
$(z^n)_{n\in\R}$ and multiplication $z^n z^m:=z^{n+m}$.
The vector space $W\set{z}$ is a module over $\C[z^{\R}]$.
There exists a natural isomorphism $a(z)\mapsto\al$
between $W\set{z}$ and the space of linear maps from $\C[z^{\R}]$ to $W$.
The map $\al:\C[z^{\R}]\to W$ is given by $p(z)\mapsto\int p(z)a(z)dz$
where $\int b(z)dz$ denotes the residue $b_{-1}$ of 
a distribution $b(z)$.

There exists a natural topology on $\C[z^{\R}]$ that induces
the $(z)$-adic topology on the polynomial ring $\C[z]$.
An $\End(V)$-valued distribution $a(z)$
is called a {\it field} if $\al:\C[z^{\R}]\to\End(V)$ is continuous
where $\End(V)$ is endowed with the weak topology.
An $\End(V)$-valued distribution $a(z)$ is a field
if and only if for any state $b$ the distribution 
$a(z)b:=\sum_{n\in\R}a_n(b)z^n$ is contained in the
$\C[z^{\R}]$-submodule of $V\set{z}$ generated by the 
subspace $V\pau{z}$ of power series $\sum_{n\in\N}a_n z^n$.

The definition and the characterization of a field
generalize directly to several variables.
If $a(z,w)$ is a distribution then $a(w,w)$ is in general {\it not}
well-defined, 
but if $a(z,w)$ is a field then $a(w,w)$ exists 
and is a field, too.

Because spacetime is two-dimensional 
the relevant distributions are of the form
$a(z_1, \bz_1, \dots, z_r, \bz_r)$ where 
$(z_i,\bz_i)$ are pairs of formal variables.
In order to simplify our notation, we write
$\vn:=(n,\bn)$ and $\vz^{\vn}:=z^n \bz^{\bn}$ for
$(n,\bn)\in\R^2$, $a(\vz):=a(z,\bz), W\set{\vz}:=W\set{z,\bz}$, etc.

\section{OPE-Algebras}
\label{S:ope alg intro}

We give the precise definition of an {\it OPE-algebra}
because in the main body of the text the definition only appears very late
despite the fact that it is elementary to state.
We motivate it by recalling the Wightman axioms.

\bigskip

The basic datum $(V,1,\cF)$ of a quantum field theory
is also the datum of an OPE-algebra.
The axioms of an OPE-algebra are analogous to the 
Wightman axioms for a quantum field theory.
The Wightman axioms are:

\begin{enumerate}
\item[$\iti$]
symmetry:
the group of spacetime symmetries acts on $V$ and 
fields in $\cF$ transform covariantly with respect to this action;

\item[$\itii$]
vacuum:
there exists a unique state $1$ that is invariant with respect to
spacetime transformations;

\item[$\itiii$]
completeness:
the states $m(\phi_1(s_1), \dots, \phi_r(s_r))1$ span the state space 
where $\phi_i\in\cF$, $s_i$ are test functions, and $m$ is a monomial;

\item[$\itiv$]
locality:
the commutator of any two fields in $\cF$ vanishes
if the positions of the fields are spacelike separated.
\end{enumerate}

\bigskip

{\bf Definition.}\:
A vector space $V$ together with a vector $1$ of $V$ and
a vector subspace $\cF$ of $\End(V)\set{\vz}$ is called
an {\bf OPE-algebra} if 

\begin{enumerate}
\item[$\iti$]
there exist two operators $T$ and $\bT$ of $V$ such that 
any distribution $a(\vz)$ in $\cF$ is {\it translation covariant}
for $T$ and $\bT$, 
i.e.~$[T,a(\vz)]=\del_z a(\vz)$ and $[\bT,a(\vz)]=\del_{\bz} a(\vz)$
where $\del_z$ and $\del_{\bz}$ denote partial derivatives;

\item[$\itii$]
the state $1$ is {\it invariant} for $T$ and $\bT$, i.e.~$T(1)=\bT(1)=0$;

\item[$\itiii$]
any distribution $a(\vz)$ in $\cF$ is {\it creative},
i.e.~$a(\vz)1\in V\pau{\vz}$, 
and the {\it field-state correspondence} $s_1:\cF\to V$, 
defined by $a(\vz)\mapsto a(\vz)1|_{\vz=0}$, is surjective;

\item[$\itiv$]
any two distributions $a(\vz)$ and $b(\vz)$ in $\cF$ are {\it mutually local}, 
i.e.~there exist fields $c^i(\vz,\vw)$ and 
$\vh_i\in\R^2$ such that $h_i-\bh_i\in\Z$ and 
\begin{equation}
\label{E:intro ope}
a(\vz)b(\vw)
\; =\;
\sum_{i=1}^r
\frac{c^i(\vz,\vw)}{(\vz-\vw)^{\vh_i}},
\qquad
b(\vw)a(\vz)
\; =\;
\sum_{i=1}^r
\frac{c^i(\vz,\vw)}{(\vz-\vw)_{w>z}^{\vh_i}}
\end{equation}
where $a(\vz)b(\vw):=\sum_{\vn, \vm\in\R^2}a_{\vn}b_{\vm}\vz^{\vn}\vw^{\vm}$, 
for $\vh\in\R^2$ we define\linebreak[0]
$(\vz-\vw)^{\vh}:=
\sum_{\vi\in\N^2}(-1)^{i+\bi}\binom{h}{i}\binom{\bh}{\bi}
\vz^{\vh-\vi}\vw^{\vi}$,
and if $h-\bh\in\Z$ we define
$(\vz-\vw)_{w>z}^{\vh}:=(-1)^{h-\bh}(\vw-\vz)^{\vh}$.
\end{enumerate}

\bigskip

If $V$ is an OPE-algebra then 
distributions in $\cF$ are actually fields, 
the operators $T$ and $\bT$ are unique, and
the field-state correspondence $s_1:\cF\to V$ is a vector space isomorphism. 
The inverse of $s_1$ is called the
{\it state-field correspondence} and is written 
$a\mapsto Y(a,\vz)\equiv a(\vz)\equiv \sum_{\vn\in\R^2}a_{(\vn)}\vz^{-\vn-1}$.

One can define the notion of an OPE-algebra also in terms of $Y$
instead of $1$ and $\cF$.
Note that to give a linear map $Y:V\to\End(V)\set{\vz}$ is equivalent
to giving an algebra structure 
$V\otimes V\to V, a\otimes b\mapsto a_{(\vn)}b$, for any $\vn\in\R^2$.
A morphism $\varphi:V\to W$ of OPE-algebras is by definition 
a linear map such that
$\varphi(a_{(\vn)}b)=\varphi(a)_{(\vn)}\varphi(b)$ and $\varphi(1)=1$.

\section{The Operator Algebra of Quantum Fields}
\label{S:intro op alg}

The original idea behind our first main result is 
the {\it operator algebra of quantum fields}.
This notion, which is due to Kadanoff and Polyakov, applies 
to conformal field theories in any dimension.

\bigskip

A salient feature of quantum fields is the fact that the
product $\phi(x)\psi(y)$ of two quantum fields is singular for $x=y$.
This is possible because quantum fields are {\it not} functions but
distributions.

The singularity of $\phi(x)\psi(y)$ at $x=y$ first appeared as one
of the divergences in quantum field theory that one had to get rid of.
Wick introduced the {\it normal ordered product} 
$\normord{\phi(x)\psi(y)}$ which is regular for $x=y$ 
and thus yields a new quantum field $\normord{\phi(x)\psi(x)}$.
In our situation, the normal ordered product of 
a {\it holomorphic} field $a(z)$, i.e.~$a(z)=\sum_{n\in\Z}a_n z^n$,
and a field $b(\vz)$ is defined as 
$\normord{a(z)b(\vw)}:=
(\sum_{n\in\N}a_n z^n)b(\vw)+b(\vw)(\sum_{n\in\Z_<}a_n z^n)$.

Wilson and Kadanoff had the idea to extract information 
from the singularity of $\phi(x)\psi(y)$ at $x=y$. 
They postulated the existence of the {\it operator product expansion} 
or {\it OPE}
\begin{equation}
\label{E:intro wilson ope}
\phi(x)\psi(y)
\; =\;
\sum_i \: f_i(x-y)\, \chi_i(y)
\end{equation}
where $f_i(x)$ are singular complex-valued functions and 
$\chi_i(x)$ are quantum fields.
In our case, 
for any holomorphic field $a(z)$ and any field $b(\vz)$
there exist fields $c^n(\vz)$ such that
for any integer $K$ there exists a field $r(z,\vw)$ such that
$$
a(z)b(\vw)
\; =\;
\sum_{n\in\N-K}\:
\frac{c^n(\vw)}{(z-w)^{n+1}}
\; +\;
(z-w)^K\, r(z,\vw).
$$
One can prove that $c^n(\vw)=\int (z-w)^n [a(z),b(\vw)]dz$ for any $n\in\Z$
under very natural assumptions.
Then $c^{-1-n}(\vz)=\:\normord{\del_z^n a(z)b(\vz)}/n!$ for any $n\in\N$.

A natural question is whether the quantum fields 
$\normord{\phi(x)\psi(x)}$ and $\chi_i(x)$ belong to $\cF$
provided that $\phi(x)$ and $\psi(x)$ do.
To answer this question,
we note that in the Wightman axioms $\cF$ only consists of 
``fundamental" fields.
Polyakov enlarges $\cF$ by including ``composite" fields and
replaces Wightman's completeness axiom by the much stronger requirement that 
$V$ is spanned by the states $\phi(\text{time}\to -\infty)1$
where $\phi(x)\in\cF$.
This is exactly axiom $\itiii$ in the definition of an OPE-algebra.
Polyakov's approach is more conceptual because different sets of
fundamental fields may define the same quantum field theory.
Polyakov claims that his completeness axiom is equivalent to
the condition that there exists an OPE for any quantum fields
$\phi(x)$ and $\psi(x)$ in $\cF$ 
such that the OPE-coefficients $\chi_i(x)$ belong to $\cF$.

Polyakov proposes to view $\cF$ as an algebra with infinitely many
products $\chi_i(x)$ for any fields $\phi(x), \psi(x)\in\cF$.
There exists a notion of associativity for the algebra $\cF$ 
that is equivalent to the well-known {\it crossing symmetry} 
of four-point functions.
The algebra $\cF$ has an identity given by the {\it identity field}, 
i.e.~the constant identity operator $1(x):=\id_V$, and 
the algebra $\cF$ is endowed with a set of derivations 
given by the spacetime derivatives $\del_{x^{\mu}}$.
The space $\cF$ together with this structure is called the 
{\it operator algebra}.

In the case of OPE-algebras,
Kapustin and Orlov show that the fields $c^n(\vz)$ are contained
in $\cF$. In fact, we have $c^n(\vz)=Y(a_{(n)}b,\vz)$ where
$a(z)=\sum_{n\in\Z}a_{(n)}z^{-n-1}$. 
Moreover, we have $1(z)=Y(1,\vz)$ and for any state $a$ we have
$\del_z a(\vz)=Y(Ta,\vz)$ and $\del_{\bz} a(\vz)=Y(\bT a,\vz)$.
This generalizes results of Li and Lian-Zuckerman for vertex algebras.

The space $\cF$ of fields of an OPE-algebra also satisfies the
following completeness property.
If $a(\vz)$ is a creative, translation covariant field
that is local to any field in $\cF$ then 
$a(\vz)$ is contained in $\cF$ and thus $a(\vz)=Y(s_1(a(\vz)),\vz)$.
This result is used to prove that the
OPE-coefficients $c^n(\vz)$, or more generally $\chi_i(x)$,
are contained in $\cF$. 
Because creativity and translation covariance of $\chi_i(x)$
are in general relatively easy to verify 
the main step in the proof is to establish a statement
of the following type.
If $\phi(x), \psi(x)$, and $\vartheta(x)$ are pairwise
mutually local fields then $\chi_i(x)$ and $\vartheta(x)$ are mutually 
local for any $i$. A statement of this type is generally called
{\it Dong's lemma}.

\section{The Operator Algebra of an OPE-Algebra}
\label{S:overview}

We give an outline of our first main result which states that 
for any fields $a(\vz)$ and $b(\vz)$ that satisfy the first
equation of \eqref{E:intro ope} 
there exist natural products $a(\vz)_{(\vn)}b(\vz)$ for any $\vn\in\R^2$ 
such that 
if $V$ is a multiply local OPE-algebra and 
$\cF$ is endowed with these products
then $Y:V\to\cF$ is an OPE-algebra isomorphism.

\bigskip

As we have explained in the previous section,
the idea of an operator algebra was realized for OPE-algebras
only as an ``operator module", i.e.~the space $\cF$ of fields is a 
module over the algebra of holomorphic fields via the products 
$a(z)\otimes b(\vz)\mapsto a(z)_{(n)}b(\vz):=c^n(\vz)$ for any integer $n$.
We present these results in sections 
\ref{S:formal distributions}--\ref{S:module q fields}.

In order to find the operator algebra structure on $\cF$
we have to ask what is the OPE of two fields $a(\vz)$ and $b(\vz)$.
It seems that the right formulation 
is provided by the first equation of \eqref{E:intro ope}.
Henceforth we call the first equation of \eqref{E:intro ope}
an {\it OPE} of $a(\vz)$ and $b(\vz)$.
In contrast to the holomorphic case, 
an OPE of two non-holomorphic fields does not need to exist and
if it exists then there is no explicit formula for the fields 
$c^i(\vz,\vw)$ in terms of $a(\vz)$ and $b(\vz)$.
In contrast to locality for vertex algebras,
in order to formulate locality for OPE-algebras one has to require
that an OPE exists.

The OPE \eqref{E:intro ope} is called {\it reduced} if the pairs 
$\vh_i\in\R^2$ lie in pairwise different cosets modulo $\Z^2$ and if 
$(\vz-\vw)^{-\vn}c^i(\vz,\vw)$ is {\it not} a field for any $i$ and
any $\vn\in\N^2\setminus\set{0}$.
If two fields have an OPE then they have a unique reduced OPE.

The OPE \eqref{E:intro ope} does {\it not}
agree with Wilson's OPE \eqref{E:intro wilson ope},
the difference being that the singular functions in \eqref{E:intro ope} 
are only powers $(\vz-\vw)^{\vh}$ and 
the fields $c^i(\vz,\vw)$ depend on both $\vz$ and $\vw$.
Because non-integral powers of $z$ and $\bz$ may appear in $c^i(\vz,\vw)$ 
it is in general {\it not} possible to express $c^i(\vz,\vw)$ 
as a finite Taylor series
$\sum_{\vn=0}^{\vN}d_{\vn}(\vw)(\vz-\vw)^{\vn}$ plus a remainder.
Still, if there exists an OPE of $a(\vz)$ and $b(\vz)$ then
by taking the formal Taylor coefficients of the fields $c^i(\vz,\vw)$ 
we may assign infinitely many new fields to $a(\vz)$ and $b(\vz)$.
A natural question is whether these fields belong to $\cF$.

If \eqref{E:intro ope} is an OPE of $a(\vz)$ and $b(\vz)$
then for any $\vn\in\R^2$ we define  
$$
a(\vw)_{(\vn)}b(\vw)
\; :=\;
\sum_{i=1}^r\:
\del_{\vz}^{(\vh_i-1-\vn)}c^i(\vz,\vw)|_{\vz=\vw}
$$
where $\del_z^{(n)}:=\del_z^n/n!$ if $n\in\N$ and $\del_z^{(n)}:=0$
otherwise and $\del_{\vz}^{(\vn)}:=\del_z^{(n)}\del_{\bz}^{(\bn)}$.
One can show that this definition does not depend on the choice of the OPE
and that if $a(z)$ is holomorphic and local to $b(\vz)$ then 
$a(z)_{(\vn)}b(\vz)=\de_{\bn,-1}\, a(z)_{(n)}b(\vz)$.
If \eqref{E:intro ope} is the reduced OPE of $a(\vz)$ and $b(\vz)$
and $\vn\in\N^2$ then the Taylor coefficient 
$\del_{\vz}^{(\vn)}c^i(\vz,\vw)|_{\vz=\vw}$ is equal to 
$a(\vw)_{(\vh_i-1-\vn)}b(\vw)$.
If $a(\vz)$ and $b(\vz)$ are creative and translation covariant
then so is $a(\vz)_{(\vn)}b(\vz)$ and 
$s_1(a(\vz)_{(\vn)}b(\vz))=a_{\vn}s_1(b(\vz))$ where 
$a(\vz)=\sum_{\vn\in\R^2}a_{\vn}\vz^{-\vn-1}$.

As we have explained at the end of the previous section,
what remains to be done in order to prove that 
$\cF$ is closed with respect to the products $a(\vz)_{(\vn)}b(\vz)$
and that $Y:V\to\cF$ is an isomorphism of OPE-algebras,
is to establish Dong's lemma.
At present, it is {\it not} known whether Dong's lemma
is satisfied for locality of OPE-algebras.
However, Dong's lemma can be proven for multiple locality.

{\it Multiple locality} is the immediate generalization
of locality from a pair of fields to a family $(a^i(\vz))_i$ of fields.
One requires that for any permutation $\si$ the distribution
$a^{\si 1}(\vz_{\si 1})\dots a^{\si r}(\vz_{\si r})$ is equal to
a sum of fractions with numerators 
some fields $c^{\io}(\vz_1, \dots, \vz_r)$ and with
denominators certain products of 
$(\vz_i-\vz_j)^{\vh_{ij}^k}$ and $(\vz_i-\vz_j)_{z_j>z_i}^{\vh_{ij}^k}$
for some $\vh_{ij}^k\in\R^2$.
See Definition \ref{SS:mult local dong} for the precise statement.
It is an open question whether a family of pairwise mutually local
fields is multiply local.

\section{Duality and Locality for Vertex Algebras}
\label{S:intro rev dual local}

Vertex algebras satisfy six basic identities and there are
various implications between them. 
We recall these results, always assuming the existence of an identity $1$
and a translation operator $T$.

\bigskip

Following Goddard, we call the vertex algebra identities
\begin{equation}
\label{E:intro duali}
(x+w)^M\, a(x+w)b(w)c
\; =\;
(w+x)^M\, (a(x)b)(w)c
\end{equation}
and 
\begin{equation}
\label{E:intro locality}
(z-w)^N\, [a(z),b(w)]
\; =\; 
0
\end{equation}
{\it duality} and {\it locality}.
Here $M$ and $N$ are sufficiently large integers. 
These identities are often called ``associativity" and ``commutativity". 

Duality and locality are probably 
the most fundamental identities for vertex algebras. 
One reason for this is provided by geometry.

Huang has shown that vertex operator algebras
can be formulated geometrically as algebras over the partial operad
of punctured Riemann surfaces of genus zero with local coordinates. 
The state-field correspondence $Y$ corresponds to a pair of pants,
i.e.~a sphere with three punctures. The axioms of a vertex algebra
translate into consistency conditions on the gluing and cutting 
of punctured spheres. Because a punctured sphere is either trivial
or can be cut into pairs of pants 
the consistency conditions reduce to the requirement that the
three ways that a four-punctured sphere can be cut into two
pairs of pants, are all equal.
If the punctures are labelled by $1, 2, 3, 4$ then the 
decompositions $(12)(34)$, $(14)(23)$, and $(13)(24)$ are 
conventionally called exchange in the 
$s$-, $t$-, and $u$-channel, respectively. Huang proved that 
equality between exchange in the $s$- and $t$-channel
is equivalent to duality and 
equality between exchange in the $s$- and $u$-channel
is equivalent to locality.

Vertex algebras are often defined in terms of the {\it Jacobi identity} 
$$
\de(z,w+x)\, (a(x)b)(w)c
=
\de(x,z-w)\, a(z)b(w)c 
-
\de(x,z-w)_{w>z}\, b(w)a(z)c
$$
where $\de(z,w):=\sum_{n\in\Z}w^n z^{-n-1}$.
Duality and locality together are equivalent to the Jacobi identity.
One may view the Jacobi identity
as an elaborate combination of duality and locality.

For a triple of integers $(r,s,t)$,
if we take the coefficients in front of $z^{-t-1}w^{-s-1}x^{-r-1}$
of both sides of the Jacobi identity we get an identity between
elements of the vertex algebra. 
For sufficiently large $t$, 
this identity coincides with duality \eqref{E:intro duali} with $M=t$.
For sufficiently large $r$, 
it coincides with locality \eqref{E:intro locality} with $N=r$.
For $t=0$,
it is the {\it associativity formula} $(a_{(r)}b)(z)=a(z)_{(r)}b(z)$.
For $r=0$, 
it is the {\it commutator formula}
$[a_{(t)},b_{(s)}]=\sum_{i\in\N}\binom{t}{i}(a_{(i)}b)_{(t+s-i)}$.
Finally, yet another specialization of the Jacobi identity
yields {\it skew-symmetry} 
$b_{(r)}a=\sum_{i\in\N}(-1)^{r+1+i}\, T^i(a_{(r+i)}b)/i!$.
Thus all the six basic identities of a vertex algebra are
special cases of the Jacobi identity. 
One can show that the associativity formula implies duality, 
the commutator formula implies locality, and
locality implies skew-symmetry.

Li proved that there are three equivalent definitions of a vertex algebra.
Vertex algebras can be defined in terms of either the Jacobi identity or 
locality or duality and skew-symmetry. 
Bakalov and Kac showed that 
vertex algebras with a Virasoro vector can be defined in terms
of the associativity formula. 

The notion of a vertex algebra has been generalized by allowing
more general fields than holomorphic ones in one variable.
The notions of a generalized vertex algebra and an intertwining algebra
use distributions in one variable with non-integral exponents.
The notion of a $G_n$-vertex algebra uses fields in $n$ variables 
with integer exponents.
The general theory of additive OPE-algebras is a straightforward
combination of the general theories of generalized vertex algebras
and $G_2$-vertex algebras.
But for arbitrary OPE-algebras this is no longer true.
Moreover, it was not known how to define products of
fields in one variable with non-integral exponents.

\section{Locality for OPE-Algebras}
\label{S:intro local}

We explain why locality for OPE-algebras seems to be 
the proper non-holomorphic generalization of locality for vertex algebras.

\bigskip

The first example of a non-holomorphic field is provided by
the Fubini-Veneziano vertex operators 
that generate a toroidal conformal field theory.
These fields satisfy
\begin{equation}
\label{E:intro addit loc}
(\vz-\vw)^{\vh}\, a(\vz)b(\vw)
\; -\;
(\vz-\vw)_{w>z}^{\vh}\, b(\vw)a(\vz)
\; =\;
0
\end{equation}
for some $\vh\in\R^2$ such that $h-\bh\in\Z$. 
This is the most obvious generalization of locality for vertex algebras.
Note that the product $(\vz-\vw)^{\vh}a(\vz)b(\vw)$ is well-defined
whereas $(\vz-\vw)^{\vh}_{w>z}a(\vz)b(\vw)$ is in general {\it not}.
Because the first term of \eqref{E:intro addit loc}
is a field in $\vw$, the second term is a field in $\vz$,
and they are equal we obtain that 
$c(\vz,\vw):=(\vz-\vw)^{\vh}a(\vz)b(\vw)$ is a field and 
\begin{equation}
\label{E:intro addit loc2}
a(\vz)b(\vw)
\; =\;
\frac{c(\vz,\vw)}{(\vz-\vw)^{\vh}},
\qquad
b(\vw)a(\vz)
\; =\;
\frac{c(\vz,\vw)}{(\vz-\vw)_{w>z}^{\vh}}.
\end{equation}
Thus $a(\vz)$ and $b(\vz)$ are mutually local in the sense of OPE-algebras,
equation \eqref{E:intro ope} is satisfied with $r=1$. 
Conversely, \eqref{E:intro addit loc2} implies \eqref{E:intro addit loc}.

We call fields $a(\vz)$ and $b(\vz)$ 
{\it additively local} if they satisfy \eqref{E:intro addit loc} 
for some $\vh$ and we call an OPE-algebra $V$ {\it additive}
if there exists a vector space gradation of $V$ such that 
$a(\vz)$ and $b(\vz)$ are additively local for any homogeneous states
$a$ and $b$. 

Arbitrary fields of an additive OPE-algebra are mutually local
but in general {\it not} additively local. 
The point is that locality is a bilinear relation whereas
additive locality is {\it not}. In general, it is {\it not} possible
to reduce \eqref{E:intro ope} to a common denominator since
$(\vz-\vw)^{\vh}c(\vz,\vw)$ is {\it not} a field if 
$\vh\notin\N^2$ and $c(\vz,\vw)$ is non-zero. 
In other words, we only can reduce the $i$-th and the $j$-th summand
of \eqref{E:intro ope} to a common denominator if $\vh_i\in\vh_j+\Z^2$.

Roughly speaking, the number $r$ of summands in the OPE \eqref{E:intro ope}
is the number of ``primary" fields into which $a(\vz)$ and $b(\vz)$
``fuse". 
To require that $r=1$ for any homogeneous fields $a(\vz)$ and $b(\vz)$
is a very strong condition 
because the fusion rules are usually {\it not} additive. 
The only additively local conformal field theories
seem to be the toroidal ones.
We prove that a family of pairwise additively local fields
is multiply local.

\section{Duality for OPE-Algebras}
\label{S:intro dual local}

We state our second and third main result which are
generalizations to OPE-algebras of 
results about the axiomatics of vertex algebras.

\bigskip

We define a notion of {\it duality} for OPE-algebras that 
generalizes duality for vertex algebras 
in the same way as Kapustin and Orlov's locality 
generalizes locality for vertex algebras,
see Definition \ref{SS:duality}.
Moreover, we show that skew-symmetry for vertex algebras has the
identity 
$$
b_{(\vr)}a
\; =\;
\sum_{\vi\in\N^2}(-1)^{r-\br+i+\bi}\, T^i\bT^{\bi}(a_{(\vr+\vi)}b)/i!\bi !
$$
as its generalization. 
Our second main result states that duality and skew-symmetry 
together imply locality. 
We also prove that locality implies skew-symmetry.
It is not clear to me whether one may hope that locality implies 
duality as in the case of vertex algebras.
I think that the proper non-holomorphic generalization of vertex 
algebras should include duality, either as part or as a consequence
of the axioms.

In analogy to the notion of additive locality,
we define the notion of {\it additive duality} and prove that
additive duality and skew-symmetry together are equivalent to
additive locality provided that there exists a translation operator. 
We define a generalization of the Jacobi identity and prove that 
it is equivalent to additive duality plus additive locality. 
The last two statements constitute our third main result. 

The generalized Jacobi identity has three new features.
First, the generalized Jacobi identity has 
{\it four} terms on its right-hand side.
If one thinks of the 
two terms on the right-hand side of the ordinary Jacobi identity 
as corresponding to $z>w$ and $w>z$
then these four terms of the generalized Jacobi identity 
correspond to $\vz>\vw, \vw>\vz, (z>w,\bw>\bz)$, and $(w>z,\bz>\bw)$. 
Second, the terms corresponding to $(z>w,\bw>\bz)$ and $(w>z,\bz>\bw)$
are only defined if $a(\vz)$ and $b(\vz)$ are additively local.

Third, the generalized Jacobi identity is indexed by 
triples of pairs lying in $\vr+\Z^2, \vs+\Z^2$, and $\vt+\Z^2$
for some $\vr, \vs, \vt\in\R^2$ depending on $a, b, c$.
Additive duality and additive locality are again special cases 
of the Jacobi identity corresponding to large
$\vt'\in\vt+\Z^2$ and large $\vr'\in\vr+\Z^2$.
But if $\vt, \vr\notin\Z^2$ then 
the specializations $\vt'=0$ and $\vr'=0$ do not exist
and hence there is no associativity formula and no commutator formula for 
additive OPE-algebras.
Our result $(a_{(\vn)}b)(\vz)=a(\vz)_{(\vn)}b(\vz)$ for 
multiply local OPE-algebras is a substitute for 
this non-existing associativity formula.
In the case of additive OPE-algebras it takes the more explicit form 
$$
Y(a_{(\vn)}b,\vw)
\; =\; 
\del_{\vz}^{(\vh-1-\vn)}((\vz-\vw)^{\vh}a(\vz)b(\vw))|_{\vz=\vw}
$$
where $a$ and $b$ are states such that \eqref{E:intro addit loc}
is satisfied.

Because of these three new features 
the generalized Jacobi identity 
becomes complicated and appears to be quite artificial.
We consider the question, whether there is a Jacobi identity
for OPE-algebras, to be an academic one.

\section{Examples of OPE-Algebras}
\label{S:intro examp}

It is a difficult and important problem to construct
conformal field theories rigorously.
In particular, it is a non-trivial problem to construct
further {\it examples of OPE-algebras}.
We sketch one approach to this problem.

\bigskip

The most direct construction of conformal field theories 
is based on the observation that 
the chiral algebra of a conformal field theory
commutes with the anti-chiral algebra of anti-holomorphic fields 
and that the conformal field theory is a module over these two algebras.

Thus one starts with two vertex algebras $V$ and $\bV$,
takes the families of irreducible modules $(M_i)$ and $(\bM_{\bi})$
over $V$ and $\bV$, resp.,
chooses a matrix $(n_{i\bi})$ of non-negative integers, and
tries to construct an OPE-algebra structure on
$\cH:=\bigoplus_{i,\bi}(M_i\otimes\bM_{\bi})^{\oplus n_{i\bi}}$
using intertwiners between the modules $M_i\otimes\bM_{\bi}$.
This problem is broken down to two.

The first one is to show that intertwiners between the modules $M_i$
satisfy some version of duality and locality.
This problem has been solved by Huang in many cases.
The formulation of duality and locality involves a braided tensor
category associated to $V$ that provides a correspondence
between intertwiners such that 
the duality and locality identities hold for 
those combinations of intertwiners that are prescribed by the 
correspondence.  

The second problem is to find families of intertwiners for the
modules $(M_i)$ and $(\bM_{\bi})$ that define a 
state-field correspondence $Y$ on $\cH$ such that 
$(\cH,Y)$ is an OPE-algebra.
The solution of the first problem reduces the second problem to
a problem about the braided tensor categories associated to
$V$ and $\bV$.
One also has to verify that 
if the chiral and anti-chiral halves are combined then
Huang's notions of duality and locality for fixed intertwiners
becomes equivalent to duality and locality for OPE-algebras.

The intertwiners of the chiral halves $(M_i)$ and $(\bM_{\bi})$ 
are multi-valued and their monodromy is encoded in the braided
tensor categories. 
One reason to combine the chiral and anti-chiral halves 
is to get single-valued operators.
Another reason is that the possible combinations of the 
chiral and anti-chiral halves often play a very important role
for applications. 

If one takes for $V$ and $\bV$ the vertex algebras associated to
the Heisenberg, Virasoro, and affine Lie algebras then
the resulting OPE-algebras $\cH$ correspond to 
toroidal conformal field theories, minimal models, and WZW-models, 
respectively.

Besides the construction of $\cH$ from $V$ and $\bV$,
there are also a number of ways of constructing 
new conformal field theories from given ones:
tensor products, cosets, orbifolds, simple current extensions, 
Gepner's $U(1)$-projection, and BRST cohomology.
These constructions should yield new examples of OPE-algebras.
We show that the tensor product of OPE-algebras is an OPE-algebra.

\section{Motivation}
\label{S:motivation}

Mathematicians may have various reasons to study 
quantum field theories, and in particular conformal field theories.
Because this is a huge program 
that is still in the process of formation and 
whose scope and impact cannot yet be fully over- and foreseen 
we formulate these reasons only in very general terms,
although it is possible to make them more concrete.

Basically, there are five reasons to study quantum field theories.
First, the proposed construction of the OPE-algebra $\cH$ 
sketched in the previous section 
shows that some interesting mathematics is involved in the 
construction of quantum field theories.
Second, the concept of a quantum field theory is 
mathematically significant since it combines in a non-trivial way
the {\it geometry} of spacetime with the {\it algebra} of operators
on Hilbert space.
Third, it is possible to associate to many algebraic and geometric
objects quantum field theories that capture non-trivial features 
of the mathematical objects.
Fourth, because there are various general principles underlying
quantum field theories the fact that
different algebraic and geometric structures map to
quantum field theories 
leads to a tremendous conceptual unification in terms of these 
principles.
Fifth, quantum field theory makes some amazing predictions in mathematics
that sometimes can be made precise and verified.
The resulting picture
probably can be best understood in terms of quantum field theory.

%\contentsline {chapter}{Preface}{V}  put this into book.toc

\chapter{The Algebra of Fields}
\label{C:fields}

\section{The $\N$-Fold Module of Holomorphic Distributions}
\label{S:formal distributions}

{\bf Summary.}\:
In sections \ref{SS:prelim} and \ref{SS:prelim sfold}
we discuss distributions and $S$-fold modules.
In section \ref{SS:n-fold module holom distr} 
we introduce the $\N$-fold module $\fg\set{\vz}$.
In sections \ref{SS:comm formu} and \ref{SS:ope}
we remark that the $n$-th products of $\fg\set{\vz}$ are the coefficients
of the commutator formula and of the operator product expansion.

In section \ref{SS:statefield corr}
we define the notion of a state-field correspondence.
In section \ref{SS:Jacobi formal distributions}
we prove that if $\fg$ is a Leibniz algebra then 
$\fg\set{\vz}$ satisfies the holomorphic Jacobi identity.
In sections \ref{SS:hol jacobi def} and \ref{SS:hol jacobi} 
we consider four special cases of the holomorphic Jacobi identity.
In section \ref{SS:pf hol jacobi}
we prove three implications between
the holomorphic Jacobi identity and its special cases.

\medskip

{\bf Conventions.}\:
Except when we discuss $S$-fold modules, 
we denote by $V, W$, and $U$ vector spaces,
$\fg$ is an algebra with multiplication
$a\otimes b\mapsto [a,b]$,
and $A$ is an associative algebra
with multiplication $a\otimes b\mapsto ab$ 
and commutator $[a,b]:=ab-\paraab\, ba$.

\subsection{Distributions}
\label{SS:prelim}

In order to keep these preliminaries about {\it distributions}
as short as possible
we only give ad hoc de\-finitions and leave out a few details.
More conceptual definitions and further details
can be found in appendix \ref{C:formal distri}.

\bigskip

A $V$-valued {\bf distribution} \index{distribution}
in the even formal variables $z_1, \dots, z_r$
is an expression of the form
$$
a(z_1, \dots, z_r)
\; =\;
\sum_{n_1, \dots, n_r\in\K}\:
a_{n_1, \dots, n_r}\,
z_1^{-n_1-1}\dots z_r^{-n_r-1}
$$
where $a_{n_1, \dots, n_r}$ are vectors of $V$.
By convention,
the coefficient of a distribution $a(z_1, \dots, z_r)$ 
in front of $z_1^{-n_1-1}\dots z_r^{-n_r-1}$
is denoted by $a_{n_1, \dots, n_r}$ and 
is called the {\bf $(n_1, \dots, n_r)$-th mode} \index{mode} of 
$a(z_1, \dots, z_r)$. This definition uses a total
order $z_1>\dotsc>z_r$ on the set of variables.
Usually, 
the variables $w, x, z$ are ordered by
$z>w>x$.
In the following we only discuss the case of one variable
if the extension to several variables 
is straightforward.

The {\bf support} \index{support} of a distribution $a(z)$ is 
the set of elements $n$ of $\K$ such that $a_n$ is non-zero.
We denote by $V\set{z}$ 
the vector space of $V$-valued distributions in one variable and 
by $V[z^{\K}]$ the subspace of distributions with finite support.
We identify $V$ with the space of distributions whose support
is contained in $\set{-1}$.
The vector space $\K[z^{\K}]$ endowed with the product $z^n z^m:=z^{n+m}$ 
is a commutative algebra. 
We define on $V\set{z}$ the structure of a $\K[z^{\K}]$-module 
by $z^m a(z):=\sum_{n\in\K}a_{n+m} z^{-n-1}$.
A morphism $V\otimes W\to U, a\otimes b\mapsto ab$, induces a morphism
\begin{align}
\notag
V\set{z}\otimes W\set{w}&\to U\set{z,w},
\\
\notag
a(z)\otimes b(w)&\mapsto a(z)b(w):=\sum_{n,m\in\K}a_n b_m\, z^{-n-1}w^{-m-1},
\end{align}
that is called {\bf juxtaposition}. \index{juxtaposition}
The operator $\del_z$ of $V\set{z}$ 
given by $a(z)\mapsto -\sum_{n\in\K} n\, a_{n-1}z^{-n-1}$
is called {\bf derivative}. \index{derivative}

We denote by $\K[z\uppm]$ the ring of 
Laurent polynomials \index{Laurent polynomial} which is obtained from the
polynomial ring $\K[z]$ by inverting $z$.
A distribution $a(z)$ is called {\bf integral} \index{integral distribution} 
if its support is contained in $\Z$.
The vector space $V\pau{z\uppm}$ of integral $V$-valued distributions 
is a module over $\K[z\uppm]$.

A distribution $a(z_1, \dots, z_r)$ is called a 
{\bf power series} \index{power series}
if its support is contained in $\Z_<^r$.
The vector space $V\pau{z_1, \dots, z_r}$ of $V$-valued power series is 
a module over $\K[z_1, \dots, z_r]$.
Its localization at the ideal generated by $z_1, \dots, z_r$ 
is a submodule of $V\pau{z_1\uppm, \dots, z_r\uppm}$
that is denoted by $V\lau{z_1, \dots, z_r}$ and 
whose elements are called {\bf Laurent series}. \index{Laurent series} 
There exists a morphism 
\begin{align}
\notag
V\lau{z,w}&\to V\lau{z}, 
\\
\notag
a(z,w)&\mapsto
a(z,w)|_{w=z}:=\sum_{n\in\Z}\bigg(\sum_{m\in\Z} a_{m,n-m-1}\bigg)\, z^{-n-1}.
\end{align}
A morphism $V\otimes W\to U$ induces a morphism
$V\lau{z}\otimes_{\K[z\uppm]} W\lau{z}\to U\lau{z},
a(z)\otimes b(z)\mapsto a(z)b(z):=a(z)b(w)|_{w=z}$.
The algebra $\K\lau{z}$ is the quotient field of $\K\pau{z}$
but the algebra $\K\lau{z,w}$ is {\it not} a field. 

Because the field $\K(z_1,\dots,z_r)$ of rational functions
is the free field over $\K$ generated by $z_1, \dots, z_r$
there exists a unique morphism of fields over $\K$ 
$$
T_{z_1,\dots,z_r}:
\;
\K(z_1,\dots,z_r)
\;\longrightarrow\;
\K\lau{z_1}\dots\lau{z_r}
$$
such that $z_i\mapsto z_i$.
For an integer $n$ and $\mu, \nu\in\K$, 
we define $(\mu z+\nu w)^n\linebreak[0]:=T_{z,w}((\mu z+\nu w)^n)$.
In this definition the order of the summands of $\mu z+\nu w$ 
determines whether we apply $T_{z,w}$ or $T_{w,z}$.
Remark \ref{SS:power series expansions} implies that
$(z-w)^{-1}$ is equal to the Taylor series
$\sum_{m\in\N}\del_w^{(m)}(z-w)^{-1}|_{w=0} \, w^m=
\sum_{m\in\N}z^{-m-1}w^m$.
Applying 
$\del_w^{(n)}$ to this equation, 
where $n$ is a non-negative integer,
we obtain
\begin{equation}
\label{E:z>w expansion}
(z-w)^{-n-1}
\; =\;
\sum_{m\in\N}\:
\binom{m}{n}\, w^{m-n} z^{-m-1}.
\end{equation}

Definitions formulated for one of the modules $V\set{z}, V[z^{\K}], \dots$,
often apply directly to others as well and will thus not be repeated.
For example, we denote by $V\pau{z^{-1}}$ the image of $V\pau{z}$ 
with respect to the involution of $V\set{z}$
given by $a(z)\mapsto a(z^{-1}):=\sum_{n\in\K}a_n z^{n+1}$.
The same notation also applies to $V\lau{z}$.

Let $z$ and $\bz$ be formal variables. 
Despite the resulting ambiguity,
we write
$V\set{\vz}:=V\set{z,\bz}$
and
$a(\vz):=a(z,\bz)$.
Distributions in $V\pau{z\uppm}$ are called 
{\bf holomorphic}, \index{holomorphic distribution}
distributions in $V\pau{\bz\uppm}$ are called 
{\bf anti-holomorphic}. \index{anti-holomorphic distribution}
If a definition or a statement consists of a holomorphic part and
an analogous anti-holomorphic part 
then the anti-holomorphic part is omitted.

In general, if $S$ and $T$ are sets and $\vs\in S\times T$ then 
we denote by $s\in S$ and $\bs\in T$ 
the first and the second component of $\vs$.
In general,
if $S_1,\dots, S_r$ are sets and 
$a_{s_1, \dots, s_r}(z)$ and $\ba_{s_1, \dots, s_r}(z)$ are 
$A$-valued distributions for any $(s_i)\in\prod_i S_i$ 
then we define
$\va_{\vs_1,\dots,\vs_r}(\vz):=
a_{s_1,\dots,s_r}(z)\ba_{\bs_1,\dots,\bs_r}(\bz)$ 
for any $(\vs_i)\in\prod_i S_i^2$.
For example,
we have $a(\vz)=\sum_{\vn\in\K^2}a_{\vn}\vz^{-\vn-1}$ and 
if $\va\in A^2$ and $\vn\in\K^2$ then $\va^{(\vn)}=a^{(n)}\ba^{(\bn)}$.

\subsection{$S$-Fold Modules}
\label{SS:prelim sfold}

We define the notions of an {\it $S$-fold module} and $S$-fold algebra.

\bigskip

For an even set $S$ and a vector space $V$,
a vector space $M$ together with 
an even morphism 
$V\otimes M\to M, a\otimes b\mapsto a_{(s)}b$, 
for any $s\in S$ is called an 
{\bf $S$-fold module} \index{Sfold module@$S$-fold module} over $V$.
The vector $a_{(s)}b$ is called the
{\bf $s$-th product} \index{sth product@$s$-th product} of $a$ and $b$.
A vector space $V$ together with
an $S$-fold module structure over $V$ 
is called an {\bf $S$-fold algebra}.
\index{Sfold algebra@$S$-fold algebra}

Notions defined for $S$-fold modules
often apply in an obvious way 
to $S$-fold algebras as well.
For this reason many definitions
will only be formulated for $S$-fold modules.

A morphism from an $S$-fold module $M$ over $V$
to an $S$-fold module $N$ over $W$ is by definition
a pair of linear maps $M\to N$ and $V\to W$
such that for any $s\in S$ we have a commutative diagram
$$
\xymatrix{
V\otimes M \ar[r]^{\;\,{}_{(s)}} \ar[d] & M  \ar[d]\\
W\otimes N \ar[r]^{\;\,{}_{(s)}}        & N.
}
$$

\subsection{The $\N$-Fold Module of Holomorphic Distributions}
\label{SS:n-fold module holom distr}

We prove that a distribution in $z^{-1}V\pau{z^{-1}}\set{\vw}$
has a unique expansion in negative powers of $z-w$
with coefficients in $V\set{\vw}$
and define on $\fg\set{\vz}$ the structure of an
$\N$-fold module over $\fg\pau{z\uppm}$
by taking the coefficients of such an expansion 
as the $n$-th products.

\bigskip

The {\bf restriction} \index{restriction!of a distribution}
of a distribution $a(z_1,\dots,z_r)$ 
to a subset $S$ of $\K^r$ is the distribution 
$$
a(z_1,\dots,z_r)|_S
\; :=\;
\sum_{(n_1,\dots,n_r)\in S}\:
a_{n_1, \dots, n_r}\,
z_1^{-n_1-1}\dots z_r^{-n_r-1}.
$$
In this definition we implicitly use a total order on the
set of variables. 
The $0$-th mode of a distribution $a(z)$ 
is called the {\bf residue} \index{residue} of $a(z)$ and
is denoted by $\res_z a(z)$.

\bigskip

{\bf Lemma.}\: {\it
For a distribution $a(z,\vw)$ in $V\pau{z\uppm}\set{\vw}$,
there exist unique $V$-valued distributions $c^n(\vz)$
for any non-negative integer $n$ such that
\begin{equation}
\label{E:commut ope}
a(z,\vw)|_{\N\times\K^2}
\; =\;
\sum_{n\in\N}
\frac{c^n(\vw)}{(z-w)^{n+1}}.
\end{equation}
Moreover, 
\eqref{E:commut ope} 
is equivalent to
\begin{equation}
\label{E:commut formula0}
a_{m,\vk}
\; =\;
\sum_{n\in\N}\:
\binom{m}{n}\,
c^n_{m+k-n,\bk}
\end{equation}
for any $m\in\N$ and $\vk\in\K^2$.
We have
$c^n(\vw)=\res_z((z-w)^n a(z,\vw))$.
}

\bigskip

\begin{pf}
Equation \eqref{E:z>w expansion} yields 
\begin{equation}
\label{E:ope lemma1}
\sum_{n\in\N}\: 
\frac{c^n(\vw)}{(z-w)^{n+1}}
\; =\; 
\sum_{m\in\N}\:
\sum_{n\in\N}
\binom{m}{n} w^{m-n}\, c^n(\vw)\; z^{-m-1}.
\end{equation}
This shows that \eqref{E:commut ope} and
\eqref{E:commut formula0} are equivalent.
The inverse of the triangular matrix
$\big(\binom{m}{i}  w ^{m-i}\big)_{m,i\in\N}$ 
is 
$\big(\binom{i}{n}(-w)^{i-n}\big)_{i,n\in\N}$
because we have 
$\binom{m}{i}\binom{i}{n}=\binom{m}{n}\binom{m-n}{i-n}$
and 
$\sum_{i=n}^m(-1)^{i-n}\linebreak[0]\binom{m-n}{i-n}=(1-1)^{m-n}=\de_{m,n}$
for any non-negative integers $m, n$, and $i$.
Thus equa\-ting \eqref{E:ope lemma1} with
$a(z,\vw)|_{\N\times\K^2}=
\sum_{m\in\N}\sum_{\vk\in\K^2}a_{m,\vk}\vw^{-\vk-1}\, z^{-m-1}$
we see that  \eqref{E:commut ope} is equivalent to
$$
c^n(\vw)
\; =\; 
\sum_{m\in\N}\:\binom{n}{m}(-w)^{n-m}\,
\sum_{\vk\in\K^2}a_{m,\vk}\vw^{-\vk-1}
\; =\; 
\res_z((z-w)^n a(z,\vw)).
$$
\end{pf}

\bigskip

We define on $\fg\set{\vz}$ the structure of an
$\N$-fold module over $\fg\pau{z\uppm}$ by 
$$
a(w)_{(n)}b(\vw)
\; :=\;
\res_z((z-w)^n[a(z),b(\vw)]).
$$
Thus we obtain functors 
$\fg\mapsto\fg\pau{z\uppm},\fg\set{\vz}$
from the category of algebras 
to the categories of $\N$-fold algebras and
$\N$-fold modules.
The $n$-th product
$a(z)_{(n)}b(\vz)$ only depends on $a_0, \dots, a_n$
and $b(\vz)$.

\subsection{Commutator Formula}
\label{SS:comm formu}

We define the {\it commutator formula} and remark that
the $n$-th products of $\fg\set{\vz}$ are the unique coefficients
of the commutator formula.

\bigskip

Let $(c^n(\vz))_{n\in\N}$ be a family of 
$\fg$-valued distributions,
$m$ be an integer, and $\vk\in\K^2$.
Assume either that $m$ is non-negative or
that $c^n(\vz)$ is zero for large $n$.
A holomorphic $\fg$-valued distribution $a(z)$ and 
a $\fg$-valued distribution $b(\vz)$ 
are said to satisfy the {\bf commutator formula}
for indices $m$ and $\vk$ with coefficients $(c^n(\vz))_n$
if 
$$
[a_m,b_{\vk}]
\; =\;
\sum_{n\in\N}\:
\binom{m}{n}\, c^n_{m+k-n,\bk}.
$$
\index{commutator formula}

Lemma \ref{SS:n-fold module holom distr} shows that the $n$-th products 
$a(z)_{(n)}b(\vz)$ are the unique $\fg$-valued distributions
$c^n(\vz)$ such that
$a(z)$ and $b(\vz)$ satisfy the commutator formula
for any indices $m\in\N$ and $\vk\in\K^2$
with coefficients $(c^n(\vz))_n$.

\subsection{Operator Product Expansion}
\label{SS:ope}

We define the normal ordered product 
and prove that the $n$-th products of $A\set{\vz}$ are the 
unique coefficients of the {\it operator product expansion}.

\bigskip

The {\bf annihilation} \index{annihilation part} and 
the {\bf creation part} \index{creation part} of 
an integral distribution $a(z)$ are defined by
$a(z)_+:=a(z)|_{\N}$ and $a(z)_-:=a(z)|_{\Z_<}$, resp.

\bigskip

{\bf Definition.} \:
The {\bf normal ordered product} \index{normal ordered product}
of a holomorphic $A$-valued distribution
$a(z)$ and an $A$-valued distribution $b(\vz)$ is defined as
$$
\normord{a(z)b(\vw)}
\; :=\;
a(z)_- b(\vw)
\: +\:
\paraab \, b(\vw)a(z)_+.
$$

\bigskip

{\bf Proposition.} \:  {\it
The $n$-th products of 
a holomorphic $A$-valued distribution $a(z)$ and 
an $A$-valued distribution $b(\vz)$
are the unique $A$-valued distributions $c^n(\vz)$ such that
\begin{equation}
\label{E:infinite OPE} 
a(z)b(\vw)
\; =\;
\sum_{n\in\N}\:
\frac{c^n(\vw)}{(z-w)^{n+1}}
\;\; +\; 
\normord{a(z)b(\vw)}.
\end{equation}
}

\begin{pf} 
This follows from Lemma \ref{SS:n-fold module holom distr} 
applied to $a(z,\vw)=[a(z),b(\vw)]$ because
\begin{equation}
\label{E:azbw decomp1}
a(z)b(\vw)
\; =\;
[a(z)_+,b(\vw)]
\; +\; 
a(z)_- b(\vw)
\; +\; 
\paraab \, b(\vw)a(z)_+.
\end{equation}
\end{pf}

\bigskip

Equation \eqref{E:infinite OPE} and its rewriting 
\begin{equation}
\label{E:tilde ope}
a(z)b(\vw)
\; \sim \;
\sum_{n\in\N}\:
\frac{c^n(\vw)}{(z-w)^{n+1}}
\end{equation}
are called the 
{\bf operator product expansion} \index{operator product expansion}
or {\bf OPE} \index{OPE} for $z>w$
of $a(z)$ and $b(\vz)$.
The distribution $c^n(\vz)=a(z)_{(n)}b(\vz)$ is called the
$(n+1)$-st {\bf OPE-coefficient}. \index{OPE-coefficient}
Because of \eqref{E:infinite OPE} the normal ordered product
$\normord{a(z)b(\vw)}$
is also called the {\bf regular part} \index{regular part of the OPE}
of the OPE.
The tilde in \eqref{E:tilde ope}
indicates that the regular part is omitted.

Taking for $A$ the enveloping algebra of a Lie algebra $\fg$,
the Proposition shows that 
the commutator formula for $a(z)\in\fg\pau{z\uppm}$ and 
$b(\vz)\in\fg\set{\vz}$ and any indices $m\in\N$ and $\vk\in\K^2$
with coefficients $(c^n(\vz))_n$
is equivalent to the OPE for $z>w$ of $a(z)$ and $b(\vz)$ with
OPE-coefficients $c^n(\vz)$.

\subsection{State-Field Correspondence}
\label{SS:statefield corr}

Distributions and $S$-fold modules are interrelated in the same way as 
operators and ordinary modules are.
First, distributions form an $\N$-fold module.
Second, to give a $\K^2$-fold module $M$ over a vector space $V$
is equivalent to giving a {\it state-field correspondence}
$V\to\End(M)\set{\vz}, a\mapsto \sum_{\vn\in\K^2}a_{(\vn)}\vz^{-\vn-1}$.

\bigskip

Let $T$ be a subset of an even set $S$.
There exists a forgetful functor from the category of
$S$-fold modules to the category of 
$T$-fold modules, which is called {\bf restriction}
\index{restriction!of an $S$-fold module} 
and is written $M\mapsto M|_T$,
and there exists a fully faithful functor from the category of
$T$-fold modules to the category of 
$S$-fold modules which is called {\bf extension by zero}
\index{extension by zero}
and is written $M\mapsto M|^S$.

We sometimes identify 
a $T$-fold module $M$ with the $S$-fold module $M|^S$.
In particular, if $T$ is a subset of $\K$ 
we identify a $T$-fold module $M$ with 
the $\K^2$-fold module $M|^{\K^2}$
where we consider $T$ as a subset of $\K^2$
by means of the injection $T\to\K^2, n\mapsto (n,-1)$.
Conversely, 
we sometimes consider $S$-fold modules as $T$-fold modules
by making implicitly use of the functor of restriction.
This practice may cause confusion but has the advantage
of unifying our discussion of conformal algebras, 
vertex algebras, and OPE-algebras.

Let $M$ be an $S$-fold module over a vector space $V$,
$a, b\in V$, $c\in M$, and $s, t\in S$.
We denote by $a_{(s)}$ the operator of $M$ given by $d\mapsto a_{(s)}d$
and we denote by $a_{(s)}b_{(t)}$ the composition of the operators 
$a_{(s)}$ and $b_{(t)}$.
Thus $a_{(s)}b_{(t)}c=a_{(s)}(b_{(t)}c)$.

A pair $T=(T,T)$ of operators of $V$ and $M$
is called a {\bf derivation} \index{derivation}
for $a$ and $c$ and the $s$-th product if 
$T(a_{(s)}c)=T(a)_{(s)}c+\zeta^{\tT\ta}\, a_{(s)}T(c)$.

Of course, to give an $S$-fold module structure 
$({}_{(s)}:V\otimes M\to M)_{s\in S}$
on a vector space $M$ is equivalent to 
giving an even morphism 
$V\to\End(M)^S, a\mapsto (a_{(s)})_{s\in S}$.

\bigskip

{\bf Definition.}\:
For vector spaces $V$ and $M$, 
an even linear map from $V$ to $\End(M)\set{\vz}$ is called 
a {\bf state-field correspondence}. \index{state-field correspondence}

\bigskip

If $Y:V\to\End(M)\set{\vz}$ is
a state-field correspondence
then the image of a vector $a$ of $V$ with respect to $Y$
is denoted by $Y(a,\vz)$ or just by $a(\vz)$.
The image of $Y$
is called the {\bf space of fields} of $M$
and is denoted by $\cF_Y$.
\index{space of fields}

If $Y:V\to\End(V)\set{\vz}$ is a state-field correspondence
then vectors of $V$ are called {\bf states} 
and $V$ is called the {\bf state space}.
\index{state}
\index{state space}

For vector spaces $V$ and $M$, 
we use the vector space isomorphism
$\End(M)\set{\vz}\to\End(M)^{\K^2},
a(\vz)\mapsto (a_{\vn})_{\vn\in\K^2}$,
to identify 
state-field correspondences $Y:V\to\End(M)\set{\vz}$
with $\K^2$-fold module structures $V\to\End(M)^{\K^2}$.
Thus for a vector $a$ of $V$ we write
$a(\vz)=\sum_{\vn\in\K^2}a_{(\vn)}\vz^{-\vn-1}$.

A $V$-valued distribution $a(\vz)$ is called
{\bf bounded} \index{bounded distribution} by $\vh_1,\dots,\vh_r\in\K^2$ 
if $a(\vz)$ is contained in 
$\sum_i \vz^{-\vh_i}V\pau{\vz}$.
An $\End(V)$-valued distribution $a(\vz)$ is called
bounded on a vector $b$ of $V$ if $a(\vz)b$ is bounded.

We denote by $V\sqbrack{z_1, \dots, z_r}$
the $\K[z_1^{\K}, \dots, z_r^{\K}]$-submodule of $V\set{z_1, \dots, z_r}$
generated by $V\pau{z_1, \dots, z_r}$.
Distributions in $V\sqbrack{z_1, \dots, z_r}$ are called {\bf vertex series}.
\index{vertex series}
A distribution $a(\vz)$ is bounded if and only if
it is a vertex series.
A linear map $V\otimes W\to U$ induces a morphism
$V\sqbrack{z}\otimes_{\K[z^{\K}]} W\sqbrack{z}\to U\sqbrack{z}$
of $\K[z^{\K}]$-modules extending the morphism
$V\lau{z}\otimes_{\K[z\uppm]} W\lau{z}\to U\lau{z}$.

Let $M$ be a $\K^2$-fold module over a vector space $V$.
A vector $a$ of $V$ is called bounded \index{bounded vector} on
a vector $b$ of $M$ if $a(\vz)$ is bounded on $b$.
We call $M$ bounded 
\index{bounded Ksquare-fold module@bounded $\K^2$-fold module}
if any pair of vectors of $V$ and $M$ is bounded.

A vector $a$ of $V$ is called {\bf holomorphic} \index{holomorphic vector}
on $M$ if $a(\vz)$ is a holomorphic distribution.
Let $a$ be a vector of $V$ that is holomorphic on $M$.
For an integer $n$, we denote by $a_{(n)}$ the operator $a_{(n,-1)}$.
If $b$ is a vector of $M$ then the least integer $N$ such that 
$a_{(n)}b=0$ for any $n\geq N$
is called the {\bf mode infimum} \index{mode infimum} of $a$ and $b$
and is denoted by $o'(a,b)\in\Z\cup\set{\pm\infty}$.
\index{oaaaaaaaaa@$o'(a,b)$}

A state of a $\K^2$-fold algebra $V$ is called holomorphic 
\index{holomorphic state}
if it is holomorphic on $V$.

\subsection{Holomorphic Jacobi Identity for Distributions}
\label{SS:Jacobi formal distributions}

Since the $n$-th products of $\fg\set{\vz}$
are defined in terms of the algebra structure of $\fg$,
it is natural to deduce properties of 
the $\N$-fold module $\fg\set{\vz}$
from pro\-perties of the algebra $\fg$.
We prove that if $\fg$ is a Leibniz algebra 
then $\fg\set{\vz}$ satisfies the {\it holomorphic Jacobi identity}.

\bigskip

{\bf Definition.}\:
Let $M$ be a $\K^2$-fold module over a $\Z$-fold algebra $V$,
$a, b\in V$, $c\in M$, $r, t\in\Z$, and $\vs\in\K^2$ 
such that $a$ is holomorphic on $M$. 
Assume either that $r$ and $t$ are non-negative or that
the pairs $(a, b), (a,c)$, and $(b,c)$ are bounded.
The {\bf holomorphic Jacobi identity} \index{holomorphic Jacobi identity}
for \index{Jacobi identity!holomorphic} $a, b, c$ and indices $r, \vs, t$
is given by
\begin{align}
\notag
&\sum_{i\in\N}\:
\binom{t}{i}\,
(a_{(r+i)}b)_{(s+t-i,\bs)}c
\\
\notag
\; =\;
&\sum_{i\in\N}\:
(-1)^i \binom{r}{i}\,
(a_{(t+r-i)}b_{(s+i,\bs)}
\; -\;
\paraab (-1)^r \, b_{(s+r-i,\bs)}a_{(t+i)})c.
\end{align}

\bigskip

We say that $a, b, c$ satisfy the holomorphic Jacobi identity if 
the pairs $(a, b), (a,c)$, and $(b,c)$ are bounded and
$a, b, c$ satisfy the holomorphic Jacobi identity for  
any indices $r\in\Z, \vs\in\K^2$, and $t\in\Z$.
We identify the index $\vs=(s,-1)\in\K^2$ with $s\in\K$.

Let $M$ be an $\N$-fold module over an $\N$-fold algebra $V$.
We say that the holomorphic Jacobi identity 
is satisfied for $a, b\in V$ and $c\in M$ if it is satisfied
for $a, b, c$ and any indices $r, s, t\in\N$.

\bigskip

{\bf Proposition.}\;{\it
If $\fg$ is a Leibniz algebra
then 
the $\N$-fold module $\fg\set{\vz}$
over
the $\N$-fold algebra $\fg\pau{z\uppm}$
satisfies the holomorphic Jacobi identity.
}

\bigskip

\begin{pf}
Let $a(z), b(z)\in\fg\pau{z\uppm}$, $c(\vz)\in\fg\set{\vz}$,
and $r, s, t\in\N$.
Expanding the $t$-th power of $z-x=(z-w)+(w-x)$ we obtain
\begin{align}
\notag
&\res_{z,w}\big(
(z-w)^r (w-x)^s (z-x)^t[[a(z),b(w)],c(\vx)]\big) 
\\
\notag
=\;
&
\sum_{i\in\N}\:
\binom{t}{i}\: 
\res_{z,w}
\big(
(z-w)^{r+i}(w-x)^{s+t-i}[[a(z),b(w)],c(\vx)]
\big)
\\
\notag 
=\;
&\sum_{i\in\N}\: 
\binom{t}{i}\: 
\res_w
\big(
(w-x)^{s+t-i}\,
[\res_z((z-w)^{r+i}[a(z),b(w)]),\, c(\vx)] 
\big) 
\\
\notag
=\;
&\sum_{i\in\N}\: \binom{t}{i}\:
(a(x)_{(r+i)}b(x))_{(s+t-i)}c(\vx). 
\end{align}

On the other hand, the Leibniz identity 
and the expansion of the $r$-th power of $z-w=(z-x)-(w-x)$ yield
\begin{align}
\notag
&\res_{z,w}\big(
(z-w)^r (w-x)^s (z-x)^t[[a(z),b(w)],c(\vx)]\big) 
\\
\notag
\; =\;
&\res_{z,w}\big(
(z-w)^r (w-x)^s (z-x)^t\,
\big(
[a(z),[b(w),c(\vx)]]
\\
\notag
&\qquad\qquad\qquad\qquad\qquad\qquad\qquad\qquad\qquad
-\,
\paraab\, [b(w),[a(z),c(\vx)]]
\big)\big)
\\
\notag
=\; 
&
\sum_{i\in\N}\: 
(-1)^i\binom{r}{i}\:
\res_{z,w}
\Big(
(w-x)^{s+i}(z-x)^{t+r-i}[a(z),[b(w),c(\vx)]]
\\
\notag
&\qquad\qquad\qquad\:
-\paraab (-1)^r
(w-x)^{s+r-i}(z-x)^{t+i}[b(w),[a(z),c(\vx)]]
\Big)
\\
\notag
=\;
&\sum_{i\in\N}\: 
(-1)^i \:
\binom{r}{i}\:
\Big(
a(x)_{(t+r-i)}b(x)_{(s+i)}
\\
\notag
&\qquad\qquad\qquad\qquad\qquad\qquad\quad\;
-\paraab (-1)^r\,  
b(x)_{(s+r-i)}a(x)_{(t+i)}
\Big)c(\vx).
\end{align}
\end{pf}

\subsection{Special Cases of the Holomorphic Jacobi Identity I}
\label{SS:hol jacobi def}

We define for $\K^2$-fold modules the identities of 
holomorphic duality and holomorphic locality,
the associativity formula, and the commutator formula.

\bigskip

{\bf Definition.}\:
Let $M$ be a $\K^2$-fold module over a $\Z$-fold algebra $V$,
$a, b\in V$, $c\in M$, $r, t\in\Z$, and $\vs\in\K^2$
such that $a$ is holomorphic on $M$. 

\begin{enumerate}
\item[($a$)]
Assume that the pairs $a, b$ and $b, c$ are bounded.
{\bf Holomorphic duality} \index{holomorphic duality}
of order at most $t$ for $a, b, c$
and indices $r$ and $\vs$ is the identity
\begin{equation}
\notag
\sum_{i\in\N}\:
\binom{t}{i}\,
(a_{(r+i)}b)_{(s+t-i,\bs)}c
\; =\;
\sum_{i\in\N}\:
(-1)^{i} \binom{r}{i}\,
a_{(t+r-i)}b_{(s+i,\bs)}c.
\end{equation}

\item[($b$)]
Assume either that
$r$ is non-negative or that 
$a$ and $b$ are bounded on $c$.
{\bf Holomorphic locality} \index{holomorphic locality}
of order at most $r$ for $a, b, c$
and indices $t$ and $\vs$ is the identity
$$
\sum_{i\in\N}\:
(-1)^{i} \binom{r}{i}\,
(a_{(t+r-i)}b_{(s+i,\bs)}
\; -\;
\paraab\, (-1)^r \, b_{(s+r-i,\bs)}a_{(t+i)})c
\; =\;
0.
$$

\item[$\itc$]
Assume either that $r$ is non-negative or that
$a$ and $b$ are bounded on $c$.
The {\bf associativity formula} \index{associativity formula}
for $a, b, c$ and indices $r$ and $\vs$
is given by
$$
(a_{(r)}b)_{(\vs)}c
=
\sum_{i\in\N}
(-1)^i \binom{r}{i}
(a_{(r-i)}b_{(s+i,\bs)}
-
\paraab (-1)^r  b_{(s+r-i,\bs)}a_{(i)})c.
$$

\item[$\itd$]
Assume either that $t$ is non-negative or that 
$a$ is bounded on $b$.
The {\bf commutator formula} \index{commutator formula}
for $a, b, c$ and indices $t$ and $\vs$
is given by
$$
[a_{(t)},b_{(\vs)}]c
\; =\;
\sum_{i\in\N}\:
\binom{t}{i}\,
(a_{(i)}b)_{(t+s-i,\bs)}c.
$$
\end{enumerate}

\bigskip

We say that $a, b, c$ satisfy holomorphic duality of order $t$ if 
the pairs $a, b$ and $b, c$ are bounded and
$a, b, c$ satisfy holomorphic duality of order $t$ for  
any indices $r\in\Z$ and $\vs\in\K^2$.
Similar terminology applies to the other three identities.

\subsection{Special Cases of the Holomorphic Jacobi Identity II}
\label{SS:hol jacobi}

We remark that the four identities defined in the previous subsection
are all special cases of the holomorphic Jacobi identity.

\bigskip

If $a(z)$ is a holomorphic $A$-valued distribution,
$b(\vz)$ is an $A$-valued distribution, and $n$ is a non-negative integer 
then
\begin{align}
\label{E:nth prod modes}
&a(z)_{(n)}b(\vz)
\\
\notag
=\;
&\sum_{\vm\in\K^2}
\sum_{i\in\N}\:
(-1)^i \binom{n}{i}
(a_{n-i}b_{m+i,\bm}
-
\paraab (-1)^n\,  b_{m+n-i,\bm}a_i)\; \vz^{-\vm-1}.
\end{align}

\bigskip

{\bf Remark.}\:
$\iti$\:
The holomorphic Jacobi identity for $a, b, c$ and indices 
\begin{enumerate}
\item[$\ita$]
$r, \vs, t$ with $t\geq o'(a,c)$
is identical with holomorphic duality of order at most $t$
for indices $r$ and $\vs$.

\item[$\itb$]
$r, \vs, t$ with $r\geq o'(a,b)$
is identical with holomorphic locality of order at most $r$
for indices $t$ and $\vs$.

\item[$\itc$]
$r, \vs, 0$
is identical with the associativity formula for indices $r$ and $\vs$.

\item[$\itd$]
$0, \vs, t$
is identical with the commutator formula for indices $t$ and $\vs$.
\end{enumerate}

$\itii$\:
Let $r$ be a non-negative integer.
Equation \eqref{E:nth prod modes} shows that 
the associativity formula is satisfied for $a, b, c$
and indices $r$ and any $\vs\in\K^2$
if and only if $(a_{(r)}b)(\vz)c=a(z)_{(r)}b(\vz)c$.

$\itiii$\:
The commutator formula is satisfied for $a, b$, and any $c\in M$ 
and indices $m$ and $\vk$
if and only if 
$a(z)$ and $b(\vz)$ satisfy the commutator formula
for indices $m$ and $\vk$ and coefficients $((a_{(n)}b)(\vz))_n$.

$\itiv$\:
The commutator formula for $a, b, c$ and indices $0$ and $\vs$ 
is equivalent to $a_{(0)}$ being a derivation
for $b$ and $c$ and the $\vs$-th product.
In particular, the commutator formula for indices $0$ and $0$
is equivalent to the Leibniz identity for the $0$-th product.

$\itv$\:
Let $r, s, t\in\Z$ and 
$a, b, c$ be states of a bounded $\Z$-fold algebra that satisfy
holomorphic duality of order $t$.
Holomorphic duality for indices $r$ and $s-t$ is
$$
(a_{(r)}b)_{(s)}c
\; =\;
\sum_{i\in\N}\: (-1)^{i} \binom{r}{i}\, a_{(t+r-i)}b_{(s-t+i)}c
\; -\;
\sum_{i>0}\: \binom{t}{i}\, (a_{(r+i)}b)_{(s-i)}c.
$$
Applying holomorphic duality for indices $r+j$ and $s-t-j$
to the right-hand side of this equation where $j=1, \dots, o'(a,b)-r-1$ 
we see that for any $N\geq o'(a,b)$ there exist integers $N_i$, 
that only depend on $r, s, t$, and $N$ but do not depend on $a, b, c$ and $V$, 
such that
$$
(a_{(r)}b)_{(s)}c
\; =\;
\sum_{i\in\N}\: N_i\, a_{(t+N-1-i)}b_{(s+r-t-N+1+i)}c.
$$

\subsection{Holomorphic Duality and Holomorphic Locality}
\label{SS:pf hol jacobi}

We prove that 
the associativity formula implies holomorphic duality,
the commutator formula implies holomorphic locality, and 
holomorphic duality and holomorphic locality together
imply the holomorphic Jacobi identity.

\bigskip

{\bf Lemma.}\: {\it
Let $M$ be a $\K^2$-fold module over a $\Z$-fold algebra $V$,
$a, b\in V$, $c\in M$, $r, t\in\Z$, and $\vs\in\K^2$ 
such that $a$ is holomorphic on $M$ and
$a, b, c$ satisfy the holomorphic Jacobi identity
for indices $r, (s+1,\bs), t$.
The vectors $a, b, c$ satisfy the holomorphic Jacobi identity
for indices $r+1, \vs, t$ if and only if 
they satisfy it for indices $r, \vs, t+1$.
}

\bigskip

\begin{pf}
Let $J^1_{r,\vs,t}$ denote the left-hand side
of the holomorphic Jacobi identity for $a, b, c$ and 
indices $r, \vs, t$ and define
$$
J^2_{r,\vs,t}
\; :=\;
\sum_{i\in\N}\: (-1)^{i} 
\binom{r}{i}\;
a_{(t+r-i)}b_{(s+i,\bs)}c
$$
and 
$$
J^3_{r,\vs,t}
\; :=\;
\sum_{i\in\N}\: (-1)^{r+i} 
\binom{r}{i}\;
b_{(s+r-i,\bs)}a_{(t+i)}c.
$$
Thus
$J^2_{r,\vs,t}-\paraab J^3_{r,\vs,t}$ is the right-hand side of 
the holomorphic Jacobi identity.
We have
$J^i_{r,\vs,t+1}=J^i_{r+1,\vs,t}+J^i_{r,(s+1,\bs),t}$
for $i=1, 2$, and $3$ because
$$
J^1_{r,\vs,t+1}
\; =\;
\sum_{i\in\N}\:
\bigg(
\binom{t}{i-1}+\binom{t}{i}
\bigg)\;
(a_{(r+i)}b)_{(s+t+1-i,\bs)}c,
$$
$$
J^2_{r,\vs,t+1}
\; =\;
\sum_{i\in\N}
\:
\bigg(
(-1)^i\binom{r+1}{i}+(-1)^{i-1}\binom{r}{i-1}
\bigg)\;
a_{(t+1+r-i)}b_{(s+i,\bs)}c,
$$
and 
$$
J^3_{r,\vs,t+1}
\; =\;
\sum_{i\in\N}\,
\bigg(
(-1)^{r+1+i}\binom{r+1}{i}+(-1)^{r+i}\binom{r}{i}
\bigg)\,
b_{(s+r-(i-1),\bs)}a_{(t+i)}c.
$$
This implies the claim.
\end{pf}

\bigskip

{\bf Proposition.}\: {\it
Let $M$ be a $\K^2$-fold module over a $\Z$-fold algebra $V$,
$a, b\in V$, and $c\in M$ such that $a$ is holomorphic on $M$.

\begin{enumerate}
\item[$\iti$]
If $a, b, c$ satisfy the associativity formula 
then they satisfy holomorphic duality of order at most $t$ for any
integer $t$ such that $t\geq o'(a,c)_+$.

\item[$\itii$]
If $a, b, c$ satisfy the commutator formula 
then they satisfy holomorphic locality of order at most $r$ for any
integer $r$ such that $r\geq o'(a,b)_+$.

\item[$\itiii$]
The vectors $a, b, c$ satisfy holomorphic duality and holomorphic locality 
if and only if they satisfy the holomorphic Jacobi identity.
\end{enumerate}
}

\bigskip

\begin{pf}
This follows from Remark \ref{SS:hol jacobi}\,$\iti$ and the Lemma.
For $\itiii$ we use that 
if the holomorphic Jacobi identity is satisfied 
for indices $r, \vs, t-1$ and $r-1, \vs, t$ then
it is satisfied for $r-1, \vs+1, t-1$ according to the Lemma.
\end{pf}

\section{Conformal Symmetry}
\label{S:vir symm}

{\bf Summary.}\:
In section \ref{SS:sfold transl op}
we define the notion of a translation operator.
In section \ref{SS:transl cov}
we discuss the notion of translation covariance.
In section \ref{SS:derivative}
we prove that the derivative is a translation operator
of $\fg\set{\vz}$
and that the translation covariant distributions form an $\N$-fold submodule.

In section \ref{SS:sfold dilat op}
we define the notion of a dilatation operator 
and show that there exists a one-to-one correspondence between 
dilatation operators and gradations.
In section \ref{SS:dilat cov}
we discuss the notion of dilatation covariance.
In section \ref{SS:grad dilat cov}
we prove that the dilatation covariant distributions
form a graded $\N$-fold module.
In section \ref{SS:witt alg}
we define the {\it Witt algebra}.
In section \ref{SS:vir alg}
we define the {\it Virasoro algebra} as an explicit central extension
of the Witt algebra and prove that it is the universal central extension.
In section \ref{SS:conf distr}
we define the notion of a {\it conformal distribution}.

\medskip

{\bf Conventions.}\:
Except when we discuss $S$-fold modules, 
we denote by $V$ a vector space,
$\fg$ is an algebra with multiplication
$a\otimes b\mapsto [a,b]$,
and $A$ is an associative algebra
with multiplication $a\otimes b\mapsto ab$ 
and commutator $[a,b]:=ab-\paraab\, ba$.

\subsection{Translation Operator}
\label{SS:sfold transl op}

We define the notions
of a translation generator, endomorphism, and ope\-rator
of an $S$-fold module where $S$ is either $\N, \Z, \K$, or $\K^2$.
In section \ref{SS:derivative} we prove that 
the pair of derivatives $(\del_z,\del_{\bz})$
is a {\it translation operator} of $\fg\set{\vz}$.
Existence of a translation generator 
is one of the three axioms in the definition of an OPE-algebra.
An OPE-algebra has a unique translation operator.

\bigskip

{\bf Definition.}\:                     
$\iti$\:                                
Let $S$ be either $\N, \Z$, or $\K$ 
and let 
$M$ be an $S$-fold module over a vector space $V$,
$a\in V$, and $b\in M$.
 
\begin{enumerate}
\item[$\ita$]
An even operator $T$ of $M$ is called
a {\bf translation generator} \index{translation generator}
for $a$ and $b$ if 
$T(a_{(n)}b)=-n\, a_{(n-1)}b+a_{(n)}T(b)$ for any $n\in S$.

\item[$\itb$]
An even operator $T$ of $V$ is called
a {\bf translation endomorphism} \index{translation endomorphism}
for $a$ and $b$ if
$T(a)_{(n)}b=-n\, a_{(n-1)}b$ for any $n\in S$.

\item[$\itc$]
A pair $T=(T,T)$ of even operators of $V$ and $M$ 
is called
a {\bf translation operator} \index{translation operator}
for $a$ and $b$ if 
$T$ is a derivation and a translation generator for $a$ and $b$.
\end{enumerate}

$\itii$\:
Let $M$ be a $\K^2$-fold module over a vector space $V$.
An even operator $T$ of $M$ is called a
{\bf holomorphic} \index{translation generator!holomorphic}
translation generator 
if $T$ is a translation generator of the $\K$-fold module 
$M|_{\K\times\set{n}}$ for any $n\in\K$.
A pair $\vT=(T,\bT)$ consisting of
a holomorphic and an anti-holomorphic 
translation generator is called a
{\bf translation generator} of $M$.
The notions of a
translation endomorphism and of a translation operator 
of $M$ are defined in the same way.

\bigskip

Note that any two of the following three statements imply the third one:
$T$ is a derivation; $T$ is a translation generator;
$T$ is a translation endomorphism.

Let $S$ be either $\N$ or $\Z$.
An even operator $T$
of an $S$-fold module $M$ is
a translation generator of $M$
if and only if 
$T$ is a translation generator of $M|^{\K}$.
The same is true for the notions of a
translation endomorphism and
translation operator.

If we identify a $\K$-fold module $M$ with 
the $\K^2$-fold module $M|^{\K^2}$
then parts $\iti$ and $\itii$ of the definition lead to 
two {\it different} notions of
a translation generator, endomorphism, and operator of $M$. 
Since the data $T$ and $\vT$ 
of these two notions are also different
this ambiguity should not cause any confusion.
 
Let $M$ be a $\K^2$-fold module over a vector space $V$, 
$a\in V$, and $b\in M$.
If $\vT$ is a translation endomorphism for $a$ and $b$
then for any $\vm, \vn\in\K^2$ we have
\begin{equation}
\label{E:transl endo}
\vT^{(\vn)}(a)_{(\vm)}b
\; =\;
(-1)^{\vn}\binom{\vm}{\vn}\, a_{(\vm-\vn)}b.
\end{equation}
We sometimes write $1$ instead of $(1,1)\in\K^2$. 
If $\vn\in\N^2$ then putting $\vm=-1$ in \eqref{E:transl endo} we obtain
\begin{equation}
\label{E:transl endo2}
\vT^{(\vn)}(a)_{(-1)}b
\; =\;
a_{(-1-\vn)}b.
\end{equation}

A pair $T=(T,T)$ of even operators of $V$ and $M$
is a derivation for $a$ and $b$ if and only if 
$[T,a(\vz)]b=(Ta)(\vz)b$.
An even operator $T$ of $M$
is a holomorphic translation generator for $a$ and $b$ if and only if
$[T,a(\vz)]b=\del_z a(\vz)b$.
An even operator $T$ of $V$ 
is a holomorphic translation endomorphism for $a$ and $b$ if and only if
$(Ta)(\vz)b=\del_z a(\vz)b$.

\subsection{Translation Covariance}
\label{SS:transl cov}

We give an infinitesimal definition of 
the notion of {\it translation covariance}
and prove that it is equivalent to 
its integrated version.

\bigskip

{\bf Definition.}\:
A $V$-valued distribution $a(\vz)$ is called
{\bf translation covariant} \index{translation covariant distribution}

\begin{enumerate}
\item[$\iti$]
for an even operator $T$ of $V$ if $Ta(\vz)=\del_z a(\vz)$;
\item[$\itii$]
for a pair $\vT$ of even operators of $V$
if $a(\vz)$ is translation covariant for $T$ and 
$\bT a(\vz)=\del_{\bz} a(\vz)$.
\end{enumerate}

\bigskip

We denote by $V\set{\vz}_T$ and $V\set{\vz}_{\vT}$ the vector spaces of 
$V$-valued distributions that are 
translation covariant for $T$ and $\vT$, resp.
Of course, the vector space $V\set{\vz}_T$ is invariant
with respect to $\del_z$ and $\del_{\bz}$.
An $A$-valued distribution $a(\vz)$ is called
translation covariant for a pair $\vT$ of even elements of $A$
if $a(\vz)$ is translation covariant for $[\vT,\;\, ]$.

Let $M$ be a $\K^2$-fold module over a vector space $V$ and $a\in V$.
A pair $\vT$ of even operators of $M$ is a translation generator for $a$ 
and any vector of $M$ if and only if
$a(\vz)$ is translation covariant for $\vT$.

If $a$ is an element of a topological algebra 
such that the sequence
$(\sum_{i=0}^n a^{(i)})_{n\in\N}$ 
has a unique limit
then this limit is called the {\bf exponential} 
of $a$ and is denoted by $e^a$.
There exists a natural morphism $V\pau{z}\to V,
a(z)\mapsto a(z)|_{z=0}\equiv a(0):=a_{-1}$.
\index{exponential}

\bigskip

{\bf Remark.}\: {\it
Let $a$ be a vector and $\vT$ be a pair of even operators of $V$.
There exists a translation covariant $V$-valued power series $b(\vz)$ 
such that $b(0)=a$ if and only if 
$T^n\bT^{\bn}a=\bT^{\bn}T^n a$ for any $\vn\in\N^2$.
If $b(\vz)$ exists then it is equal to $e^{\vz\vT}a$.
}

\bigskip

\begin{pf}
If $b(\vz)$ exists then
$T^n\bT^{\bn}a=\del_z^n \del_{\bz}^{\bn} b(\vz)|_{\vz=0}
=\del_{\bz}^{\bn}\del_z^n b(\vz)|_{\vz=0}=
\bT^{\bn}T^n a$ for any $\vn\in\N^2$ 
and
$b(\vz)=e^{\vz\del_{\vw}}b(\vw)|_{\vw=0}
=e^{\vz\vT}b(\vw)|_{\vw=0}
=e^{\vz\vT}a$.
Conversely, if $T^n\bT^{\bn}a=\bT^{\bn}T^n a$ for any $\vn\in\N^2$
then $e^{\vz\vT}a$ is obviously a translation covariant power series 
such that $e^{\vz\vT}a|_{\vz=0}=a$.
\end{pf}

\bigskip

If $a(z)$ is a $V$-valued distribution then we define
$a(z+w):=e^{w\del_z} a(z)\in V\set{z}\pau{w}$
and $a(z+w)_{w>z}:=a(w+z)$.
Remark \ref{SS:power series expansions} shows that 
for any integer $n$ we have 
$(z+w)^n
=\sum_{m\in\N}\del_w^{(m)}(z+w)^n|_{w=0}w^m
=\sum_{m\in\N}\binom{n}{m}z^{n-m}w^m
=e^{w\del_z}z^n$.
This implies that 
for any integral $V$-valued distribution $a(z)$ we have 
$a(z+w)=\sum_{n\in\Z}a_n (z+w)^{-n-1}$.

\bigskip

{\bf Proposition.}\: {\it
An $A$-valued distribution $a(\vz)$ is translation covariant for 
a pair $\vT$ of even elements of $A$ if 
\begin{equation}
\label{E:transl conj distr}
e^{\vw\vT}\,
a(\vz)\,
e^{-\vw\vT}
\; =\;
a(\vz+\vw).
\end{equation}
The converse is true if $T$ and $\bT$ commute.
}

\bigskip

\begin{pf}
Considering the coefficients of $w$ and $\bw$ in \eqref{E:transl conj distr}
we obtain $[T,a(\vz)]=\del_z a(\vz)$ and $[\bT,a(\vz)]=\del_{\bz} a(\vz)$.
This proves the first claim.

Suppose that $a(\vz)$ is translation covariant and that
$T$ and $\bT$ commute.
Let $L$ and $R$ denote the left-hand and the
right-hand side of \eqref{E:transl conj distr}.
Thus $L$ and $R$ are power series in $\vw$.
The constant terms of $L$ and $R$ are both equal to $a(\vz)$.
We have $\del_w L=[T,L]$ and
$\del_w R=\del_w e^{\vw\del_{\vz}}a(\vz)=
\del_z e^{\vw\del_{\vz}}a(\vz)=[T,R]$.
Similarly, we have $\del_{\bw}L=[\bT,L]$ and 
$\del_{\bw}R=[\bT,R]$.
Thus \eqref{E:transl conj distr} follows from the Remark.
\end{pf}

\subsection{Derivatives}
\label{SS:derivative}

We prove that 
the pair of derivatives $(\del_z,\del_{\bz})$ 
is a translation operator of $\fg\set{\vz}$
and that 
the translation covariant distributions 
form an $\N$-fold submodule.

\bigskip

Because for any $V$-valued distribution $a(z)$
the residue of $\del_z a(z)$ is zero
the product formula implies the
{\bf integration-by-parts} formula
$$
\res_z(\del_z a(z)b(z))
\; =\;
-\res_z(a(z)\del_z b(z))
$$
for any $a(z)\in V\set{z}$ and $b(z)\in\K[z^{\K}]$.

Let $S$ be an even set, 
$M$ be an $S$-fold module over a vector space $V$,
and $W$ and $N$ be subsets of $V$ and $M$, resp.
For a subset $T$ of $S$,
we denote by $W_{(T)}N$ the span of the $s$-th products $a_{(s)}b$ 
where $a\in W, b\in N$, and $s\in T$.
A vector subspace $N$ of $M$ is called an {\bf $S$-fold submodule} 
\index{Sfold submodule@$S$-fold submodule}
over $W$ if $W_{(S)}N$ is contained in $N$.

We denote by $\del_{\vz}$
the pair $(\del_z,\del_{\bz})$ of derivatives of $V\set{\vz}$.

\bigskip

{\bf Proposition.}\: {\it
$\iti$\:
The pair of derivatives $\del_{\vz}$ 
is a translation operator of the
$\N$-fold module $\fg\set{\vz}$ over $\fg\pau{z\uppm}$.

$\itii$\:
If $T$ is a derivation of the algebra $\fg$
then $T$ defines a derivation of the $\N$-fold module $\fg\set{\vz}$.

$\itiii$\:
If $\vT$ is a pair of even derivations of the algebra $\fg$
then the subspace $\fg\set{\vz}_{\vT}$ of $\fg\set{\vz}$
is an $\N$-fold submodule over $\fg\pau{z\uppm}_{\vT}$.
}

\bigskip

\begin{pf}
$\iti$\:
The integration-by-parts formula yields
$$
\res_z((z-w)^n[\del_z a(z),b(\vw)])
\; =\;
-n\, 
\res_z((z-w)^{n-1} [a(z),b(\vw)]).
$$
Thus $\del_z$ is a holomorphic translation endomorphism.
The product formula implies
\begin{align}
\notag
\del_w 
&\res_z((z-w)^n [a(z),b(\vw)])
\\
\notag
=\;
-n\,
&\res_z((z-w)^{n-1} [a(z),b(\vw)])
\: +\:
\res_z((z-w)^n [a(z),\del_w b(\vw)]).
\end{align}
This shows that $\del_z$ is a 
holomorphic translation generator.
From 
$$
\del_{\bz}a(z)_{(n)}b(\vz)
\; =\;
0
\; =\;
a(z)_{(n,-2)}b(\vz)
$$
and
$$
\del_{\bw}
\res_z((z-w)^n [a(z),b(\vw)])
\; =\;
\res_z((z-w)^n [a(z),\del_{\bw}b(\vw)])
$$
we see that $\del_{\bz}$ is an anti-holomorphic
translation operator.

$\itii$\:
This follows from  
$$
T\res_z((z-w)^n [a(z),b(\vw)])
\; =\;
\res_z((z-w)^n T[a(z),b(\vw)]).
$$

$\itiii$\:
Parts $\iti$ and $\itii$ show that
$\del_{\vz}$ and $\vT$ 
are even derivations of the $\N$-fold module $\fg\set{\vz}$.
This implies the claim.
\end{pf}

\subsection{Dilatation Operator}
\label{SS:sfold dilat op}

We define the notion of a {\it dilatation operator} 
and show that there exists a one-to-one correspondence between 
dilatation operators and gradations.

\bigskip

Let $V$ be a vector space and $S$ be an even set.
Recall that an {\bf $S$-gradation} \index{gradation of a vector space} 
of $V$ is a family $(V_s)_{s\in S}$ of subspaces of $V$ such that
$V=\bigoplus_{s\in S} V_s$.
The subspaces $V_s$ are called the {\bf homogeneous} subspaces 
\index{homogeneous subspaces} of the {\bf $S$-graded} vector space $V$.
\index{graded vector space}

If $\vH$ is a pair of even commuting diagonalizable operators of $V$
then the simultaneous eigenspace decomposition 
$V=\bigoplus_{\vh\in\K^2}V_{\vh}$, where $V_{\vh}$ consists of simultaneous 
eigenvectors of $H$ and $\bH$ with eigenvalues $h$ and $\bh$, resp.,
defines a $\K^2$-gradation of $V$.
Thus we may identify pairs of even commuting diagonalizable operators of $V$
with $\K^2$-gradations of $V$.

Let $V$ be a $\K^2$-graded vector space. 
A vector $a$ of $V_{\vh}$ is called {\bf homogeneous} 
\index{homogeneous vector}
of {\bf weight} \index{weight} $\vh$ and we write $\vh(a):=\vh$.
If $a$ is homogeneous of weight $\vh$ then
$h-\bh$ is called the {\bf spin} \index{spin} and 
$h+\bh$ is called the {\bf scaling dimension} \index{scaling dimension} 
of $a$.

\bigskip

{\bf Definition.}\:                     
$\iti$\:                                
Let $M$ be a $\K$-fold module over a vector space $V$.
A pair $H=(H,H)$ of even diagonalizable operators of $V$ and $M$ 
is called a {\bf dilatation operator} \index{dilatation operator} if
$H(a_{(n)}b)=H(a)_{(n)}b+a_{(n)}H(b)-(n+1)\, a_{(n)}b$
for any $a\in V, b\in M$, and $n\in\K$.

$\itii$\:
Let $M$ be a $\K^2$-fold module over a vector space $V$.
A pair $H=(H,H)$ of even diagonalizable operators of $V$ and $M$ is called a
{\bf holomorphic} dilatation operator 
\index{dilatation operator!holomorphic}
if $H$ is a dilatation operator of the $\K$-fold module 
$M|_{\K\times\set{n}}$ for any $n\in\K$.
A pair $\vH=(H,\bH)$ consisting of a holomorphic and an anti-holomorphic
dilatation operator is called a {\bf dilatation operator} if 
$H$ and $\bH$ commute.

\bigskip

Let $M$ be a $\K$-fold module over a vector space $V$.
A pair $H=(H,H)$ of even diagonalizable operators of $V$ and $M$
is a dilatation operator if and only if 
$[H,a(\vz)]b=(Ha)(\vz)b+z\del_z a(\vz)b$ for any $a\in V$ and $b\in M$.

Let $M$ be a $\K^2$-fold module over a vector space $V$.
A pair $\vH$ of pairs of even commuting diagonalizable operators of $V$ and $M$
is a dilatation operator if and only if 
$V_{\vh}{}_{(\vn)}M_{\vh'}$ is contained in $M_{\vh+\vh'-\vn-1}$
for any $\vh, \vh', \vn\in\K^2$.
If the vector space $V$ is $\K^2$-graded then a $\K^2$-gradation of $M$ 
is called a {\bf gradation} of the $\K^2$-fold module $M$
if the corresponding pair $\vH$ is a dilatation operator of $M$.
A $\K^2$-fold module together with a dilatation operator is called
a {\bf graded} $\K^2$-fold module. \index{graded $\K^2$-fold module}
If $\vH$ is a dilatation operator of a $\K^2$-fold module $M$ then
the operator $H+\bH$ is called the {\bf Hamiltonian} \index{Hamiltonian}
of $M$.

Let $V$ be a graded $\K^2$-fold algebra and $a$ be a homogeneous state.
The weight $\vh:=\vh(a)$ of $a$ 
is also called the {\bf conformal weight} \index{conformal weight} of $a$.
For $\vn\in\K^2$, we define
$a_{\vn}:=a_{(\vn-1+\vh)}$ so that 
$Y(a,\vz)=\sum_{\vn\in\K^2}a_{\vn} \vz^{-\vn-\vh}$.
The operator $a_{\vn}$ is homogeneous of weight $-\vn$. 

Let $b$ be a state.
From $H a_{\vn}(b)=a_{\vn}(Hb)-n\, a_{\vn}(b)$ we get
$[H,a_{\vn}]=-n\, a_{\vn}$.
If $T$ is a translation generator
then $(Ta)_{\vn}=-(n+h(a))a_{\vn}$.
If $a$ is holomorphic then the commutator formula 
for the new indexing of modes reads
$$
[a_n,b_{\vm}]
\; =\;
\sum_{i\in\N}\:
\binom{n+h(a)-1}{i}\, (a_{(i)}b)_{n+m,\bm}.
$$

\bigskip

{\bf Remark.} \: {\it
Let $L$ be an even holomorphic state of a bounded $\K^2$-fold algebra 
such that $L$ and any two states satisfy the commutator formula.
If $L_{(0)}$ is a translation generator and
$L_{(1)}$ is diagonalizable then $L_{(1)}$ is a dilatation operator.
}

\bigskip

\begin{pf} 
For any state $a$ and $\vn\in\K^2$, the commutator formula yields
\begin{equation}
\notag
[L_{(1)},a_{(\vn)}]
\; =\;
(L_{(0)}a)_{(n+1,\bn)}
\, +\,
(L_{(1)}a)_{(\vn)}
\; =\;
-(n+1)\, a_{(\vn)}
\, +\,
(L_{(1)}a)_{(\vn)}.
\end{equation}
\end{pf}

\subsection{Dilatation Covariance}
\label{SS:dilat cov}

We give an infinitesimal definition of 
the notion of {\it dilatation covariance}
and prove that it is equivalent to its integrated version.

\bigskip

{\bf Definition.} \:
Let $\vH$ be a pair of even operators of $V$ and $\vh\in\K^2$.
A $V$-valued distribution $a(\vz)$ is called
{\bf dilatation covariant} \index{dilatation covariant} 

\begin{enumerate}
\item[$\iti$]
of weight $h$ for $H$ if $Ha(\vz)=h\, a(\vz)+z\del_z a(\vz)$;

\item[$\itii$]
of weight $\vh$ for $\vH$ if 
$a(\vz)$ is dilatation covariant for $H$ of weight $h$ and
$\bH a(\vz)=\bh\, a(\vz)+\bz\del_{\bz} a(\vz)$.
\end{enumerate}

\bigskip

A distribution $a(\vz)$ is dilatation covariant of weight $h$ if and only if
for any $\vn\in\K^2$ the mode $a_{\vn}$ is an eigenvector of $H$
with eigenvalue $h-n-1$.

Let $V$ be a vector space endowed with a pair $\vH$ of even operators.
We denote by $V\set{\vz}_{\vh}$ the vector space of 
$V$-valued distributions that are dilatation covariant of weight $\vh$ 
and we denote by $V\set{\vz}_{\vast}$ the direct sum of these spaces.
We denote by $V\pau{z\uppm}_h$ the vector space of holomorphic
$V$-valued distributions that are dilatation covariant of weight 
$(h,0)$ and we denote by $V\pau{z\uppm}_{\ast}$ their direct sum.
An $V$-valued distribution is called dilatation covariant for $\vH$
if it is contained in $V\set{\vz}_{\vast}$.
The derivative $\del_z$ maps $V\set{\vz}_{\vh}$ to $V\set{\vz}_{h+1,\bh}$.

An $A$-valued distribution $a(\vz)$ is called
dilatation covariant for a pair $\vH$ of even elements of $A$
if $a(\vz)$ is dilatation covariant for $[\vH,\;\, ]$.

Let $\vH$ be a pair of commuting diagonalizable even operators of $V$.
We denote by $\End(V)_{\vh}$ the space of operators of $V$ that map
$V_{\vh'}$ to $V_{\vh'+\vh}$ for any $\vh'\in\K^2$.
An $\End(V)$-valued distribution $a(\vz)$ is 
dilatation covariant of weight $\vh$ if and only if 
$a_{\vn}$ is contained in $\End(V)_{\vh-\vn-1}$
for any $\vn\in\K^2$.

Let $M$ be a $\K^2$-fold module over a vector space $V$.
A pair $H=(H,H)$ of diagonalizable even operators of $V$ and $M$
is a dilatation operator if and only if for any $a\in V$
the distribution $a(\vz)$ is dilatation covariant for $H$.
A pair $\vH$ of pairs of commuting diagonalizable even operators of $V$ and $M$
is a dilatation operator if and only if 
the state-field correspondence is a morphism of 
$\K^2$-graded vector spaces from $V$ to $\End(M)\set{\vz}_{\vast}$.

If $a(z)$ is a $V$-valued distribution then we define 
$a(wz):=\sum_{n\in\K}a_n z^{-n-1}\linebreak[0]w^{-n-1}$.
Let $H$ be a diagonalizable even operator of a vector space $V$.
We define an operator $w^H\equiv e^{H\ln w}$ of $V[w^{\K}]\set{z}$
by $a(z,w)\mapsto \sum_{n,m,h\in\K}a_{n,m}^h z^{-n-1}\linebreak[0]w^{-m-1+h}$
where $a_{n,m}=\sum_h a_{n,m}^h$ is the decomposition into
eigenvectors so that $H a_{n,m}^h=h a_{n,m}^h$.

\bigskip

{\bf Proposition.}\: {\it
$\iti$\:
Let $\vH$ be a pair of commuting diagonalizable even operators of $V$ 
and $\vh\in\K^2$.
An $\End(V)$-valued distribution $a(\vz)$ is dilatation covariant
of weight $\vh$ if and only if 
\begin{equation}
\label{E:dilat cov}
\vw^{\vH}\, a(\vz)\, \vw^{-\vH}
\; =\;
\vw^{\vh}\, a(\vw\vz).
\end{equation}

$\itii$\:
Let $M$ be a $\K^2$-fold module over a vector space $V$.
A pair $\vH$ of pairs of commuting diagonalizable even operators of $V$ and $M$
is a dilatation operator if and only if for any $a\in V$ we have
$$
\vw^{\vH}\, a(\vz)\, \vw^{-\vH}
\; =\;
(\vw^{\vH}a)(\vw\vz).
$$
}

% no \bigskip

\begin{pf}
$\iti$\:
If $b$ is a homogeneous vector of $V$ of weight $\vh'$ and 
$\sum_{\vk\in\K^2}(a_{\vn}b)_{\vk}$ is the decomposition
of $a_{\vn}b$ into eigenvectors then \eqref{E:dilat cov}
applied to $b$ is the identity
$$
\sum_{\vn, \vk\in\K^2}\:
(a_{\vn}b)_{\vk}\, \vz^{-\vn-1}\vw^{\vk-\vh'}
\; =\;
\sum_{\vn\in\K^2}\:
a_{\vn}b\, \vz^{-\vn-1}\vw^{\vh-\vn-1}.
$$
This is equivalent to $(a_{\vn}b)_{\vk}=\de_{\vk,\vh+\vh'-\vn-1}a_{\vn}b$.

$\itii$\:
This follows from $\iti$.
\end{pf}

\subsection{Gradation for Dilatation Covariant Distributions}
\label{SS:grad dilat cov}

We prove that the space of dilatation covariant distributions
is a graded $\N$-fold module.

\bigskip

For a subset $S$ of $\K$, 
we identify $S$-graded vector spaces with $\K^2$-graded vector spaces 
via the inclusion $S\to\K^2, h\mapsto (h,0)$.

\bigskip

{\bf Proposition.}\: {\it
If $\vH$ is a pair of even derivations of $\fg$ then
the space $\fg\set{\vz}_{\vast}$ is a graded $\N$-fold module
over $\fg\pau{z\uppm}_{\ast}$.
}

\bigskip

\begin{pf}
For $a(z)\in\fg\pau{z\uppm}_h$, $b(\vz)\in\fg\set{\vz}_{\vh'}$,
and $n\in\N$, we have
\begin{align}
\notag
&H(a(w)_{(n)}b(\vw))
\\
\notag
=\; 
&\res_z((z-w)^n\, H[a(z),b(\vw)])
\\
\notag
=\;
&(h+h')\, a(w)_{(n)}b(\vw)
\; +\;
\res_z((z-w)^n\, (z\del_z+w\del_w)\, [a(z),b(\vw)]).
\end{align}
Thus the claim follows from
\begin{align}
\notag
&\res_z((z-w)^n\, z\, [\del_z a(z),b(\vw)])
\\
\notag
=\;
&\res_z(((z-w)^{n+1}+(z-w)^n w)\, [\del_z a(z),b(\vw)])
\\
\notag
=\;
&\del_w a(w)_{(n+1)}b(\vw)
\; +\;
w\del_w a(w)_{(n)}b(\vw)
\end{align}
and from the fact that $\del_z$ is a translation operator.
\end{pf}

\subsection{The Witt Algebra}
\label{SS:witt alg}

We define the {\it Witt algebra}.

\bigskip

The general linear group $\GL(2,\C)$ acts on the complex projective line
$\CP^1$ by M\"obius transformations \index{M\"obius transformation}
$$
\begin{pmatrix}
a & b \\
c & d
\end{pmatrix}
\;\mapsto\;
\bigg(
z\mapsto
\frac{az+b}{cz+d}
\bigg).
$$
Thus one obtains an isomorphism between
$\PSL(2,\C):=\SL(2,\C)/\set{\pm 1}$ and the automorphism group of $\CP^1$.
The corresponding isomorphism of Lie algebras from 
$\fsl_2$ to the Lie algebra of global holomorphic vector fields on $\CP^1$
is given by 
$$
\begin{pmatrix}
0 & 1 \\
0 & 0
\end{pmatrix}
\mapsto -\del_z,\qquad
\frac{1}{2}
\begin{pmatrix}
1 & 0 \\
0 & -1
\end{pmatrix}
\mapsto -z\del_z,\qquad
\begin{pmatrix}
0 & 0 \\
-1 & 0
\end{pmatrix}
\mapsto -z^2\del_z.
$$
The vector field $-\del_z$ generates translations,
$-z\del_z$ generates dilatations, and
$-z^2\del_z$ generates {\it special conformal transformations}.
\index{special conformal transformations}

The Lie algebra of global meromorphic vector fields on $\CP^1$
with poles only at zero and infinity is called the
{\bf Witt algebra} \index{Witt algebra} and is denoted by $\Witt$.
The vector fields $l_n:=-z^{n+1}\del_z$ for $n\in\Z$ form a 
basis of $\Witt$. We have $[l_n,l_m]=(n-m)l_{n+m}$
for any integers $n$ and $m$.

\subsection{The Virasoro Algebra}
\label{SS:vir alg}

We define the {\it Virasoro algebra} as an explicit central extension
of the Witt algebra and prove that it is the universal central extension.

\bigskip

{\bf Remark.}\: {\it
The map $\varep:\bigwedge^2\Witt\to\K, 
l_n\wedge l_m\mapsto (n^3-n)\de_{n+m,0}/12$, is 
a Chevalley-Eilenberg 2-cocycle on $\Witt$, 
i.e.~for any elements $a, b, c$ of $\Witt$ we have
$$
\varep([a,b],c)
\; -\;
\varep([a,c],b)
\; +\;
\varep([b,c],a)
\; =\;
0.
$$
}

% no bigskip

\begin{pf}
We have to show that for any integers $n, m, k$ the number
\begin{align}
\notag
\de_{n+m+k,0}\Big((n -m)\big( (n+m)^3-(n+m)\big)
\; -\;
&(n-k)\big( (n+k)^3-(n+k)\big)
\\
\notag
+\;
&(m-k)\big( (m+k)^3-(m+k)\big)\Big)
\end{align}
vanishes. This is clear if $n+m+k$ is non-zero.
If $n+m+k$ is zero then this number is equal to $1/6$ times 
\begin{align}
\notag
&(n-m)((n+m)^3-(n+m))
\: +\:
(2n+m)(m^3-m)
\: -\:
(2m+n)(n^3-n)
\\
\notag
=\,
&n^4{}_a +
3n^3m\, {}_b +
3n^2m^2{}_c+
nm^3{}_d 
\; -\;
n^3m\, {}_b -
3n^2m^2{}_c -
3nm^3{}_d -
m^4{}_e
\\
\notag
&-
n^2{}_f +
m^2{}_g
\\
\notag 
&+
2nm^3{}_d -
2nm\, {}_h +
m^4{}_e -
m^2{}_g
\; -\;
2mn^3{}_b +
2mn\, {}_h -
n^4{}_a +
n^2{}_f
\\
\notag
=\,
&0
\end{align}
where we have indicated which terms cancel against each other.
\end{pf}

\bigskip

The central extension of the Witt algebra corresponding to the 
2-cocycle $\varep$ is called the {\bf Virasoro algebra} \index{Virasoro algebra}
and is denoted by $\Vir=\Witt\oplus\K\hc$.

\bigskip

{\bf Proposition.}\: {\it
The Virasoro algebra is the universal central extension
of the Witt algebra.
}

\bigskip

\begin{pf}
Let $\fg$ be a central extension of the Witt algebra by a subalgebra $\fh$. 
Choose inverse images $L_n$ of $l_n$ in $\fg$ and define the elements 
$c_{n,m}$ of $\fh$ by 
$$
[L_n,L_m]
\; =\;
(n-m)\, L_{n+m}
\; +\;
c_{n,m}.
$$
Skew-symmetry implies that $c_{n,m}=-c_{m,n}$. 
The vector space $\fg$ is the direct sum of $\fh$ and 
the span of the vectors $L_n$ for $n\in\Z$.
The component in $\fh$ of $[L_k,[L_n,L_m]]$ is equal to $(n-m)c_{k,n+m}$.
Taking the components in $\fh$ of both sides of the Jacobi identity
$$
[L_0,[L_n,L_m]]
\; =\;
[[L_0,L_n],L_m]
\; +\;
[L_n,[L_0,L_m]]
$$
we obtain $(n-m)c_{0,n+m}=-(n+m)c_{n,m}$.
Defining $L'_n:=L_n-c_{0,n}/n$ if $n$ is non-zero and
$L'_0:=L_0+c_{1,-1}/2$ we get
$$
[L'_n,L'_m]
\; =\;
(n-m)\, L'_{n+m}
\; +\;
\de_{n+m,0}\, c_n
$$
where $c_n:=c_{n,-n}-n\, c_{1,-1}$. We have $c_{\pm 1}=c_0=0$.
Taking the components in $\fh$ of both sides of the Jacobi identity
$$
[L'_1,[L'_n,L'_{-n-1}]]
\; =\;
[[L'_1,L'_n],L'_{-n-1}]
\; +\;
[L'_n,[L'_1,L'_{-n-1}]]
$$
yields $0=-(n-1)c_{n+1}+(n+2)c_n$.
Dividing this last equation by $(n-1)n(n+1)(n+2)$ we see that
$c:=2c_n/\binom{n+1}{3}$ for $|n|\geq 2$ does not depend on $n$.
Thus the map $\Vir\to\fg, L_n\mapsto L'_n, \hc\mapsto c$, 
is a morphism of Lie algebras that is compatible with the morphisms 
to $\Witt$.
\end{pf}

\subsection{Conformal Distribution}
\label{SS:conf distr}

We define the notion of a {\it conformal distribution}
and express translation and dilatation covariance in terms of an OPE.

\bigskip

If we denote a holomorphic distribution by $L(z)$ then
we write $L_{n+1}$ for the $n$-th mode of $L(z)$ so that
$L(z)=\sum_{n\in\Z}L_n z^{-n-2}$.
This deviation from our standard notation should not cause confusion
and will be justified later.

\bigskip

{\bf Definition.} \:
$\iti$\:
For $c\in\K$, a holomorphic $\End(V)$-valued distribution $L(z)$
is called a {\bf Virasoro distribution} \index{Virasoro distribution}
of {\bf central charge} $c$
if the map $\Vir\to\End(V), L_n\mapsto L_n, \hc\mapsto c$, 
is a representation of the Virasoro algebra.

$\itii$\:
Let $\cF$ be a subset of $\End(V)\set{\vz}$ and $c\in\K$.
A holomorphic distribution $L(z)$ of $\cF$ 
is called a {\bf conformal distribution} \index{conformal distribution}
of $\cF$ of central charge $c$
if $L(z)$ is a Virasoro distribution of central charge $c$ and
$\cF$ is translation covariant for $L_{-1}$ and dilatation covariant
for $L_0$.

\bigskip

For $\vc\in\K^2$,
a pair $\vL(\vz)=(L(z),\bL(\bz))$ of a holomorphic and an anti-holomorphic 
distribution of $\cF$ is called a pair of conformal distributions 
\index{pair of conformal distributions}
of central charges $\vc$ if $L(z)$ and $\bL(\bz)$ are Virasoro distributions 
of central charges $c$ and $\bc$, resp., and 
$\cF$ is translation covariant for $(L_{-1},\bL_{-1})$ and 
dilatation covariant for $(L_0,\bL_0)$.

Let $V$ be a $\K^2$-fold algebra.
A state $L$ is called a {\bf Virasoro vector} \index{Virasoro vector}
if $L(z)$ is a Virasoro distribution.
A state $L$ is called a {\bf conformal vector} \index{conformal vector}
if $L(z)$ is a conformal distribution of $\cF_Y$.

\bigskip

{\bf Remark.}\: {\it
Let $L(z)$ be a holomorphic $\End(V)$-valued distribution and $h\in\K$.
An $\End(V)$-valued distribution $a(\vz)$ is translation covariant for 
$L_{-1}$ and dilatation covariant of weight $h$ for $L_0$ if and only if
we have
$$
L(z)a(\vw)
\; \sim \;
\dots
\; +\;
\frac{h\, a(\vw)}{(z-w)^2}
\; +\; 
\frac{\del_w a(\vw)}{z-w}.
$$
}

% no \bigskip

\begin{pf}
This follows from $L(z)_{(0)}a(\vz)=[L_{-1},a(\vz)]$ and 
$$
L(z)_{(1)}a(\vz)
\; =\;
-z[L_{-1},a(\vz)]
\; +\;
[L_0,a(\vz)].
$$
\end{pf}

\section{Holomorphic Locality}
\label{S:locality distr}

{\bf Summary.}\:
In section \ref{SS:delta distr}
we study the delta distribution.
In section \ref{SS:derivative taylor}
we prove Taylor's formula.

In section \ref{SS:locality}
we give three equivalent formulations of locality
for a holomorphic distribution.
In section \ref{SS:locality of formal distributions}
we define the notion of holomorphic locality 
and relate it to the commutator formula.
In section \ref{SS:locality and ope}
we show that locality is equivalent to 
the existence of finite OPEs for $z>w$ and $w>z$. 

In section \ref{SS:skew-symmetry distr}
we prove that holomorphically local integral distributions 
satisfy holomorphic skew-symmetry.
In section \ref{SS:Dong}
we prove Dong's lemma for holomorphic locality
which states that the $n$-th products of $\fg\set{\vz}$
preserve holomorphic locality.
In section \ref{SS:vir ope}
we express the property of being a Virasoro distribution
in terms of an OPE.

\medskip

{\bf Conventions.}\:
Except when we discuss $S$-fold modules, 
we denote by $V$ a vector space,
$\fg$ is an algebra with multiplication
$a\otimes b\mapsto [a,b]$,
and $A$ is an associative algebra
with multiplication $a\otimes b\mapsto ab$ 
and commutator $[a,b]:=ab-\paraab ba$.

\subsection{The Delta Distribution}
\label{SS:delta distr}

We define the {\it delta distribution}
and prove an equation that relates 
the delta distribution to $(z-w)^n$ and
$T_{w,z}((z-w)^n)$.

\bigskip

There exists a na\-tural isomorphism
$V\pau{z\uppm}\tra \Vect(\K[z\uppm],V), a(z)\mapsto\al$, where 
$\al:p(z)\mapsto\res_z(p(z)a(z))$.
The distribution $a(z)$ is called the {\bf kernel} 
\index{kernel distribution} of $\al$.
We have $a(z)=\sum_{n\in\Z}\al(z^n)z^{-n-1}$.

The kernel of the identity of $\K[z\uppm]$
is called the {\bf delta distribution} \index{delta distribution}
and is denoted by $\de(z,w)=\sum_{n\in\Z}w^n z^{-n-1}$.
In appendix \ref{SS:kernel delta}
a delta distribution $\de_C$ is defined
for any commutative algebra $C$ 
and five properties of $\de_C$ are proven.
We recall these properties in the case that
$C=\K[z\uppm]$ where $\de_C=\de(z,w)$.

For any $a(z,w)\in\K[z\uppm,w\uppm]$ we have
\begin{equation}
\label{E:delta is delta}
\de(z,w)a(z,w)
\; =\; 
\de(z,w)a(w,w).
\end{equation}
An integral distribution $a(z)$ with finite support
is called a {\bf Laurent polynomial}.
The vector space of $V$-valued Laurent polynomials
is denoted by $V[z\uppm]$.
By projecting onto a basis of $V$ 
we see that \eqref{E:delta is delta}
also holds for any 
$a(z,w)\in V[z\uppm,w\uppm]$.
\index{Laurent polynomial}

If $c(z)\in V\set{z}$ and $n\in\K$ 
then $c(w)\del_w^{(n)}\de(z,w)$ is the kernel of
the morphism 
$c(z)\del_z^{(n)}:\K[z\uppm]\to V\set{z}$.
We have
\begin{equation}
\label{E:delta distr symmetric}
\del_w \de(z,w)
\; =\;
-\del_z \de(z,w),
\end{equation}
\begin{equation}
\label{E:delta properties}
(z-w)\del_w^{(n)}\de(z,w)
\; =\;
\del_w^{(n-1)}\de(z,w),
\end{equation}
and $\de(z,w)=\de(w,z)$.
Equation \eqref{E:delta distr symmetric} implies
\begin{equation}
\label{E:delta distr homogeneous}
\de(z-x,w)
\; =\;
e^{-x\del_z}\de(z,w)
\; =\;
e^{x\del_w}\de(z,w)
\; =\;
\de(z,w+x).
\end{equation}

For an integer $n$ and $\mu, \nu\in\K$, 
we define
$(\mu z+\nu w)_{w>z}^n:=(\nu w+\mu z)^n\equiv
T_{w,z}((\mu z+\nu w)^n)$.

\bigskip

{\bf Remark.}\: {\it
For any integer $n$ we have
\begin{equation}
\label{E:delta and z-w}
\del_w^{(n)}\de(z,w)
\; =\;
(z-w)^{-n-1}
\; -\;
(z-w)^{-n-1}_{w>z}.
\end{equation}
}

% no bigskip

\begin{pf}
Equation \eqref{E:delta and z-w} holds true if $n$ is negative
because $(z-w)^m=(z-w)_{w>z}^m$ for any non-negative integer $m$.
It holds true for $n=0$ because of \eqref{E:z>w expansion}
and $(z-w)_{w>z}^{-1}=-(w-z)^{-1}$.
If we apply $\del_w$ to \eqref{E:delta and z-w} for $n=m\in\N$
then we obtain \eqref{E:delta and z-w} for $n=m+1$ times $m+1$.
Thus the claim follows by induction.
\end{pf}

\bigskip

From
$\de(z,w+x)=\sum_{n\in\Z}\del_w^{(n)}\de(z,w)x^n$
and 
$$
\de(x,z-w)
-
\de(x,z-w)_{w>z}
\; =\;
\sum_{n\in\Z}\:
((z-w)^{-n-1}-(z-w)^{-n-1}_{w>z})\,
x^n
$$
we see that
\eqref{E:delta and z-w} is equivalent to 
\begin{equation}
\label{E:jacobi for id}
\de(z,w+x)
\; =\;
\de(x,z-w)
\; -\;
\de(x,z-w)_{w>z}.
\end{equation}
In Remark \ref{SS:h Jac id distr}\,$\iti$
we will see that \eqref{E:jacobi for id} is 
the holomorphic Jacobi identity for $1, 1, 1$
where $1$ is the identity of a unital $\K^2$-fold algebra.

Because 
$(z-w)^{-1}\in z^{-1}\Z[z^{-1}]\pau{w}$
and
$(z-w)_{w>z}^{-1}\in w^{-1}\Z[w^{-1}]\pau{z}$
and
because $T_{z,w}$ and $T_{w,z}$ are algebra morphisms
\eqref{E:delta and z-w} implies that for any non-negative integer $n$ we have
\begin{equation}
\label{E:delta pm z-w ><}
\del_w^{(n)}\de(z,w)|_{\N\times\Z}
=
(z-w)^{-n-1},
\quad
\del_w^{(n)}\de(z,w)|_{\Z_<\times\Z}
=
-(z-w)^{-n-1}_{w>z}.\;\;
\end{equation}

\subsection{Taylor's Formula}
\label{SS:derivative taylor}

We prove that a vertex series can be expanded
with respect to an integral variable
as a finite Taylor series plus a remainder.

\bigskip

We denote by $V[(z/w)^{\K}]\set{z}$
the vector space of distributions $a(z,w)$
such that for any $n\in\K$
the distribution $\sum_{m\in\K}a_{m,n-m}x^m$ has finite support.
The morphism 
$$
V[(z/w)^{\K}]\set{z}
\;\to\;
V\set{z},
\quad
a(z,w)
\;\mapsto\;
\sum_{n\in\K}\:
\Big(\sum_{m\in\K} a_{m,n-m-1}\Big)\: z^{-n-1},
$$
is written $a(z,w)\mapsto a(z,z)\equiv a(z,w)|_{w=z}$.
We define
$V(\!(z,w\rangle :=V\set{z}\sqbrack{w}\cap V\set{w}\lau{z}$.

\bigskip

{\bf Lemma.}\: {\it
If $a(z,w)$ is a vertex series in $V(\!(z,w\rangle$
then there exists a unique vertex series $r(z,w)$
in $V(\!(z,w\rangle$ such that
$$
a(z,w)
\; =\;
a(w,w)
\; +\;
(z-w)r(z,w).
$$
Moreover, 
if $a(z,w)$ is a Laurent polynomial then
so is $r(z,w)$.
}

\bigskip

\begin{pf}
By projecting onto a basis of $V$ we may assume that $V=\K$.
Furthermore, by considering the restriction of 
$w^h a(z,w)$ to $\Z^2$ for any $h\in\K$ we may assume that 
$a(z,w)\in\K\lau{z,w}$. Thus we have to show that if 
$b(z,w)\in\K\lau{z,w}$ vanishes for $z=w$ then $b(z,w)$
is divisible by $z-w$. 
For this we may assume that $b(z,w)\in\K\pau{z,w}$. 
Write $b(z,w)=\sum_{n\in\N}b_n(z,w)$ where 
$b_n(z,w)\in\K[z,w]$ is the homogeneous component of
$b(z,w)$ of degree $n$. 
Then $b_n(w,w)=0$ for any $n$. 
Thus $b_n(z,w)$ is divisible by $z-w$ and hence so is
$b(z,w)$.

The second claim follows from the arguments we already gave.
\end{pf}

\bigskip

{\bf Proposition.} (Taylor's Formula)  \: {\it
If $a(z,\vw)$ is a vertex series in $V(\!(z,\vw\rangle$
and $N$ is a non-negative integer 
then there exists a unique vertex series $r(z,\vw)$
in $V(\!(z,\vw\rangle$ such that
$$
a(z,\vw)
\; =\;
\sum_{n=0}^{N-1}\:
\del_z^{(n)}a(z,\vw)|_{z=w} (z-w)^n
\; +\;
(z-w)^N\, r(z,\vw).
$$
Moreover, 
if $a(z,\vw)$ is a Laurent polynomial then
so is $r(z,\vw)$.
}

\bigskip

\begin{pf}
We may assume that $a(z,\vw)\in V(\!(z,w\rangle$.
By applying the Lemma $N$ times we obtain
$c_i(z)\in V\sqbrack{z}$
and $r(z,w)\in V(\!(z,w\rangle$ such that
$$
a(z,w)
\; =\;
\sum_{n=0}^{N-1}\:
c_n(w)\, (z-w)^n
\; +\;
(z-w)^N r(z,w).
$$
Acting with $\del_z^{(n)}$ on this equation
and setting $z=w$ we get $\del_z^{(n)}a(z,w)|_{z=w}=c_n(w)$.
Because $V\lau{z}\sqbrack{w}$ is a vector space over $\K(z,w)$
the distribution $r(z,w)$ is unique.
\end{pf}

\subsection{Locality for a Holomorphic Distribution}
\label{SS:locality}

We define the notion of {\it locality for a holomorphic distribution}
and prove that $a(z,\vw)$ is local if and only if 
$a(z,\vw)$ is the kernel of a linear differential operator.
This is an algebraic analogue
of the theorem of Peetre 
\cite{peetre.thm}
which states that 
a linear operator acting on functions is local 
if and only if 
it is a differential operator.

\bigskip

{\bf Definition.}\:
For a non-negative integer $N$,
a distribution $a(z,\vw)$ in $V\pau{z\uppm}\set{\vw}$ 
is called {\bf local} of order at most $N$
\index{local holomorphic distribution}
if 
$$
(z-w)^N a(z,\vw)
\; =\;
0.
$$

\bigskip

If $a(z,\vw)$ is local of order $N$ then 
$\del_z a(z,\vw)$ and $\del_w a(z,\vw)$ are local 
of order at most $N+1$.

\bigskip

{\bf Proposition.\;}{\it
For a non-negative integer $N$,
a distribution $a(z,\vw)$ in $V\pau{z\uppm}\set{\vw}$ 
is local of order $N$
if and only if
there exist $V$-valued distributions 
$c^0(\vz), \dots, c^{N-1}(\vz)$ 
such that $c^{N-1}(\vz)$ is non-zero and 
\begin{equation}
\label{E:sum of del de}
a(z,\vw)
\; =\;
\sum_{n=0}^{N-1}\: 
c^n(\vw)\, \del_w^{(n)}\de(z,w).
\end{equation}
Moreover, \eqref{E:sum of del de} is equivalent to
\begin{equation}
\label{E:Lie of modes}
a_{m,\vk}
\; =\;
\sum_{n=0}^{N-1}\: 
\binom{m}{n}\, c^n_{m+k-n,\bk}
\end{equation}
for any $m\in\Z$ and $\vk\in\K^2$.
Equation \eqref{E:sum of del de} implies that
$c^n(\vw)$ is equal to $\res_z((z-w)^n a(z,\vw))$.
}

\bigskip

\begin{pf}
Suppose that $a(z,\vw)$ is local of order $N$.
Let $\al:\K[z\uppm,w\uppm]\to V\set{\vw}$
be the morphism of $\K[w\uppm]$-modules 
whose kernel is $a(z,\vw)$.
Because of Taylor's formula
there exists for any Laurent polynomial $p(z,w)$
a Laurent polynomial $r(z,w)$ such that
\begin{align}
\notag
&\al(p(z,w))
\\
\notag
=\;
&\res_z\Bigg( \Bigg(
\sum_{n=0}^{N-1}\:
\del_z^{(n)}p(z,w)|_{z=w}(z-w)^n
\; +\;
(z-w)^N r(z,w)\Bigg) \: a(z,\vw)\Bigg) 
\\
\notag
=\;
&\sum_{n=0}^{N-1}\:
\res_z((z-w)^n a(z,\vw))\;\;
\del_z^{(n)}p(z,w)|_{z=w}.
\end{align}
This shows that the kernel $a(z,\vw)$ of $\al$
is given by \eqref{E:sum of del de} with
$c^n(\vw)=\res_z((z-w)^n a(z,\vw))$.
We have $c^{N-1}(\vz)\ne 0$ because otherwise $a(z,\vw)$
would be local of order at most $N-1$ by \eqref{E:delta properties}.

The converse statement
follows from \eqref{E:delta properties}.
Equations
\eqref{E:sum of del de} and \eqref{E:Lie of modes}
are equivalent because 
$\del_w^{(n)}\de(z,w)=
\sum_{m\in\Z}\binom{m}{n}w^{m-n} z^{-m-1}$
for any integer $n$.
The last statement follows from Lemma \ref{SS:n-fold module holom distr}
by restricting \eqref{E:sum of del de} to $\N\times\K^2$
and making use of \eqref{E:delta pm z-w ><}.
\end{pf}

\bigskip

{\bf Corollary.}\: {\it
The only local distribution in $V\set{\vw}\lau{z}$  
is the zero distribution.
The same is true for the spaces 
$V\set{\vw}\lau{z^{-1}}, V\pau{z\uppm}\sqbrack{\vw}$,
and 
$V\pau{z\uppm}\sqbrack{\vw^{-1}}$.
}

\bigskip

\begin{pf}
Let $a(z,\vw)$ be a local distribution in
$V\set{\vw}\lau{z}$ or one of the other three spaces.
We may assume that 
$a(z,\vw)$ is local of order one.
The Proposition implies that
$a(z,\vw)=c(\vw)\de(z,w)$
for some $c(\vz)\in V\set{\vz}$.
Because 
$c(\vw)\de(z,w)=\sum_{n\in\Z}c(\vw)w^n z^{-n-1}$
we obtain $c(\vz)=0$ and hence $a(z,\vw)=0$.
${}_{}$
\end{pf}

\subsection{Holomorphic Locality}
\label{SS:locality of formal distributions}

We define the notion of {\it holomorphic locality},
relate this notion to the identity of holomorphic locality,
remark that for $\fg$-valued distributions the commutator formula
is equivalent to holomorphic locality,
and observe again that for $\K^2$-fold modules 
the commutator formula implies holomorphic locality.

\bigskip

{\bf Definition.}\:
For a non-negative integer $N$,
a holomorphic $\fg$-valued distribution $a(z)$
and a $\fg$-valued distribution $b(\vz)$
are called 
{\bf holomorphically local} \index{holomorphically local}
of order $N$ if the commutator $[a(z),b(\vw)]$ is local of order $N$.

\bigskip

If $a(z)$ and $b(\vz)$ 
are holomorphically local of order $N$ then
$a(z)$ is bounded on $b(\vz)$ by $N$.

For a $\K^2$-fold module $M$ over a vector space $V$,
vectors $a$ and $b$ of $V$ are called 
holomorphically local on $M$
if $a$ is holomorphic on $M$ and
$a(z)$ and $b(\vz)$ are holomorphically local.
Two states of a $\K^2$-fold algebra $V$ 
are called holomorphically local if 
they are holomorphically local on $V$.

\bigskip

{\bf Remark.}\:
$\iti$\:
Proposition \ref{SS:locality} shows that
a holomorphic $\fg$-valued distribution $a(z)$
and a $\fg$-valued distribution $b(\vz)$
are holomorphically local of order $N$ if and only if 
$a(z)$ and $b(\vz)$ satisfy the commutator formula
for some coefficients $(c^n(\vz))_n$ such that
$c^n(\vz)$ is zero for $n\geq N$ and 
$c^{N-1}(\vz)$ is non-zero.
In this case we have
$c^n(\vz)=a(z)_{(n)}b(\vz)$ for any $n$.

$\itii$\:
Let $M$ be a $\K^2$-fold module over a $\Z$-fold algebra $V$,
$a$ and $b$ be states of $V$ such that $a$ is holomorphic on $M$,
$c$ be a vector of $M$,
and $r$ be a non-negative integer.  
The vectors $a, b, c$ satisfy holomorphic locality of order at most $r$ 
if and only if 
the distribution $[a(z),b(\vw)]c$ is local of order at most $r$.
Thus holomorphic locality of order at most $r$ is satisfied
for $a, b$, and any $d\in M$
if and only if
$a$ and $b$ are holomorphically local on $M$
of order at most $r$.
In particular, if $V$ and $M$ are bounded then
the holomorphic Jacobi identity for $M$ implies that
$V$ is holomorphically local on $M$
because of Remark \ref{SS:hol jacobi}\,$\iti\itb$ and 
because $o'(a,b)$ does not depend on $c$.

Moreover, 
Proposition \ref{SS:locality} implies the following result.
Assume that $a$ is bounded on $b$.
The vectors $a, b, c$ satisfy the commutator formula
if and only if 
$[a(z),b(\vw)]c$ is a local distribution  
and $(a_{(n)}b)(\vz)c=a(z)_{(n)}b(\vz)c$
for any non-negative integer $n$.
Thus if $a, b, c$ satisfy the commutator formula then 
they satisfy holomorphic locality of order at most $o'(a,b)_+$.
This result was already obtained in 
Proposition \ref{SS:pf hol jacobi}\,$\itii$.

\subsection{Holomorphic Locality and OPE}
\label{SS:locality and ope}
We prove that holomorphic locality is equivalent to 
the existence of two finite OPEs for $z>w$ and $w>z$
that have the same OPE-coefficients and the same remainder.
Moreover, we show that
the normal ordered product is uniquely
determined as the remainder of these two OPEs.

\bigskip

{\bf Proposition.}\: {\it
For a non-negative integer $N$,
a holomorphic $A$-valued distribution $a(z)$
and an $A$-valued distribution $b(\vz)$
are holomorphically local of order $N$
if and only if there exist 
$A$-valued distributions $c^n(\vz)$ and $r(z,\vw)$
such that
$c^{N-1}(\vz)$ is non-zero and we have
\begin{equation}
\label{E:OPEzw}
\quad\hspace{2mm}
a(z)b(\vw)
\; =\;
\sum_{n=0}^{N-1}\: 
\frac{c^n(\vw)}{(z-w)^{n+1}}
\: +\:
r(z,\vw)
\end{equation}
and
\begin{equation}
\label{E:OPEwz}
\paraab\, b(\vw)a(z)
\; =\;
\sum_{n=0}^{N-1}\: 
\frac{c^n(\vw)}{(z-w)_{w>z}^{n+1}}
\: +\:
r(z,\vw).
\end{equation}
If this is the case then  
$c^n(\vz)=a(z)_{(n)}b(\vz)$ 
and $r(z,\vw)=\, \normord{a(z)b(\vw)}$.
}

\bigskip

\begin{pf}
``$\Rightarrow$"\:
Proposition \ref{SS:locality} implies  
\begin{equation}
\label{E:commut deldelta}
[a(z),b(\vw)]
\; =\;
\sum_{n=0}^{N-1}\:
c^n(\vw)\,
\del_w^{(n)}\de(z,w)
\end{equation}
where $c^n(\vz)=a(z)_{(n)}b(\vz)$.
If we restrict \eqref{E:commut deldelta}
to $\N\times\K^2$ and $\Z_<\times\K^2$
then using \eqref{E:delta pm z-w ><} we get
\begin{equation}
\notag
[a(z)_+,b(\vw)]
\, =\,
 \sum_{n=0}^{N-1}
\frac{c^n(\vw)}{(z-w)^{n+1}},
\qquad
[a(z)_-,b(\vw)]
\, =\,
-\sum_{n=0}^{N-1}
\frac{c^n(\vw)}{(z-w)_{w>z}^{n+1}}.
\end{equation}
The OPEs \eqref{E:OPEzw} and \eqref{E:OPEwz} 
with $r(z,\vw)=\,\normord{a(z)b(\vw)}$
thus follow from 
\eqref{E:azbw decomp1} and
$$
\paraab\, b(\vw)a(z)
\; =\;
-[a(z)_-,b(\vw)]
\: +\: a(z)_- b(\vw)
\: +\:
\paraab\, b(\vw)a(z)_+.
$$

``$\Leftarrow$"\:
Taking the difference of 
\eqref{E:OPEzw} and \eqref{E:OPEwz}
we obtain
\eqref{E:commut deldelta}
because of \eqref{E:delta and z-w}.
Thus by Proposition \ref{SS:locality}
the distributions 
$a(z)$ and $b(\vz)$ are holomorphically local of order $N$
and $c^n(\vz)=a(z)_{(n)}b(\vz)$.
This implies that 
$a(z)_{(n)}b(\vz)$ vanishes for $n\geq N$
so that 
$r(z,\vw)=\, \normord{a(z)b(\vw)}$
follows from Proposition \ref{SS:ope}.
\end{pf}

\subsection{Holomorphic Skew-Symmetry}
\label{SS:skew-symmetry distr}
In section \ref{SS:Jacobi formal distributions}
we gave a first answer to the question,
how properties of the algebra $\fg$ imply
properties of the $\N$-fold module $\fg\set{\vz}$,
when we showed that the Leibniz identity for $\fg$
implies the holomorphic Jacobi identity for $\fg\set{\vz}$.
Here we prove that if $\fg$ is skew-symmetric then 
holomorphically local distributions in $\fg\pau{z\uppm}$ satisfy 
{\it holomorphic skew-symmetry}.

\bigskip

{\bf Definition.}\:
Let $V$ be a $\K^2$-fold algebra, $T$ be an even operator,
and $n$ be an integer.
A holomorphic state $a$ 
and a state $b$ are said to satisfy
{\bf holomorphic skew-symmetry} 
for $T$ and the $n$-th product if
$a$ is bounded on $b$ and
$$
\paraab\, b_{(n,-1)}a
\; =\;
\sum_{i\in\N}\:
(-1)^{n+1+i}\, T^{(i)}(a_{(n+i)}b).
$$
\index{holomorphic skew-symmetry!for a pair of states}

\bigskip

States $a$ and $b$ of an $\N$-fold algebra 
are said to satisfy holomorphic skew-symmetry
if they satisfy holomorphic skew-symmetry
for the $n$-th product for any non-negative integer $n$.

\bigskip

{\bf Lemma.}\: {\it
If $a(z,w)$ is a distribution in $V\pau{z\uppm}\set{w}$
such that $a(w,w)$ is well-defined 
then
the product $\de(z,w)a(z,w)$ is well-defined
and 
$$
\de(z,w)a(z,w)
\; =\;
\de(z,w)a(w,w).
$$
}

% no bigskip

\begin{pf}
Because $a(w,w)$ is well-defined 
we have 
$a(z,w)\in V[(z/w)^{\K}]\set{w}$.
Thus $a(z,w)\in V\pau{z\uppm}\set{w}$ implies
$a(z,w)\in V[(z/w)\uppm]\set{w}$.
Hence we may consider $a(z,w)$ as the kernel of
a morphism $\K[w^{\K}]\to V[(z/w)\uppm]$.
Since $\de(z,w)\in z^{-1}\Z\pau{(z/w)\uppm}$
this shows that $\de(z,w)a(z,w)$ is well-defined
and together with \eqref{E:delta is delta}
it implies the equation.
\end{pf}

\bigskip

The Lemma implies that
for any integral $V$-valued distribution $a(z)$ 
we have
\begin{equation}
\label{E:delta test}
\res_z(\de(z,w)a(z))
\; =\;
\res_z(\de(z,w)a(w))
\; =\;
a(w).
\end{equation}

\bigskip

{\bf Proposition.}\: {\it
If $\fg$ is skew-symmetric then
holomorphically local distributions of 
the $\N$-fold algebra $\fg\pau{z\uppm}$
satisfy holomorphic skew-symmetry for $\del_z$.
}

\bigskip

\begin{pf}
Let $a(z)$ and $b(z)$ be 
holomorphically local distributions in $\fg\pau{z\uppm}$
and $n$ be a non-negative integer.
From
the definition of the $n$-th products,
Proposition \ref{SS:locality},
\eqref{E:delta properties},
the integration-by-parts formula, and
\eqref{E:delta test}
we obtain 
\begin{align}
\notag
\paraab\, b(w)_{(n)}a(w)
\; &=\;
-\res_z((z-w)^n [a(w),b(z)])
\\
\notag
&=\;
-\sum_{i\in\N}\:
\res_z((z-w)^n\, a(z)_{(i)}b(z)\, \del_z^{(i)}\de(z,w))
\\
\notag
&=\;
\sum_{i\in\N}\:
(-1)^{n+1}\,
\res_z(a(z)_{(n+i)}b(z)\, \del_z^{(i)}\de(z,w))
\\
\notag
&=\;
\sum_{i\in\N}\:
(-1)^{n+1+i}\,
\res_z(\del_z^{(i)}( a(z)_{(n+i)}b(z))\, \de(z,w))
\\
\notag
&=\;
\sum_{i\in\N}\:
(-1)^{n+1+i}\, \del_w^{(i)}( a(w)_{(n+i)}b(w)).
\end{align}
\end{pf}

\subsection{Dong's Lemma for Holomorphic Locality}
\label{SS:Dong}
We prove {\it Dong's lemma for holomorphic locality} 
which states that 
the $n$-th pro\-ducts 
of the $\N$-fold module of holomorphic distributions
preserve the relation of holomorphic locality.

\bigskip

{\bf Lemma.}\: {\it
Let $\fg$ be a Lie algebra
and 
$a(z), b(z)$, and $c(z)$ be integral $\fg$-valued distributions
that are pairwise holomorphically local of 
orders $N_{ab}, N_{bc}$, and $N_{ac}$.
If
$n_{ab}, n_{bc}$, and $n_{ac}$
are non-negative integers such that
$n_{ab}+n_{bc}+n_{ac}\geq N_{ab}+N_{bc}+N_{ac}-1$
then
$$
(z-w)^{n_{ab}}
(w-x)^{n_{bc}}
(z-x)^{n_{ac}}\,
[[a(z),b(w)],c(x)]
\; =\;
0.
$$
}

\begin{pf}
The claim is obvious if $n_{ab}\geq N_{ab}$.
Assume that $n_{ab}<N_{ab}$ and thus
$n_{bc}+n_{ac}\geq N_{bc}+N_{ac}$.
If $n_{bc}\geq N_{bc}$ and $n_{ac}\geq N_{ac}$ then
the claim follows by applying the Jacobi identity.
Assume that $n_{bc}<N_{bc}$ so that $n_{ac}>N_{ac}$.
The alternative case $n_{ac}<N_{ac}$
is proven in the same way.

Expanding the
$(n_{ac}-N_{ac})$-th power of $z-x=(z-w)+(w-x)$
we obtain
$$
(z-x)^{n_{ac}}
\; =\;
\sum_{i=0}^{n_{ac}-N_{ac}}
\binom{n_{ac}-N_{ac}}{i}
(z-w)^i\, (w-x)^{n_{ac}-N_{ac}-i}\, (z-x)^{N_{ac}}.
$$
If $n_{ab}+i\geq N_{ab}$ then 
$[[a(z),b(w)],c(x)]$
is annihilated by $(z-w)^{n_{ab}}(w-x)^{n_{bc}}$ times
the $i$-th summand.
If $n_{ab}+i<N_{ab}$ then $n_{bc}+n_{ac}\geq N_{bc}+N_{ac}+i$
and thus $n_{bc}+(n_{ac}-N_{ac}-i)\geq N_{bc}$.
From the Jacobi identity follows that 
$[[a(z),b(w)],c(x)]$
is annihilated by $(z-w)^{n_{ab}}(w-x)^{n_{bc}}$ times
the $i$-th summand.
\end{pf}

\bigskip

{\bf Dong's Lemma.}\: {\it
Let $\fg$ be a Lie algebra, 
$a(z)$ and $b(z)$ be holomorphic $\fg$-valued distributions,
$c(\vz)$ be a $\fg$-valued distibution,
and $n$ be a non-negative integer.
If $a(z), b(z)$, and $c(\vz)$ 
are pairwise holomorphically local of 
orders $N_{ab}, N_{bc}$, and $N_{ac}$
then
the pairs 
$a(z)_{(n)}b(z), c(\vz)$ and 
$a(z), b(z)_{(n)}c(\vz)$ are both holomorphically local
of order at most
$(N_{ab}+N_{bc}+N_{ac}-n-1)_+$.
}

\bigskip

\begin{pf}
This is a direct consequence of the Lemma
because
we may assume that $c(\vz)=c(z)\in\fg\pau{z\uppm}$ 
and 
because we have
$$
(z-w)^m[a(z)_{(n)}b(z),c(w)]
\; =\;
\res_x(
(z-w)^m(x-z)^n\, [[a(x),b(z)],c(w)]).
$$
\end{pf}

\subsection{The Virasoro OPE}
\label{SS:vir ope}

We express the property of being a Virasoro distribution
in terms of an OPE.

\bigskip

{\bf Remark.}\: {\it
Let $c\in\K$. The following statements 
about a holomorphic $\End(V)$-valued distribution $L(z)$ are equivalent:

\begin{enumerate}
\item[$\iti$]
$L(z)$ is a Virasoro distribution of central charge $c$;

\item[$\itii$]
$L(z)$ is holomorphically local with OPE
\begin{equation}
\label{E:vir ope}
L(z)L(w)
\; \sim \;
\frac{c/2}{(z-w)^4}
\; +\;
\frac{2L(w)}{(z-w)^2}
\; +\;
\frac{\del_w L(w)}{z-w};
\end{equation}

\item[$\itiii$]
$L(z)$ is holomorphically local of order at most four
and the second and fourth OPE-coefficient of $L(z)L(w)$
is $2L(z)$ and $c/2$, respectively.
\end{enumerate}
}

\bigskip

\begin{pf}
$\iti\Leftrightarrow\itii$\:
Proposition \ref{SS:locality} shows that $L(z)$ is holomorphically local with 
OPE-coefficients $c^i(z)$ given by \eqref{E:vir ope} if and only if 
for any integers $n$ and $m$ the operator $[L_n,L_m]$ of $V$
is equal to 
\begin{align}
\notag
&\sum_{i\in\N}\:\binom{n+1}{i} c^i_{n+m+2-i}
\\
\notag
=\;
&(\del_z L(z))_{n+m+2}
\; +\;
(n+1)(2L(z))_{n+m+1}
\; +\;
\binom{n}{3}\, c\, \de_{n+m-1,-1}/2
\\
\notag
=\;
&(-(n+m+2)+2(n+1))L_{n+m}
\; +\;
(n^3-n)\, c\, \de_{n+m,0}/12.
\end{align}

$\itii\Rightarrow\itiii$\:
The $n$-th product $L(z)_{(n)}L(z)$ is zero for $n>3$ because
$L(z)$ is holomorphically local of order at most four.
Holomorphic skew-symmetry implies that
$$
L(z)_{(2)}L(z)
\; =\;
-L(z)_{(2)}L(z)
\; +\;
\del_z(L(z)_{(3)}L(z)),
$$
thus $L(z)_{(2)}L(z)=0$, and
$$
L(z)_{(0)}L(z)
\; =\;
-L(z)_{(0)}L(z)
\; +\;
\del_z(2L(z))
\; +\;
\del_z^{(3)}(L(z)_{(3)}L(z)),
$$
hence $L(z)_{(0)}L(z)=\del_z L(z)$.

$\itiii\Rightarrow\itii$\:
This is obvious.
\end{pf}

\section{The $\Z$-Fold Module of Holomorphic Fields}
\label{S:module q fields}

{\bf Summary.}\:
In section \ref{SS:nth products} 
we introduce the $\Z$-fold module $\cF(V)$.
In section \ref{SS:normord}
we prove that the $n$-th products for negative $n$ 
are the Taylor coefficients of the normal ordered product.
In section \ref{SS:holom loc implies itself}
we prove that holomorphically local holomorphic fields 
satisfy holomorphic locality.

In section \ref{SS:sfold id}
we define various notions of an identity.
In section \ref{SS:creat}
we prove that the field-state correspondence 
is a morphism from the $\Z$-fold submodule of creative fields 
to the $\Z$-fold module $V$.
In section \ref{SS:id field}
we prove that the identity field is a full identity of $\cF(V)$.

\medskip

{\bf Conventions.}\:
Except when we discuss $S$-fold modules, 
we denote by $V$ a vector space,
$1$ is an even vector, and $\vT$ is a pair of even operators of $V$.
Unless stated otherwise,
all distributions take values in $\End(V)$.

\subsection{The $\Z$-Fold Module of Holomorphic Fields}
\label{SS:nth products}

We define the notion of a {\it field}
and using the concept of radial ordering
we define on the vector space of fields
the structure of a $\Z$-fold module over the
$\Z$-fold algebra of holomorphic fields.

\bigskip

{\bf Definition.} 
An $\End(V)$-valued distribution $a(\vz_1,\dots,\vz_r)$
is called a {\bf field} on $V$
if $a(\vz_1,\dots,\vz_r)$ is bounded on any vector of $V$.
\index{field}

\bigskip

We denote by $\cF_r(V)$ 
the vector space of fields $a(\vz_1,\dots,\vz_r)$ on $V$  
and we write $\cF(V):=\cF_1(V)$. We define  
$$
\cF\sqbrack{\vz}(V)
\; :=\;
\Vect(V,V\sqbrack{\vz}).
$$
The vector spaces 
$\cF\sqbrack{\vz}\sqbrack{\vw}(V), \cF\lau{z}(V), 
\cF(\!(z,\vw\rangle(V), \dots$, 
are defined in the same way.
Using the identification $\End(V)\set{\vz}=\Vect(V,V\set{\vz})$
we have $\cF(V)=\cF\sqbrack{\vz}(V)$.

Let $a(\vz)$ and $b(\vz)$ be fields.
The distribution
$a(\vz)b(\vw)$ is contained in
$\cF\sqbrack{\vz}\sqbrack{\vw}(V)$
which is a module over 
$\K\sqbrack{\vz}\sqbrack{\vw}$.
Thus for any integer $n$
the product $(z-w)^n a(\vz)b(\vw)$ is well-defined.
For the same reason the product
$$
(z-w)^n b(\vw)a(\vz)
\; :=\;
(z-w)_{w>z}^n b(\vw)a(\vz)
$$
is well-defined. 
This short-hand notation should not lead to confusion
because the product of $(z-w)_{z>w}^n$ and $b(\vw)a(\vz)$
does in general not exist.
Using these conventions,
we define
$$
(z-w)^n[a(\vz),b(\vw)]
\; :=\;
(z-w)^n a(\vz)b(\vw)
\: -\: 
\paraab\, (z-w)^n b(\vw)a(\vz).
$$

Similarly, 
the product of a rational function
$f(z_1,\dots,z_r)$
with an expression involving 
juxtapositions and Lie brackets of fields
$a^1(\vz), \dots, a^r(\vz)$
is defined by first replacing the Lie brackets by
commutators, so that one obtains a sum of 
juxtapositions 
$A_{\si}:=
\pm a^{\si 1}(\vz_{\si 1})\dots a^{\si r}(\vz_{\si r})$
where $\si\in\bbS_r$,
and then
multiplying the summand $A_{\si}$
by 
$T_{z_{\si 1},\dots,z_{\si r}}(f)$
for any $\si\in\bbS_r$.

This relation between 
the power series expansion
and the order of fields 
corresponds to the notion of 
radial ordering in physics. 
Therefore we will refer to 
the above notational conventions as
{\bf radial ordering}.
\index{radial ordering}

A field $a(\vz)$ is called {\bf holomorphic}
\index{holomorphic field} if \index{field!holomorphic}
it is a holomorphic distribution.
We denote by $\cF_z(V)$ the vector space of holomorphic fields on $V$.
Thus $\cF_z(V)=\cF\lau{z}(V)$.

We define on $\cF(V)$ 
the structure of a $\Z$-fold module over $\cF_z(V)$ 
by 
$$
a(w)_{(n)}b(\vw)
\; :=\;
\res_z((z-w)^n [a(z),b(\vw)]).
$$
Note that $a(z)_{(n)}b(\vz)$ is indeed a field because
$\res_z$ maps 
$\cF\sqbrack{z}\sqbrack{\vw}(V)$ and 
$\cF\sqbrack{\vw}\sqbrack{z}(V)$ to
$\cF\sqbrack{\vw}(V)$.

Of course,
the restriction of the $\Z$-fold module $\cF(V)$ to $\N$
is an $\N$-fold submodule of the $\N$-fold module
$\End(V)\set{\vz}$ over $\cF_z(V)$.
The subspace $\cF_z(V)$ is a $\Z$-fold submodule of $\cF(V)$
and is thus a $\Z$-fold algebra.

The identity
\begin{equation}
\label{E:nth prod fields}
a(w)_{(n)}b(\vw)
\; =\;
\sum_{i\in\N}\:
(-1)^i \binom{n}{i} \,
(a_{n-i}b(\vw)w^i
\: -\:
\paraab\, (-1)^n w^{n-i} b(\vw)a_i)
\end{equation}
yields the first part of the following claim.

\bigskip

{\bf Remark.}\:
Let $M$ be a $\K^2$-fold module over a $\Z$-fold algebra $V$,
$a$ and $b$ be states of $V$ such that $a$ is holomorphic on $M$,
$c$ be a vector of $M$, and $r$ and $t$ be integers. 
Assume that $a$ and $b$ are bounded on $M$.
The vectors $a, b, c$
satisfy the associativity formula for indices $r$ and any $\vs\in\K^2$
if and only if 
$(a_{(r)}b)(\vz)c=a(z)_{(r)}b(\vz)c$.
The vectors $a, b, c$
satisfy holomorphic locality of order at most $t$
if and only if $(z-w)^t[a(z),b(\vw)]c$ vanishes.

\subsection{Normal Ordered Product}
\label{SS:normord}

We prove that the $n$-th products $a(z)_{(n)}b(\vz)$
of the $\Z$-fold module of holomorphic fields for negative $n$ 
are the Taylor coefficients of the normal ordered product
$\normord{a(z)b(\vw)}$.

\bigskip

If $a(z)$ is a holomorphic field and $b(\vz)$ is a field
then \eqref{E:nth prod fields} for $n=-1$ yields
\begin{equation}
\label{E:-1st=nop}
a(z)_{(-1)}b(\vz)
\; =\; 
a(z)_- b(\vz)\: +\:
\paraab\, b(\vz)a(z)_+
\; =\;
\normord{a(z)b(\vz)}.
\end{equation}
The field $\normord{a(z)b(\vz)}$
is called like $\normord{a(z)b(\vw)}$
the {\bf normal ordered product} of $a(z)$ and $b(\vz)$.
\index{normal ordered product}

\bigskip

{\bf Remark.}\: {\it
$\iti$\:
The derivative $\del_{\vz}$ is a translation operator
of the $\Z$-fold module $\cF(V)$ over $\cF_z(V)$.
If $a(z)$ is a holomorphic field,
$b(\vz)$ is a field, and 
$n$ is a non-negative integer
then
\begin{equation}
\label{E:-1th normord}
a(z)_{(-1-n)}b(\vz)
\; =\; 
\normord{\del_z^{(n)}a(z)b(\vz)}.
\end{equation}

$\itii$\:
If $a(z)$ is a holomorphic field and
$b(\vz)$ is a field then
the normal ordered product
$\normord{a(z)b(\vw)}$ is a field in 
$\cF(\!(z,\vw\rangle(V)$.
}

\bigskip

\begin{pf}
$\iti$\:
The first claim is proven in the same way as 
the corrresponding
statement for the $\N$-fold module $\fg\set{\vz}$,
see Remark \ref{SS:derivative}.
The second claim is a consequence of the first claim
because of \eqref{E:transl endo2} and 
\eqref{E:-1st=nop}.

$\itii$\:
From $a(z)_-\in\End(V)\pau{z}$
follows that 
$a(z)_-b(\vw)\in\cF\pau{z}\sqbrack{\vw}(V)$.
If $c\in V$
then $a(z)_+c\in V[z\uppm]$ 
and thus
$b(\vw)a(z)_+c\in V[z\uppm,\vw\rangle$.
\end{pf}

\bigskip

{\bf Proposition.}\: {\it
If $a(z)$ is a holomorphic field,
$b(\vz)$ is a field, 
and 
$N$ is an integer 
then there exists a unique field
$r(z,\vw)$ in $\cF(\!(z,\vw\rangle(V)$
such that 
\begin{equation}
\notag
a(z)b(\vw)
\; =\;
\sum_{n\in\N-N}\:
\frac{a(w)_{(n)}b(\vw)}{(z-w)^{n+1}}
\: +\:
(z-w)^N\, r(z,\vw).
\end{equation}
}

\bigskip

\begin{pf}
From
$\del_z(a(z)_{\pm})=(\del_z a(z))_{\pm}$
and \eqref{E:-1th normord} we obtain
\begin{equation}
\label{E:-n-1 th =del normord}
\del_z^{(n)}\normord{a(z)b(\vw)}|_{z=w} 
\; =\; 
\normord{\del_w^{(n)} a(w)b(\vw)}
\; =\;
a(w)_{(-1-n)}b(\vw)
\end{equation}
for any $n\in\N$.
Thus existence of $r(z,\vw)$ follows from 
Proposition \ref{SS:ope},
part $\itii$ of the Remark,
and Taylor's formula applied to $\normord{a(z)b(\vw)}$.
The field $r(z,\vw)$ is unique
because 
$\cF\sqbrack{\vz}\sqbrack{\vw}(V)$ is 
a vector space over $\K\lau{z}\lau{w}$.
\end{pf}

\subsection{Holomorphic Locality of Holomorphically Local Fields}
\label{SS:holom loc implies itself}

We prove that holomorphically local holomorphic fields 
satisfy holomorphic locality for the $\Z$-fold module of
holomorphic fields.

\bigskip

{\bf Proposition.}\: {\it
For $N\in\N$, 
if $a(z)$ and $b(z)$ are holomorphic fields
that are holomorphically local of order $N$ 
then the $\End(\cF(V))$-valued distributions 
$a(z)(Z)$ and $b(z)(Z)$ are holomorphically local
of order at most $N$.
}

\bigskip

\begin{pf}
For any field $c(\vz)$, we have
\begin{align}
\notag
a(w)(Z)\, c(\vw)
\; &=\;
\sum_{n\in\Z}\:\res_z((z-w)^n [a(z),c(\vw)])\; Z^{-n-1}
\\
\notag
&=\;
\res_z(\de(Z,z-w)\, [a(z),c(\vw)]).
\end{align}
Thus
$$
a(x)(Z)\, b(x)(W)\, c(\vx)
\; =\;
\res_{z,w}(\de(Z,z-x)\, \de(W,w-x)\, [a(z),[b(w),c(\vx)]]).
$$
By \eqref{E:delta is delta} we have 
\begin{align}
\notag
&(Z-W)^N\, \de(Z,z-x)\, \de(W,w-x)
\\
\notag
=\;
&((z-x)-(w-x))^N\, \de(Z,z-x)\, \de(W,w-x)
\\
\notag
=\;
&(z-w)^N\, \de(Z,z-x)\, \de(W,w-x).
\end{align} 
Together with the Leibniz identity we obtain
\begin{align}
\notag
&(Z-W)^N\, [a(x)(Z),b(x)(W)]c(\vx)
\\
\notag
=\;
&\res_{z,w}((z-w)^N\, \de(Z,z-x)\, \de(W,w-x)\, [[a(z),b(w)],c(\vx)])
\\
\notag
=\;
&0.
\end{align}
\end{pf}

\subsection{Identity}
\label{SS:sfold id}

We define various notions of an {\it identity}
for $\Z$-fold and $\K^2$-fold modules.
In section \ref{SS:id field} we prove that
the identity field is a full identity of $\cF(V)$.
Exis\-tence of an invariant right identity   
is one of the three axioms in the definition 
of an OPE-algebra.
An OPE-algebra has a full identity.

\bigskip

For a state $1$ of a $\K^2$-fold algebra $V$,
we denote by $\vT_1$ the pair of operators of $V$
that are given by $a\mapsto a_{(-2,-1)}1$ and $a\mapsto a_{(-1,-2)}1$.

A vector $1$ of a vector space $V$
is called {\bf invariant} \index{invariant vector} for 
a family $(T_i)$ of operators of $V$
if $T_i$ annihilates $1$ for any $i$.

Let $V$ be a $\K^2$-fold algebra and 
$1$ be an even vector of $V$.
A state $a$ of $V$ is called {\bf creative} \index{creative state}
for $1$ if $a_{(\vn)}1=\de_{\vn,-1}a$ for any 
$\vn\in\K^2\setminus(\Z_<^2\setminus\set{-1})$.
For a pair $\vT$ of even operators of $V$,
a state $a$ of $V$ is called
{\bf strongly creative} \index{strongly creative state}
for $1$ and $\vT$ if $a_{(\vn)}1=\vT^{(-1-\vn)}(a)$
for any $\vn\in\K^2$.

\bigskip

{\bf Definition.}\:
Let $V$ be a vector space, $1$ be an even vector of $V$,
and $\vT$ be a pair of even operators of $V$.

$\iti$\:
If $M$ is a $\K^2$-fold module over $V$ then 
the vector $1$ is called a {\bf left identity} \index{left identity}
of $M$ if $1_{(\vn)}a=\de_{\vn,-1}a$ for any $a\in M$ and $\vn\in\K^2$.

$\itii$\:
If $V$ is endowed with a $\K^2$-fold algebra structure then
the state $1$ is called 
\begin{enumerate}
\item[$\ita$]
a {\bf right identity} \index{right identity} if
any state of $V$ is creative for $1$;
\item[$\itb$]
a {\bf strong right identity} \index{strong right identity} for $\vT$ if
any state of $V$ is strongly creative for $1$ and $\vT$.
\end{enumerate}

$\itiii$\:
If $M$ is a $\K^2$-fold module over $V$ 
such that $V$ is a vector subspace and a $\K^2$-fold subalgebra of $M$ then
the state $1$ is called 
\begin{enumerate}
\item[$\ita$]
an {\bf identity} \index{identity}
if $1$ is a left and a right identity;

\item[$\itb$]
a {\bf strong identity} \index{strong identity} for $\vT$ if 
$1$ is a left identity and a strong right identity for $\vT$;

\item[$\itc$]
a {\bf full identity} \index{full identity} if
$1$ is a strong identity for $\vT_1$ and 
$\vT_1$ is a translation operator of $V$.
\end{enumerate}

\bigskip

For a $\Z$-fold module $M$ over a vector space $V$,
an even vector $1$ of $V$ is called a left identity of $M$ 
if $1$ is a left identity of $M|^{\K^2}$.
For a $\Z$-fold algebra $V$ and an even operator $T$ of $V$,
an even vector $1$ of $V$ is called a strong right identity for $T$ of $V$ 
if $1$ is a strong right identity for $\vT=(T,0)$ of $V|^{\K^2}$.
The notions of a right identity, an identity, a strong identity, and
a full identity are defined in the same way.

\medskip

Let $V$ be a $\K^2$-fold algebra and $\vT$ be a pair of even operators of $V$.
A right identity which is invariant for $\vT$ 
will just be called an invariant right identity.
\index{invariant right identity}
Of course,
if $1$ is a strong right identity for $\vT$ then $\vT=\vT_1$.
An identity is unique.

\bigskip

{\bf Remark.}\: {\it 
Let $V$ be a $\K^2$-fold algebra with a right identity $1$.
If $\vT$ is a 
translation endomorphism for $a$ and $1$ for any $a\in V$
then the endomorphisms $T$ and $\bT$ commute
and $1$ is a strong right identity for $\vT$.
The same conclusion holds if 
$\vT$ is a translation generator for $a$ and $1$
for any $a\in V$ and $1$ is invariant.
}

\bigskip

In Proposition \ref{SS:creat}\,$\itii$
we prove a stronger version of the claim
that if $\vT$ is a translation generator and 
$1$ is an invariant right identity
then $1$ is a strong right identity.

\bigskip

\begin{pf}
If $\vT$ is a translation endomorphism for $a$ and $1$ then
$T\bT(a)=T\bT(a)_{(-1,-1)}1=a_{(-2,-2)}1$
by \eqref{E:transl endo2}.
On the other hand, we have
$\bT T(a)=\bT T(a)_{(-1,-1)}1=T(a)_{(-1,-2)}1=a_{(-2,-2)}1$.
Thus $T$ and $\bT$ commute.
Moreover,
taking $b=1$ in \eqref{E:transl endo2}
shows that $1$ is a strong right identity for $\vT$.

If $\vT$ is a translation generator for $a$ and $1$
and $\vT(1)=0$ then
for any $\vn\in\K^2$ we have
$T(a_{(\vn)}1)=
-n\, a_{(n-1,\bn)}1+a_{(\vn)}T(1)
=-n\, a_{(n-1,\bn)}1$ and similarly
$\bT(a_{(\vn)}1)=-\bn\, a_{(n,\bn-1)}1$.
Thus the claim follows from the same arguments as those
that were used in the case of a translation endomorphism.
${}^{}$
\end{pf}

\bigskip

A $\K^2$-fold algebra $V$ 
is called {\bf unital} \index{unital $\K^2$-fold algebra}
if there exists an identity of $V$.
The notions of a strong unital \index{strong unital $\K^2$-fold algebra}
and a fully unital $\K^2$-fold algebra \index{fully unital $\K^2$-fold algebra}
are defined in the same way.
A morphism $V\to W$ of unital $\K^2$-fold algebras
is by definition a morphism of $\K^2$-fold algebras $V\to W$
such that the identity of $V$ is mapped to the identity of $W$.
For a unital $\K^2$-fold algebra $V$,
a $\K^2$-fold module $M$ over the vector space $V$
is called {\bf unitary} \index{unitary $\K^2$-fold module} if 
the identity of $V$ is a left identity of $M$.

If $1$ is an identity of a graded $\K^2$-fold algebra $V$
then $1$ is homogeneous of weight $0$. 
Indeed, let $1_{\vh}\in V_{\vh}$ such that $1=\sum_{\vh} 1_{\vh}$.
For any $\vh\in\K^2$, we have
$1_{\vh}=1_{\vh}{}_{(-1)}1=\sum_{\vh'} 1_{\vh}{}_{(-1)}1_{\vh'}$.
Because $1_{\vh}{}_{(-1)}1_{\vh'}\in V_{\vh+\vh'}$ the state 
$1_{\vh}{}_{(-1)}1_{\vh'}$ is zero if $\vh'\ne 0$.
Hence $1_{\vh'}=1_{(-1)}1_{\vh'}=\sum_{\vh} 1_{\vh}{}_{(-1)}1_{\vh'}=0$
if $\vh'\ne 0$.  
From $\vh(1)=0$ follows that $1_{\vn}=\de_{\vn,0}$ for any $\vn\in\K^2$.

\subsection{Field-State Correspondence and Creativity}
\label{SS:creat}

We define the notions of a {\it field-state correspondence}
and of {\it creativity}
and prove that the field-state correspondence 
is a morphism from the $\Z$-fold submodule of creative fields 
to the $\Z$-fold module $V$.
Moreover, we prove that if $1$ is invariant then
translation covariance implies strong creativity
for distributions that are bounded on $1$.

\bigskip

The morphism $s_1:\End(V)\set{\vz}\to V$ defined by
$a(\vz)\mapsto a_{-1}(1)$ is called the
{\bf field-state correspondence} \index{field-state correspondence}
of $V$ with respect to $1$.

\bigskip

{\bf Definition.}\:
A distribution $a(\vz_1,\dots,\vz_r)$ is called
{\bf creative} \index{creative distribution} for $1$ if
$a(\vz_1,\linebreak[0]\dots,\vz_r)1$ is a power series.
A distribution $a(\vz)$ is called
{\bf strongly creative} \index{strongly creative distribution} 
for $1$ and $\vT$ if 
$$
a(\vz)1
\; =\;
e^{\vz\vT}s_1(a(\vz)).
$$

\bigskip

We denote by $\End(V)\set{\vz}_1$ the space of 
creative distributions.
The spaces $\cF(V)_1$ and $\cF_z(V)_1$ are defined similarly.

\bigskip

{\bf Remark.}\: 
Let $V$ be a $\K^2$-fold algebra and $1$ be an even state.
A state $a$ is creative if and only if 
$a(\vz)$ is creative and $s_1(a(\vz))=a$.
Thus if $1$ is a right identity then $\cF_Y$ is creative and
$Y:V\to\cF_Y$ is a vector space isomorphism
with inverse $s_1:\cF_Y\to V$.

Conversely, if $V$ is a vector space,
$1$ is an even vector,
and $\cF$ is a vector subspace of $\End(V)\set{\vz}_1$
such that $s_1:\cF\to V$ is a vector space isomorphism 
then the inverse of $s_1:\cF\to V$
is the unique $\K^2$-fold algebra structure on $V$
such that $1$ is a right identity and 
$\cF$ is the space of fields.

Let $V$ be a $\K^2$-fold algebra, $1$ be an even state,
and $\vT$ be a pair of even operators of $V$.
A state $a$ is strongly creative if and only if 
$a(\vz)$ is strongly creative and $s_1(a(\vz))=a$.

\bigskip

{\bf Proposition.}\: {\it
$\iti$\:
A creative distribution $a(\vz)$ is strongly creative 
if and only if 
$s_1(\del_{\vz}^{(\vn)}a(\vz))=\vT^{(\vn)}s_1(a(\vz))$
for any $\vn\in\K^2$.

$\itii$\:
If $1$ is invariant for $\vT$ then
a translation covariant distribution $a(\vz)$ is strongly creative
if and only if $a(\vz)$ is bounded on $1$.

$\itiii$\:
The subspace $\cF(V)_1$ of $\cF(V)$ is a 
$\Z$-fold submodule over $\cF_z(V)_1$
and the field-state correspondence 
$s_1:\cF(V)_1\to V$ is a morphism of $\Z$-fold modules.
}

\bigskip

\begin{pf}
$\iti$\:
This is clear.

$\itii$\:
Suppose that $a(\vz)$ is bounded on $1$.
We may assume that $a(\vz)1$ is non-zero.
For any 
$\vh\in\K^2$ such that $a_{\vh}1$ is non-zero
there exists $\vn\in\vh+\Z^2$ 
such that
$a_{\vn}1$ is non-zero and
$a_{\vn+\vi}1$ vanishes for any non-zero $\vi\in\N^2$.
From 
$T(a_{n+1,\bn}1)=
-(n+1)a_{\vn}1+a_{n+1,\bn}T1$
we get 
$(n+1)a_{\vn}1=0$ and thus $n=-1$.
Similarly, we have $\bn=-1$. 
This shows that $a(\vz)$ is creative.
Because $Ta(\vz)1=[T,a(\vz)]1=\del_z a(\vz)1$ 
and similarly $\bT a(\vz)1=\del_{\bz}a(\vz)1$
the claim follows from Remark \ref{SS:transl cov}.
The converse statement is trivial.

$\itiii$\:
This follows from the explicit expression
\eqref{E:nth prod fields}
for the $n$-th products 
which yields for 
any creative holomorphic field $a(z)$ and 
any creative field $b(\vz)$ the identity 
$$
a(z)_{(n)}b(\vz)1
\; =\;
\sum_{i\in\N}\:
(-1)^i \binom{n}{i}\, z^i\,
a_{n-i}b(\vz)1.
$$
\end{pf}

\subsection{Identity Field}
\label{SS:id field}

We prove that 
the {\it identity field} is a full identity of $\cF(V)$.

\bigskip

The field $id_V$ is called the {\bf identity field} \index{identity field}
on $V$ and is denoted by $1(z)$.
The identity field is translation covariant for
any $\vT$ and creative for any $1$.
It is strongly creative for $1$ and $\vT$ if and only if
$1$ is invariant for $\vT$.

For a $\K^2$-fold module $M$ over a vector space $V$,
a vector $1$ of $V$ is a left identity if and only if 
$Y(1,\vz)$ is the identity field on $M$.

\bigskip

{\bf Proposition.}\: {\it
The identity field is a full identity of the
$\Z$-fold module $\cF(V)$ over $\cF_z(V)$.
}

\bigskip

\begin{pf}
The commutator $[1(z),a(\vw)]$ vanishes for any distribution $a(\vz)$.
Thus $1(z)_{(n)}a(\vz)$ and $b(z)_{(n)}1(z)$ are zero 
for any non-negative integer $n$.
Equation \eqref{E:-1th normord} yields that
$1(z)_{(-1-n)}a(\vz)=\;\normord{\del_z^{(n)}1(z)a(\vz)}=\de_{n,0}a(\vz)$ 
for any non-negative integer $n$.
This shows that $1(z)$ is a identity.
From Remark \ref{SS:nth products}\,$\iti$ we know that
$\del_z$ is a translation operator of $\cF_z(V)$.
Thus Remark \ref{SS:sfold id} implies that
$T_1=\del_z$ and $1$ is a full identity.
\end{pf}

\section{The Partial $\K^2$-Fold Algebra of Fields}
\label{S:algebra q fields}

{\bf Summary.}\: 
In sections \ref{SS:expansions powers} and
\ref{SS:canonical expansions powers}
we characterize the non-canonical and the 
canonical expansions of 
$(\mu z+\nu w)^h$ and $(\mu\vz+\nu\vw)^{\vh}$
where $\vh\in\K^2$ and $\mu, \nu\in\set{\pm 1}$.

In section \ref{SS:ope finite}
we define the notion of OPE-finiteness and 
prove a uniqueness result about the summands of an OPE.
In section \ref{SS:reduced ope}
we define the notion of a reduced OPE and prove that
OPE-finite distributions have a unique reduced OPE.
In section \ref{SS:ope holom distr}
we prove that the OPE of a holomorphic distribution
is holomorphic.

In section \ref{SS:nth prod ope fin}
we introduce a partial $\K^2$-fold algebra structure on $\End(V)\set{\vz}$. 
In section \ref{SS:nth prod holom}
we compare it with the
$\Z$-fold module structure on the subspace $\cF(V)$.
In section \ref{SS:id field derivatives}
we prove that the identity field is a full identity
of the partial $\K^2$-fold algebra $\cF(V)$,
that $\del_{\vz}$ is a translation operator,
and that the $\vn$-th products preserve translation covariance.
In section \ref{SS:grad dilat cov fields}
we prove that the space of dilatation covariant distributions
is a graded partial $\K^2$-fold algebra.

\medskip

{\bf Conventions.}\: 
We denote by $V$ a vector space,
$1$ is an even vector, and $\vT$ is a pair of even operators of $V$.
Unless stated otherwise,
all distributions take values in $\End(V)$.

\subsection{Non-Canonical Expansions of Non-Integral Powers}
\label{SS:expansions powers}

We define in terms of two characterizing properties 
the {\it non-integral powers}
$(\mu z+\nu w)_b^h$ where
$\mu, \nu\in\set{\pm 1}$, $h\in\K$, and $b$ is a power map.

\bigskip

Let $C$ be a commutative algebra, 
$S$ be a subset of $C\setminus\set{0}$,
$M$ be an abelian monoid,
and $1$ be an element of $M$.
A {\bf power map} \index{power map} is a map 
$M\times S\to C, (m,s)\mapsto s^m$, 
such that
$s^m s^n=s^{m+n}, s^0=1, s^1=s$, and 
if $1\in S$ then $1^m=1$.
A power map $p:\K\times S\to C$ is called {\bf compatible}
\index{power map!compatible with a differential}
with a derivation $\del$ of $C$ if
$\del(s^h)=hs^{h-1}\del(s)$ for any $s\in S$ and $h\in\K$.

Let $\K'$ be a commutative algebra 
and $\ka$ be an involution of $\K'$ that leaves $\K$ invariant.
A power map 
$b:\K\times\set{\pm 1}\to\K', (h,\mu)\mapsto (\mu)_b^h$,
is called {\bf compatible} \index{power map!compatible with an involution}
with $\ka$ if $\ka((-1)_b^h)=(-1)_b^{-\ka h}$ for any $h\in\K$.

For example, if $\K=\C$ then there exists a natural series 
of power maps $b_n:\C\times\set{\pm 1}\to\C$
indexed by $n\in\Z$ which are  
given by $(-1)_{b_n}^h:=e^{(2n+1)\pi ih}$.
The power maps $b_n$ are compatible with complex conjugation of $\C$.
This example is also the reason for putting
the minus sign in the definition of compatibility.

For an arbitrary field $\K$,
the algebra $\K':=\K[z^{\K}]/(z+1)$
has a power map $b:\K\times\set{\pm 1}\to\K'$
where $(-1)_b^h$ is given by the residue class of $z^h$ in $\K'$.
This power map is compatible with the involution of $\K'$
induced by $z^h\mapsto z^{-h}$.

Fix a power map $b:\K\times\set{\pm 1}\to\K'$
that is compatible with $\ka$
and let $z>w$ be two formal variables together with a total order. 
Define 
$\cL:=\set{\mu z+\nu w\mid\mu,\nu\in\set{\pm 1}}$.
There exists a unique power map
$p:\K\times\cL\to\K'\sqbrack{z}\pau{w}, (h,l)\mapsto (l)_b^h$,
such that
\begin{enumerate}
\item[$\iti$]\: 
$p$ compatible with the derivation $\del_w$ and
\item[$\itii$]\:
we have $(\mu z+\nu w)_b^h|_{w=0}=(\mu)_b^h\, z^h$.
\end{enumerate}
Indeed, uniqueness follows from
\begin{align}
\notag
(\mu z+\nu w)_b^h
\; &=\;
\sum_{i\in\N}\:
\del_w^{(i)}(\mu z+\nu w)_b^h|_{w=0}\, w^i
\\
\notag
&=\;
\sum_{i\in\N}\:
\binom{h}{i}(\mu z+\nu w)_b^{h-i}|_{w=0}\, \nu^i w^i
\\
\label{E:linear powers}
&=\;
\sum_{i\in\N}\:
\binom{h}{i}\, (\mu)_b^{h+i}\nu^i\, z^{h-i}w^i.
\end{align}
A direct calculation using the identity
$\binom{h+h'}{i}=\sum_{j=0}^i\binom{h}{j}\binom{h'}{i-j}$
shows that \eqref{E:linear powers} defines a power map
that satisfies $\iti$ and $\itii$.

In the same way one constructs a power map
$p:\K\times\cL'\to\K\sqbrack{z}\pau{w},
(h,l)\linebreak[0]\mapsto l_{z>w}^h$,
where 
$\cL':=\set{z+\la w, 1+\la z^{-1}w\mid\la\in\K}$.
This power map does {\it not} depend on $b$.
We define
$(z+\la w)^h:=(z+\la w)_{z>w}^h$
and
$(1+\la z^{-1}w)^h:=(1+\la z^{-1}w)_{z>w}^h$,
i.e.~the order of the variables $z$ and $w$ is determined by
the order of the summands of $z+\la w$.
If $a(z)$ is a $V$-valued distribution
then $a(z+w)=\sum_{n\in\K}a_n(z+w)^{-n-1}$.

\subsection{Canonical Expansions of Non-Integral Powers}
\label{SS:canonical expansions powers}

We define the non-integral powers $(\mu\bz+\nu\bw)_b^h$
and show that the product $(\mu\vz+\nu\vw)_b^{\vh}$
is independent of $b$ if $h-\bh$ is an integer.

\bigskip

Let $\K'$ be a commutative algebra with an involution $\ka$ 
that leaves $\K$ invariant and 
$b:\K\times\set{\pm 1}\to\K'$ be a power map that is compatible with $\ka$.
We define an extension of $\ka$ to 
an involution $\ka$ of $\K'\set{\vz}$ by
$$
\ka(a(\vz))
\; :=\;
\sum_{\vn\in\K^2}\:
\ka(a_{\vn})\:
\bz^{-\ka(n)-1}\, z^{-\ka(\bn)-1}.
$$
This definition is motivated
by the identity $\overline{a^b}=\ba^{\bb}$ where
$a\in\C^{\times}, b\in\C, 
a^b:=e^{b\ln a}, \ba^{\bb}:=e^{\bb\,\overline{\ln a}}$,
and $\ln a$ is a logarithm of $a$.

Define $\check{\cL}:=\cL\cup\ka\cL$.
There exists a unique extension of the power map
$p:\K\times\cL\to\K'\sqbrack{z}\pau{w}$
to a power map
$p:\K\times\check{\cL}\to\K'\sqbrack{\vz}\pau{\vw}, 
(h,l)\mapsto (l)_b^h$,
that commutes with $\ka$, i.e.~such that
$\ka((l)_b^h)=(\ka l)_b^{\ka h}$ for any $h\in\K$
and $l\in\check{\cL}$.
For $\mu,\nu\in\set{\pm 1}$ and $\vh\in\K^2$, 
we have
\begin{align}
\notag
(\mu\vz+\nu\vw)_b^{\vh}
\; &=\;
(\mu z+\nu w)_b^h\: \ka((\mu z+\nu w)_b^{\ka\bh})
\\
\notag
&=\;
\sum_{\vi\in\N^2}\:
\binom{\vh}{\vi}\,
(\mu)_b^{h+i}\, \ka((\mu)_b^{\ka\bh+\bi})\,
\nu^{\vi}\, \vz^{\vh-\vi}\, \vw^{\vi}
\\
\label{E:check z-w expan}
&=\;
\sum_{\vi\in\N^2}\:
\binom{\vh}{\vi}\,
(\mu)_b^{h-\bh+i+\bi}\, \nu^{\vi}\, \vz^{\vh-\vi}\, \vw^{\vi}.
\end{align}

We denote by $\vbbK$
the abelian subgroup of $\K^2$ consisting of $\vh$
such that $h-\bh$ is an integer.
Equation \eqref{E:check z-w expan} shows that
if $\vh\in\vbbK$ then 
$(\mu\vz+\nu\vw)_b^{\vh}$
does {\it not} depend on $b$ and is contained in
$\K\sqbrack{\vz}\pau{\vw}$. 
Moreover, we have
$(-\mu\vz-\nu\vw)_b^{\vh}=(-1)^{h-\bh}(\mu\vz+\nu\vw)_b^{\vh}$.
\index{Kaaa@$\vbbK$}

If the variables are ordered by $z>w$
then we define
$(\mu\vz+\nu\vw)_{z>w}^{\vh}:=(\mu\vz+\nu\vw)_b^{\vh}$
for any $\vh\in\vbbK$.
For unordered variables $z$ and $w$ we define
$(\mu\vz+\nu\vw)^{\vh}:=(\mu\vz+\nu\vw)_{z>w}^{\vh}$,
i.e.~the order of the variables $z$ and $w$ is determined by
the order of the summands of $\mu z+\nu w$.

\subsection{OPE-Finite Distributions}
\label{SS:ope finite}

We define the notion of {\it OPE-finiteness}
and prove that the space of
OPE-finite distributions
is the direct sum over cosets $H\in\K^2/\Z^2$ of
spaces of distributions of the form 
$c(\vz,\vw)(\vz-\vw)^{\vh}$
where $c(\vz,\vw)$ is a field and $\vh\in H$.
This result is used in the proof of many 
other statements.

\bigskip

A family of fields $c^i(\vz,\vw)$
together with $\vh_i\in\K^2$ where $i=1, \dots, r$
are called an {\bf OPE} \index{OPE} of a distribution $a(\vz,\vw)$ if 
\begin{equation}
\label{E:check ope}
a(\vz,\vw)
\; =\;
\sum_{i=1}^r\:
\frac{c^i(\vz,\vw)}{(\vz-\vw)^{\vh_i}}.
\end{equation}
We call $r$ the {\bf length}, \index{length of an OPE}
$c^i(\vz,\vw)$ the {\bf numerators}, 
\index{numerator of an OPE} and 
$\vh_i$ the {\bf pole orders} \index{pole order of an OPE} of the OPE.
For the sake of brevity,
we call \eqref{E:check ope} an OPE of $a(\vz,\vw)$.

\bigskip

{\bf Definition.}\:
For $\vh_i\in\K^2$,
a distribution $a(\vz,\vw)$ is called {\bf OPE-finite} 
\index{OPE-finite distribution}
of orders at most $\vh_1, \dots, \vh_r$
if there exists an OPE of $a(\vz,\vw)$ with pole orders
$\vh_1, \dots, \vh_r$.

\bigskip

The notions of an OPE and of OPE-finiteness 
for a distribution $a(\vz,\vw)$ 
depend implicitly on a total order on the set of variables $z$ and $w$.
OPE-finite distributions $a(\vz,\vw)$ are contained in
$\cF\sqbrack{\vz}\sqbrack{\vw}(V)$.
Distributions $a(\vz)$ and $b(\vz)$ are called OPE-finite 
if their juxtaposition $a(\vz)b(\vw)$ is OPE-finite.

\bigskip

{\bf Lemma.}\: {\it
If $K$ is a set of representatives of $\K /\Z$
then the power series $(1-x)^n$ for $n\in K$
are linearly independent vectors of 
the vector space $\K\lau{x}$ over $\K(x)$.
}

\bigskip

\begin{pf}
Suppose that 
$\sum_{i=1}^r p_i(x)(1-x)^{n_i}=0$ where $p_i(x)\in\K(x)$
and $n_i\in K$.
We may assume that $p_i(x)\in\K[x]$.
Writing $p_i(x)=\sum_{j=0}^s a_{ij} (1-x)^j$ we see that
it suffices to prove that the power series
$(1-x)^n$ in $\K\pau{x}$ for $n\in\K$ 
are linearly independent over $\K$.
This follows from the fact that
$(1-x)^n$ is an eigenvector of the
operator $(1-x)\del_x$ with eigenvalue $-n$.
\end{pf}

\bigskip

{\bf Proposition.}\: {\it
$\iti$\: 
If $c^i(z,w)$ are $V$-valued vertex series 
and $h_i\in\K$ 
such that $h_i\notin h_j+\Z$ for any $i\ne j$
and
$$
\sum_{i=1}^r\:
c^i(z,w)\, (z-w)^{h_i}
\; =\;
0
$$
then $c^i(z,w)$ is zero for any $i$.

$\itii$\:
If $c^i(\vz,\vw)$ are fields
and $\vh_i\in\K^2$
such that $\vh_i\notin \vh_j+\Z^2$ for any $i\ne j$
and
\begin{equation}
\label{E:ope relation check}
\sum_{i=1}^r\:
c^i(\vz,\vw)\, (\vz-\vw)^{\vh_i}
\; =\;
0
\end{equation}
then $c^i(\vz,\vw)$ is zero for any $i$.
}

\bigskip

\begin{pf}
$\iti$\:
By projecting onto a basis of $V$ we may assume $V=\K$.
Let $K$ be a set of representatives of $\K /\Z$.
There exists an isomorphism
$\io:\K\set{z,w}\to\K\pau{x\uppm}^{\K\times K}$
given by
$a(z,w)\mapsto\big(\sum_{i\in\Z}a_{n+i,m-i}x^i\big)_{n\in\K,m\in K}$.
We have 
$\io:\K\sqbrack{z,w}\to\K[x\uppm]^{\K\times K}$.
From 
$$
a(z,w)
\; =\;
\sum_{n\in\K, m\in K}\:
\bigg(
\sum_{i\in\Z}\: a_{n+i,m-i} \bigg( \frac{w}{z}\bigg)^i\, 
\bigg)
\:
z^{-n-1} w^{-m-1}
$$
we see that
\begin{align}
\notag
\io\sum_{i=1}^r\:
c^i(z,w)\, (z-w)^{h_i}
\; &=\; 
\io\sum_{i=1}^r\:
c^i(z,w)\, 
z^{h_i}\,\bigg( 1-\frac{w}{z} \bigg)^{h_i}
\\
\notag
&=\;
\sum_{i=1}^r\:
\io(c^i(z,w)z^{h_i})\,
(1-x)^{h_i}.
\end{align}
Therefore the claim follows from the Lemma.

$\itii$\:
By applying 
\eqref{E:ope relation check} to vectors of $V$
we see that it suffices to prove the claim 
for $V$-valued vertex series $c^i(\vz,\vw)$.
Part $\iti$ applied to
$(\bz-\bw)^{\bh_i} c^i(\vz,\vw)\in
V\set{\bz,\bw}\sqbrack{z,w}$ 
shows that for any coset $H\in\K/\Z$ we have
$\sum_{i:\, h_i\in H} (\bz-\bw)^{\bh_i} c^i(\vz,\vw)=0$.
Applying $\iti$ again the claim follows.
\end{pf}

\subsection{Reduced OPE}
\label{SS:reduced ope}

We define the notion of a {\it reduced OPE}
and prove that for any OPE-finite distribution 
there exists a unique reduced OPE.

\bigskip

{\bf Definition.}\:
An OPE \eqref{E:check ope} is called {\bf reduced} \index{reduced OPE} 
if $\vh_i\notin \vh_j+\Z^2$ for any $i\ne j$ 
and $(\vz-\vw)^{-\vn} c^i(\vz,\vw)$ is {\it not} a field
for any $i$ and $\vn\in\N^2\setminus\set{0}$.

\bigskip

{\bf Lemma.}\: {\it
$\iti$\:
If $d(\vz,\vw)$ is a non-zero $V$-valued Laurent series then the set $S$, 
consisting of $\vn\in\N^2$ such that
$(\vz-\vw)^{-\vn}d(\vz,\vw)$
lies in $V\lau{\vz,\vw}$, 
contains a supremum $\vN(d(\vz,\vw))$.
Moreover, 
the set $S$ coincides with the set 
consisting of $\vn\in\N^2$ such that
$(\vz-\vw)_{w>z}^{-\vn}d(\vz,\vw)$
lies in $V\lau{\vz,\vw}$.

$\itii$\:
If $c(\vz,\vw)$ is a non-zero field 
then the set, consisting of $\vn\in\N^2$ such that
$(\vz-\vw)^{-\vn}c(\vz,\vw)$ is a field, 
contains a supremum $\vN$.
}

\bigskip

\begin{pf}
$\iti$\:
Denote by $d^i(\vz,\vw)\in\K\lau{\vz,\vw}$ 
the projections of $d(\vz,\vw)$ 
with respect to a basis of $V$.
If 
$(\vz-\vw)^{-\vn}d(\vz,\vw)\in V\lau{\vz,\vw}$
then
$(\vz-\vw)^{-\vn}d^i(\vz,\vw)\in\K\lau{\vz,\vw}$
for any $i$.
We will prove the converse statement,
i.e.~if $(\vz-\vw)^{-\vn}d^i(\vz,\vw)\in\K\lau{\vz,\vw}$ for any $i$ then
there exists a family $((\vn_i,\vm_i))_i$ in $\Z^4$
that is bounded from below and such that
$(\vz-\vw)^{-\vn}d^i(\vz,\vw)\in
\vz^{\vn_i}\vw^{\vm_i}\K\pau{\vz,\vw}$.
This shows that $S$ is the set consisting of 
$\vn\in\N^2$ such that 
$d^i(\vz,\vw)$ is divisible by $(\vz-\vw)^{\vn}$ for any $i$.

We may assume 
$\vn=(1,0)$ and 
$d(\vz,\vw)\in V\lau{z,w}$.
If $D:=(z-w)^{-1}d(z,w)\notin V\lau{z,w}$
then $D\notin V\pau{w\uppm}\lau{z}$ and for any $n\in\Z$ 
there exists $n'\geq n$ and $m\in\Z$ such that $D_{n',m}\ne 0$. 
Then there exists $i$ such that $D^i_{n',m}\ne 0$ 
and by assumption we may choose $n'\geq n$ such that
$D^i_{n'+j,m-j}=0$ for $j>0$.
We get $d^i_{n',m-1}=D^i_{n'+1,m-1}-D^i_{n',m}\ne 0$
and hence $d_{n',m-1}\ne 0$.
This contradicts $d(z,w)\in V\lau{z,w}$.

The ring $\K\pau{\vz,\vw}$ is a unique factorization domain,
see e.g.~\cite{bourbaki.commutative.algebra.1-7},
Chap.~VII, \S 3.9, Prop.~8, hence so is $\K\lau{\vz,\vw}$.
Because $z-w$ and $\bz-\bw$ are non-associated 
prime elements of $\K\lau{\vz,\vw}$ 
there exists for any
$e(\vz,\vw)\in\K\lau{\vz,\vw}$
the supremum $\vN(e(\vz,\vw))$
of the set consisting of $\vn\in\N^2$
such that 
$e(\vz,\vw)$ is divisible by $(\vz-\vw)^{\vn}$.
The supremum $\vN(d(\vz,\vw))$
is the minimum of the family 
$(\vN(d^i(\vz,\vw)))_i$.

The second claim also follows from the 
description of $S$ 
in terms of divisibility of $d^i(\vz,\vw)$.

$\itii$\:
This follows from $\iti$
because $\vN$ is the minimum of the elements
$\vN(\vz^{\vh}\vw^{\vh'} c(\vz,\vw)e|_{\Z^4})$ of $\N^2$
where $e\in V$ and $\vh, \vh'\in\K^2$.
\end{pf}

\bigskip

Proposition \ref{SS:ope finite}
and the Lemma yield the following claim.

\bigskip

{\bf Proposition.}\: {\it
Any OPE-finite distribution $a(\vz,\vw)$ has a unique reduced OPE.
\hfill $\square$
}

\bigskip

Let $a(\vz,\vw)$ be an OPE-finite distribution.
The length of the reduced OPE of $a(\vz,\vw)$
is called the {\bf OPE-length} \index{OPE-length} of $a(\vz,\vw)$.
The pole orders of the reduced OPE of $a(\vz,\vw)$
are called the pole orders \index{pole orders of a distribution}
of $a(\vz,\vw)$.
If \eqref{E:check ope} is the reduced OPE of $a(\vz,\vw)$
and $\vh\in\K^2$ then 
the {\bf $\vh$-th numerator} of $a(\vz,\vw)$
is defined by 
$N_{\vh}(a(\vz,\vw)):=c^i(\vz,\vw)$ if $\vh=\vh_i$ for some $i$
and $N_{\vh}(a(\vz,\vw)):=0$ otherwise.
By definition, we have
$$
a(\vz,\vw)
\; =\;
\sum_{\vh\in\K^2}\:
\frac{N_{\vh}(a(\vz,\vw))}{(\vz-\vw)^{\vh}}.
$$

\index{numerator!of an OPE}

\bigskip

{\bf Remark.}\: {\it
Let $a(\vz,\vw)$ be an OPE-finite distribution,
$d$ be a vector of $V$,
and $S$ be the support in $\vw$ of $a(\vz,\vw)d$.
The support in $\vw$ of $N_{\vh}(a(\vz,\vw))c$
is contained in $S+\Z^2$ for any $\vh\in\K^2$.
}

\bigskip

\begin{pf}
Define $T:=\K^2\times(\K^2\setminus (S+\Z^2))$
and let \eqref{E:check ope} be the reduced OPE of $a(\vz,\vw)$.
Because $(\vz-\vw)^{\vh}\in\K\set{\vz}\pau{\vw}$ for any $\vh\in\K^2$
we have
$$
\sum_{i=1}^r\:
\frac{c^i(\vz,\vw)d|_T}{(\vz-\vw)^{\vh_i}}
\; =\;
\Bigg(\sum_{i=1}^r\:
\frac{c^i(\vz,\vw)d}{(\vz-\vw)^{\vh_i}}\Bigg)
\Bigg|_T
\; =\;
a(\vz,\vw)d|_T
\; =\;
0.
$$
By Proposition \ref{SS:ope finite} this implies that
$c^i(\vz,\vw)d|_T$ is zero for any $i$. That is the claim.
\end{pf}

\subsection{OPE of a Holomorphic Distribution}
\label{SS:ope holom distr}

We prove that an OPE-finite distribution $a(z,\vw)$
whose pole orders are contained in $\vbbK$
has an OPE of the form $c(z,\vw)(z-w)^N$ where $N$ is an integer.

\bigskip

{\bf Proposition}.\: {\it
Let $a(z,\vw)$ be an OPE-finite distribution in $\End(V)\set{z,\vw}$.
If the pole orders $\vh_1, \dots, \vh_r$ of $a(z,\vw)$
are contained in $\vbbK$ then $r\leq 1$ and 
if $r=1$ then $\vh_1$ is contained in $\Z\times\set{0}$.
}

\bigskip

\begin{pf}
Let \eqref{E:check ope} be the reduced OPE of $a(z,\vw)$.
Applying $(\bz-\bw)\del_{\bz}$ to \eqref{E:check ope}
we obtain
$$
0
\; =\;
\sum_{i=1}^r\:
\frac{(-\bh_i+(\bz-\bw)\del_{\bz})\, c^i(\vz,\vw)}{(\vz-\vw)^{\vh_i}}.
$$
Proposition \ref{SS:ope finite}
implies that $(\bz\del_{\bz}-\bh_i)c^i(\vz,\vw)=
\bw\del_{\bz}c^i(\vz,\vw)$ for any $i$.
In terms of modes this reads
\begin{equation}
\label{E:ope local}
(-\bn-1-\bh_i)c^i_{\vn,\vm}
\; =\;
-\bn\, c^i_{n,\bn-1,m,\bm+1}.
\end{equation}
Fix $i$.
Assume that there exist 
$d\in V$ and $(\vn_0,\vm_0)\in\K^4$ such that
$c^i_{\vn_0,\vm_0}(d)$ is non-zero.
We may assume that 
$c^i_{n_0,\bn_0-j,m_0,\bm_0+j}(d)=0$ for any positive integer $j$.
Then \eqref{E:ope local} implies that $\bn_0=-1-\bh_i$ and
$$
c^i_{n_0,\bn_0+j,m_0,\bm_0-j}(d)
\; =\;
\binom{\bn_0+j}{j}\, c^i_{\vn_0,\vm_0}(d).
$$
Because $\binom{\bn_0+j}{j}=(-1)^j\binom{\bh_i}{j}$
this shows that if $\bh_i\notin\N$ then
$c^i_{n_0,\bn_0+j,m_0,\bm_0-j}(d)$ is non-zero for any $j\in\N$.
But this contradicts that 
$c^i(\vz,\vw)(d)$ is a vertex series.
Thus $r\leq 1$ and if $r=1$ then 
$\vh_1\in\Z\times\N$.
Because
$(\bz-\bw)\del_{\bz}c^1(\vz,\vw)=\bh_1 c^1(\vz,\vw)$
and because \eqref{E:check ope} is the reduced OPE
of $a(z,\vw)$ we obtain $\bh_1=0$.
${}_{}$
\end{pf}

\subsection{The Partial $\K^2$-Fold Algebra of Fields}
\label{SS:nth prod ope fin}

We define the structure of a partial $\K^2$-fold algebra on 
$\End(V)\set{\vz}$
and show that the $\vn$-th products can be expressed in terms of an
arbitrary OPE.

\bigskip

{\bf Definition.}\:
For $\vn\in\K^2$,
the {\bf $\vn$-th OPE-coefficient}
of an OPE-finite distribution $a(\vz,\vw)$ is defined as
$$
c_{\vn}(a(\vz,\vw))(\vw)
\; :=\;
\sum_{\vh\in\K^2}\:
\del_{\vz}^{(\vh-\vn)}N_{\vh}(a(\vz,\vw))|_{\vz=\vw}.
$$

\bigskip

By definition, the OPE-coefficients are fields.
If \eqref{E:check ope} is the reduced OPE of $a(\vz,\vw)$
then the $\vh_i$-th OPE-coefficient
$c_{\vh_i}(a(\vz,\vw))(\vz)=c^i(\vz,\vz)$
is called the {\bf leading term} \index{leading term}
of order $\vh_i$ of $a(\vz,\vw)$.

If $a(\vz,\vw)$ is a field then
for any $\vn\in\K^2$ and $\vi\in\N^2$ we have
\begin{equation}
\label{E:taylor shift}
\del_{\vz}^{(\vn+\vi)}
((\vz-\vw)^{\vi}a(\vz,\vw))|_{\vz=\vw}
\; =\;
\del_{\vz}^{(\vn)}
a(\vz,\vw)|_{\vz=\vw}.
\end{equation}
In fact, if $\vn+\vi\notin\N^2$ then
$\vn\notin\N^2$ and 
\eqref{E:taylor shift} is trivial.
If $\vn+\vi\in\N^2$ then the Leibniz identity implies that 
both sides of \eqref{E:taylor shift} are equal to
$$
\sum_{\vj=0}^{\vn+\vi}\:
\del_{\vz}^{(\vj)}\, (\vz-\vw)^{\vi}|_{\vz=\vw}\;
\del_{\vz}^{(\vn+\vi-\vj)}a(\vz,\vw)|_{\vz=\vw}.
$$

From Proposition \ref{SS:ope finite} and \eqref{E:taylor shift}
follows that if \eqref{E:check ope} is an OPE of a distribution
$a(\vz,\vw)$ then 
$$
c_{\vn}(a(\vz,\vw))(\vw)
\; =\;
\sum_{i=1}^r\:
\del_{\vz}^{(\vh_i-\vn)}c^i(\vz,\vw)|_{\vz=\vw}.
$$

For an even set $S$, a vector space $V$ together with a subset $D$
of $V\times V$ and an even map ${}_{(s)}:D\to V$ for any $s\in S$
is called a {\bf partial} \index{partial $S$-fold algebra}
$S$-fold algebra if 
$(a,c), (b,c)\in D$ and $\la\in\K$ implies that 
$(a+\la b,c)\in D$ and $(a+\la b)_{(s)}c=a_{(s)}c+\la(b_{(s)}c)$
for any $s\in S$ and the same is true for the second component.
A morphism $\varphi:V\to W$ of partial $S$-fold algebras
is a linear map such that $(a,b)\in D$ implies $(\varphi a,\varphi b)\in D$
and $\varphi(a_{(s)}b)=\varphi(a)_{(s)}\varphi(b)$ for any $s\in S$.

An operator $T$ of a partial $S$-fold algebra $V$ is called a
derivation if $(a,b)\in D$ implies that $(Ta,b), (a,Tb)\in D$ and
$T(a_{(s)}b)=T(a)_{(s)}b+\zeta^{\tT\ta}\, a_{(s)}T(b)$.
An even state $1$ of a partial $\K^2$-fold algebra $V$ is called a
left identity if $(1,a)\in D$ and $1_{(\vn)}a=\de_{\vn,-1}a$ 
for any $a\in V$ and $\vn\in\K^2$. 
In the same way other notions defined for $S$-fold algebras
can be generalized to partial $S$-fold algebras.

We endow the module $\End(V)\set{\vz}$ with 
the structure of a partial $\K^2$-fold algebra by
letting $D$ be the set of pairs $(a(\vz),b(\vz))$ of 
OPE-finite distributions and defining
$a(\vz)_{(\vn)}b(\vz):=c_{\vn+1}(a(\vz)b(\vw))(\vz)$.

The subspace $\cF(V)$ is a partial $\K^2$-fold subalgebra of 
$\End(V)\set{\vz}$.
The definition of the $\vn$-th products shows that 
if $a(\vz)$ and $b(\vz)$ are OPE-finite distributions then
$a(\vz)$ is bounded on $b(\vz)$ by the pole orders of $a(\vz)b(\vw)$.
In particular, 
$\End(V)\set{\vz}$ is a bounded partial $\K^2$-fold algebra.

\subsection{The Partial $\K^2$-Fold Algebra of Holomorphic Fields}
\label{SS:nth prod holom}

We prove that the partial $\K^2$-fold algebra structure on $\cF(V)$
agrees for a holomorphic field $a(z)$ and a field $b(\vz)$ 
with the $\Z$-fold module structure 
if and only if $a(z)$ is bounded on $b(\vz)$ 
with respect to the $\Z$-fold module structure.

\bigskip

The $(n,-1)$-st products of the partial $\K^2$-fold algebra $\cF(V)$ 
in general do {\it not} coincide with the $n$-th products 
of the $\Z$-fold module $\cF(V)$ over $\cF_z(V)$
because the partial $\K^2$-fold algebra $\cF(V)$ is bounded
whereas there are OPE-finite holomorphic fields that
are {\it not} bounded with respect to the $\Z$-fold module structure. 
For example, if $a(z), b(z)\in z^{-1}\End(V)\pau{z}$ 
such that $[a_0,b(z)]$ is non-zero then
$a(z)b(w)$ is a field and hence OPE-finite but
the $n$-th products $a(z)_{(n)}b(z)=[a_0,b(z)](-z)^n$
of the $\Z$-fold module are non-zero for any non-negative integer $n$.

\bigskip

{\bf Proposition.}\: {\it
Let $a(z)$ be a holomorphic field, $b(\vz)$ be a field,
and $N$ be the mode infimum of $a(z)$ and $b(\vz)$
defined with respect to the $\Z$-fold module structure on $\cF(V)$.
We have $N<\infty$ if and only if 
$a(z)$ and $b(\vz)$ are OPE-finite and 
for any $\vn\in\K^2$ 
the $\vn$-th product of the partial $\K^2$-fold algebra coincides 
for $a(z)$ and $b(\vz)$
with the $\vn$-th product $c^{\vn}(\vz)$ of the $\Z$-fold module.
Moreover, in this case the following is true:

\begin{enumerate}
\item[$\iti$]
we have $N=-\infty$ if and only if $a(z)b(\vw)$ is zero;

\item[$\itii$]
we have $N>-\infty$ if and only if
the OPE-length of $a(z)b(\vw)$ is one;
in this case the pole order is $(N,0)$.
\end{enumerate}
}

\bigskip

\begin{pf}
The ``if" part of the first statement follows from the boundedness of
the partial $\K^2$-fold algebra $\cF(V)$.

Assume that $N=-\infty$
so that by definition $c^{\vn}(\vz)=0$ for any $\vn$.
Proposition \ref{SS:nth products}
shows that $a(z)b(\vz)$ is a field that is
divisible as a field by $(z-w)^n$ for any 
non-negative integer $n$.
Thus $a(z)b(\vz)$ vanishes by Lemma \ref{SS:reduced ope}
and the claim follows in this case.
Conversely, assume that $N<\infty$ and $a(z)b(\vw)$ is zero.
If $N$ were positive then Proposition \ref{SS:ope} yields
$$
\frac{c^{(N-1,-1)}(\vw)}{(z-w)^N}
\; =\;
-\sum_{n=0}^{N-2}\:
\frac{c^{(n,-1)}(\vw)}{(z-w)^{n+1}}
\; -\;
\normord{a(z)b(\vw)}.
$$
Thus multiplying both sides by $(z-w)^N$ and setting $z=w$
we arrive at the contradiction $c^{(N-1,-1)}(\vz)=0$.
Therefore $N\leq 0$ and $\normord{a(z)b(\vw)}=a(z)b(\vw)=0$ so that
$c^{\vn}(\vz)$ is zero for any $\vn$ by \eqref{E:-n-1 th =del normord}.

Now assume that $N\in\Z$.
Define $c(z,\vw):=(z-w)^N a(z)b(\vw)$.
Proposition \ref{SS:nth products} shows that for any integer $M$
there exists a field $r(z,\vw)$ such that
$$
c(z,\vw)
\; =\;
\sum_{n=M}^{N-1}\:
c^{(n,-1)}(\vw)\, (z-w)^{N-n-1}
\; +\;
(z-w)^{N-M}r(z,\vw).
$$
By taking $M=N$ we see that $c(z,\vw)$ is a field and
thus $a(z)b(\vw)=c(z,\vw)(z-w)^{-N}$ is an OPE.
Hence $c^{\vn}(\vw)=\del_{\vz}^{((N,0)-1-\vn)}c(z,\vw)|_{\vz=\vw}$
for any $\vn$.
Because $c(z,\vz)=c^{(N-1,0)}(\vz)$ is non-zero
the field $c(z,\vw)$ is not divisible as a field
by $z-w$ or $\bz-\bw$.
Thus $(N,0)$ is the pole order of $a(z)b(\vw)$.
\end{pf}

\subsection{Derivatives and Identity Field}
\label{SS:id field derivatives}

We prove that 
$\del_{\vz}$ is a translation operator of $\End(V)\set{\vz}$,
the identity field is a full identity of $\cF(V)$,
and translation covariant distributions form a partial $\K^2$-fold subalgebra.

\bigskip

If $a(\vz,\vw)$ is a field then
\begin{equation}
\label{E:del and evaluate}
\del_w a(\vw,\vw)
\; =\;
(\del_z+\del_w)a(\vz,\vw)|_{\vz=\vw}
\end{equation}
because both sides are equal to 
$\sum_{\vn,\vm\in\K^2}(-n-m-2)\, a_{\vn,\vm}\vw^{-\vn-\vm-(3,2)}$.

\bigskip

{\bf Proposition.}\: {\it
$\iti$\:
The derivative $\del_{\vz}$ is a translation operator 
of the partial $\K^2$-fold algebra $\End(V)\set{\vz}$.

$\itii$\:
The identity field is a full identity 
of the partial $\K^2$-fold algebra $\cF(V)$.

$\itiii$\:
The space $\End(V)\set{\vz}_{\vT}$ 
is a partial $\K^2$-fold subalgebra of $\End(V)\set{\vz}$.
}

\bigskip

\begin{pf}
$\iti$\:
If \eqref{E:check ope} is an OPE of the distributions $a(\vz)$ and $b(\vz)$
then
\begin{equation}
\label{E:del a  b ope}
\del_z a(\vz)b(\vw)
\; =\;
\sum_{i=1}^r\:
\frac{-h_i\, c^i(\vz,\vw)\; +\; (z-w)\, \del_z c^i(\vz,\vw)}
{(\vz-\vw)^{\vh_i+(1,0)}}
\end{equation}
and
$$
a(\vz)\del_w b(\vw)
\; =\;
\sum_{i=1}^r\:
\frac{h_i\, c^i(\vz,\vw)\; +\; (z-w)\, \del_w c^i(\vz,\vw)}
{(\vz-\vw)^{\vh_i+(1,0)}}
$$
are OPEs. 
On the other hand, from \eqref{E:del and evaluate} we get
$$
\del_w(a(\vw)_{(\vn)}b(\vw))
\; =\;
\sum_{i=1}^r\:
(\del_z+\del_w)\del_{\vz}^{(\vh_i-1-\vn)}c^i(\vz,\vw)|_{\vz=\vw}.
$$
This shows that $\del_z$ is a derivation.
From \eqref{E:del a  b ope} we obtain
\begin{align}
\notag
\del_w a(\vw)_{(\vn)}b(\vw)
\; =\;
&\sum_{i=1}^r\:
(\del_{\vz}^{(\vh_i-1-\vn)}\del_z-h_i\del_{\vz}^{(\vh_i+(1,0)-1-\vn)})\,
c^i(\vz,\vw)|_{\vz=\vw}
\\
\notag
=\;
&-n\, a(\vw)_{(n-1,\bn)}b(\vw).
\end{align}

$\itii$
If $a(\vz)$ is a field
then
the identities $1(z)a(\vw)=a(\vw)$ and $a(\vz)1(w)=a(\vz)$
are OPEs. Thus we have
$1(w)_{(\vn)}a(\vw)=
\del_{\vz}^{(-1-\vn)}a(\vw)|_{\vz=\vw}
=\de_{\vn,-1}a(\vw)$
and $a(\vw)_{(\vn)}1(\vw)
=\del_{\vw}^{(-1-\vn)}a(\vw)$ for any $\vn\in\K^2$.

$\itiii$\:
If \eqref{E:check ope} is an OPE of 
the translation covariant distributions $a(\vz)$ and $b(\vz)$
then
\begin{align}
\notag
[T,a(\vw)_{(\vn)}b(\vw)]
\; &=\;
\sum_{i=1}^r\: 
[T,\del_{\vz}^{(\vh_i-1-\vn)}c^i(\vz,\vw)|_{\vz=\vw}]
\\
\notag
\; &=\;
\sum_{i=1}^r\: 
\del_{\vz}^{(\vh_i-1-\vn)}[T,c^i(\vz,\vw)]|_{\vz=\vw}.
\end{align}
Thus $[T,a(\vz)_{(\vn)}b(\vz)]$ is the $(\vn+1)$-st OPE-coefficient
of $[T,a(\vz)b(\vw)]$.
Because $[T,a(\vz)b(\vw)]=[T,a(\vz)]b(\vw)+a(\vz)[T,b(\vw)]=
\del_z a(\vz)b(\vw)+a(\vz)\del_w b(\vw)$ 
the claim follows from $\itii$.
\end{pf}

\subsection{Gradation for Dilatation Covariant Fields}
\label{SS:grad dilat cov fields}

We prove that the space of dilatation covariant distributions
is a graded partial $\K^2$-fold algebra.

\bigskip

{\bf Proposition.}\: {\it
Let $V$ be a vector space endowed with a pair $\vH$ of even ope\-rators.
The subspace $\End(V)\set{\vz}_{\vast}$ is a graded partial $\K^2$-fold 
subalgebra of $\End(V)\set{\vz}$.
}

\bigskip

\begin{pf}
Let $a(\vz)$ and $b(\vz)$ be dilatation covariant distributions
of weights $\vh$ and $\vh'$, resp., $\vn\in\K^2$,
and \eqref{E:check ope} be an OPE of $a(\vz)b(\vw)$.
From the proof of Proposition \ref{SS:id field derivatives}\,$\itiii$
we know that $[H,a(\vz)_{(\vn)}b(\vz)]$ is the 
$(\vn+1)$-st OPE-coefficient of $[H,a(\vz)b(\vw)]$.
Because 
$$
[H,a(\vz)b(\vw)]
\; =\;
(h+h')\, a(\vz)b(\vw)
\; +\;
z\del_z a(\vz)b(\vw)
\; +\;
w\, a(\vz)\del_w b(\vw)
$$
and 
\begin{align}
\notag
&\del_{\vz}^{(\vh_i-1-\vn)}(z\del_z c^i(\vz,\vw))|_{\vz=\vw}
\\
\notag
=\;
&w\del_{\vz}^{(\vh_i-1-\vn)}\del_z c^i(\vz,\vw)|_{\vz=\vw}
\; +\;
\del_{\vz}^{(\vh_i-1-(n+1,\bn))}\del_z c^i(\vz,\vw)|_{\vz=\vw}
\end{align}
the claim follows from the fact that $\del_z$ is a translation operator.
\end{pf}

\section{Locality}
\label{S:locality}

{\bf Summary.}\: 
In section \ref{SS:ope commut}
we define the notion of locality
and prove that in the case of fields it reduces to the notion
of holomorphic locality if one field is holomorphic.
In section \ref{SS:vacuum skew}
we prove that the field-state correspondence 
is a morphism of partial $\K^2$-fold modules from 
creative and mutually local distributions to $V$.
In section \ref{SS:skew symm}
we define the notion of skew-symmetry and 
prove that mutually local distributions satisfy skew-symmetry.
In section \ref{SS:goddard}
we prove that if the field-state correspondence 
from a local set of strongly creative distributions is surjective
then it is bijective.

In section \ref{SS:mult local dong}
we define the notion of multiple locality and prove that
the $\vn$-th products of $\End(V)\set{\vz}$ preserve the
relation of multiple locality.
In section \ref{SS:locality of fields}
we define the notion of additive locality and prove that
it is equivalent to locality of OPE-length at most one.
In section \ref{SS:addi local dong}
we prove that if two mutually local fields
are additively local to a third field 
then the $\vn$-th product of any two of them
is mutually local to the third.

\medskip

{\bf Conventions.}\: 
Except when we discuss $S$-fold modules, 
we denote by $V$ a vector space,
$1$ is an even vector, and $\vT$ is a pair of even operators of $V$.
Unless stated otherwise, all distributions take values in $\End(V)$.

\subsection{Locality}
\label{SS:ope commut}
We define the notion of {\it locality} 
for a pair of distributions $a(\vz)$ and $b(\vz)$ 
and prove that 
if $a(z)$ is a holomorphic field and $b(\vz)$ is a field 
then locality is equivalent to holomorphic locality.

\bigskip

{\bf Definition.}\:
Distributions $a(\vz)$ and $b(\vz)$ are called {\bf mutually local} 
\index{local fields} if \index{mutually local fields} 
$a(\vz)$ and $b(\vz)$ are OPE-finite,
the pole orders of $a(\vz)b(\vw)$ are contained in $\vbbK$, and we have
$$
\paraab\, b(\vw)a(\vz)
\; =\;
\sum_{\vh\in\vbbK}\:
\frac{N_{\vh}(a(\vz)b(\vw))}{(\vz-\vw)_{w>z}^{\vh}}.
$$

\bigskip

From Proposition \ref{SS:ope finite} follows that 
distributions $a(\vz)$ and $b(\vz)$ are mutually local 
if and only if 
there exist fields $c^i(\vz,\vw)$ and $\vh_i\in\vbbK$
such that \eqref{E:check ope} and
$$
\paraab\, b(\vw)a(\vz)
\; =\;
\sum_{i=1}^r\:
\frac{c^i(\vz,\vw)}{(\vz-\vw)_{w>z}^{\vh_i}}
$$
are satisfied.
Because 
$(\vz-\vw)^{\vh}=(-1)^{h-\bh}(\vw-\vz)_{z>w}^{\vh}$ 
for any $\vh\in\vbbK$
we see that
locality is a symmetric relation,
the pole orders are independent of the order of 
$a(\vz)$ and $b(\vz)$, 
and for any $\vh\in\vbbK$ we have
\begin{equation}
\label{E:N(ba)=N(ab)}
\paraab\, N_{\vh}(b(\vw)a(\vz))
\; =\;
(-1)^{h-\bh}\, N_{\vh}(a(\vz)b(\vw)).
\end{equation}

\bigskip

{\bf Proposition.}\: {\it
$\iti$\:
If a holomorphic distribution $a(z)$ 
and a distribution $b(\vz)$ are mutually local 
then they are holomorphically local.

$\itii$\:
If a holomorphic field $a(z)$ 
and a field $b(\vz)$ are holomorphically local 
then they are mutually local.
Moreover, 
if $a(z)b(\vw)$ is non-zero then
$a(z)b(\vw)$ has OPE-length one,
the mode infimum $N$ of $a(z)$ and $b(\vz)$
is the least integer such that 
$(z-w)^N[a(z),b(\vw)]=0$, 
and the pole order of $a(z)b(\vw)$ is equal to $(N,0)$.
}

\bigskip

\begin{pf}
$\iti$\:
Proposition \ref{SS:ope holom distr}
implies that there exist OPEs
$a(z)b(\vw)=c(z,\vw)(z-w)^{-h}$ and
$\paraab b(\vw)a(z)=c(z,\vw)(z-w)_{w>z}^{-h}$ 
where $c(z,\vw)$ is a field and $h$ is a non-negative integer.
Thus $(z-w)^h[a(z),b(\vw)]=c(z,\vw)-c(z,\vw)=0$.

$\itii$\:
Suppose that $(z-w)^N[a(z),b(\vw)]=0$ for some integer $N$.
The distribution $c(z,\vw):=(z-w)^N a(z)b(\vw)$ is contained in 
$\cF\set{z}\sqbrack{\vw}(V)$.
From 
$c(z,\vw)=\paraab (z-w)^N b(\vw)a(z)$ 
follows that
$c(z,\vw)$ is also contained in 
$\cF\set{\vw}\lau{z}(V)$.
Thus $c(z,\vw)$ is a field and we obtain
$a(z)b(\vw)=c(z,\vw)(z-w)^{-N}$ and
$\paraab b(\vw)a(z)=c(z,\vw)(z-w)_{w>z}^{-N}$.
Together with Proposition \ref{SS:nth prod ope fin}
this shows everything.
\end{pf}

\bigskip

If a holomorphic field $a(z)$ and a field $b(\vz)$ are 
mutually local with OPE-length one and pole order $\vh\in\Z\times\set{0}$
then we say that $a(z)$ and $b(\vz)$ are mutually local of order $h$
\index{mutually local of order $h$} and \index{local of order $h$}
we call $h$ the pole order \index{pole order} of $a(z)b(\vw)$.
If $a(z)$ and $b(\vz)$ are mutually local with OPE-length zero
then we define the pole order of $a(z)b(\vw)$ to be $-\infty$.

For example, 
$1(z)$ is local to any field and if $a(\vz)$ is a non-zero field
then the pole order of $1(z)a(\vw)=a(\vw)$
is zero because $(z-w)^n a(\vw)$ is not a field
for any negative integer $n$.

\bigskip

{\bf Lemma.}\: {\it
Let $a(\vz)$ and $b(\vz)$ be mutually local distributions
with pole orders $\vh_{ab}^1, \dots, \vh_{ab}^r$,
$c$ be a vector such that
$a(\vz)$ and $b(\vz)$ are bounded on $c$ by 
$(\vh_{ac}^i)$ and $(\vh_{bc}^j)$, resp.,
and $\vh\in\vbbK$.

\begin{enumerate}
\item[$\iti$]
The support in $\vz$ of $N_{\vh}(a(\vz)b(\vw))c$
is contained in $\sum_i \vh_{ac}^i+\Z^2$.

\item[$\itii$]
The support in $\vz$ of $a(\vz)b(\vw)c$ 
is contained in
$\sum_{k,j} \vh_{ab}^k+\vh_{ac}^j+\Z^2$.

\item[$\itiii$]
If the sets 
$\vh_{ab}^k+\set{\vh_{ac}^i\mid i}+\Z^2$
are pairwise disjoint for $k=1,\dots,r$
then
$N_{\vh}(a(\vz)b(\vw))c$ 
is bounded in $\vw$ by $(\vh_{bc}^j)$.
\end{enumerate}
}

\bigskip

\begin{pf}
$\iti$\:
This follows from Remark \ref{SS:reduced ope} because
$N_{\vh}(a(\vz)b(\vw))$ is equal to
$\paraab\,(-1)^{\vh}N_{\vh}(b(\vw)a(\vz))$.

$\itii$\:
This follows from $\iti$ and the OPE of $a(\vz)b(\vw)$.

$\itiii$\:
Define $S_l:=(\vh_{ab}^l+\set{\vh_{ac}^i\mid i}+\Z^2)\times\K^2$
for $l=1, \dots, r$.
Part $\iti$ and the assumption imply that for any $l$ we have
$$
a(\vz)b(\vw)c|_{S_l}
=
\Bigg(\sum_{k=1}^r
\frac{N_{\vh_{ab}^k}(a(\vz)b(\vw))c}{(\vz-\vw)^{\vh_{ab}^k}}\Bigg)
\Bigg|_{S_l}
=
\frac{N_{\vh_{ab}^l}(a(\vz)b(\vw))c}
{(\vz-\vw)^{\vh_{ab}^l}}.
$$
Because the left-hand side is bounded in $\vw$ by $(\vh_{bc}^j)$
the claim follows. 
\end{pf}

\subsection{Creativity}
\label{SS:vacuum skew}

We prove that the field-state correspondence is a
morphism of partial $\K^2$-fold modules 
for creative and mutually local distributions
and that the $\vn$-th products of mutually local distributions
$a(\vz)$ and $b(\vz)$ are strongly creative if 
$a(\vz)$ is translation covariant and $b(\vz)$ is strongly creative.

\bigskip

{\bf Proposition.}\: {\it
If $a(\vz)$ and $b(\vz)$ 
are creative and mutually local distributions
and $\vn\in\K^2$ then

\begin{enumerate}
\item[$\iti$]
the fields $N_{\vn}(a(\vz)b(\vw))$ 
and $a(\vz)_{(\vn)}b(\vz)$ are creative,

\item[$\itii$]
the distribution $a(\vz)$ is bounded on $s_1(b(\vz))$ by
the pole orders of $a(\vz)b(\vw)$, 

\item[$\itiii$]
we have $s_1(a(\vz)_{(\vn)}b(\vz))\, =\, a_{\vn}s_1(b(\vz))$, and

\item[$\itiv$]
if $a(\vz)$ is translation covariant and 
$b(\vz)$ is strongly creative then the field
$a(\vz)_{(\vn)}b(\vz)$ is strongly creative.
\end{enumerate}
}

\bigskip

\begin{pf}
$\iti$\:
Lemma \ref{SS:ope commut}\,$\itiii$ shows that 
$N_{\vn}(a(\vz)b(\vw))$ is creative.
From the definition of the
$\vn$-th products follows that
$a(\vz)_{(\vn)}b(\vz)$ is creative, too.

$\itii$\:
Let $\vh_1, \dots, \vh_r$ be the pole orders
of $a(\vz)b(\vw)$.
The fields $N_{\vn}(a(\vz)b(\vw))$ are creative.
Thus if we apply both sides of the reduced OPE of $a(\vz)b(\vw)$
to $1$ and put $\vw=0$ we get
\begin{equation}
\label{E:a on b}
a(\vz)s_1(b(\vz))
\; = \;
\sum_{\vh\in\vbbK}\:
\frac{N_{\vh}(a(\vz)b(\vw))1|_{\vw=0}}{\vz^{\vh}}
\; \in \;
\sum_{i=1}^r\:
\vz^{-\vh_i}V\pau{\vz}.
\end{equation}

$\itiii$\:
If $\vh$ is a pole order of $a(\vz)b(\vw)$ and
$\vn\in\N^2$ then \eqref{E:a on b} implies 
$$
a(\vw)_{(\vh-1-\vn)}b(\vw)1|_{\vw=0}
\; =\;
\del_{\vz}^{(\vn)} N_{\vh}(a(\vz)b(\vw))1|_{\vz=\vw=0}
\; =\;
a_{\vh-1-\vn}s_1(b(\vz)).
$$

$\itiv$\:
Let $\vm\in\N^2$. 
From $\itiii$ and 
the fact that $\del_{\vz}$ is a translation generator 
we obtain
\begin{align}
\notag
s_1(\del_{\vz}^{(\vm)}(a(\vz)_{(\vn)}b(\vz)))
&=
\sum_{\vi=0}^{\vm}
s_1\big(
(-1)^{\vi}\binom{\vn}{\vi} a(\vz)_{(\vn-\vi)}b(\vz)
+
a(\vz)_{(\vn)}\del_{\vz}^{(\vm-\vi)}b(\vz)
\big)
\\
\notag
&=
\sum_{\vi=0}^{\vm}
(-1)^{\vi}\binom{\vn}{\vi} a_{\vn-\vi}s_1(b(\vz))
+
a_{\vn}s_1(\del_{\vz}^{(\vm-\vi)}b(\vz)).
\end{align}
On the other hand,
using $\itiii$ and the assumption that 
$a(\vz)$ is translation covariant
we get
\begin{align}
\notag
\vT^{(\vm)}s_1(a(\vz)_{(\vn)}b(\vz))
\; &=\;
\vT^{(\vm)}a_{\vn}s_1(b(\vz))
\\
\notag
&=\;
\sum_{\vi=0}^{\vm}\:
(-1)^{\vi}\binom{\vn}{\vi}\, a_{\vn-\vi}s_1(b(\vz))
\, +\,
a_{\vn}\vT^{(\vm-\vi)}s_1(b(\vz)).
\end{align}
Thus Proposition \ref{SS:creat}\,$\iti$ shows 
that the claim follows from 
$\iti$ and strong creativity of $b(\vz)$.
\end{pf}

\subsection{Skew-Symmetry for Local Fields}
\label{SS:skew symm}

We define the notion of {\it skew-symmetry} and 
prove that mutually local distributions satisfy skew-symmetry.

\bigskip

A $V$-valued distribution $a(\vz)$ is called {\bf statistical}
\index{statistical distribution} if its support is contained in $\vbbK$.

\bigskip

{\bf Definition.}\:
Let $V$ be a partial $\K^2$-fold algebra,
$\vT$ be a pair of even operators of $V$, and $\vn\in\vbbK$.
States $a$ and $b$ are said to satisfy {\bf skew-symmetry} 
\index{skew-symmetry}
for $\vT$ and the $\vn$-th product if
$(a,b)\in D$, $a(\vz)b$ is a statistical vertex series, and
$$
\paraab\, b_{(\vn)}a
\; =\; 
\sum_{\vi\in\N^2}\:
(-1)^{n-\bn+i+\bi}\: \vT^{(\vi)}(a_{(\vn+\vi)}b).
$$

\bigskip

Let $b:\K\times\set{\pm 1}\to\K', (h,\mu)\mapsto (\mu)_b^h$, 
be a power map
that is compatible with an involution $\ka$ of $\K'$.
We define automorphisms $\io_z, \io_{\bz}$, and $\io$ of 
the vector space of $(V\otimes\K')$-valued distributions $a(\vz)$
by $\io_za(\vz)\equiv a(-z,\bz):=
\sum_{\vn\in\K^2}(-1)_b^{-n+1}\, a_{\vn}\vz^{-\vn-1},
\io_{\bz}:=\ka\circ\io_z\circ\ka$, and
$\io:=\io_z\circ\io_{\bz}$.
Because $\io_{\bz}(\vz^{\vn})=
\ka\io_z(\bz^{\ka n} z^{\ka\bn})=
\ka((-1)_b^{\ka\bn} \bz^{\ka n} z^{\ka\bn})=
(-1)_b^{-\bn}\vz^{\vn}$
we obtain
$$
\io:\;
a(\vz)
\;\mapsto\; 
a(-\vz)
\; :=\;
\sum_{\vn\in\K^2}\:
(-1)_b^{-n+\bn}\, a_{\vn}\vz^{-\vn-1}.
$$
This shows that $\io$ induces an involution
$a(\vz)\mapsto a(-\vz)$ of the vector space
of statistical $V$-valued distributions
that does not depend on $b$.

\bigskip

{\bf Remark.}\: 
Let $V$ be a $\K^2$-fold algebra and $\vT$ be a pair of even operators.
States $a$ and $b$ satisfy skew-symmetry for all $\vn$-th products
if and only if $a(\vz)b$ is a statistical vertex series and
$$
\paraab\, 
b(\vz)a
\; =\;
e^{\vz\vT}a(-\vz)b.
$$
Assuming that $T$ and $\bT$ commute, 
this shows that $a$ and $b$ satisfy skew-symmetry 
if and only if $b$ and $a$ do.
If $1$ is an even state then any two of the following three statements 
imply the third one:
$1$ is a left identity; 
$1$ is a strong right identity for $\vT$; 
$1$ and any state of $V$ satisfy skew-symmetry for $\vT$.

\bigskip

{\bf Proposition.}\: {\it
Mutually local distributions $a(\vz)$ and $b(\vz)$ 
of the partial $\K^2$-fold algebra $\End(V)\set{\vz}$
satisfy skew-symmetry for $\del_{\vz}$.
}

\bigskip

\begin{pf}
Let $\vh$ be a pole order of $a(\vz)b(\vw)$ and $\vn\in\N^2$. 
From the definition of the $\vn$-th products
and \eqref{E:del and evaluate}
we obtain
\begin{align}
\notag
&\sum_{\vi\in\N^2}\:
(-1)^{h-\bh+n+\bn+i+\bi}\,
\del_{\vw}^{(\vi)}(a(\vw)_{(\vh-1-\vn+\vi)}b(\vw))
\\
\notag
\; =\;
&\sum_{\vi=0}^{\vn}\:
(-1)^{h-\bh+n+\bn+i+\bi}\,
\del_{\vw}^{(\vi)}(\del_{\vz}^{(\vn-\vi)}N_{\vh}(a(\vz)b(\vw))|_{\vz=\vw})
\\
\label{E:pf skew symm}
\; =\;
&\sum_{\vi=0}^{\vn}\:
(-1)^{h-\bh+n+\bn+i+\bi}\,
(\del_{\vz}+\del_{\vw})^{(\vi)}\del_{\vz}^{(\vn-\vi)}
N_{\vh}(a(\vz)b(\vw))|_{\vz=\vw}.
\end{align}
In general, 
if $\al$ and $\be$ are even elements of a commutative algebra 
and $n$ is a non-negative integer then 
$\sum_{i=0}^n
(-1)^i(\al+\be)^{(i)} \al^{(n-i)}
=(\al-(\al+\be))^{(n)}
=(-1)^n \be^{(n)}$.
Together with \eqref{E:N(ba)=N(ab)} this shows 
that the right-hand side of \eqref{E:pf skew symm}
is equal to
$$
(-1)^{h-\bh}\,\del_{\vw}^{(\vn)}N_{\vh}(a(\vz)b(\vw))|_{\vz=\vw}
\; =\;
\paraab\,
b(\vw)_{(\vh-1-\vn)}a(\vw).
$$
\end{pf}

\subsection{Goddard's Uniqueness Theorem}
\label{SS:goddard}
We prove {\it Goddard's uniqueness theorem}
which states that
if $\cF$ is a local set of creative distributions
such that $s_1(\cF)=V$
then distributions in $\cF$
are determined by their action on $1$.

\bigskip

A set $\cF$ of creative distributions in $\End(V)\set{\vz}$
is called {\bf complete} \index{complete set of distributions} 
if $s_1(\cF)=V$.

\bigskip

{\bf Goddard's Uniqueness Theorem.}\: {\it
Let $\cF$ be a complete local set of creative distributions
and $a(\vz)$ be a creative distribution that is local to $\cF$.

\begin{enumerate}
\item[$\iti$]
If there exists a distribution $b(\vz)$ in $\cF$
such that $a(\vz)1=b(\vz)1$ then $a(\vz)=b(\vz)$.

\item[$\itii$]
If $\cF$ and $a(\vz)$ are strongly creative then
$s_1:\cF\to V$ is a vector space isomorphism 
and $a(\vz)=(s_1|_{\cF})^{-1}(s_1(a(\vz)))$.
\end{enumerate}
}

\bigskip

\begin{pf}
$\iti$\:
We define $c(\vz):=a(\vz)-b(\vz)$ and
prove that $c(\vz)d=0$ for any vector $d$ of $V$.
Because $s_1(\cF)=V$
there exists a distribution $d(\vz)$ in $\cF$ such that $s_1(d(\vz))=d$.
Locality yields
$$
\sum_{\vh\in\vbbK}\:
\frac{N_{\vh}(c(\vz)d(\vw))1}{(\vz-\vw)_{w>z}^{\vh}}
\; =\;
\zeta^{\tc\td}\,
d(\vw)c(\vz)1
\; =\;
0.
$$
Proposition \ref{SS:ope finite}
implies that $N_{\vh}(c(\vz)d(\vw))1=0$ for any $\vh$. 
Thus
$$
c(\vz)d
\; =\;
c(\vz)d(\vw)1|_{\vw=0}
\; =\;
\Bigg(
\sum_{\vh\in\vbbK}\:
\frac{N_{\vh}(c(\vz)d(\vw))1}{(\vz-\vw)^{\vh}}
\Bigg) \Bigg|_{\vw=0}
\; =\;
0.
$$

$\itii$\:
This follows immediately from $\iti$.
\end{pf}

\subsection{Dong's Lemma for Multiple Locality}
\label{SS:mult local dong}

We define the notion of multiple locality and prove 
{\it Dong's lemma for multiple locality}
which states that 
the $\vn$-th products of the partial $\K^2$-fold algebra of fields
preserve the relation of multiple locality.
It is an open question
whether any family of pairwise mutually local fields
is multiply local.
Proposition \ref{SS:addi local dong}
gives a positive answer to this question
in a special case.

\bigskip

If $\si$ is a permutation on $r$ letters, 
$i, j\in\set{1, \dots, r}$, and 
$\vn\in\vbbK$
then we define
$$
(\vz_i-\vz_j)_{\si}^{\vn}
\; :=\;
\begin{cases}
\qquad (\vz_i-\vz_j)^{\vn}&
\text{if $\si^{-1}(i)<\si^{-1}(j)$ and}\\
\qquad (\vz_i-\vz_j)_{z_j>z_i}^{\vn}&
\text{otherwise.}
\end{cases}
$$

\bigskip

{\bf Definition.}\: 
Let $r_{ij}$ be non-negative integers and
$\vh_{ij}^k\in\vbbK$ for any
$i,j=1,\dots,r$ and $k=1,\dots,r_{ij}$.
Distributions $a^1(\vz),\dots, a^r(\vz)$ are called
{\bf multiply local} of orders at most $(\vh_{ij}^k)_{i,j,k}$
if there exist fields $c^{\al}(\vz_1,\dots,\vz_r)$ 
for any
$\al\in\prod_{1\leq i<j\leq r}\set{1,\dots,r_{ij}}$
such that for any permutation $\si$ we have
\begin{equation}
\label{E:mult loc}
\zeta'\,
a^{\si 1}(\vz_{\si 1})\dots  a^{\si r}(\vz_{\si r})
\; =\;
\sum_{\al}\:
\frac{c^{\al}(\vz_1,\dots,\vz_r)}
{\prod_{i<j}(\vz_i-\vz_j)_{\si}^{\vh_{ij}^{\al(i,j)}}}
\end{equation}
where $\zeta'$ is the sign of the permutation $\si$
restricted to the set of $i$ such that $a^i(\vz)$ is odd.

\bigskip

Of course, a pair of distributions is multiply local 
if and only if it is mutually local.

\bigskip

{\bf Dong's Lemma.}\: {\it
If $a^1(\vz), \dots, a^r(\vz)$
are multiply local distributions of orders at most 
$(\vh_{ij}^k)$ and
$a^1(\vz)$ and $a^2(\vz)$ are OPE-finite
of orders at most $\vh_{12}^1, \dots, \vh_{12}^{r_{12}}$
then for any $k\in\set{1, \dots, r_{12}}$ and
$\vn\in\vh_{12}^k+\Z^2$
the distributions 
$a^1(\vz)_{(\vn)}a^2(\vz), a^3(\vz), \dots, a^r(\vz)$
are multiply local of orders at most
$(\vh_{ij}^{k'}+
\de_{i,2}(\vh_{12}^k+\vh_{1j}^{k''}-\vn-1))_{i,j,k',k''}$
where $1<i<j\leq r$, $k'\in\set{1, \dots, r_{ij}}$, and 
$k''\in\set{1, \dots, r_{1j}}$.
}

\bigskip

\begin{pf}
Let $\si$ be a permutation on $r$ letters
such that $\si^{-1}(2)=\si^{-1}(1)+1$.
Because $a^1(\vz)$ and $a^2(\vz)$ are OPE-finite
of orders at most $\vh_{12}^1, \dots, \vh_{12}^{r_{12}}$
there exist fields $c^k(\vz,\vw)$ such that
$a^1(\vz)a^2(\vw)=
\sum_{k=1}^{r_{12}}(\vz-\vw)^{-\vh_{12}^k}\, c^k(\vz,\vw)$.
Thus we have
$$
a^{\si 1}(\vz_{\si 1})\dots  a^{\si r}(\vz_{\si r})
\; =\;
\sum_{k=1}^{r_{12}}
\frac{a^{\si 1}(\vz_{\si 1})\dots c^k(\vz_1,\vz_2)\dots a^{\si r}(\vz_{\si r})}
{(\vz_1-\vz_2)^{\vh_{12}^k}}.
$$
We may assume $h_{12}^i+\Z\ne h_{12}^j+\Z$ 
for any $i\ne j$. Fix $k$.
Equation \eqref{E:mult loc} and
Proposition \ref{SS:ope finite} imply that
\begin{align}
\label{E:mult dong}
\zeta'\,
a^{\si 1}(\vz_{\si 1})\dots c^k(&\vz_1,\vz_2)\dots a^{\si r}(\vz_{\si r})
\\
\notag
\; =\;
&\sum_{\al:\: \al(1,2)=k}\:
\frac{c^{\al}(\vz_1,\dots,\vz_r)}
{\prod_{(i,j)\ne (1,2)}(\vz_i-\vz_j)_{\si}^{\vh_{ij}^{\al(i,j)}}}.
\end{align}

We have
$a^1(\vw)_{(\vh_{12}^k-1-\vn)}a^2(\vw)=
\del_{\vz}^{(\vn)}c^k(\vz,\vw)|_{\vz=\vw}$
for any $\vn\in\N^2$.
Thus by acting with $\del_{\vz_1}^{(\vn)}$ on
\eqref{E:mult dong} and putting 
$\vz_1=\vz_2$ we obtain
\begin{align}
\notag
\zeta'\,
a^{\si 1}(\vz_{\si 1})\dots 
a^1(&\vz_2)_{(\vh_{12}^k-1-\vn)}a^2(\vz_2)\dots a^{\si r}(\vz_{\si r})
\\
\notag
\; =\;
&\sum_{\al:\: \al(1,2)=k}\:
\frac{\tc^{\al}(\vz_2,\dots,\vz_r)}
{\prod_{1<i<j\leq r}(\vz_i-\vz_j)_{\si}
^{\vh_{ij}^{\al(i,j)}+\de_{i,2}(\vh_{1j}^{\al(1,j)}+\vn)}}
\end{align}
for some fields $\tc^{\al}$.
This establishes the claim.
\end{pf}

\subsection{Tensor Product of Fields}
\label{SS:tensor prod fields}

We define the {\it tensor product of fields} and prove
that the tensor products of pairwise mutually local fields are
mutually local.

\bigskip

Let $V$ and $V'$ be vector spaces and $a(\vz)$ and $a'(\vz)$ be
fields on $V$ and $V'$, resp.
Recall that there exists a canonical morphism
$V\sqbrack{\vz}\otimes V'\sqbrack{\vz}\to (V\otimes V')\sqbrack{\vz}$.
The field on $V\otimes V'$ that is defined by
$c\otimes c'\mapsto a(\vz)c\otimes a'(\vz)c'$
is called the {\bf tensor product} of $a(\vz)$ and $a'(\vz)$
and is denoted by $a(\vz)\otimes a'(\vz)$.

\bigskip

{\bf Proposition.}\: {\it
Let $a(\vz)$ and $b(\vz)$ be fields on a vector space $V$ and 
$a'(\vz)$ and $b'(\vz)$ be fields on a vector space $V'$.

\begin{enumerate}
\item[$\iti$]
If the pairs $a(\vz), b(\vz)$ and $a'(\vz), b'(\vz)$ are OPE-finite
with pole orders $\vh_1, \dots, \vh_r$ and $\vh'_1, \dots, \vh'_{r'}$, resp., 
then the fields $a(\vz)\otimes a'(\vz)$ and $b(\vz)\otimes b'(\vz)$
are OPE-finite of orders at most $(\vh_i+\vh'_j)_{i,j}$.
Moreover, for $\vn\in\K^2$ we have
$$
(a(\vz)\otimes a'(\vz))_{(\vn)}(b(\vz)\otimes b'(\vz))
=
\sum_{\vm\in\K^2}
a(\vz)_{(\vm)}b(\vz)\otimes a'(\vz)_{(\vn-1-\vm)}b'(\vz).
$$

\item[$\itii$]
If the pairs $a(\vz), b(\vz)$ and $a'(\vz), b'(\vz)$ are both mutually local
then the fields $a(\vz)\otimes a'(\vz)$ and $b(\vz)\otimes b'(\vz)$
are mutually local as well.

\item[$\itiii$]
If $a(\vz)$ and $a'(\vz)$ are translation covariant for 
the pairs of operators $\vT$ and $\vT'$ on $V$ and $V'$, resp.,
then $a(\vz)\otimes a'(\vz)$ is translation covariant for
$\vT\otimes\id_{V'} +\id_V\otimes\vT'$.

\item[$\itiv$]
If $a(\vz)$ and $a'(\vz)$ are dilatation covariant of weights 
$\vh$ and $\vh'$ for the pairs of even operators $\vH$ and $\vH'$ on 
$V$ and $V'$, resp.,
then $a(\vz)\otimes a'(\vz)$ is dilatation covariant of weight $\vh+\vh'$
for $\vH\otimes\id_{V'} +\id_V\otimes\vH'$.
\end{enumerate}
}

\bigskip

\begin{pf}
$\iti$\:
We have
\begin{align}
\notag
&(a(\vz)\otimes a'(\vz))(b(\vw)\otimes b'(\vw))
\\
\notag
=\; 
&a(\vz)b(\vw)\otimes a'(\vz)b'(\vw)
\\
\label{E:tensor ope}
=\;
&\sum_{\vh, \vh'\in\K^2}\:
\frac{N_{\vh}(a(\vz)b(\vw))\otimes N_{\vh'}(a'(\vz)b(\vw))}
{(\vz-\vw)^{\vh+\vh'}}.
\end{align}
This proves the first claim. 

For fields $c(\vz,\vw)$ and $c'(\vz,\vw)$ on $V$ and $V'$, resp.,
and $\vk\in\N^2$ we have 
$$
\del_{\vz}^{(\vk)}(c(\vz,\vw)\otimes c'(\vz,\vw))|_{\vz=\vw}
\; =\;
\sum_{\vi=0}^{\vk}\:
\del_{\vz}^{(\vi)}c(\vz,\vw)|_{\vz=\vw}
\: \otimes\:
\del_{\vz}^{(\vk-\vi)}c'(\vz,\vw)|_{\vz=\vw}.
$$
Thus the second claim follows from \eqref{E:tensor ope} and the identity
$(\vh+\vh'-1-\vk)-1=(\vh-1-\vi)+(\vh'-1-(\vk-\vi))$.

$\itii$\:
This follows from \eqref{E:tensor ope}.

$\itiii$\:
This follows from 
\begin{equation}
\label{E:del a otimes a'}
\del_z(a(\vz)\otimes a'(\vz))
\; =\;
\del_z a(\vz)\otimes a'(\vz)
\; +\;
a(\vz)\otimes \del_z a'(\vz).
\end{equation}

$\itiv$\:
This also follows from \eqref{E:del a otimes a'}.
\end{pf}

\subsection{Additive Locality}
\label{SS:locality of fields}

We define the notion of {\it additive locality} 
and
prove that additive locality is equivalent
to locality with OPE-length at most one.

\bigskip

{\bf Definition.}\:
For $\vh\in\vbbK$,
fields $a(\vz)$ and $b(\vz)$ are called {\bf additively local}
\index{additively local fields} of order at most $\vh$ if 
$$
(\vz-\vw)^{\vh}[a(\vz),b(\vw)]
\; =\; 
0.
$$

\bigskip

Note that the definition of additive locality involves radial ordering.

Fields $a(\vz)$ and $b(\vz)$ are called additively local of order $\vh$
if $\vh$ is the infimum of the set of $\vh'\in\vh+\Z^2$ 
such that $a(\vz)$ and $b(\vz)$ 
are additively local of order at most $\vh'$.

\bigskip

{\bf Remark.}\: {\it
For fields $a(\vz)$ and $b(\vz)$ and $\vh\in\vbbK$,
the following statements are equivalent:
\begin{enumerate}
\item[$\iti$]
$a(\vz)$ and $b(\vz)$ are additively local of order at most $\vh$;
\item[$\itii$]
$a(\vz)$ and $b(\vz)$
are mutually local with OPE-length at most one and
if the OPE-length is one and
$\vh_0$ is the pole order of $a(\vz)b(\vw)$
then $\vh\geq\vh_0$;
\end{enumerate}
the following statements are equivalent:
\begin{enumerate}
\item[$\itip$]
$a(\vz)$ and $b(\vz)$ are additively 
local of order $\vh+\vn$ for any $\vn\in\Z^2$;
\item[$\itiip$]
$a(\vz)$ and $b(\vz)$
are mutually local with OPE-length zero;
\item[$\itiiip$]
$a(\vz)b(\vz)=b(\vz)a(\vz)=0$;
\end{enumerate}
and the following statements are equivalent:
\begin{enumerate}
\item[$\itipp$]
$a(\vz)$ and $b(\vz)$ are additively local of order $\vh$;
\item[$\itiipp$]
$a(\vz)$ and $b(\vz)$
are mutually local with OPE-length one and pole order $\vh$.
\end{enumerate}
}

\bigskip

\begin{pf}
$\iti\Rightarrow\itii$\:
The distribution
$c(\vz,\vw):=(\vz-\vw)^{\vh}a(\vz)b(\vw)$
is contained in $\cF\set{\vz}\sqbrack{\vw}(V)$.
From 
$c(\vz,\vw)=\paraab(\vz-\vw)^{\vh}b(\vw)a(\vz)$ 
follows that $c(\vz,\vw)$ is also contained in
$\cF\set{\vw}\sqbrack{\vz}(V)$.
Thus $c(\vz,\vw)$ is a field.
Because
$a(\vz)b(\vw)$ and $c(\vz,\vw)$
are contained in the $\K\sqbrack{\vz}\sqbrack{\vw}$-module
$\cF\sqbrack{\vz}\sqbrack{\vw}(V)$
and
$b(\vw)a(\vz)$ and $c(\vz,\vw)$
are contained in the $\K\sqbrack{\vw}\sqbrack{\vz}$-module
$\cF\sqbrack{\vw}\sqbrack{\vz}(V)$
we obtain
$$
a(\vz)b(\vw)
\; =\;
\frac{c(\vz,\vw)}{(\vz-\vw)^{\vh}}
\qquad\text{and}\qquad
\paraab\, b(\vw)a(\vz)
\; =\;
\frac{c(\vz,\vw)}{(\vz-\vw)_{w>z}^{\vh}}.
$$
The second claim of $\itii$ follows from
the uniqueness of the reduced OPE.

$\itii\Rightarrow\iti$\:
Of course,
if the OPE-length is zero then $a(\vz)$ and $b(\vz)$ are 
additively local of order at most $\vh'$ for any 
$\vh'\in\vbbK$.
If the OPE-length is one and $\vh=\vh_0+\vn$ where $\vn\in\N^2$
then
$$
(\vz-\vw)^{\vh}a(\vz)b(\vw)
\; =\;
(\vz-\vw)^{\vn}N_{\vh}(a(\vz)b(\vw))
\; =\;
\paraab(\vz-\vw)^{\vh}b(\vw)a(\vz).
$$
Thus $a(\vz)$ and $b(\vz)$ are additively local of order
at most $\vh$.

$\itip\Rightarrow\itiiip$\:
This follows from Lemma \ref{SS:reduced ope}\,$\itii$.

$\itiip\Leftrightarrow\itiiip\Rightarrow\itip$\:
This is trivial.

$\itipp\Leftrightarrow\itiipp$\:
This follows from the equivalence of $\iti$ and $\itii$.
\end{pf}

\bigskip

If $a(\vz)$ and $b(\vz)$ are fields 
that are additively local of order at most $\vh$ then
$$
a(\vw)_{(\vn)}b(\vw)
\; =\;
\del_{\vz}^{(\vh-1-\vn)}
((\vz-\vw)^{\vh}a(\vz)b(\vw))|_{\vz=\vw}.
$$
Proposition \ref{SS:ope commut} shows that
a holomorphic field $a(z)$ and a field $b(\vz)$
are additively local of order $\vh$
if and only if they are mutually local of order $h$ and $\bh=0$.

\subsection{Dong's Lemma for Additive Locality}
\label{SS:addi local dong}

We prove {\it Dong's lemma for additive locality}
which states that if 
two mutually local fields
are additively local to a third field 
then the $\vn$-th product of any two of them
is mutually local to the third.

\bigskip

{\bf Proposition.}\: {\it
Let $a^1(\vz), \dots, a^r(\vz)$ be pairwise mutually local
fields and 
$\vh_{ij}^1, \dots, \vh_{ij}^{r_{ij}}$ 
be the pole orders of $a^i(\vz)a^j(\vw)$.
If $a^i(\vz)$ and $a^j(\vz)$ are additively local
for any $i<j$ such that $(i,j)\ne (1,2)$ 
then $a^1(\vz), \dots, a^r(\vz)$ are multiply local
of orders at most $(\vh_{ij}^k)$.
}

\bigskip

\begin{pf}
For a permutation $\si$ on $r$ letters,
define $F_{\si}$ as the product of 
$(\vz_i-\vz_j)_{\si}^{\vh_{ij}^1}$
over any $i,j=1, \dots, r$ such that
$i<j, (i,j)\ne (1,2)$, and $r_{ij}=1$.
Let \eqref{E:check ope} be the reduced OPE of 
$a^1(\vz)a^2(\vw)$.
Denote by $S_{\pm}$ the set of
$\si$ such that $\si^{-1}(2)=\si^{-1}(1)\pm 1$.
For $\si\in S_+\cup S_-$, we define
$d_{\si}^i(\vz_1,\dots,\vz_r)$ to be the
distribution obtain from
$\zeta' F_{\si}a^{\si 1}(\vz_{\si 1})\dots a^{\si r}(\vz_{\si r})$
by replacing
$a^1(\vz_1)a^2(\vz_2)$ or $a^2(\vz_2)a^1(\vz_1)$ by
$c^i(\vz_1,\vz_2)$;
in other words, we define
\begin{align}
\notag
&d_{\si}^i(\vz_1,\dots,\vz_r)
\\
\notag
:=\, 
\zeta'&\,F_{\si}\,
a^{\si 1}(\vz_{\si 1})\dots 
a^{\si k}(\vz_{\si k})\, c^i(\vz_1,\vz_2)\, a^{\si (k+3)}(\vz_{\si (k+3)})
\dots a^{\si r}(\vz_{\si r})
\end{align}
where $k:=\min(\si^{-1}(1),\si^{-1}(2))-1$
and 
$\zeta'$ is defined in Definition \ref{SS:mult local dong}.
Let $\io:S_+\to S_-$ be the bijection 
such that
$\si^{-1}(2)=(\io\si)^{-1}(1)$ and
$\si^{-1}(i)=(\io\si)^{-1}(i)$ for any $i>2$.
By definition, we have $d_{\si}^i=d_{\io\si}^i$.

From the OPE \eqref{E:check ope} follows that
for any $\si\in S_+\cup S_-$ we have
\begin{equation}
\label{E:pf dong additive}
\zeta'\, F_{\si}\,
a^{\si 1}(\vz_{\si 1})\dots a^{\si r}(\vz_{\si r})
\; =\;
\sum_{i=1}^{r_{12}}\:
\frac{d_{\si}^i(\vz_1,\dots,\vz_r)}
{(\vz_1-\vz_2)_{\si}^{\vh_{12}^i}}.
\end{equation}
The left-hand side of \eqref{E:pf dong additive}
only depends on whether $\si$ lies in $S_+$ or $S_-$.
Thus Proposition \ref{SS:ope finite} implies that 
likewise 
$d_{\si}^i$ only depends on whether $\si$ lies in $S_+$ or $S_-$.
From $d_{\si}^i=d_{\io\si}^i$ we get that
$d^i:=d_{\si}^i$ is independent of $\si\in S_+\cup S_-$.
Because for any $j=2, \dots, r$
there exists $\si\in S_+$ such that 
$\si r=j$ and 
because
$a^1(\vz), \dots, a^r(\vz)$ and $c^i(\vz,\vw)$ are fields
the distributions $d^i$ are fields, too.

For an arbitrary permutation $\si$ in $\bbS_r$,
the left-hand side of \eqref{E:pf dong additive}
only depends on whether $\si^{-1}(1)$
is smaller or bigger than $\si^{-1}(2)$.
Therefore \eqref{E:pf dong additive} holds for any $\si$
if we define $d_{\si}^i:=d^i$.
Dividing \eqref{E:pf dong additive} by $F_{\si}$
we obtain multiple locality. 
\end{pf}

\bigskip

The Proposition combined with 
Dong's lemma for multiple locality 
yields the following result.

\bigskip

{\bf Dong's Lemma.}\: {\it
Let $a(\vz)$ and $b(\vz)$ be mutually local fields
with pole orders $\vh_1, \dots, \vh_r$
and
$c(\vz)$ be a field that is additively local to both
$a(\vz)$ and $b(\vz)$
of orders at most $\vh$ and $\vh'$.
Then
for any $k\in\set{1, \dots, r}$ and $\vn\in\vh_k+\Z^2$
the fields $a(\vz)_{(\vn)}b(\vz)$ and $c(\vz)$
are additively local of order at most 
$\vh_k+\vh+\vh'-\vn-1$.
Moreover,
for any $\vn\in\K^2$
the pairs of fields $a(\vz), b(\vz)_{(\vn)}c(\vz)$
and $a(\vz), c(\vz)_{(\vn)}b(\vz)$ 
are both mutually local of orders at most
$\vh_1+\vh+\vh'-\vn-1, \dots, \vh_r+\vh+\vh'-\vn-1$.
\hfill $\square$
}

%\contentsline {chapter}{Preface}{V}  put this into book.toc

\chapter{OPE-Algebras}
\label{C:opea}

\section{Duality}
\label{S:duality}

{\bf Summary.}\:
In section \ref{SS:skew symm states}
we prove that creative, mutually local states
$a$ and $b$ satisfy skew-symmetry 
if $a(\vz)$ is translation covariant and $b$ is strongly creative.
In section \ref{SS:h Jac id distr}
we rewrite the holomorphic Jacobi identity and holomorphic duality
in terms of distributions.
In section \ref{SS:duality}
we define the notions of duality and additive duality
and prove that duality of dual length at most one
is equivalent to additive duality.
In section \ref{SS:dual and local}
we prove that if there exists a translation endomorphism
then duality for $a, b$ and $b, a$ 
together with skew-symmetry for $a$ and $b$ 
implies locality for $a$ and $b$.

\subsection{Skew-Symmetry for Local States}
\label{SS:skew symm states}

We prove that creative, mutually local states
$a$ and $b$ of a $\K^2$-fold algebra satisfy skew-symmetry
if $a(\vz)$ is translation covariant and $b$ is strongly creative.

\bigskip

For a $\K^2$-fold module $M$ over a vector space $V$,
vectors $a$ and $b$ of $V$ are called 
mutually local \index{local states} on $M$
if $a(\vz)$ and $b(\vz)$ are mutually local.
Two states of a $\K^2$-fold algebra 
are called mutually local if 
they are mutually local on $V$.
Similar terminology applies to the notion of multiple locality.
\index{multiply local states}

\bigskip

{\bf Proposition.}\: {\it
Let $V$ be a bounded $\vbbK$-fold algebra, $1$ be a state, and
$\vT$ be a pair of even operators of $V$.
If $a$ and $b$ are creative, mutually local states,
$a(\vz)$ is translation covariant, and
$b$ is strongly creative then 
$a$ and $b$ satisfy skew-symmetry.
}

\bigskip

\begin{pf}
Proposition \ref{SS:skew symm} shows that
$a(\vz)$ and $b(\vz)$
satisfy skew-symmetry for $\del_{\vz}$.
From Proposition \ref{SS:vacuum skew}\,$\itiii$, $\itiv$
and Proposition \ref{SS:creat}\,$\iti$
follows that for any $\vn\in\vbbK$ we have
\begin{align}
\notag
\paraab\, b_{(\vn)}a
\; &=\;
\paraab\, s_1(b(\vz)_{(\vn)}a(\vz))
\\
\notag
&=\;
\sum_{\vi\in\N^2}\:
(-1)^{n-\bn+i+\bi}\, s_1(\del_{\vz}^{(\vi)}(a(\vz)_{(\vn+\vi)}b(\vz)))
\\
\notag
&=\;
\sum_{\vi\in\N^2}\:
(-1)^{n-\bn+i+\bi}\, \vT^{(\vi)}(a_{(\vn+\vi)}b).
\end{align}
\end{pf}

\subsection{Holomorphic Jacobi Identity in terms of Distributions}
\label{SS:h Jac id distr}

We rewrite the holomorphic Jacobi identity and holomorphic duality
in terms of distributions.

\bigskip

{\bf Remark.}\: {\it
Let $M$ be a bounded $\K^2$-fold module over a 
bounded $\Z$-fold algebra $V$, 
$a, b\in V$, $c\in M$, $r, t\in\Z$, and $\vs\in\K^2$
such that $a$ is holomorphic on $M$. 

$\iti$\:
The vectors $a, b, c$ satisfy the holomorphic Jacobi identity
for indices $r, \vs, t$ if and only if 
the $(t,\vs,r)$-th modes of
$$
\de(z,w+x)(a(x)b)(\vw)c 
\qquad \text{and} \qquad
\de(x,z-w)[a(z),b(\vw)]c
$$
are equal.

$\itii$\:
The vectors $a, b, c$ satisfy holomorphic duality of order at most $t$
for indices $r$ and $\vs$ if and only if 
the $(\vs,r)$-th modes of
$$
(w+x)^t(a(x)b)(\vw)c
\qquad \text{and} \qquad
(x+w)^t a(x+w)b(\vw)c
$$ 
are equal.
}

\bigskip

\begin{pf}
$\iti$\:
We have
\begin{align}
\notag
&\de(z,w+x)\, (a(x)b)(\vw)c
\\
\notag
=\;
&\sum_{t\in\Z}\: (w+x)^t\, (a(x)b)(\vw)c\; z^{-t-1}
\\
\notag
=\;
&\sum_{t,r\in\Z}\,\sum_{\vs\in\vbbK}
\bigg(
\sum_{i\in\N}\:
\binom{t}{i}\, (a_{(r+i)}b)_{(s+t-i,\bs)}c
\bigg)
\;
z^{-t-1} \vw^{-\vs-1} x^{-r-1}
\end{align}
and
\begin{align}
\notag
&
\de(x,z-w)[a(z),b(\vw)]c
\\
\notag
=\; 
&\sum_{r\in\Z}\:
(z-w)^r [a(z),b(\vw)]c\; x^{-r-1}
\\
\notag
=\;
&\sum_{r,t\in\Z}\,\sum_{\vs\in\K^2}
\bigg(
\sum_{i\in\N}\:
(-1)^i\binom{r}{i}\,
(a_{(t+r-i)} b_{(s+i,\bs)}
-\paraab (-1)^r b_{(s+r-i,\bs)}a_{(t+i)})c
\bigg) 
\\
\notag
&\qquad\qquad\qquad\qquad\qquad\qquad\qquad\qquad\qquad\qquad\qquad\quad
z^{-t-1} \vw^{-\vs-1} x^{-r-1}.
\end{align}

$\itii$\:
It is clear that 
the left-hand side of holomorphic duality 
is equal to the 
$(\vs,r)$-th mode of $(w+x)^t(a(x)b)(\vw)c$.
The right-hand side of holomorphic duality 
is equal to the 
$(\vs,r)$-th mode of 
$(x+w)^t a(x+w)b(\vw)c$ because
\begin{align}
\notag
&(x+w)^t a(x+w)b(\vw)c
\\
\notag
=\;
&\sum_{r\in\Z, \vs\in\K^2}
a_{(t+r)}b_{(\vs)}c\, (x+w)^{-r-1}\vw^{-\vs-1}
\\
\notag
=\;
&\sum_{\vr\in\Z, \vs\in\K^2}
\bigg( \sum_{i\in\N}\:
\binom{-(r-i)-1}{i}\,
a_{(t+r-i)}b_{(s+i,\bs)}c\bigg)\; 
x^{-r-1}\vw^{-\vs-1}
\end{align}
and because
$\binom{-r-1+i}{i}=(-1)^i\binom{r}{i}$.
\end{pf}

\bigskip

States $a$ and $b$ of $V$ are called
{\bf holomorphically dual} of order at most $t$
on a vector $c$ of $M$ if $a$ is holomorphic on $M$
and 
$a, b, c$ satisfy holomorphic duality of order $t$.
\index{holomorphically dual states}

\subsection{Duality and Additive Duality}
\label{SS:duality}

We define the notions of 
{\it duality} and {\it additive duality}
and prove that duality of dual length at most one
is equivalent to additive duality.

\bigskip

Let $M$ be a bounded $\vbbK$-fold module over a bounded
$\vbbK$-fold algebra $V$, $a, b\in V$, and $c\in M$.
A family of $M$-valued vertex series $d^j(\vz,\vw)$ together with
$\vk_j\in\K^2$ where $j=1, \dots, s$ 
are called a {\bf dual OPE}  \index{dual OPE} of $a, b, c$ if
\begin{equation}
\label{E:duality1}
a(\vx+\vw)b(\vw)c
\; =\;
\sum_{j=1}^s\:
\frac{d^j(\vx,\vw)}{(\vx+\vw)^{\vk_j}},
\qquad
(a(\vx)b)(\vw)c
\; =\;
\sum_{j=1}^s\:
\frac{d^j(\vx,\vw)}{(\vx+\vz)^{\vk_j}}.
\end{equation}
We call $s$ the {\bf length}, \index{length!of dual OPE}
$d^j(\vz,\vw)$ the {\bf numerators}, \index{numerator!of dual OPE}
and $\vk_j$ the {\bf dual pole orders} \index{dual pole order!of dual OPE}
of the dual OPE.
For the sake of brevity, 
we call \eqref{E:duality1} a dual OPE of $a, b, c$ and 
we call the first and second equation of \eqref{E:duality1}
a dual OPE for $z-w>w$ \index{dual OPE!for $z-w>w$}
and for $w>z-w$, \index{dual OPE!for $w>z-w$}
resp.

\bigskip

{\bf Definition.}\:
Let $M$ be a bounded $\vbbK$-fold module over a bounded
$\vbbK$-fold algebra $V$, $a, b\in V$, $c\in M$, and 
$\vh_{ab}^i, \vh_{ac}^j, \vt\in\vbbK$.

$\iti$\:
The states $a$ and $b$ are called 
{\bf mutually dual} \index{mutually dual states}
on \index{dual states} $c$ 
of orders at most $\vh_{ab}^1, \dots, \vh_{ab}^r$
and of dual orders at most $\vh_{ac}^1, \dots, \vh_{ac}^s$
if there exists a dual OPE of $a, b, c$ 
with dual pole orders $\vh_{ac}^1, \dots, \vh_{ac}^s$
such that the numerators $d^j(\vz,\vw)$ are bounded in $\vz$ by
$\vh_{ab}^1, \dots, \vh_{ab}^r$.

$\itii$\:
The states $a$ and $b$ are called 
{\bf additively dual} \index{additively dual states}
on $c$ of order at most $\vt$ if 
$$
(\vw+\vx)^{\vt}\, (a(\vx)b)(\vw)c
\; =\;
(\vx+\vw)^{\vt}\, a(\vx+\vw)b(\vw)c.
$$

\bigskip

States $a$ and $b$ are called mutually dual on $M$
of orders at most $\vh_{ab}^1, \dots, \vh_{ab}^r$
if they are mutually dual 
of orders at most $\vh_{ab}^1, \dots, \vh_{ab}^r$ 
on any vector of $M$.
States $a$ and $b$ are called dual on $M$
with {\bf pole orders} $\vh_{ab}^1, \dots, \vh_{ab}^r$
in the case that for any $\vn_i\N^2$ 
the states $a$ and $b$ are dual 
of orders at most $\vh_{ab}^1+\vn_1, \dots, \vh_{ab}^r+\vn_r$
if and only if $\vn_i=0$ for any $i$.

Assume that the underlying vector space of $M$ is endowed
with a gradation. 
States $a$ and $b$ are called additively dual on $M$
if they are additively dual on any homogeneous vector of $M$.

From Proposition \ref{SS:ope finite} and
Lemma \ref{SS:reduced ope} follows that 
there exists a unique {\bf reduced} dual OPE, \index{reduced dual OPE}
i.e.~such that $\vh_{ac}^j\notin \vh_{ac}^k+\Z^2$ for any $j\ne k$
and $(\vz+\vw)^{-\vn}d^j(\vz,\vw)$ is {\it not} a vertex series
for any $j$ and $\vn\in\N^2\setminus\set{0}$.

If \eqref{E:duality1} is the reduced dual OPE of $a, b, c$  
then we call $s$ the {\bf dual length}, \index{dual length}
we call $\vh_{ac}^1, \dots, \vh_{ac}^s$ the {\bf dual pole order},
\index{dual pole order!of vectors}
and we define the $\vh$-th {\bf dual numerator} \index{dual numerator}
by $M_{\vh}(a(\vz)b(\vw)c):=d^j(\vz,\vw)$
if $\vh=\vh_{ac}^j$ and
$M_{\vh}(a(\vz)b(\vw)c):=0$ otherwise.

\bigskip

{\bf Remark.}\: {\it
Let $M$ be a bounded $\vbbK$-fold module over a bounded
$\vbbK$-fold algebra $V$ and 
let the vector space $M$ be endowed with a gradation.
For states $a$ and $b$ of $V$, the following statements are equivalent:

\begin{enumerate}
\item[$\iti$]
$a$ and $b$ are additively dual on $M$;

\item[$\itii$]
$a$ and $b$ are mutually dual on $M$ and 
if $c$ is a homogeneous vector of $M$ then
the dual length of $a, b, c$ is at most one.
\end{enumerate}
}

\bigskip

\begin{pf}
Except for one additional argument
this is proven in the same way as the 
equivalence of $\iti$ and $\itii$ of Remark \ref{SS:locality of fields}.
Namely, in the proof that $\iti$ implies $\itii$
we also need to show that there exist 
$\vh_1, \dots, \vh_r\in\vbbK$ 
that do not depend on $c$ and such that the numerator
$(\vw+\vx)^{\vt} (a(\vx)b)(\vw)c$ is bounded in $\vx$
by $\vh_1, \dots, \vh_r$.
This follows from the assumption that $a$ is bounded on $b$.
\end{pf}

\bigskip

{\bf Lemma.}\: {\it
Let $M$ be a bounded $\vbbK$-fold module over a bounded
$\vbbK$-fold algebra $V$,
$a, b\in V$, $c\in M$, and
$\vh_{ab}^i$ and $\vh_{bc}^j$ be the bounds of $a(\vz)b$ and $b(\vz)c$, resp.
If $a$ and $b$ are dual 
then for any $\vh\in\vbbK$ the vertex series
$M_{\vh}(a(\vz)b(\vw)c)$ lies in
$\sum_{i,j}\vz^{-\vh_{ab}^i}\vw^{-\vh_{bc}^j}V\lau{\vz,\vw}$.
}

\bigskip

\begin{pf}
This is proven in the same way as 
Remark \ref{SS:reduced ope} and Lemma \ref{SS:ope commut}.
\end{pf}

\subsection{Duality and Skew-Symmetry}
\label{SS:dual and local}

We prove that if there exists a translation endomorphism
then duality for $a, b$ and $b, a$ 
together with skew-symmetry for $a$ and $b$ 
implies locality for $a$ and $b$.

\bigskip

{\bf Proposition.}\: {\it
Let $M$ be a bounded $\vbbK$-fold module over a bounded
$\vbbK$-fold algebra $V$,
$\vT$ be a pair of even operators of $V$,
and $a$ and $b$ be states of $V$ such that 
$a$ and $b$ are dual with pole orders
$\vh_{ab}^1,\dots, \vh_{ab}^r$ and 
$b$ and $a$ are dual with pole order
$\vh_{ba}^1,\dots, \vh_{ba}^s$.
If $a$ and $b$ satisfy skew-symmetry 
and $\vT$ is a translation endomorphism of $M$ 
then $a$ and $b$ are mutually local with
pole orders $\vh_{ab}^1, \dots, \vh_{ab}^r$.
Moreover, we have $r=s$ and 
there exists a permutation $\si$ such that
$\vh_{ba}^i=\vh_{ab}^{\si i}$.
}

\bigskip

\begin{pf}
Let $c\in M$. 
We may choose dual OPEs of $a, b, c$ and $b, a, c$ with 
numerators $d^j(\vz,\vw)$ and $e^k(\vz,\vw)$ and
dual pole orders $\vh_{ac}^j$ and $\vh_{bc}^k$ in $\K^2$, resp.,
such that 
\begin{equation}
\label{E:duality coeff ab}
d^j(\vz,\vw)
\; =\;
\sum_{i,k}\:
\vz^{-\vh_{ab}^i}\,\vw^{-\vh_{bc}^k}\:
p_{ijk}^{ab}(\vz,\vw)
\end{equation}
for any $j$ and
\begin{equation}
\label{E:duality coeff ba}
e^k(\vz,\vw)
\; =\;
\sum_{i,j}\:
\vz^{-\vh_{ba}^i}\,\vw^{-\vh_{ac}^j}\:
p_{ijk}^{ba}(\vz,\vw)
\end{equation}
for any $k$ where 
$p_{ijk}^{ab}(\vz,\vw)$ and $p_{ijk}^{ba}(\vz,\vw)$ are 
$M$-valued power series.

Inserting \eqref{E:duality coeff ab}
into the dual OPE of $a, b, c$ for $w>z-w$ we get
$$
(a(\vz)b)(\vw)c
\; =\;
\sum_{i,j,k}\:
(\vw+\vz)^{-\vh_{ac}^j}\vz^{-\vh_{ab}^i}\vw^{-\vh_{bc}^k}
p_{ijk}^{ab}(\vz,\vw).
$$
On the other hand,
inserting \eqref{E:duality coeff ba}
into the dual OPE of $b, a, c$ for $w>z-w$ 
and using skew-symmetry and the assumption that 
$\vT$ is a translation endomorphism we obtain
\begin{align}
\notag
(a(\vz)b)(\vw)c
\; =\;
&
\paraab\,
(e^{\vz\vT}b(-\vz)a)(\vw)c
\\
\notag
\; =\;
&
\paraab\,
e^{\vz\del_{\vw}}(b(-\vz)a)(\vw)c
\\
\notag
\; =\;
&
\paraab\,
e^{\vz\del_{\vw}}
\sum_{i,j,k}\:
(\vw-\vz)^{-\vh_{bc}^k}(-\vz)^{-\vh_{ba}^i}\vw^{-\vh_{ac}^j}
p_{ijk}^{ba}(-\vz,\vw)
\\
\notag
\; =\;
&
\paraab
\sum_{i,j,k}\:
\vw^{-\vh_{bc}^k}(-\vz)^{-\vh_{ba}^i}(\vw+\vz)^{-\vh_{ac}^j}
p_{ijk}^{ba}(-\vz,\vw+\vz).
\end{align}
We may assume that
$\vh_{ac}^j\notin\vh_{ac}^k+\Z^2$ for any $j\ne k$.
Proposition \ref{SS:ope finite} implies that
if $p_{ijk}^{ab}$ is non-zero then there exist unique
$\si i$ and $\vn_{ijk}\in\Z^2$ such that
$\vh_{ab}^i=\vh_{ba}^{\si i}+\vn_{ijk}$ and
\begin{equation}
\label{E:pijkab}
\vz^{-\vn_{ijk}}p_{ijk}^{ab}(\vz,\vw)
\; =\;
\paraab\,
(-1)^{\vh_{ba}^{\si i}}p_{\si ijk}^{ba}(-\vz,\vw+\vz).
\end{equation}
Moreover, $\si i$ does not depend on $j$ or $k$.

Inserting \eqref{E:duality coeff ab}
into the dual OPE of $a, b, c$ for $z-w>w$ we get
\begin{align}
\notag
a(\vz)b(\vw)c
\; =\;
&e^{-\vw\del_{\vz}}
a(\vz+\vw)b(\vw)c
\\
\notag
\; =\;
&
e^{-\vw\del_{\vz}}
\sum_{i,j,k}\:
(\vz+\vw)^{-\vh_{ac}^j}\vz^{-\vh_{ab}^i}\vw^{-\vh_{bc}^k}\,
p_{ijk}^{ab}(\vz,\vw)
\\
\notag
\; =\;
&
\sum_{i,j,k}\:
\vz^{-\vh_{ac}^j}(\vz-\vw)^{-\vh_{ab}^i}\vw^{-\vh_{bc}^k}\,
p_{ijk}^{ab}(\vz-\vw,\vw).
\end{align}
On the other hand,
inserting \eqref{E:duality coeff ba}
into the dual OPE of $b, a, c$ for $z-w>w$ we obtain
$$
b(\vw)a(\vz)c
\; =\;
\sum_{i,j,k}\:
\vw^{-\vh_{bc}^k}(\vw-\vz)^{-\vh_{ba}^i}\vz^{-\vh_{ac}^j}\,
p_{ijk}^{ba}(\vw-\vz,\vz).
$$
Because $p_{ijk}^{ab}$ and $p_{ijk}^{ba}$ are power series 
\eqref{E:pijkab} implies
\begin{align}
\notag
(\vz-\vw)^{-\vn_{ijk}}p_{ijk}^{ab}(\vz-\vw,\vw)
\; =\;
&\paraab\,
(-1)^{\vh_{ba}^{\si i}}p_{\si ijk}^{ba}(-\vz+\vw,\vw+\vz-\vw)
\\
\notag
\; =\;
&\paraab\,
(-1)^{\vh_{ba}^{\si i}}p_{\si ijk}^{ba}(\vw-\vz,\vz).
\end{align}
Thus locality holds.
Reversing the argument, we see that if the above OPE of $a(\vz)b(\vw)$
were {\it not} reduced then $\vh_{ab}^1,\dots, \vh_{ab}^r$
would not be the pole orders of the dual OPE of $a$ and $b$.
Because the statement of the Proposition is symmetric in $a$ and $b$
the numbers $\vh_{ba}^1,\dots, \vh_{ba}^s$ are also
the pole orders of $a(\vz)b(\vw)$ and hence must be equal to
$\vh_{ab}^1,\dots, \vh_{ab}^r$.
\end{pf}

\section{OPE-Algebras}
\label{S:ope algebras}

{\bf Summary.}\:
In section \ref{SS:ko vertex}
we define the notion of an OPE-algebra
and prove that for multiply local OPE-algebras 
the state-field correspondence is an OPE-algebra isomorphism.
In section \ref{SS:chiral alg}
we prove that the subspaces of 
holomorphic and anti-holomorphic states 
of an OPE-algebra are vertex subalgebras.
In section \ref{SS:space fields}
we show that the notion of an OPE-algebra
can be equivalently defined in terms of 
a space of fields
instead of a state-field correspondence.
In section \ref{SS:exist thm}
we prove the existence theorem for OPE-algebras.
In section \ref{SS:crossing algs}
we define the notion of a crossing algebra
and give an equivalent reformulation of it.
In section \ref{SS:cft type}
we discuss the notion of an OPE-algebra of CFT-type.

\subsection{OPE-Algebras}
\label{SS:ko vertex}

We define the notion of an {\it OPE-algebra}
and prove that for multiply local OPE-algebras  
the state-field correspondence is an OPE-algebra isomorphism.

\bigskip

{\bf Remark.}\: {\it
Let $V$ be a $\K^2$-fold algebra and $1$ be a state.
If $a$ and $b$ are creative, mutually local states
then $a$ is bounded on $b$ by the pole orders of $a(\vz)b(\vw)$.
}

\bigskip

\begin{pf}
This is a special case of Proposition \ref{SS:vacuum skew}\,$\itii$.
\end{pf}

\bigskip

{\bf Definition.}\:
A local $\K^2$-fold algebra is called an {\bf OPE-algebra} \index{OPE-algebra}
if there exist a translation generator and an invariant right identity.

\bigskip

{\bf Proposition.}\: {\it
Let $V$ be an OPE-algebra with translation generator $\vT$ and
invariant right identity $1$.

\begin{enumerate}
\item[$\iti$]
The $\K^2$-fold algebra $V$ is a bounded $\vbbK$-fold algebra.

\item[$\itii$]
If $a(\vz)$ is a distribution that is strongly creative and 
local to $\cF_Y$ then $a(\vz)=Y(s_1(a(\vz)),\vz)$.

\item[$\itiii$]
The state $1$ is a full identity.

\item[$\itiv$]
The $\K^2$-fold algebra $V$ satisfies skew-symmetry for $\vT$.

\item[$\itv$]
If $a$ and $b$ are states and $\vn\in\K^2$ such that
$a(\vz)_{(\vn)}b(\vz)$ is local to $\cF_Y$ then
$a(\vz)_{(\vn)}b(\vz)=(a_{(\vn)}b)(\vz)$.
In particular, if $V$ is multiply local then $\cF_Y$ is 
a unital $\K^2$-fold subalgebra of $\cF(V)$ and 
the state-field correspondence $Y:V\to\cF_Y$
is an isomorphism of OPE-algebras.
\end{enumerate}
}

\bigskip

By $\itiii$ OPE-algebras are unital $\K^2$-fold algebras.
A morphism $V\to W$ of OPE-algebras \index{morphism of OPE-algebras}
is by definition a morphism $V\to W$ of unital $\K^2$-fold algebras.

\bigskip

\begin{pf}
$\iti$\:
This follows from the Remark because, by definition, the pole orders
of mutually local distributions are contained in $\vbbK$.

$\itii$\:
Proposition \ref{SS:creat}\,$\itii$ implies that 
$\cF_Y$ is strongly creative.
Moreover, we have $s_1(\cF_Y)=V$.
Thus the claim follows from Goddard's uniqueness theorem.

$\itiii$\:
Proposition \ref{SS:creat}\,$\itii$ implies that 
$1$ is a strong right identity.

The identity field $1(z)$ is strongly creative and mutually local 
to $\cF_Y$.
By $\itii$ we get $1(z)=Y(s_1(1(z)),z)=Y(1,z)$.
Thus $1$ is a left identity.

For a state $a$ of $V$, the distribution $\del_z a(\vz)$ is strongly creative
and mutually local to $\cF_Y$.
By $\itii$ we get $\del_z a(\vz)=Y(s_1(\del_z a(\vz)),\vz)=Y(T a,\vz)$.
Similarly, we have $\del_{\bz} a(\vz)=Y(\bT a,\vz)$.
Thus $\vT$ is a translation operator and $1$ is a full identity.

$\itiv$\:
This follows from Proposition \ref{SS:skew symm states}.

$\itv$\:
Proposition \ref{SS:vacuum skew}\,$\itiv$ implies that
$a(\vz)_{(\vn)}b(\vz)$ is strongly creative.
By $\itii$ and Proposition \ref{SS:vacuum skew}\,$\itiii$
we get 
$a(\vz)_{(\vn)}b(\vz)=
Y(s_1(a(\vz)_{(\vn)}b(\vz)),\vz)=
Y(a_{(\vn)}b,\vz)$.
The second claim follows from Dong's lemma for multiple locality.
${}_{}$
\end{pf}

\bigskip

The category of unital $\K^2$-fold algebras is complete.
We denote by $\OPEA$ the subcategory of 
the category of unital $\K^2$-fold algebras consisting of OPE-algebras.
The category $\OPEA$ contains finite limits. 
Indeed, if $F$ is a functor from a category with finitely many objects
to $\OPEA$ then the limit unital $\K^2$-fold algebra of $F$ is an OPE-algebra.

Let $V$ be an OPE-algebra.
A vector subspace $I$ of $V$ is called an {\bf ideal} \index{ideal}
if $V_{(\K^2)}I$ is contained in $I$.
Let $I$ be an ideal of $V$.
Skew-symmetry implies that $I_{(\K^2)}V$ is also contained in $I$.
Thus there exists an induced $\K^2$-fold algebra structure on $V/I$.
The full identity $1$ of $V$ induces a full identity of $V/I$.
If $a, b\in V$ and $c\in I$ then 
$a(\vz)b(\vw)c=0\in V/I$ and hence 
$N_{\vh}(a(\vz)b(\vw))c=0\in V/I$ for any $\vh\in\K^2$ by
Proposition \ref{SS:ope finite}.
Thus we see that
$a+I$ and $b+I$ are mutually local states of $V/I$
and $V/I$ is an OPE-algebra.

If $V$ and $V'$ are OPE-algebras then the vector space $V\otimes V'$
\index{tensor product of OPE-algebras}
with state-field correspondence $a\otimes a'\mapsto a(\vz)\otimes a'(\vz)$
is an OPE-algebra.
Indeed, if $1$ and $1'$ are the identities of $V$ and $V'$, resp., then
$\vT_1\otimes\id_{V'}+\id_V\otimes\vT_{1'}$ is a translation generator
of $V\otimes V'$ by Proposition \ref{SS:tensor prod fields}\,$\itiii$,
the state $1\otimes 1'$ is obviously an invariant right identity, and 
$V\otimes V'$ is local by Proposition \ref{SS:tensor prod fields}\,$\itii$.

\subsection{Vertex Algebras and Chiral Algebras}
\label{SS:chiral alg}

We define the notion of a {\it vertex algebra}
and prove that the subspaces of holomorphic and anti-holomorphic states 
of an OPE-algebra are vertex subalgebras.
These subalgebras are called the {\it chiral} and the anti-chiral algebra.

\bigskip

{\bf Definition.}\:
A holomorphic OPE-algebra is called a {\bf vertex algebra}.
\index{vertex algebra}

\bigskip

The category $\VertA$ of vertex algebras 
is a subcategory of $\OPEA$ that is closed with respect to
finite limits and the tensor product.

For a $\K^2$-fold algebra $V$, we denote by
$V_z$ and $V_{\bz}$ the vector spaces of 
holomorphic and anti-holomorphic states of $V$.
By definition, 
an OPE-algebra $V$ is a vertex algebra if and only if $V=V_z$.

Let $S$ be an even set and $M$ be an $S$-fold module over a vector space.
For subset $T$ of $S$, a $T$-fold submodule $N$ of $M|_T$ is called 
a {\bf $T$-fold submodule} \index{Tfold submodule@$T$-fold submodule} of $M$ 
if $N|^S$ is an $S$-fold submodule of $M$.

A holomorphic state $a$ and a state $b$ of a $\K^2$-fold algebra 
are said to {\bf commute} \index{commuting states}
if $a(z)$ and $b(\vz)$ are holomorphically local of order zero.
Of course, this is equivalent to the condition that the modes
$a_{(n)}$ and $b_{(\vm)}$ commute for any $n\in\Z$ and $\vm\in\K^2$.

\bigskip

{\bf Proposition.}\: {\it
Let $V$ be a $\vbbK$-fold algebra. 

\begin{enumerate}
\item[$\iti$]
If $Y$ is injective and 
$\vT$ is a translation endomorphism of $V$ 
then $V_z=\ker(\bT)$ and $V_{\bz}=\ker(T)$.

\item[$\itii$]
If $V$ has a full identity $1$
then $V_z$ and $V_{\bz}$ 
are unital $\Z$-fold subalgebras;
the $\Z$-fold subalgebra structure on $V_{\bz}$ 
is induced by the injection $\Z\to\K^2, n\mapsto (-1,n)$.

\item[$\itiii$]
A holomorphic state $a$ and an anti-holomorphic state $b$
are holomorphically local if and only if $a$ and $b$ commute.

\item[$\itiv$]
If $V$ is an OPE-algebra
then $V_z$ and $V_{\bz}$ are commuting vertex subalgebras.
\end{enumerate}
}

\bigskip

\begin{pf}
$\iti$\:
Because  
$\del_{\bz}Y(a,\vz)=Y(\bT a,\vz)$ we have
$Y(a,\vz)\in\End(V)\set{z}$ if and only if
$\bT a=0$.
Since $Y(a,\vz)$ is statistical
$Y(a,\vz)\in\End(V)\set{z}$ 
is equivalent to 
$Y(a,\vz)\in\End(V)\pau{z\uppm}$.

$\itii$\:
Because $\bT_1$ is a derivation of the 
$\vbbK$-fold algebra $V$ 
the subspace $V_z=\ker(\bT_1)$ is a
$\vbbK$-fold subalgebra.
From $Y(V_z)\subset\End(V)\pau{z\uppm}$
follows that 
$V_z$ is a $\Z$-fold subalgebra.
Since the identity field is holomorphic
the identity $1$ is contained in $V_z$.

$\itiii$\:
Corollary \ref{SS:locality} shows that
the distribution $[a(z),b(\bw)]$ is local 
if and only if it is zero because $[a(z),b(\bw)]$ is bounded in $w$.

$\itiv$\:
This follows from $\itii$ and $\itiii$.
\end{pf}

\bigskip

{\bf Definition.}\: 
For a $\vbbK$-fold algebra $V$ with full identity,
the $\Z$-fold subal\-gebras $V_z$ and $V_{\bz}$
are called the {\bf chiral} \index{chiral algebra}
and the {\bf anti-chiral algebra} \index{anti-chiral algebra} of $V$.

\bigskip

Part $\iti$ of the Proposition shows that 
if a $\vbbK$-fold algebra has a full identity 
then holomorphic states satisfy skew-symmetry
if and only if they satisfy holomorphic skew-symmetry.

\bigskip

{\bf Remark.}\: {\it
Let $V$ be a $\K^2$-fold algebra, 
$a$ be a holomorphic state, and $1$ and $b$ be states such that
$a$ and $b$ are creative and mutually local of order $h$.
Then $h\geq o'(a,b)$. 
Moreover, if $Y:V\to\End(V)\set{\vz}$ is a morphism of $\K^2$-fold algebras
then $h=o'(a,b)$.
}

\bigskip

\begin{pf}
The first claim follows from Remark \ref{SS:ko vertex}.
Let us prove the second claim.
We may assume that $h:=o(a,b)>-\infty$.
Remark \ref{SS:ope finite} shows that the field
$c(z,\vw):=(z-w)^h a(z)b(\vw)$ is not divisible as a field by $z-w$.
Thus $Y(a_{(h-1)}b,\vz)=c(z,\vz)\ne 0$ and hence $o'(a,b)\geq h$.
\end{pf}

\subsection{Space of Fields}
\label{SS:space fields}

We show that the notion of an OPE-algebra
can be equivalently defined in terms of
a {\it space of fields}
instead of a state-field correspondence.

\bigskip

{\bf Proposition.}\: {\it
If $V$ is an OPE-algebra with identity $1$  
then the vector space $V$
together with the even vector $1$ and the vector subspace
$\cF:=\cF_Y$ of $\End(V)\set{\vz}$ satisfy the following properties:

\begin{enumerate}
\item[$\iti$]
there exists a pair $\vT$ of even endomorphisms of $V$
such that $1$ is invariant and
$\cF$ is creative, translation covariant, complete, and local;
\end{enumerate}

Moreover, if the OPE-algebra $V$ is multiply local then 
the vector space $V$ together with $1$ and $\cF$ satisfy 
the following properties:

\begin{enumerate}
\item[$\itii$]
$\cF$ is a $\K^2$-fold subalgebra of $\End(V)\set{\vz}$
that is creative and local and 
the map $s_1:\cF\to V$ is a vector space isomorphism.
\end{enumerate}

Conversely, if $V$ is a vector space together with
an even vector $1$ and 
a vector subspace $\cF$ of $\End(V)\set{\vz}$
such that either $\iti$ or $\itii$ are satisfied
then 
$s_1:\cF\to V$ is a vector space isomorphism and
its inverse is the unique OPE-algebra structure
on $V$ such that $1$ is a right identity and 
$\cF$ is the space of fields. 
Moreover, if $\itii$ is satisfied then $Y$ is an isomorphism 
of OPE-algebras. 
}

\bigskip

\begin{pf}
``$\Rightarrow$"\:
This follows from Proposition \ref{SS:ko vertex}.

``$\Leftarrow$"\,$\iti$\:
Goddard's uniqueness theorem shows that $s_1|_{\cF}$ is invertible.
It is clear that
$(V,(s_1|_{\cF})^{-1})$ is an OPE-algebra 
with right identity $1$ and space of fields $\cF$.
Uniqueness follows from Remark \ref{SS:creat}.

$\itii$\:
Let $Y:V\to\cF$ denote the inverse of $s_1:\cF\to V$.
Proposition \ref{SS:vacuum skew}\,$\itiii$ 
and the fact that $\cF$ is a $\K^2$-fold subalgebra imply that
$s_1:\cF\to (V,Y)$ is an isomorphism of $\K^2$-fold algebras.
Together with creativity of $\cF$ this shows that 
$Y(1,\vz)$ is a right identity of $\cF$.

Proposition \ref{SS:id field derivatives}\,$\itii$ says
that $\del_{\vz}$ is a translation operator for OPE-finite distributions. 
Because $\cF$ is local we get
$\del_z a(\vz)=\del_z a(\vz)_{(-1)}Y(1,\vz)=a(\vz)_{(-2,-1)}Y(1,\vz)$
for any $a(\vz)\in\cF$.
In the same way we get $\del_{\bz}a(\vz)=a(\vz)_{(-1,-2)}Y(1,\vz)$.
This shows that $\cF$ is invariant with respect to $\del_{\vz}$.
From Remark \ref{SS:sfold id} follows that
$Y(1,\vz)$ is a strong right identity for $\del_{\vz}$.

Proposition \ref{SS:vacuum skew} and locality of $\cF$ imply that
any pair of distributions in $\cF$ satisfies skew-symmetry for $\del_{\vz}$.
Thus Remark \ref{SS:skew symm} shows that
$Y(1,\vz)$ is a left identity of $\cF$.
Because $s_1$ is an isomorphism of $\K^2$-fold algebras
the state $s_1(Y(1,\vz))=1$ is a left identity of $V$, 
i.e.~$Y(1,\vz)$ is the identity field.

From $s_1(\cF)=V$ and Proposition \ref{SS:vacuum skew}\,$\itii$
follows that the distributions in $\cF$ are fields.
Proposition \ref{SS:id field derivatives}\,$\iti$ and $\itii$
imply that the identity field is a full identity of $\cF$.
Hence $s_1(1(z))=1$ is a full identity of $V$.
This shows that $(V,Y)$ is an OPE-algebra.
Uniqueness follows from Remark \ref{SS:creat}. 
\end{pf}

\subsection{Existence Theorem}
\label{SS:exist thm}

We prove the {\it existence theorem} for OPE-algebras which
states that a vector space together with a complete
set of creative, translation covariant, and 
multiply local fields 
carries a unique OPE-algebra structure.

\bigskip

For a vector space $V$ and a vector $1$ of $V$,
a  subset $S$ of $\End(V)\set{\vz}$ is called
{\bf generating} \index{generating set of distributions}
if
$$
V
\; =\;
\rspan\set{ a^1_{\vn_1}\dots a^r_{\vn_r}1\mid
a^i(\vz)\in S,\, \vn_i\in\K^2,\, r\in\N }.
$$

\bigskip

For a $\K^2$-fold algebra $V$ and a subset $S$ of $V$,
we denote by $\sqbrack{S}$ the $\K^2$-fold subalgebra of $V$
that is generated by $S$.

\bigskip

{\bf Existence Theorem.} \: {\it
Let $V$ be a vector space and 
$S$ be a subset of $\End(V)\set{\vz}$
such that there exist an even vector $1$ and
a pair $\vT$ of even endomorphisms of $V$ such that 
$1$ is invariant and
$S$ is creative, translation covariant, generating, and multiply local.
Then there exists a unique OPE-algebra structure $Y$
on $V$ such that $1$ is a right identity and
$Y(a_{-1}(1),\vz)=a(\vz)$ for any $a(\vz)\in S$.
Moreover, we have $\cF_Y=\sqbrack{S\cup\set{1(z)}}$
and $Y$ is an isomorphism of OPE-algebras.
In particular, if $S$ consists of holomorphic distributions
then $(V,Y)$ is a vertex algebra.
}

\bigskip

\begin{pf}
Define $\cF:=\sqbrack{S\cup\set{1(z)}}$.
Proposition \ref{SS:vacuum skew}\,$\iti$,
Proposition \ref{SS:id field derivatives}\,$\itiii$,
and Dong's lemma for multiple locality 
show that $\cF$ is creative, translation covariant, and multiply local.

Proposition \ref{SS:vacuum skew}\,$\itiii$ implies that 
for any $a(\vz), b(\vz)\in\cF$ and $\vn\in\K^2$
we have
$s_1(a(\vz)_{(\vn)}b(\vz))=a_{\vn}s_1(b(\vz))$.
In particular, we have
$s_1(a(\vz)_{(\vn)}1(z))=a_{\vn}(1)$
for any $a(\vz)\in\cF$ and $\vn\in\K^2$.
Because $S$ is complete we obtain that
$s_1:\cF\to V$ is surjective.
Thus the claim follows from Proposition \ref{SS:space fields}
and Proposition \ref{SS:ko vertex}\,$\itv$.
\end{pf}

\subsection{Crossing Algebras}
\label{SS:crossing algs}

We define the notion of a {\it crossing algebra}
and give an equivalent reformulation of it.

\bigskip

{\bf Definition.}\: 
An OPE-algebra $V$ is called a {\bf crossing algebra} if $V$ is dual.

\bigskip

{\bf Remark.}\: {\it
A bounded $\vbbK$-fold algebra $V$ is 
a crossing algebra if and only if
there exist a right identity and a translation operator $\vT$
and $V$ is dual and satisfies skew-symmetry for $\vT$.
}

\bigskip

\begin{pf}
``$\Rightarrow$"\:
This follows from the definition of a
crossing algebra and Proposition \ref{SS:skew symm states}.

``$\Leftarrow$"\:
From Proposition \ref{SS:duality} follows that $V$ is local.
By Remark \ref{SS:sfold id} a right identity of $V$ 
is a strong right identity for $\vT$.
By Remark \ref{SS:skew symm states}
a strong right identity $1$ is a left identity
and hence $\vT_1(1)=0$.
\end{pf}

\subsection{OPE-Algebras of CFT-Type}
\label{SS:cft type}

We define the notions of a conformal OPE-algebra and of an 
{\it OPE-algebra of CFT-type}.

\bigskip

{\bf Definition.}\:
For $\vc\in\K^2$, an OPE-algebra together with a pair of conformal 
vectors of central charges $\vc$  
is called a {\bf conformal} OPE-algebra \index{conformal OPE-algebra}
of central charges $\vc$.

\bigskip

{\bf Remark.}\: {\it
Let $V$ be an OPE-algebra and $L$ and $\bL$ be a holomorphic and
an anti-holomorphic even state. 

\begin{enumerate}
\item[$\iti$]
If $L_{(1)}$ is a holomorphic dilatation operator of $V$ then
any anti-holomorphic state is homogeneous of weight zero for $L_{(1)}$.

\item[$\itii$]
Assume that $\cF_Y$ is translation covariant for $(L_{(0)},\bL_{(0)})$
and dilatation covariant for $(L_{(1)},\bL_{(1)})$.
If $L_{(3)}L, \bL_{(3)}\bL\in\K 1$ and 
$L(z)$ and $\bL(\bz)$ are both local to themselves of order at most four and 
dilatation covariant of weights $\vh_1$ and $\vh_2$ then
$\vL$ is a pair of conformal vectors,
the pair $(L_{(0)},\bL_{(0)})$ is the translation operator of $V$,
the pair $(L_{(1)},\bL_{(1)})$ is a dilatation operator,
and $\vh_1=(2,0)$ and $\vh_2=(0,2)$.
\end{enumerate}
}

\bigskip

\begin{pf}
$\iti$\:
Let $a$ be an anti-holomorphic state 
that is homogeneous of weight $h$ for $L_{(1)}$.
Dilatation covariance yields 
$[L_{(1)},a(\bz)]=h\, a(\bz)+z\del_z a(\bz)=h\, a(\bz)$.
By Remark \ref{SS:chiral alg}\,$\itiii$
the states $L$ and $a$ commute. Thus $h=0$.

$\itii$\:
The identity $1$ of $V$ is invariant for $L_{(0)}$ and $L_{(1)}$
because $L(z)$ is creative for $1$.
Remark \ref{SS:sfold id} shows that $L_{(0)}$ is equal to the
holomorphic translation operator $T_1$ of $V$.
If $a$ is a state such that $a(\vz)$ is dilatation covariant for
$L_{(1)}$ of weight $h$ then setting $\vz=0$ in the equation
$$
[L_{(1)},a(\vz)]1
\; =\;
h\, a(\vz)1
\; +\;
z\del_z a(\vz)1
$$
yields $L_{(1)}a=h\, a$. By assumption, $\cF_Y$ is a graded subspace
of $\End(V)\set{\vz}_{\vast}$. Since $Y$ is a vector space isomorphism the 
operator $L_{(1)}$ of $V$ is diagonalizable.
Remark \ref{SS:sfold dilat op} implies that $L_{(1)}$ is a 
holomorphic dilatation operator. Since $L_{(1)}$ and $\bL_{(1)}$ commute 
the pair $(L_{(1)},\bL_{(1)})$ is a dilatation operator.
Because $L_{(0)}=T_1$ we have 
$\vh_1-(0,-1)-1=(1,0)$ and hence $\vh_1=(2,0)$.
Remarks \ref{SS:vir ope} and \ref{SS:conf distr} show
that $L(z)$ is a Virasoro distribution.
\end{pf}

\bigskip

For an even set $S$, an $S$-graded vector space is called
{\bf locally finite} \index{locally finite gradation}
if its homogeneous subspaces are finite-dimensional.

Let $V$ be a $\K^2$-graded vector space and $\vH$ be the corresponding
pair of even commuting diagonalizable operators. 
The {\bf spectrum} \index{spectrum of a graded vector space} of $V$ 
is by definition the spectrum of $\vH$, i.e.~the set of
elements $\vh$ of $\K^2$ such that $V_{\vh}$ is non-zero.
The $\K^2$-graded vector space $V$ is called {\bf statistical} 
\index{statistical graded vector space}
if $V\even=\bigoplus_{h-\bh\in\Z}V_{\vh}$ and
$V\odd=\bigoplus_{h-\bh\in 1/2+\Z}V_{\vh}$.

\bigskip

{\bf Definition.}\:
Let $\K=\C$. A conformal OPE-algebra $V$ is called an
OPE-algebra {\bf of CFT-type} \index{OPE-algebra!of CFT-type}
if its gradation is locally finite
and statistical, the spectrum is contained in $\R_{\geq}^2$, and $V_0=\K 1$.

\section{Additive OPE-Algebras}
\label{S:addi opea}

{\bf Summary.}\:
In sections \ref{SS:ope alg} and \ref{SS:addi dual and skew}
we define the notion of an additive OPE-algebra
and prove that this notion can be equivalently defined
in terms of additive duality and skew-symmetry.
In section \ref{SS:jacobi modes}
we define the Jacobi identity in terms of modes and
consider nine special cases of it.
In section \ref{SS:jacobi pretex}
we rewrite the Jacobi identity and two of its special cases
in terms of distributions.
In section \ref{SS:dual local pretex}
we prove that the Jacobi identity is equivalent to 
additive duality and additive locality.

\subsection{Additive OPE-Algebras}
\label{SS:ope alg}

We define the notion of an {\it additive OPE-algebra}
and prove that additive duality implies that
any strong identity is a full identity.

\bigskip

{\bf Definition.}\:
A vector space gradation of a bounded $\K^2$-fold algebra 
is called {\bf local} \index{local gradation}
if any two homogeneous states are additively local.
An OPE-algebra is called {\bf additive} \index{additive OPE-algebra}
if there exists a local gradation of it.

\bigskip

If a bounded $\K^2$-fold algebra $V$ possesses a local gradation
then $V$ is local
because additive locality implies locality 
and the relation of locality is bilinear.
Thus an additive OPE-algebra is a 
bounded $\K^2$-fold algebra such that 
there exists a translation generator, an invariant right identity,
and a local gradation.

Because the identity field is additively local
to any field
there exists for any additive OPE-algebra $V$
a local gradation such that
the identity of $V$ is homogeneous.

\bigskip

{\bf Lemma.}\: {\it
Let $V$ be a bounded $\vbbK$-fold algebra with a vector space gradation,
$1$ be a strong right identity, $\vT:=\vT_1$, and $a$ and $c$ be states.

\begin{enumerate}
\item[$\iti$]
If for any homogeneous state $b$
the states $a, b, 1$ satisfy additive duality
then $a(\vz)$ is translation covariant.

\item[$\itii$]
If $1$ is a strong identity and 
the states $a, 1, c$ satisfy additive duality of order at most $\vt$
then $a$ is bounded on $c$ by $\vt$ and
$\vT$ is a translation endomorphism for $a$ and $c$.
\end{enumerate}
}

\bigskip

\begin{pf}
$\iti$\:
Acting with $e^{-\bw\bT}e^{-wT}$ on both sides of 
the additive duality identity
$$
(\vx+\vw)^{\vt}\, a(\vx+\vw)b(\vw)1
\; =\;
(\vw+\vx)^{\vt}\, (a(\vx)b)(\vw)1
$$
and using strong creativity we obtain
\begin{equation}
\label{E:transl gen pf}
(\vx+\vw)^{\vt}\, e^{-\bw\bT}e^{-wT}\, a(\vx+\vw)e^{\vw\vT}b
\; =\;
(\vw+\vx)^{\vt}\, a(\vx)b.
\end{equation}

Assume that $a(\vz)b$ is non-zero. 
The left-hand side of \eqref{E:transl gen pf}
is contained in $V\sqbrack{\vx}\pau{\vw}$
and the right-hand side
is non-zero and contained in
$\vw^{\vt}V\sqbrack{\vx}\pau{\vw^{-1}}$.
This implies $\vt\in\N^2$.
Since $e^{-\bw\bT}e^{-wT}\, a(\vx+\vw)e^{\vw\vT}b$
and $a(\vx)b$ are both contained in 
$V\sqbrack{\vx}\pau{\vw}$
we may cancel
$(\vx+\vw)^{\vt}=(\vw+\vx)^{\vt}$ on both sides of 
\eqref{E:transl gen pf}. 

If $a(\vz)b$ is zero then additive duality shows that
$a(\vx+\vw)b(\vw)1=0$ and
hence additive duality of order $\vt=0$ is satisfied.
Thus for any homogeneous $b$ we have
$e^{-\bw\bT}e^{-wT}a(\vx+\vw)e^{\vw\vT}b=a(\vx)b$.
The proof of Proposition \ref{SS:transl cov} shows that this implies 
that $a(\vz)$ is translation covariant.

$\itii$\:
Additive duality 
$(\vx+\vw)^{\vt} a(\vx+\vw)1(\vw)c=(\vw+\vx)^{\vt}(a(\vx)1)(\vw)c$
reduces to
\begin{equation}
\label{E:transl endo pf}
(\vx+\vw)^{\vt}\, a(\vx+\vw)c
\; =\;
(\vw+\vx)^{\vt}\, (e^{\vx\vT}a)(\vw)c.
\end{equation}
The right-hand side is contained in $V\sqbrack{\vw}\pau{\vx}$.
The left-hand side is contained in $V\sqbrack{\vx}\pau{\vw}$
and its constant term is $\vx^{\vt}a(\vx)c$.
Thus $\vx^{\vt}a(\vx)c\in V\pau{\vx}$ and
$(\vx+\vw)^{\vt}\, a(\vx+\vw)c=(\vw+\vx)^{\vt}\, a(\vw+\vx)c$.
From \eqref{E:transl endo pf} we get
$a(\vw+\vx)c=(e^{\vx\vT}a)(\vw)c$.
Taking the coefficients in front of $x$ and $\bx$ of this last equation
the claim follows.
\end{pf}

\subsection{Additive Duality and Skew-Symmetry}
\label{SS:addi dual and skew}

We prove that additive OPE-algebras can be equivalently
defined in terms of additive duality and skew-symmetry.

\bigskip

{\bf Proposition.}\: {\it
Let $V$ be a bounded $\vbbK$-fold algebra,
$\vT$ be a pair of even commuting operators, 
and $a, b, c$ be states such that 
$a(\vz)$ is translation covariant
and the pairs $b, c$ and $a_{(\vn)}b, c$ satisfy 
skew-symmetry for any $\vn\in\vbbK$.
The states $a, b, c$ satisfy additive duality of order $\vt$
if and only if $a, c, b$ satisfy additive locality of order $\vt$.
}

\bigskip

\begin{pf}
Skew-symmetry and translation covariance imply that
\begin{align}
\notag
(\vx+\vw)^{\vt}\, a(\vx+\vw)b(\vw)c
\; &=\;
\zeta^{\tb\tc}\, (\vx+\vw)^{\vt}\, a(\vx+\vw)e^{\vw\vT}c(-\vw)b
\\
\notag
&=\;
\zeta^{\tb\tc}\, (\vx+\vw)^{\vt}\, e^{\vw\vT}a(\vx)c(-\vw)b
\end{align}
and
$$
(\vw+\vx)^{\vt}\, (a(\vx)b)(\vw)c
\; =\;
\zeta^{(\ta+\tb)\tc}\, (\vw+\vx)^{\vt}\, e^{\vw\vT}c(-\vw)a(\vx)b.
$$
The left-hand sides of these equations
are equal
if and only if
$a, b, c$ satisfy additive duality of order $\vt$.
The right-hand sides of these equations
are equal
if and only if
$a, c, b$ satisfy additive locality of order $\vt$.
\end{pf}

\bigskip

{\bf Corollary.}\: {\it
For a bounded $\K^2$-fold algebra $V$
together with a gradation of the vector space $V$,
the following statements are equivalent:

\begin{enumerate}
\item[\iti]
$V$ is an additive OPE-algebra, the gradation is local, and 
the identity is homogeneous;

\item[\itii]
there exist a translation generator and 
a homogeneous invariant right identity 
and the gradation is local;

\item[\itiii]
$V$ is a $\vbbK$-fold algebra,
there exists a homogeneous strong right identity $1$,
homogeneous states are additively dual, and
skew-symmetry is satisfied for $\vT_1$.
\end{enumerate}
\hfill $\square$
}

\bigskip

\begin{pf}
It is clear that $\iti$ and $\itii$ are equivalent.
The Proposition shows that $\iti$ implies $\itiii$.

Suppose that $\itiii$ is satisfied.
By Remark \ref{SS:skew symm} the state $1$ is a left identity. 
By Lemma \ref{SS:ope alg}\,$\itii$ 
for any homogeneous states $a$ and $b$ there exists $\vh\in\vbbK$
such that $a$ is bounded on $b$ by $\vh$ and
$\vT_1$ is a translation endomorphism.
The proof of Remark \ref{SS:duality} shows that
$a$ and $b$ are mutually dual of order at most $\vh$.
By Proposition \ref{SS:dual and local} this implies that
$a$ and $b$ are additively local of order at most $\vh$.
\end{pf}

\subsection{Jacobi Identity in terms of Modes}
\label{SS:jacobi modes}

We define the {\it Jacobi identity} in terms of modes,
consider two special cases of it,
and observe that the holomorphic Jacobi identity
and skew-symmetry can also be obtained 
as specializations of it.

\bigskip

Let $\K'$ be a commutative algebra with an involution $\ka$ 
that leaves $\K$ invariant and
$p:\K\times\set{\pm 1}\to\K'$ be a power map that is compatible with $\ka$.
For a vector space $V$, define $V_p:=V\otimes\K'$.

Let $a(\vz)$ and $b(\vz)$ be fields on $V$ which
are additively local of order $\vh$.
Because
$$
(z-w)^h a(\vz)b(\vw)
\; =\;
(\bz-\bw)^{-\bh}((\vz-\vw)^{\vh}a(\vz)b(\vw))
$$
is contained in 
$\cF\sqbrack{\bz}\sqbrack{z,\vw}(V)$
the product
$$
a(\vz)b(\vw)_{w>z}
\; :=\;
(z-w)_{w>z}^{-h}
(
(z-w)^h a(\vz)b(\vw)
)
\; \in \;
\cF\sqbrack{w,\bz}\sqbrack{z,\bw}(V_p)
$$
exists. Similarly, we define
$$
a(\vz)b(\vw)_{\bw >\bz}
\; :=\;
(\bz -\bw)_{\bw >\bz}^{-\bh}
(
(\bz -\bw)^{\bh} a(\vz)b(\vw)
)
\; \in \;
\cF\sqbrack{z,\bw}\sqbrack{w,\bz}(V_p).
$$
We have
\begin{align}
\notag
a(\vz)b(\vw)_{w>z}
\; &=\;
(z-w)_{w>z}^{-h}
(\bz -\bw)_{\bz >\bw}^{-\bh}
(
(\vz -\vw)^{\vh} a(\vz)b(\vw)
)
\\
\notag
&=\;
\paraab\,
(w-z)_{w>z}^{-h}
(\bw -\bz)_{\bz >\bw}^{-\bh}
(
(\vw -\vz)^{\vh} b(\vw)a(\vz)
)
\\
\label{E:ab ba zw wz}
&=\;
\paraab\,
b(\vw)a(\vz)_{\bz >\bw}
\end{align}
Interchanging $a(\vz)$ and $b(\vz)$
in \eqref{E:ab ba zw wz} yields
$a(\vz)b(\vw)_{\bw >\bz}=\paraab\, b(\vw)a(\vz)_{z>w}$.
If $a(z)$ is holomorphic
then $\bh=0$ by Proposition \ref{SS:ope commut}
and thus $a(z)b(\vw)_{\bw >\bz}=a(z)b(\vw)$.
Similarly,
we have
$a(z)b(\vw)_{w>z}=\paraab\, b(\vw)a(z)$.

Let $a_{[n,\bn}b_{m],\bm}\in\End(V_p)$ denote the modes of
$a(\vz)b(\vw)_{w>z}$ so that
$$
a(\vz)b(\vw)_{w>z}
\; =\;
\sum_{\vn,\vm\in\K^2}\:
a_{[n,\bn}b_{m],\bm}\,
\vz^{-\vn-1}\vw^{-\vm-1}
$$
and let $a_{n,[\bn}b_{m,\bm ]}$ denote the modes of 
$a(\vz)b(\vw)_{\bw >\bz}$.
We have
\begin{equation}
\label{E:a[b] modes}
a_{[n,\bn}b_{m],\bm}
\; =\;
(-1)_p^{-h}
\sum_{j\in\N}
\bigg(
\sum_{i\in\N}\:
(-1)^{i+j}\binom{h}{i}\binom{-h}{j}\,
a_{n+h-i+j,\bn}b_{m+i-h-j,\bm}
\bigg)
\end{equation}
and 
\begin{equation}
\label{E:[a]b modes}
a_{n,[\bn}b_{m,\bm ]}
\; =\;
(-1)_p^{\bh}
\sum_{j\in\N}
\bigg(
\sum_{i\in\N}\:
(-1)^{i+j}\binom{\bh}{i}\binom{-\bh}{j}\,
a_{n,\bn+\bh-i+j}b_{m,\bm+i-\bh-j}
\bigg).
\end{equation}
Equation \eqref{E:ab ba zw wz} in terms of modes reads
\begin{equation}
\label{E:a[b]=[a]b}
a_{[n,\bn}b_{m],\bm}
\; =\;
\paraab\,
b_{m,[\bm}a_{n,\bn ]}.
\end{equation}
If $a(z)$ is holomorphic then 
$a_{[n,\bn}b_{m],\bm}=\paraab\, b_{\vm}a_{\vn}$
and
$a_{n,[\bn}b_{m,\bm ]}=a_{\vn}b_{\vm}$.

\bigskip

{\bf Definition.}\:
Let $M$ be a bounded $\vbbK$-fold module over a
bounded $\vbbK$-fold algebra $V$,
$a, b\in V$, $c\in M$, and $\vr, \vs, \vt\in\vbbK$. 

\begin{enumerate}
\item[($J$)]
The vectors $a, b, c$ are said to satisfy
the {\bf Jacobi identity} 
\index{Jacobi identity!for bounded $\vbbK$-fold algebras}
for indices $\vr, \vs, \vt$
if $a$ and $b$ are additively local and we have
\begin{align}
\label{E:jacobi id pretex}
&\sum_{\vi\in\N^2}\:
\binom{\vt}{\vi}\,
(a_{(\vr+\vi)}b)_{(\vs+\vt-\vi)}c
\\
\notag
\; =\;
&\sum_{\vi\in\N^2}\:
(-1)^{\vi} \binom{\vr}{\vi}
\Big(
a_{(\vt+\vr-\vi)}b_{(\vs+\vi)}
\, -\,
(-1)_p^r\,
a_{([t+i,\bt+\br-\bi)}b_{(s+r-i],\bs+\bi)}
\\
\notag
&\;\;
-(-1)_p^{-\br}\,
a_{(t+r-i,[\bt+\bi)}b_{(s+i,\bs+\br-\bi])}
\, +\,
\paraab (-1)^{r-\br}\, 
b_{(\vs+\vr-\vi)}a_{(\vt+\vi)}
\Big)c.
\end{align}

\item[$\ita$]
{\bf Additive duality} 
of order $\vt$ for $a, b, c$ and indices $\vr, \vs$
\index{additive duality}
is the identity
\begin{equation}
\label{E:duality pretex}
\sum_{\vi\in\N^2}\:
\binom{\vt}{\vi}\,
(a_{(\vr+\vi)}b)_{(\vs+\vt-\vi)}c
\; =\;
\sum_{\vi\in\N^2}\:
(-1)^{\vi} \binom{\vr}{\vi}\,
a_{(\vt+\vr-\vi)}b_{(\vs+\vi)}c.
\end{equation}

\item[$\itb$]
{\bf Additive locality} 
of order $\vr$ for $a, b, c$ and indices $\vt, \vs$
\index{additive locality}
is the identity
\begin{equation}
\notag
\sum_{\vi\in\N^2}\:
(-1)^{\vi} \binom{\vr}{\vi}\,
\big(
a_{(\vt+\vr-\vi)}b_{(\vs+\vi)}
\; -\;
\paraab\, (-1)^{r-\br} \, b_{(\vs+\vr-\vi)}a_{(\vt+\vi)}
\big)c
\, =\,
0.
\end{equation}
\end{enumerate}

We say that the Jacobi identity of orders $\vr$ and $\vt$
holds for $a, b, c$ 
if the Jacobi identity holds for $a, b, c$ and any
indices $\vn\in\vr+\Z^2, \vs\in\vbbK$, and $\vm\in\vt+\Z^2$.
We say that $a, b, c$ satisfy the Jacobi identity 
if there exist $\vr, \vt\in\vbbK$ such that
the Jacobi identity of orders $\vr$ and $\vt$
is satisfied for $a, b, c$.
The same terminology applies to the special cases of
the Jacobi identity.

\bigskip

If $a$ and $b$ are additively local of order $\vh$
and $\vr\in\vh+\Z^2$ then the right-hand side of the 
Jacobi identity is contained
in $V$ because of \eqref{E:a[b] modes} and \eqref{E:[a]b modes}.

From \eqref{E:a[b] modes}, \eqref{E:[a]b modes}, and \eqref{E:a[b]=[a]b}
we obtain the following statements.

\bigskip

{\bf Remark.}\:
($J_z$)\:
If $a$ is holomorphic then
the Jacobi identity for $a, b, c$ and 
indices $(r,-1), \vs, \linebreak[0](t,-1)$
is identical with the holomorphic Jacobi identity for 
$a, b, c$ and indices $r, \vs-(0,1), t$.

$\ita$\:
If $a$ is bounded on $c$ by $\vt$ then
the Jacobi identity for indices $\vr, \vs, \vt$ 
is identical with additive duality of order $\vt$ for indices $\vr$ and $\vs$.

($a_z$)\:
If $a$ is holomorphic then
additive duality of order $(t,0)$ for $a, b, c$ and indices $(r,-1)$ and $\vs$
is identical with holomorphic duality of order $t$ for $a, b, c$
and indices $r$ and $\vs$.

$\itb$\:
If $a$ is bounded on $b$ by $\vr$ then
the Jacobi identity for indices $\vr, \vs, \vt$ 
is identical with additive locality of order $\vr$ for indices $\vt$ and $\vs$.

($b_z$)\:
If $a$ is holomorphic then
additive locality of order $(r,0)$ for $a, b, c$ and indices $(t,-1)$ and $\vs$
is identical with holomorphic locality of order $r$ for $a, b, c$
and indices $t$ and $\vs$.

$\ite$\:
If $1$ is a full identity then
the Jacobi identity for $a, b, 1$ and indices 
$\vr, 0, -1$ 
is identical with skew-symmetry for 
$a, b, \vT_1$, and the $\vr$-th product. 

($e_z$)\:
If $1$ is a full identity and $a$ is holomorphic then
the holomorphic Jacobi identity for $a, b, 1$ and indices 
$r, 0, -1$ 
is identical with holomorphic skew-symmetry for 
$a, b, T_1$, and the $r$-th product. 

$\iti$\:
If $1$ is a right identity then
the holomorphic Jacobi identity for $1, 1, 1$ and indices
$-1, -1, -1$ is identical with $T_1(1)=T_1(1)+T_1(1)$.

$\itii$\:
If $1$ is a right identity then
the Jacobi identity for $a, b, 1$ and indices $\vr, (-2,-1), 0$
is identical with 
$T_1(a_{(\vr)}b)=a_{(\vr)}T_1(b)\, -\, r\, a_{(r-1,\br)}b$.

\subsection{Jacobi Identity in terms of Distributions}
\label{SS:jacobi pretex}

We rewrite the Jacobi identity and its special cases
in terms of distributions.

\bigskip

For $h\in\K$,
we define the {\bf $h$-shifted} delta distribution by
\index{shifted delta distribution}
\index{delta distribution!shifted}
$$
\de_h(z,w)
\; :=\; 
z^{-h}w^h\de(z,w)
\; =\;
\sum_{n\in h+\Z}\: w^n z^{-n-1}.
$$
Thus $\de_0(z,w)=\de(z,w)$.
We have
$\de_{h+n}(z,w)=\de_h(z,w)$ for any integer $n$ and
$\de_h(z,w)=\de_{-h}(w,z)$.
From \eqref{E:delta distr symmetric},
the identity
$w\del_w(z^{-h}w^h)=-z\del_z(z^{-h}w^h)$, and
\eqref{E:delta is delta} we obtain
\begin{equation}
\label{E:shifted delta distr symmetric}
\del_w \de_h(z,w)
\; =\;
-\del_z \de_h(z,w).
\end{equation}
As in \eqref{E:delta distr homogeneous}
this implies 
\begin{equation}
\label{E:shifted delta distr homogeneous}
\de_h(z-x,w)
\; =\;
\de_h(z,w+x).
\end{equation}

\bigskip

{\bf Lemma.}\: {\it
For $h, h'\in\K$,
if $a(z,w)$ is a distribution in $z^{-h'}V\pau{z\uppm}\set{w}$
such that $a(w,w)$ is well-defined 
then
the product $\de_h(z,w)a(z,w)$ is well-defined
and 
\begin{equation}
\label{E:shifted delta distr is distr}
\de_h(z,w)a(z,w)
\; =\;
\de_{h+h'}(z,w)a(w,w).
\end{equation}
}

% no bigskip

\begin{pf}
We may assume $h=0$.
From $\de(z,w)a(z,w)=\de(z,w)z^{-h'}(z^{h'}a(z,w))$
and $\de_{h'}(z,w)a(w,w)=z^{-h'}w^{h'}\de(z,w)a(w,w)$
we see that we may also assume $h'=0$.
But the claim for $h=h'=0$ is just
Proposition \ref{SS:skew-symmetry distr}.
\end{pf}

\bigskip

If $a(\vz)$ and $b(\vz)$ are additively local fields on a vector space
and 
$\vr\in\vbbK$ then we define
$$
(\vz-\vw)^{\vr}a(\vz)b(\vw)_{w>z}
\; :=\;
(z-w)_{w>z}^r (\bz-\bw)_{\bz >\bw}^{\br}\,
a(\vz)b(\vw)_{w>z}
$$
and 
$$
(\vz-\vw)^{\vr}a(\vz)b(\vw)_{\bw >\bz}
\; :=\;
(z-w)_{z>w}^r (\bz-\bw)_{\bw >\bz}^{\br}\,
a(\vz)b(\vw)_{\bw >\bz}.
$$

\bigskip

{\bf Remark.}\:
Let $M$ be a bounded $\vbbK$-fold module over a
bounded $\vbbK$-fold algebra $V$,
$a, b\in V$, $c\in M$, and $\vr, \vs, \vt\in\vbbK$. 

($J$)\:
If $a$ and $b$ are additively local of order $\vr$,
$a(\vz)b$ is bounded by $\vr$, and 
$a(\vz)c$ is bounded by $\vt$
then 
the Jacobi identity for $a, b, c$ and any indices 
$\vn\in\vr+\Z^2, \vs\in\vbbK$, and $\vm\in\vt+\Z^2$ 
is equivalent to
\begin{align}
\label{E:Jacobi in fields}
&
\de_{\vt}(\vz,\vw+\vx)(a(\vx)b)(\vw)c
\\
\notag
\; =\;
&
\de_{\vr}(\vx,\vz-\vw)a(\vz)b(\vw)c
\; -\;
\de_{\vr}(\vx,\vz-\vw)a(\vz)b(\vw)_{w>z}c
\\
\notag
&\, -\de_{\vr}(\vx,\vz-\vw)a(\vz)b(\vw)_{\bw >\bz}c
\; +\;
\paraab\,
\de_{\vr}(\vx,\vz-\vw)b(\vw)a(\vz)c.
\end{align}
Here we use radial ordering and
write $\de_{\vr}(\vx,\vz-\vw)b(\vw)a(\vz)$ instead of
$\de_{\vr}(\vx,\vz-\vw)_{\vw>\vz}b(\vw)a(\vz)$.

($J_z$)\:
If $a$ is holomorphic then $a, b, c$ satisfy the
holomorphic Jacobi identity if and only if they satisfy
the Jacobi identity.

$\ita$\:
The vectors $a, b, c$ satisfy 
additive duality of order $\vt$ for 
indices $\vr$ and $\vs$ if and only if
the $(\vs,\vr)$-th modes of 
$(\vw+\vx)^{\vt}(a(\vx)b)(\vw)c$
and of 
$(\vx+\vw)^{\vt} a(\vx+\vw)b(\vw)c$
are equal.

$\itb$\:
The vectors $a, b, c$ satisfy 
additive locality of order $\vr$
for indices $\vt$ and $\vs$ if and only if 
the $(\vt,\vs)$-th mode of
$(\vz-\vw)^{\vr}[a(\vz),b(\vw)]c$ vanishes.

\bigskip

\begin{pf}
($J$)\:
We have
\begin{align}
\label{E:Jacobi in fields left}
&\de_{\vt}(\vz,\vw+\vx)\, (a(\vx)b)(\vw)c
\\
\notag
=\;
&\sum_{\vm\in\vt+\Z^2}\:
(\vw+\vx)^{\vm}\, (a(\vx)b)(\vw)c\; \vz^{-\vm-1}
\\
\notag
=\;
&\sum_{\vm\in\vt+\Z^2}\sum_{\vn, \vs\in\vbbK}
\bigg(
\sum_{\vi\in\N^2}\:
\binom{\vm}{\vi}\, (a_{(\vn+\vi)}b)_{(\vs+\vm-\vi)}c
\bigg)
\;
\vz^{-\vm-1} \vw^{-\vs-1} \vx^{-\vn-1}.
\end{align}
Since $a(\vz)b$ is bounded by $\vr$
the $(\vm,\vs,\vn)$-th mode
of the left-hand side of \eqref{E:Jacobi in fields left}
is equal to the left-hand side of the
Jacobi identity for $a, b, c$ and indices $\vn, \vs, \vm$
if $\vn\in\vr+\Z^2$ and $\vm\in\vt+\Z^2$
and it is zero otherwise.

On the other hand, we have
\begin{align}
\label{E:Jacobi in fields right}
&
\de_{\vr}(\vx,\vz-\vw)a(\vz)b(\vw)c
\; -\;
\de_{\vr}(\vx,\vz-\vw)a(\vz)b(\vw)_{w>z}c
\\
\notag
&\qquad\qquad
-\de_{\vr}(\vx,\vz-\vw)a(\vz)b(\vw)_{\bw >\bz}c
\; +\;
\paraab\,
\de_{\vr}(\vx,\vz-\vw)b(\vw)a(\vz)c
\\
\notag
\; = \; 
&\sum_{\vn\in\vr+\Z^2}\:
\Big(
(\vz-\vw)^{\vn}a(\vz)b(\vw)
\; -\;
(\vz-\vw)^{\vn}a(\vz)b(\vw)_{w>z}
\\
\notag
&\qquad\qquad
-(\vz-\vw)^{\vn}a(\vz)b(\vw)_{\bw >\bz}
\; +\;
\paraab\,
(\vz-\vw)^{\vn}b(\vw)a(\vz)
\Big) c\;\; \vx^{-\vn-1}
\\
\notag
= \;
&\sum_{\vn\in\vr+\Z^2}\,\sum_{\vs, \vm\in\K^2}\,\sum_{\vi\in\N^2}\:
(-1)^{\vi}\binom{\vn}{\vi}\,
\Big(
a_{(\vm+\vn-\vi)} b_{(\vs+\vi)}
\\
\notag
&-
(-1)^n
a_{([m+i,\bm+\bn-\bi)}b_{(s+n-i],\bs+\bi)}
-
(-1)^{-\bn}
a_{(m+n-i,[\bm+\bi)}b_{(s+i,\bs+\bn-\bi])}
\\
\notag
&\qquad\qquad\qquad\quad +\,
\paraab\, (-1)^{n-\bn} \, 
b_{(\vs+\vn-\vi)}a_{(\vm+\vi)}
\Big) c
\;\;
\vz^{-\vm-1} \vw^{-\vs-1} \vx^{-\vn-1}.
\end{align}
Because of the locality and boundedness assumptions
Lemma \ref{SS:ope commut}\,$\itii$,
\eqref{E:a[b] modes}, and \eqref{E:[a]b modes}
show that 
the $(\vm,\vs,\vn)$-th mode
of the left-hand side of \eqref{E:Jacobi in fields right}
is zero if $\vm\notin\vt+\Z^2$.
If $\vm\in\vt+\Z^2$ then it is equal to the right-hand side of the
Jacobi identity for $a, b, c$ and indices $\vn, \vs, \vm$.

($J_z$)\:
From \eqref{E:ab ba zw wz} we see that 
the Jacobi identity \eqref{E:Jacobi in fields}
is obtained from the holomorphic Jacobi identity
$$
\de(z,w+x)(a(x)b)(\vw)c
\; =\;
\de(x,z-w)[a(z),b(\vw)]c
$$
by multiplying the left-hand side with
$\de(\bz,\bw+\bx)$ and
multiplying the right-hand side with
$\de(\bx,\bz-\bw)-\de(\bx,\bz-\bw)_{\bw>\bz}$.
These two factors are equal by \eqref{E:jacobi for id}.
Conversely,
we know from Remark \ref{SS:jacobi modes}\,($J_z$)
that the Jacobi identity implies the
holomorphic Jacobi identity.

$\ita$\:
This is proven in the same way as
Remark \ref{SS:h Jac id distr}\,$\ita$.

$\itb$\:
This is clear.
\end{pf}

\subsection{Additive Duality and Additive Locality}
\label{SS:dual local pretex}

We prove that the Jacobi identity is equivalent to 
additive duality and additive locality.

\bigskip

{\bf Proposition.}\: {\it
Let $M$ be a bounded $\vbbK$-fold module over a
bounded $\vbbK$-fold algebra $V$, 
$a, b\in V$, $c\in M$, and $\vr, \vt\in\vbbK$
such that $a$ is bounded on $b$ and $c$ by $\vr$ and $\vt$, resp.
If $a, b, c$ satisfy additive duality of order $\vt$
and additive locality of order $\vr$
then 
they satisfy the Jacobi identity of orders $\vr$ and $\vt$.
}

\bigskip

\begin{pf}
Define 
$d(\vz,\vw,\vx):=\vx^{-\vr}(\vz-\vw)^{\vr}a(\vz)b(\vw)c$.
Additive locality and 
$a(\vz)c\in \vz^{-\vt}V\pau{\vz}$ implies that 
$d(\vz,\vw,\vx)$ is contained in
$\vz^{-\vt}\vx^{-\vr}V[\! [\vz,\vw\rangle$.
Together with additive duality this yields
\begin{align}
\notag
(\vw+\vx)^{\vt}d(\vw+\vx,\vw,\vx)
\; &=\;
(\vx+\vw)^{\vt}d(\vx+\vw,\vw,\vx)
\\
\notag
&=\;
(\vx+\vw)^{\vt}\vx^{-\vr}((\vx+\vw)-\vw)^{\vr}a(\vx+\vw)b(\vw)c
\\
\notag
&=\;
(\vx+\vw)^{\vt}a(\vx+\vw)b(\vw)c
\\
\notag
&=\;
(\vw+\vx)^{\vt}(a(\vx)b)(\vw)c.
\end{align}
Since $d(\vw+\vx,\vw,\vx)$ and $(a(\vx)b)(\vw)c$
are both contained in 
$V\sqbrack{\vw}\sqbrack{\vx}$ 
we can multiply the above equation by $(\vw+\vx)^{-\vt}$
and get that these two distributions are equal.
Lemma \ref{SS:jacobi pretex} implies
$$
\de(\vz,\vw+\vx)d(\vz,\vw,\vx)
=
\de_{\vt}(\vz,\vw+\vx)d(\vw+\vx,\vw,\vx)
=
\de_{\vt}(\vz,\vw+\vx)(a(\vx)b)(\vw)c.
$$
This is the left-hand side of the Jacobi identity 
\eqref{E:Jacobi in fields}. In order to obtain the
right-hand side we use the identity
$$
\de(\vz,\vw+\vx)
=
(\de(x,z-w)-\de(x,z-w)_{w>z})
(\de(\bx,\bz-\bw)-\de(\bx,\bz-\bw)_{\bw>\bz})
$$
which follows from \eqref{E:jacobi for id}.
Since $d(\vz,\vw,\vx)\in \vx^{-\vr}V\sqbrack{\vz,\vw}$ we 
can apply Lemma \ref{SS:jacobi pretex} again and
obtain that $\de(\vz,\vw+\vx)d(\vz,\vw,\vx)$ is also equal
to the right-hand side  of the Jacobi identity 
\eqref{E:Jacobi in fields} because
$$
\de(\vx,\vz-\vw)d(\vz,\vw,\vx)
=
\de_{\vr}(\vx,\vz-\vw)d(\vz,\vw,\vz-\vw)
=
\de_{\vr}(\vx,\vz-\vw)a(\vz)b(\vw)c,
$$
$$
\de(x,z-w)_{w>z}\de(\bx,\bz-\bw)d(\vz,\vw,\vx)
\; =\;
(\de_{\vr}(\vx,\vz-\vw)a(\vz)b(\vw)c)_{w>z},
$$
and similarly for the other two terms of \eqref{E:Jacobi in fields}.
\end{pf}

\bigskip

{\bf Corollary.}\: {\it
For a bounded $\vbbK$-fold algebra $V$ 
together with a gradation of the vector space $V$,
the following statements are equivalent:

\begin{enumerate}
\item[\iti]
$V$ is an additive OPE-algebra, the gradation is local,
and the identity is homogeneous;

\item[\itii]
there exists a homogeneous holomorphic right identity $1$
and homogeneous states satisfy additive duality and are
additively local;

\item[\itiii]
there exists a homogeneous holomorphic right identity $1$
and homogeneous states satisfy the Jacobi identity.
\end{enumerate}
}

\bigskip

\begin{pf}
$\iti\Rightarrow\itii$\:
This follows from Proposition \ref{SS:ope alg}.

$\itii\Rightarrow\itiii$\:
For any homogeneous states $a$ and $b$ 
there exists $\vh_{ab}\in\vbbK$ such that 
$a$ and $b$ are additively local of order $\vh_{ab}$.
Proposition \ref{SS:vacuum skew}\,$\itii$
implies that $a$ is bounded on $b$ by $\vh_{ab}$.
The left-hand side of the additive duality identity
$$
(\vw+\vx)^{\vt}(a(\vx)b)(\vw)c)
\; =\;
(\vx+\vw)^{\vt} a(\vx+\vw)b(\vw)c)
$$
is contained in $\vx^{-\vh_{ab}}V\sqbrack{\vw}\pau{\vx}$.
Applying $e^{-\vw\del_{\vx}}$ to the right-hand side
we obtain that $\vx^{\vt}a(\vx)b(\vw)c$ is contained in
$\vx^{-\vh_{ab}}V\lau{\vx}\sqbrack{\vw}$.
On the other hand, Lemma \ref{SS:ope commut}\,$\itii$
implies that 
$a(\vx)b(\vw)c$ is contained in 
$\vx^{-\vh_{ab}-\vh_{ac}}V\lau{\vx}\sqbrack{\vw}$.
Thus $\vt\in\vh_{ac}+\Z^2$ if $a(\vx)b(\vw)c$ is non-zero.
Hence the claim follows from the Proposition.

$\itiii\Rightarrow\iti$\:
Remark \ref{SS:jacobi modes}\,($J_z$), $\iti$, $\itii$
shows that $1$ is invariant for $\vT_1$
and that $\vT_1$ is a translation generator of $V$.
\end{pf}

\section{Vertex Algebras}
\label{S:vertex algebras}

{\bf Summary.}\:
In section \ref{SS:vertex algebras}
we give eight equivalent formulations of the notion of a vertex algebra. 
In section \ref{SS:modules}
we define the notion of a module over a vertex algebra and
prove that the space of fields of a module is local.
In section \ref{SS:modules conf vas}
we define four classes of modules over a conformal vertex algebra.
In section \ref{SS:commu vertex}
we prove that the category of commutative vertex algebras
is equivalent to the category of commutative algebras with a derivation.

\subsection{Vertex Algebras}
\label{SS:vertex algebras}

We give eight equivalent formulations of the notion of a vertex algebra. 
The three most important ones are in terms of
translation covariance and locality, in terms of
duality and skew-symmetry, and in terms of the Jacobi identity.

\bigskip

{\bf Proposition.}\: {\it
For a $\Z$-fold algebra $V$,
the following statements are equivalent:

\begin{enumerate}
\item[\iti]
$V$ is a vertex algebra;

\item[\itii]
there exist a translation generator and 
an invariant right identity and $V$ is holomorphically local;

\item[\itiii]
$V$ is bounded, there exists a left identity, and
$V$ satisfies the associativity formula and holomorphic skew-symmetry;

\item[\itiv]
$V$ is bounded, there exists a strong right identity, and
$V$ satisfies holomorphic duality and holomorphic skew-symmetry;

\item[\itv]
$V$ is bounded, there exists a right identity, 
and $V$ satisfies holomorphic duality and holomorphic locality;

\item[\itvi]
$V$ is bounded, there exists a right identity, 
and $V$ satisfies the holomorphic Jacobi identity;

\item[\itvii]
$V$ is bounded, $Y$ is injective,
there exists a left identity, 
and $V$ satisfies the holomorphic Jacobi identity;

\item[\itviii]
$V$ is bounded, $Y:V\to\cF_z(V)$ is a monomorphism of $\Z$-fold algebras,
and $\cF_Y$ contains the identity field and is holomorphically local.
\end{enumerate}
}

\bigskip

\begin{pf}
\iti$\Leftrightarrow$\itii\:
Let $1$ be an invariant right identity.
If $a$ and $b$ are holomorphically local states of $V$ 
of order $N$ then
$(z-w)^N a(z)b(w)1=\paraab(z-w)^N b(w)a(z)1$ is contained 
in $V\pau{z,w}$ and setting $w=0$ we obtain
$z^N a(z)b\in V\pau{z}$, 
i.e.~$a$ is bounded on $b$.
Proposition \ref{SS:ope commut} shows that
locality for $V$ implies holomorphic locality and 
the converse is true if $V$ is bounded.
Thus the claim follows.

\iti$\Rightarrow$\itiii\:
Proposition \ref{SS:ko vertex}\,$\itv$ 
and Dong's lemma for additive locality
show that
$Y:V\to\cF_z(V)$ is a morphism of $\Z$-fold algebras.
By Remark \ref{SS:nth products} this implies
the associativity formula.
Holomorphic skew-symmetry follows from
Proposition \ref{SS:ko vertex}\,$\itiii$.

\itiii$\Rightarrow$\itiv\:
Remark \ref{SS:skew symm} shows that
a left identity of $V$ is a strong identity.
Proposition \ref{SS:pf hol jacobi}\,$\iti$ shows that
$V$ satisfies holomorphic duality.

\itiv$\Rightarrow$\iti\:
This follows from Corollary \ref{SS:ope alg}
and the proof of the equivalence of $\iti$ and $\itii$.

\iti$\Leftrightarrow$\itv$\Leftrightarrow$\itvi\:
This follows from Corollary \ref{SS:dual local pretex}
by taking the trivial gradation $V=V$
because the Jacobi identity for holomorphic $a$
is identical with the holomorphic Jacobi identity
and additive duality and additive locality reduce
to holomorphic duality and holomorphic locality, resp.,
according to 
Remark \ref{SS:jacobi modes}\,($J_z$), ($a_z$), ($b_z$).

\iti$\Rightarrow$\itvii\:
This follows from the equivalence of $\iti$ and $\itv$.

\itvii$\Rightarrow$\itviii\:
This follows from the fact that the associativity
formula and holomorphic locality are special cases 
of the holomorphic Jacobi identity. 

\itviii$\Rightarrow$\iti\:
Because the identity field is a stable identity of $\cF_z(V)$
its inverse $Y^{-1}(1(z))$ is a stable identity of $V$. 
Thus $V$ is a vertex algebra.
\end{pf}

\subsection{Modules and Vertex Algebras of Fields}
\label{SS:modules}

We define the notion of a {\it module} over a vertex algebra and
prove that the space of fields of a module is local.

\bigskip

{\bf Definition.}\:
For a unital $\Z$-fold algebra $V$,
a vector space $M$ together with a morphism
$Y:V\to\cF_z(M)$ of unital $\Z$-fold algebras
is called a {\bf module} \index{module over a $\Z$-fold algebra}
over $V$.

\bigskip

Thus a module over a unital $\Z$-fold algebra $V$ 
is in particular a unitary bounded $\Z$-fold module over 
the vector space $V$.

\bigskip

{\bf Proposition.}\: {\it
If $M$ is a module over a unital $\Z$-fold algebra $V$ 
and $a$ and $b$ are states of $V$ 
that satisfy skew-symmetry for $T_1$
then $a$ and $b$ are mutually local on $M$.
In particular,
if $V$ is a vertex algebra 
then $V$ is local on $M$.
}

\bigskip

\begin{pf}
From the assumption that $Y:V\to\cF_z(M)$ is a morphism of $\Z$-fold algebras
follows that for any vector $c$ of $M$
the vectors $a, b, c$ and $b, a, c$ satisfy the associativity formula.
By Proposition \ref{SS:pf hol jacobi}\,$\iti$ this implies that 
the vectors $a, b$ and $b, a$ satisfy holomorphic duality.
The endomorphism $T_1$ of $V$ is a translation endomorphism of $M$
because
$(T_1 d)(z)=(d_{(-2)}1)(z)=d(z)_{(-2)}1(z)=\del_z d(z)$ for any $d\in V$.
Thus the claim follows from Proposition \ref{SS:dual and local}.
\end{pf}

\bigskip

{\bf Remark.}\: {\it
If $V$ is a vector space and
$\cF$ is a unital $\Z$-fold subalgebra of $\cF_z(V)$
consisting of mutually local fields
then $\cF$ is a vertex algebra. 
In particular, 
if $S$ is a local subset of $\cF_z(V)$
then the $\Z$-fold subalgebra $\sqbrack{S\cup\set{1(z)}}$ 
is a vertex algebra.

Conversely, if $\cF$ is a unital vertex subalgebra of $\cF_z(V)$
then the set $\cF$ of fields is local. 
}

\bigskip

\begin{pf}
Proposition \ref{SS:id field} states that 
the identity field is a full identity of $\cF_z(V)$.
By Proposition \ref{SS:holom loc implies itself} 
locality of the set $\cF$ implies that 
the $\Z$-fold algebra $\cF$ is local.
This yields the first claim.
The second claim follows from first claim and Dong's lemma.

If $\cF$ is a unital vertex subalgebra of $\cF_z(V)$
then $V$ together with the inclusion $\cF\to\cF_z(V)$
is a module over $\cF$. 
The Proposition implies that $\cF$ is local on $V$.
\end{pf}

\subsection{Modules over Conformal Vertex Algebras}
\label{SS:modules conf vas} 

We define four classes of modules over a conformal vertex algebra.

\bigskip

For $c\in\K$, a vertex algebra together with a conformal vector of 
central charge $c$  
is called a {\bf conformal} vertex algebra \index{conformal vertex algebra}
of central charge $c$. 
We identify a conformal vertex algebra $(V,L)$ with the
holomorphic conformal OPE-algebra $(V,(L,0))$.

A subset $S$ of $\K^2$ is called {\bf integrally bounded}
if for any $\vh\in\K^2$ there exists $\vN\in\Z^2$ such that
$S\cap(\vh+\Z^2)\subset \vh+\vN+\N^2$.
A $\K^2$-graded vector space $V$ is called integrally bounded
if its spectrum is integrally bounded.

A conformal vertex algebra is called a 
{\bf vertex operator algebra} \index{vertex operator algebra}
if its gradation is integrally bounded, locally finite, and statistical.
A vertex operator algebra is called a 
{\bf vertex algebra of CFT-type} \index{vertex algebra of CFT-type}
if the spectrum is contained in $(1/2)\N$ and $V_0=\K 1$.

\bigskip

{\bf Definition.}\:
Let $V$ be a conformal vertex algebra.

\begin{enumerate}
\item[$\iti$]
A $\Z$-fold module $M$ over $V$ is called
\begin{enumerate}
\item[$\ita$]
a {\bf weak module} if $M$ is a module over the vertex algebra $V$;
\item[$\itb$]
a {\bf generalized module} if $M$ is a weak module and
the operator $L_0$ of $M$ is diagonalizable;
\item[$\itc$]
a {\bf module} if $M$ is a generalized module and the 
spectrum of the operator $L_0$ of $M$ 
is locally finite and integrally bounded.
\end{enumerate}

\item[$\itii$]
A {\bf level gradation} of a weak $V$-module $M$ is an
$(1/2)\N$-gradation $(M_n)$ of $M$ such that
$a_{(n)}M_m$ is contained in $M_{m+h(a)-n-1}$ for any $a\in V$ and
$m\in(1/2)\N$.
A weak module is called 
{\bf gradable} if there exists a level gradation on $M$.
A weak module together with a level gradation is called
an {\bf admissible module}.
\end{enumerate}

\subsection{Commutative Vertex Algebras}
\label{SS:commu vertex}

We prove that the category of {\it commutative} vertex algebras
is equivalent to the category of commutative algebras with a derivation.

\bigskip

{\bf Definition.}\:
A $\Z$-fold algebra $V$ is called {\bf commutative} 
\index{commutative $\Z$-fold algebra}
if $V_{(\N)}V=\set{0}$.

\bigskip

A $\Z$-fold algebra $V$ is commutative if and only if
$\cF_Y$ is contained in $\End(V)\pau{z}$.
Let $V$ be a $\Z$-fold algebra.
If $V$ has a translation generator $T$ and
$a(z)V\subset V\lau{z}$ for any $a\in V$ 
then $V$ is commutative as follows by induction from
$-(n+1)a_{(n)}b=T(a_{(n+1)}b)-a_{(n+1)}Tb$.
Similarly, if $V$ has a translation endomorphism $T$ and
$\cF_Yb\subset V\lau{z}$ for any $b\in V$ 
then $V$ is commutative as follows by induction from
$-(n+1)a_{(n)}b=Ta_{(n+1)}b$.
In particular, any finite-dimensional vertex algebra
is commutative.

Let $A$ be a unital algebra with an even derivation $T$.
Define on $A$ the structure of a $\Z$-fold algebra by
$a_{(n)}b:=T^{(-1-n)}(a)b$. In other words, we have $a(z)b=(e^{zT}a)b$.
By definition, we have $A_{(\N)}A=0$ and $a_{(-1)}b=ab$. 
It is clear that the identity $1$ of the algebra $A$ 
is a strong identity of the $\Z$-fold algebra $A$.
Because $T^{(-1-n)}T=-n T^{(-1-(n-1))}$ the operator $T$
is a translation endomorphism of $A$. 
Since $T$ is clearly a derivation of the $\Z$-fold algebra $A$ 
the operator $T$ is a translation operator and $1$ is a full identity.

Conversely, let $V$ be a commutative fully unital $\Z$-fold algebra.
Define on $V$ the structure of an algebra by $ab:=a_{(-1)}b$. 
The identity $1$ of the vertex algebra $V$
is an identity of the algebra $V$.
The translation operator $T$ of $V$ is a derivation
of the algebra $V$.
Commutativity of the vertex algebra $V$ and \eqref{E:transl endo2} 
imply that $a_{(n)}b=T^{(-1-n)}(a)b$ for any integer $n$.
Thus the category of unital algebras with an even derivation
can be identified with 
the category of commutative fully unital $\Z$-fold algebras.

\bigskip

{\bf Lemma.}\: {\it
If $V$ and $W$ are $\Z$-fold algebras with a translation endomorphisms $T$
and $\phi:V\to W$ is a vector space morphism that is compatible with $T$
and such that $\phi:V|_{\set{-1}}\to W|_{\set{-1}}$ is an algebra morphism then
$\phi:V|_{\Z_<}\to W|_{\Z_<}$ is a $\Z_<$-fold algebra morphism.
}

\bigskip

\begin{pf}
For states $a$ and $b$ of $V$ and a negative integer $n$, we have
by induction $\phi(a_{(n-1)}b)=-n^{-1}\phi(Ta_{(n)}b)=
-n^{-1}T(\phi a)_{(n)}\phi b=\phi a_{(n-1)}\phi b$.
\end{pf}

\bigskip

{\bf Proposition.}\: {\it
$\iti$\:
A unital algebra $A$ with an even derivation is associative if and only if
the $\Z$-fold algebra $A$ satisfies the associativity formula.

$\itii$\:
The following statements about a vertex algebra $V$ are equivalent:

\begin{enumerate}
\item[$\ita$]
$V$ is commutative;
\item[$\itb$]
any two states of $V$ commute;
\item[$\itc$]
$V$ is a commutative algebra with an even derivation.
\end{enumerate}
}

\bigskip

\begin{pf}
$\iti$\:
Suppose that $A$ is associative. Let $a$ and $b$ be states. We have
$(a_{(-1)}b)(z)=(ab)(z)=e^{zT}(ab)=e^{zT}(a)e^{zT}(b)=a(z)_-b(z)=
a(z)_{(-1)}b(z)$.
From $\cF_Y\subset\End(A)\pau{z}$ follows that the $n$-th product
of $a(z)$ and $b(z)$ vanishes for any non-negative integer $n$.
Thus the associativity formula follows from the Lemma.
Conversely, the constant terms of $(a_{(-1)}b)(z)c=(e^{zT}(ab))c$ 
and $a(z)_{(-1)}b(z)c=e^{zT}(a)(e^{zT}(b)c)$ are $(ab)c$ and $a(bc)$.

$\itii$\:
The equivalence of $\ita$ and $\itb$ follows directly from
Remark \ref{SS:chiral alg}.

Suppose that $V$ is commutative.
By Corollary \ref{SS:locality} 
from $\cF_Y\subset\End(V)\pau{z}$ follows that 
$[a(z),b(w)]c=0$ for any states $a, b$, and $c$,
in particular $a(bc)=\paraab\, b(ac)$.
Taking $c=1$ shows that the algebra $V$ is commutative.
Moreover, we have 
$(ab)c=\zeta^{(\ta+\tb)\tc}c(ab)=\zeta^{\tb\tc}a(cb)=a(bc)$.
Thus $V$ is associative.

Conversely, it is clear that if $V$ is a commutative algebra
with an even derivation then the $\Z$-fold algebra $V$ is local
and hence is a commutative vertex algebra.
\end{pf}

\bigskip

The Proposition shows that there exist bounded fully unital
$\Z$-fold algebras that satisfy the associativity formula
but which are {\it not} vertex algebras.
Part $\itii$ of the Proposition implies that a vertex algebra has all 
$n$-th products equal to zero except for one
if and only if it is a commutative algebra.

%\contentsline {chapter}{Preface}{V}  put this into book.toc

\chapter{Structure of Vertex Algebras}
\label{C:str va}

\medskip

{\bf Conventions.}\:
Because we only consider $\N$-fold and $\Z$-fold algebras
in this chapter
the holomorphic Jacobi identity will just be called
the Jacobi identity.
Similarly, holomorphic duality, holomorphic locality, and 
holomorphic skew-symmetry will just be called
duality, locality, and skew-symmetry.
Integral distributions will just be called distributions.

\section{Conformal Algebras and Local Lie Algebras}
\label{S:conformal algebras}

{\bf Summary.}\:
In sections \ref{SS:conformal algebras}, \ref{SS:locmode}, and \ref{SS:loclie}
we define the notions of a conformal algebra, local mode algebra, 
and local Lie algebra
and construct a functor $\fg\mapsto R(\fg):=\bF_{\fg}$ 
from local Lie algebras to conformal algebras.
In section \ref{SS:affiniz} we construct local Lie algebras
by superaffinization of Lie algebras.

In section \ref{SS:conf affinization}
we prove that the affinization of a conformal algebra is 
a conformal algebra. 
In section \ref{SS:mod alg of n fold alg}
we define the mode algebra $\fg(R)$ of 
a bounded $\N$-fold algebra $R$ with a translation operator. 
In section \ref{SS:loc mod alg of n fold alg}
we prove that $\fg(R)$ is a local mode algebra such that $R=F_{\fg(R)}$.
In section \ref{SS:Lie of ConfA}
we remark that the functor $\ConfA\to\locLie, R\mapsto\fg(R)$,
is the left adjoint of $\fg\mapsto R(\fg)$ 
and prove that $\ConfA$ is a localization of $\locLie$.

In section \ref{SS:free conf} we prove the existence of 
the free conformal algebra gene\-rated by a set with a pole order bound.
In section \ref{SS:vertex reps}
we study the relation between vertex modules over
conformal algebras and modules over local Lie algebras.

\subsection{Conformal Algebras}
\label{SS:conformal algebras}

We define the notion of a {\it conformal algebra} and prove that 
if $\fg$ is a Lie algebra then any local subset of $\fg\pau{z\uppm}$ 
generates a conformal algebra.

\bigskip

{\bf Remark.}\: {\it
The Jacobi identity, the commutator formula, and 
the associativity formula are all equivalent
for elements $a, b, c$ of an $\N$-fold algebra.
}

\bigskip

\begin{pf}
This follows from Remark \ref{SS:hol jacobi}\,$\iti$ and
Lemma \ref{SS:pf hol jacobi}.
\end{pf}

\bigskip

{\bf Definition.}\:
A bounded $\N$-fold algebra $R$ together with a translation operator 
is called a {\bf conformal algebra} \index{conformal algebra}
if $R$ satisfies the Jacobi identity and skew-symmetry.

\bigskip

Let $R$ be an $\N$-fold algebra.
For elements $a$ and $b$ of $R$, we call $o_R(a,b):=o'(a,b)_+$ 
\index{oaaaaaaaab@$o_R(a,b)$}
the mode infimum \index{mode infimum!for $\N$-fold algebras} of $a$ and $b$.
If $R$ is endowed with an operator $T$ and
$S$ is a subset of $R$ then
we denote by $\bar{S}$ \index{faaaaaaaa@$\bF$}
the smallest $\N$-fold subalgebra of $R$
that contains $S$ and that is invariant with respect to $T$.

\bigskip

{\bf Proposition.}\: {\it
If $\fg$ is a Lie algebra and 
$R$ is a local $\N$-fold subalgebra of $\fg\pau{z\uppm}$
that is invariant with respect to $\del_z$
then $R$ together with $\del_z$ is a conformal algebra.
In particular,
if $F$ is a local subset of $\fg\pau{z\uppm}$
then $\bF$ is a conformal algebra.
}

\bigskip

\begin{pf}
Propositions \ref{SS:derivative}\,$\iti$
and \ref{SS:Jacobi formal distributions}
state that $\del_z$ is a translation operator of $\fg\pau{z\uppm}$ and 
that $\fg\pau{z\uppm}$ satisfies the Jacobi identity.
By Proposition \ref{SS:skew-symmetry distr}
locality of $R$ implies that $R$ satisfies skew-symmetry.
This yields the first claim.
The second claim follows from first claim and Dong's lemma.
\end{pf}

\bigskip

Conversely, 
a conformal subalgebra of $\fg\pau{z\uppm}$ need {\it not} be local,
e.g. the conformal subalgebra $\fg$ of $\fg\pau{z\uppm}$
is non-local if $\fg$ is non-abelian.

Let $R$ be an $\N$-fold algebra with an operator $T$ and $k\in\K$.
An even element $\hk$ of $R$ is called {\bf central}
\index{central element of an $\N$-fold algebra}
if $R_{(\N)}\hk=\hk_{(\N)}R=\set{0}$ and $\hk$ is invariant for $T$.
A conformal algebra together with a central element $\hk$
is called a {\bf centered} conformal algebra
\index{centered conformal algebra} 
of {\bf level} \index{level!of a centered conformal algebra}
$k$ if $k\ne 0$ implies $\hk\ne 0$.

We denote by $\ConfA$ \index{ConfA@$\ConfA$}
the category of conformal algebras and by 
$\ConfA_k$ \index{ConfAk@$\ConfA_k$}
the category of centered conformal algebras of level $k$.
The forgetful functor $f:\ConfA_k\to\ConfA$
has the functor $R\mapsto (R\oplus\K\hk,\hk)$ as left adjoint.
There is a forgetful functor $f_k$ from non-zero vertex algebras to
$\ConfA_k$ that maps a non-zero vertex algebra $V$ with identity $1$ 
to $(V|_{\N},T_1,k1)$.

A vector subspace $I$ of an $\N$-fold algebra $R$ with an operator $T$
is called an {\bf ideal} \index{ideal of an $\N$-fold algebra}
if $R_{(\N)}I, I_{(\N)}R$, and $T(I)$ are contained in $I$.

\subsection{Local Mode Algebras}
\label{SS:locmode}

We define the notion of a {\it local mode algebra}
and make some remarks about it.

\bigskip

{\bf Definition.}\:
An algebra $\fg$ together with 
a local subset $F\equiv F_{\fg}$ of $\fg\pau{z\uppm}$
is called a {\bf local mode algebra} \index{local mode algebra}
if $\fg$ is spanned by the vectors $a_n$ where $a(z)\in F$ and $n\in\Z$
and if there exists an even derivation $T$ of $\fg$ such that
$F$ is translation covariant for $T$.

\bigskip

A morphism $\fg\to\fg'$ of local mode algebras 
is a morphism $\fg\to\fg'$ of algebras such that $F\to\bF'$. 
The even derivation of a local mode algebra $\fg$,
for which $F$ is translation covariant, 
is unique and is called {\it the} derivation of $\fg$.
A morphism $\phi:\fg\to\fg'$ of local mode algebras 
is compatible with the derivations of $\fg$ and $\fg'$
because by Proposition \ref{SS:derivative}\,$\itiii$
the set $\bF'$ is translation covariant for the derivation of $\fg'$ 
and hence $\phi(Ta_n)=-n\phi(a_{n-1})=T\phi(a_n)$ for any 
$a(z)\in F$ and $n\in\Z$.

Let $\fg$ be a local mode algebra.
We denote by $\fg_{\geq}$ and $\bfg_{\geq}$
the spans of the vectors $a_n$ of $\fg$ where 
$a(z)\in F$, resp., $a(z)\in\bF$ and $n\in\N$.
The vector subspaces $\fg_<$ and $\bfg_<$ 
are defined in the same way with $\N$ replaced by $\Z_<$.

We have $\fg=\fg_<+\fg_{\geq}$.
Because $Ta_n=-n\, a_{n-1}$ 
we have $T\fg_{\geq}=\fg_{\geq}$ and $T\bfg_{\geq}=\bfg_{\geq}$
and the subspaces $\fg_<$ and $\bfg_<$ are invariant with respect to $T$.
The subspace $\bfg_{\geq}$ is the subalgebra of $\fg$
that is generated by $\fg_{\geq}$. Indeed,
the commutator formula
\begin{equation}
\label{E:comm form}
[a_n,b_m]
\; =\;
\sum_{i\in\N}\:
\binom{n}{i}\, (a(z)_{(i)}b(z))_{n+m-i},
\end{equation}
that holds for any $n\in\N$ and $m\in\Z$,
shows that $\bfg_{\geq}$ is a subalgebra. 
Conversely, 
$\bfg_{\geq}$ is contained in the subalgebra generated by $\fg_{\geq}$
because the $m$-th mode of $a(z)_{(n)}b(z)$ is equal to
$\sum_{i\in\N}(-1)^{n+i}\binom{n}{i}[a_i,b_{m+n-i}]$.

An even element $\hk$ of a local mode algebra $\fg$ 
is called an {\bf central}
if the constant distribution $\hk$ is contained in $F$ 
and if $\hk\ne 0$ implies that $\hk\notin\bfg_{\geq}$.
Let $\hk$ be a central element of $\fg$.
By Corollary \ref{SS:locality} 
locality of $F$ implies that $[\hk,\fg]=0$.
Translation covariance of $\hk$
implies that $\hk$ is invariant for $T$.

\bigskip

{\bf Lemma.}\: {\it 
Let $\fg$ be a local mode algebra, $\fg'$ be an algebra, 
and $\phi:\fg\to\fg'$ be a morphism of vector spaces.
If $\phi F$ is local and 
$\phi:\bF\to\fg'\pau{z\uppm}$ is a morphism of $\N$-fold algebras
then $\phi:\fg\to\fg'$ is an algebra morphism.
}

\bigskip

\begin{pf}
Let $a(z), b(z)\in F$ and $n, m\in\Z$.
Proposition \ref{SS:locality} shows that locality of 
$F$ and $\phi F$ implies that $a(z), b(z)$ and 
$\phi a(z), \phi b(z)$ satisfy the commutator formula.
Because $\phi:\bF\to\fg'\pau{z\uppm}$ is a morphism of $\N$-fold algebras
we have 
\begin{align}
\notag
\phi[a_n,b_m]
\; =\;
&\phi\sum_{i\in\N}\binom{n}{i}(a(z)_{(i)}b(z))_{n+m-i}
\\
\notag
\; =\;
&\sum_{i\in\N}\binom{n}{i}(\phi a(z)_{(i)}\phi b(z))_{n+m-i}
\\
\notag
\; =\;
&[\phi a_n,\phi b_m].
\end{align}
\end{pf}

\subsection{Local Lie Algebras}
\label{SS:loclie}

We define the notion of a {\it local Lie algebra}
and construct a functor $\fg\mapsto R(\fg)$ 
from local Lie algebras to conformal algebras.

\bigskip

{\bf Definition.}\:
A local mode algebra $\fg$ is called a 
{\bf local Lie algebra} \index{local Lie algebra}
if $\fg$ is a Lie algebra.

\bigskip

Proposition \ref{SS:conformal algebras} shows that 
there exists a functor
$\fg\mapsto R(\fg):=(\bF,\del_z)$ 
from the category of local Lie algebras to 
the category of conformal algebras.
The conformal algebra $R(\fg)$ is called the conformal algebra of $\fg$.

If $\fg$ is a local Lie algebra then 
the subspaces $\bfg_{\geq}$ and $\bfg_<$ are subalgebras of $\fg$
because by Remark \ref{SS:locality of formal distributions}\,$\iti$
locality of $\bF$ implies that the commutator formula
\eqref{E:comm form} is satisfied 
for any $a(z), b(z)\in\bF$ and $n, m\in\Z$.

Let $k\in\K$.
A local Lie algebra together with a central element $\hk$ is called
a {\bf centered} local Lie algebra \index{centered local Lie algebra}
of level $k$ if $k\ne 0$ implies that $\hk\ne 0$.
We denote by $\locLie$ the category of local Lie algebras
and by $\locLie_k$ the category of centered local Lie algebras
of level $k$. 
The forgetful functor $f:\locLie_k\to\locLie$ has 
$\fg\mapsto \fg\oplus\K\hk=
(\fg\oplus\K\, \hk,F\cup\set{\hk},\hk)$ as left adjoint.

In general, if $M$ is a module over a Lie algebra $\fg$
then we denote by $\rho$ the corresponding representation $\fg\to\End(M)$.
A {\bf module} \index{module!over a local Lie algebra}
over a local Lie algebra $\fg$ is 
a module $M$ over the Lie algebra $\fg$ such that
$\rho F$ is contained in $\cF_z(M)$.
For $k\in\K$,
a module $M$ over a centered local Lie algebra $\fg$ is said to be
of level $k$ \index{level!of a module}
if $\hk$ acts by multiplication with $k$ on $M$.

\bigskip

{\bf Remark.}\: {\it 
Let $\fg$ be a local Lie algebra and 
$M$ be a module over the Lie algebra $\fg$.
If there exists a subspace $W$ of $M$ such that
$F$ is bounded on $W$ and $U(\fg)W=M$ then 
$M$ is a module over the local Lie algebra $\fg$.
}

\bigskip

\begin{pf}
The canonical filtration of $U(\fg)$ induces an exhaustive increasing
filtration $(M_r)_{r\in\N}$ of $M$ such that $M_0=W$.
We prove by induction on $r$ that $F$ is bounded on $M_r$.

Let $a(z), b(z)\in F$ and $c\in M_r$.
Because $a(z)$ and $b(z)$ are mutually local on $M$
there exists a non-negative integer $n$ such that 
$(z-w)^n a(z)b(w)=\paraab (z-w)^n b(w)a(z)$. 
Since $a(z)$ and $b(z)$ are bounded on $c$ this implies that 
$d(z,w):=(z-w)^n a(z)b(w)c\in M\lau{z,w}$. 
The distribution $(z-w)^{-n}d(z,w)-a(z)b(w)c$ is bounded in $w$
and is annihilated by $(z-w)^n$.
By Corollary \ref{SS:locality} we get $(z-w)^{-n}d(z,w)=a(z)b(w)c$
and hence $a(z)b(w)c\in M\lau{z}\lau{w}$. 
This shows that $a(z)$ is bounded on $b_n c$ for any integer $n$.
\end{pf}

\subsection{Superaffinization of Lie Algebras}
\label{SS:affiniz}

We construct centered local Lie algebras 
by {\it superaffinization} of Lie algebras.

\bigskip

Let $\fg$ be a Lie algebra.
If $C$ is a commutative algebra together with an even derivation $\del$
then $\fg\otimes C$ with Lie bracket
$$
[a\otimes x,b\otimes y]
\; :=\;
\zeta^{\tilde{x}\tb}\,
[a,b]\otimes xy
$$
is a Lie algebra with an even derivation 
$T: a\otimes x\mapsto a\otimes \del x$.
The Lie algebra $\fg$ is a Lie subalgebra of $\fg\otimes C$
via $a\mapsto a\otimes 1$.

Let $\K[t,\theta]$ be the polynomial algebra
in an even and an odd variable $t$ and $\theta$
and let $\K[t\uppm,\theta]$ be the algebra 
obtained from $\K[t,\theta]$ by inverting $t$.
The Lie algebra 
$\tfg^s:=\fg\otimes\K[t\uppm,\theta]$
is called the {\bf superloop Lie algebra} \index{superloop Lie algebra}
of $\fg$.
The derivation $-\del_t$ of $\K[t\uppm,\theta]$ 
induces a derivation $T$ of $\tfg^s$.

The superloop Lie algebra is the semidirect product
$\tfg\ltimes\tfg^c$ of the {\bf loop Lie algebra} \index{loop Lie algebra}
$\tfg:=\fg\otimes\K[t\uppm]$ and the abelian Lie algebra 
$\tfg^c:=\fg\otimes\K[t\uppm]\theta$.

\bigskip

Define $a_n:=a\otimes t^n$ and $\ba_m:=a\otimes t^{m-1/2}\theta$
for $a\in\fg, n\in\Z$, and $m\in 1/2+\Z$.
An even symmetric pairing $(\; ,\;)$ on $\fg$ 
defines a Chevalley-Eilenberg 2-cochain
$\varep:\bigwedge^2\tfg^s\to\K$ on $\tfg^s$ by
$$
a_n\wedge b_m
\;\mapsto\; 
n\, (a,b)\, \de_{n+m,0},
\qquad\quad
\ba_n\wedge\bb_m
\;\mapsto\; 
(b,a)\, \de_{n+m,0}, 
$$
and $a_n\wedge\bb_m\mapsto 0$.
Recall that an even symmetric pairing $(\; ,\;)$ on $\fg$ 
is called {\bf invariant} if 
$([a,b],c)=(a,[b,c])$ for any $a, b, c\in\fg$.

\bigskip

{\bf Remark.}\: {\it
If $(\; ,\;)$ is an even invariant symmetric pairing on $\fg$
then $\varep$ is a Chevalley-Eilenberg 2-cocycle on $\tfg^s$, 
i.e.~for any $a, b, c\in\tfg^s$ we have
$$
\varep([a,b],c)
\; -\;
\zeta^{\tb\tc}  \,\varep([a,c],b)
\; +\;
\zeta^{\ta(\tb+\tc)}  \,\varep([b,c],a)
\; =\;
0.
$$
}

% no bigskip

\begin{pf}
In the case that $a, b, c$ are equal to $a_n, b_m, c_k\in\tfg$
the claim follows from vanishing of 
\begin{align}
\notag
&\big(
(n+m)([a,b],c)
-
\zeta^{\tb\tc}  (n+k) ([a,c],b)
+
\zeta^{\ta(\tb+\tc)}  (m+k) ([b,c],a)
\big)
\de_{n+m+k,0}
\\
\notag
&=\;
2(n+m+k)\, ([a,b],c)\, \de_{n+m+k,0}.
\end{align}
If $a, b, c$ are equal to $a_n, b_m, \bc_k$ then the claim
follows from $\varep(a_n,\bb_m)=0$.
If $a, b, c$ are equal to $a_n, \bb_m, \bc_k$ then the claim
follows from vanishing of 
\begin{align}
\notag
&\varep(\overline{[a,b]}_{n+m},\bc_k)
\; -\;
\zeta^{(\tb+1)(\tc+1)}\, \varep(\overline{[a,c]}_{n+k},\bb_m)
\\
\notag
=\;
&\big(
\zeta^{(\ta+\tb)\tc}
\, +\,
\zeta^{(\ta+\tc)\tb+\tb+tc+1}
\big)\:
([a,b],c)\,\de_{n+m+k,0}.
\end{align}
Here we use that $([a,b],c)$ vanishes unless $\ta+\tb+\tc=0$
because $(\; ,\;)$ and $\K$ are even.
Finally, if $a, b, c$ are equal to $\ba_n, \bb_m, \bc_k$
then the claim is trivial because $\tfg^c$ is abelian.
\end{pf}

\bigskip

Recall that a Chevalley-Eilenberg 2-cocycle 
$\varep:\bigwedge^2\tfg^s\to\K$
defines a
one-dimensional central extension 
$\hfg^s:=\tfg^s\oplus\K\hk$ of $\tfg^s$
with Lie bracket
$[a,b]_{\hfg^s}=[a,b]_{\tfg^s}+\varep(a,b)\hk$
for any $a, b\in\tfg^s$.
If we take for $\varep$ the 2-cocycle
induced by an even invariant symmetric pairing on $\fg$
then $\hfg^s$ is called the {\bf superaffinization} \index{superaffinization}
of $\fg$.
The derivation $T$ of $\tfg^s$ is extended to 
a derivation $T$ of $\hfg^s$ by $T\hk=0$.

By definition, the Lie brackets of $\hfg^s$ are given by
$$
[a_n,b_m]
\; =\;
[a,b]_{n+m}
\; +\;
n\, (a,b)\,\de_{n+m,0}\hk,
$$
$$
[a_n,\bb_m]
\; =\;
\overline{[a,b]}_{n+m},
\qquad\text{and}\qquad
[\ba_n,\bb_m]
\; =\;
(b,a)\,\de_{n+m,0}\hk.
$$
%For an element $a$ of $\fg$,
For $a\in\fg$, we define the $\hfg^s$-valued distributions
$a(z):=\sum_{n\in\Z}a_n z^{-n-1}$ and
$\ba(z):=\sum_{n\in\Z}\ba_{n+1/2} z^{-n-1}$.
Proposition \ref{SS:locality} and the above commutators of $\hfg^s$
show that the set $F:=\set{a(z),\ba(z),\hk\mid a\in\fg}$ is local and 
we have
\begin{align}
\label{E:lie opes}
a(z)b(w)
\;&\sim\;
\frac{(a,b)\hk}{(z-w)^2}
\; +\;
\frac{[a,b](w)}{z-w},
\\
\notag
\ba(z)\bb(w)
\;\sim\;
\frac{(b,a)\hk}{z-w},
&\qquad\text{and}\qquad
a(z)\bb(w)
\;\sim\;
\frac{\overline{[a,b]}(w)}{z-w}.
\end{align}
Moreover, $F$ is translation covariant because
$Ta_n=T(a\otimes t^n)=-n\, a_{n-1}$.
Thus $\hfg^s$ together with $F$ and $\hk$ is a centered local Lie algebra.

The local Lie subalgebra $\hfg:=\tfg\oplus\K\hk$ of $\hfg^s$
is called the {\bf affinization} \index{affinization!of a Lie algebra}
of $\fg$.
If $V$ is a vector space with an even symmetric pairing
$(\;, \;)$ 
then the local Lie subalgebra $\hV^c:=\tilde{V}^c\oplus\K\hk$
of $\hat{V}^s$ is called the {\bf Clifford affinization}
\index{Clifford affinization}
of $V$ where we consider $V$ as an abelian Lie algebra
endowed with the even symmetric pairing $a\otimes b\mapsto (b,a)$.

\subsection{Affinization of Conformal Algebras}
\label{SS:conf affinization}

We construct an endofunctor of the category of conformal algebras
that is called {\it affinization}.

\bigskip

{\bf Proposition.}\:{\it
Let $R$ be a bounded $\N$-fold algebra together with an even operator $T$ 
and $C$ be a commutative algebra with an even derivation $\del$.
Define on $R\otimes C$ the structure of a
bounded $\N$-fold algebra by
$$
(a\otimes f)_{(n)}(b\otimes g)
\; :=\;
\sum_{i\in\N}\:
(a_{(n+i)}b) \otimes (\del^{(i)}(f)g) 
$$
and define an operator $\tilde{T}$ of $R\otimes C$
by $a\otimes f\mapsto T a\otimes f+a\otimes \del f$.
\begin{enumerate}
\item[\iti]
If $T$ is a translation generator so is $\tilde{T}$.
\item[\itii]
If $T$ is a derivation so is $\tilde{T}$.
\item[\itiii]
If $R$ satisfies the Jacobi identity so does $R\otimes C$.
\item[\itiv]
If $R$ satisfies skew-symmetry for $T$ then
$R\otimes C$ satisfies skew-symmetry for $\tilde{T}$.
\item[\itv]
If $\hk$ is a central element of $R$ then 
the element $\hk\otimes 1$ of $R\otimes C$ is central.
\item[\itvi]
There exists a bifunctor $(R,C)\mapsto R\otimes C$
from conformal algebras and 
commutative algebras with an even derivation to conformal algebras.
\end{enumerate}
}

\bigskip

We will show in Remark \ref{SS:vertex envelope} that
any conformal algebra $R$ can be embedded into a vertex algebra $V(R)$.
If we already had this result at our disposal
we would get a more conceptual proof of part $\itvi$ of
the Proposition because $R\otimes C\to V(R)\otimes C$ is a monomorphism 
of $\N$-fold algebras with an operator and $V(R)\otimes C$ is a 
vertex algebra by Proposition \ref{SS:commu vertex}\,$\itii$.

\bigskip

\begin{pf}
For $a\in R$ and $f\in C$, we write $af:=a\otimes f$.
We omit the supersigns in the proof.

$\iti$\:
This follows from
\begin{align}
\notag
\tilde{T}(af)_{(n)}bg
\; &=\;
(T(a)f+a\del f)_{(n)}bg
\\
\notag
&=\;
\sum_{i\in\N}\:
(Ta_{(n+i)}b)\, \del^{(i)}(f)g
\, +\,
(a_{(n+i-1)}b)\, \del^{(i-1)}(\del f)g
\\
\notag
&=\;
\sum_{i\in\N}\:
(-(n+i)+i)(a_{(n-1+i)}b)\, \del^{(i)}(f)g
\\
\notag
&=\;
-n\, af_{(n-1)}bg.
\end{align}

$\itii$\:
This is clear.

$\itiii$\: 
Due to Remark \ref{SS:conformal algebras}
it suffices to establish the commutator formula for $R\otimes C$.
We have
\begin{align}
\notag
&\sum_{i\in\N}\,
\binom{n}{i}
(af_{(i)}bg)_{(m+n-i)}ch
\\
\notag
=\;
&\sum_{i,j,k\in\N}\,
\binom{n}{i}
((a_{(i+j)}b)_{(m+n-i+k)}c)
\, \del^{(k)}(\del^{(j)}(f)g)h
\\
\notag
=\;
&\sum_{i,j,k,l\in\N}\:
\binom{n}{i}
((a_{(i+j)}b)_{(m+n-i+k)}c)
\del^{(l)}\del^{(j)}(f)\, \del^{(k-l)}(g)h
\\
\notag
=\;
&\sum_{i,j,k,l\in\N}\:
\binom{n}{i-j}
((a_{(i)}b)_{(m+n-i+k+l)}c)
\, \del^{(l-j)}\del^{(j)}(f)\, \del^{(k)}(g)h.
\end{align}
Because $\sum_j\binom{n}{i-j}\binom{l}{j}=\binom{n+l}{i}$
the last expression is equal to 
\begin{align}
\notag
&af_{(n)}(bg_{(m)}ch)
\, -\,
bg_{(m)}(af_{(n)}ch)
\\
\notag
=\;
&\sum_{l,k\in\N}\,
\big(
a_{(n+l)}b_{(m+k)}c\, -\, b_{(m+k)}a_{(n+l)}c
\big)
\,
\del^{(l)}(f)\del^{(k)}(g)h
\\
\notag
=\;
&\sum_{i,l,k\in\N}\:
\binom{n+l}{i}((a_{(i)}b)_{(n+l+m+k-i)}c)\,
\del^{(l)}(f)\, \del^{(k)}(g)h.
\end{align}

$\itiv$\:
We have
\begin{align}
\notag
&\sum_{j\in\N}\,
(-1)^{n+1+j}\tilde{T}^{(j)}(af_{(n+j)}bg)
\\
\notag
=\;
&\sum_{i,j,k\in\N}
(-1)^{n+1+j}
T^{(k)}(a_{(n+j+i)}b)\del^{(j-k)}(\del^{(i)}(f)g)
\\
\notag
=\;
&\sum_{i,j,k,l\in\N}
(-1)^{n+1+j}
T^{(k)}(a_{(n+j+i)}b)\del^{(l)}\del^{(i)}(f)\del^{(j-k-l)}(g)
\\
\notag
=\;
&\sum_{i,j,k,l\in\N}
(-1)^{n+1+j+k+i}
T^{(k)}(a_{(n+j+k)}b)\del^{(l)}\del^{(i)}(f)\del^{(j-(l+i))}(g).
\end{align}
Because 
$\sum_{l+i=m}(-1)^i\, \del^{(l)}\del^{(i)}=(\del-\del)^{(m)}=0$ 
for $m>0$ the last expression is equal to
\begin{align}
\notag
bg_{(n)}af
\; &=\;
\sum_{j\in\N}\,
(b_{(n+j)}a) \, \del^{(j)}(g)f
\\
\notag
&=\;
\sum_{j,k\in\N}
(-1)^{n+j+1+k}\, T^{(k)}(a_{(n+j+k)}b) \, f\del^{(j)}(g).
\end{align}

$\itv$\:
This is clear.

$\itvi$\:
This follows from $\iti$--$\itiv$.
\end{pf}

\bigskip

{\bf Definition.}\:
The endofunctor $R\mapsto R\otimes (\K[t\uppm],\del_t)$
of the category of bounded $\N$-fold algebras with a translation operator 
is called {\bf affinization} \index{affinization!of an $\N$-fold algebra}
and is written $R\mapsto \tR$.

\subsection{The Mode Algebra of an $\N$-Fold Algebra}
\label{SS:mod alg of n fold alg}

We prove that if $R$ is a conformal algebra then $R/T(R)$ is a Lie
algebra and we define the {\it mode algebra} $\fg(R)$ of a 
bounded $\N$-fold algebra $R$ with a translation operator.

\bigskip

{\bf Remark.}\: {\it
Let $R$ be an $\N$-fold algebra with a translation operator $T$.

\begin{enumerate}
\item[$\iti$]
The $0$-th product induces an algebra structure on $R/T(R)$.

\item[$\itii$]
If $R$ satisfies the Jacobi identity then $R/T(R)$ is a Leibniz algebra.

\item[$\itiii$]
If $R$ is a conformal algebra then $R/T(R)$ is a Lie algebra.
\end{enumerate}
}

\bigskip

The multiplication on $R/T(R)$, that is induced by the $0$-th product,
is written $a\otimes b\mapsto [a,b]$.

\bigskip

\begin{pf}
$\iti$\:
This follows from $Ta_{(0)}b=0$ and $a_{(0)}Tb=T(a_{(0)}b)$.

$\itii$\:
The Jacobi identity for indices $0, 0, 0$
is just the Leibniz identity
$(a_{(0)}b)_{(0)}c=a_{(0)}(b_{(0)}c)-\paraab\, b_{(0)}(a_{(0)}c)$.

$\itiii$\:
This follows from $\itii$ and 
$\paraab\,b_{(0)}a=\sum_{i\in\N}(-1)^{1+i}T^{(i)}(a_{(i)}b)$.
\end{pf}

\bigskip

If $R$ is a bounded $\N$-fold algebra with a translation operator $T$
then the algebra $\tR /\tT(\tR)$ is called the 
{\bf mode algebra} \index{mode algebra!of an $\N$-fold algebra} 
of $R$ and is denoted by $\fg(R)$.
Thus there exists a functor $R\mapsto\fg(R)$ from
the category of bounded $\N$-fold algebras
with a translation operator to the category of algebras.
Proposition \ref{SS:conf affinization} and the Remark
show that $R\mapsto\fg(R)$ defines a functor from
conformal algebras to Lie algebras.

\subsection{The Local Mode Algebra of an $\N$-Fold Algebra}
\label{SS:loc mod alg of n fold alg}

We show that the mode algebra $\fg(R)$ of a 
bounded $\N$-fold algebra with a translation operator
has the structure of a local mode algebra and 
prove that there exists a natural isomorphism $Y:R\to R(\fg(R))$.

\bigskip

Let $R$ be bounded $\N$-fold algebra
with a translation operator $T$ and $\fg$ be the mode algebra of $R$.
For an integer $n$,
we define a map $R\to\fg$ by $a\mapsto a_n:=\io(a\otimes t^n)$ 
where $\io:\tR\to\fg$ is the quotient map.
These maps define a map 
$Y:R\to\fg\pau{z\uppm}, a\mapsto a(z):=\sum_{n\in\Z}a_n z^{-n-1}$.

The operator $T\otimes\id_{\K[t\uppm]}$ of $\tR$ 
is a derivation of $\tR$ and induces a derivation $T$ of $\fg$
because $T\otimes\id$ and $\tT$ commute.

For an even set $S$, we denote by $\SAlg{S}$ \index{SAlg@$\SAlg{S}$}
the category of $S$-fold algebras and even morphisms.
We denote by $\SAlg{\N}_T$ the category of $\N$-fold algebras
with a translation operator and even morphisms.

\bigskip

{\bf Proposition.}\: {\it
Let $R$ be a bounded $\N$-fold algebra 
with a translation ope\-rator $T$ and $\fg$ be the mode algebra of $R$.

\begin{enumerate}
\item[$\iti$]
The map
$Y:R\to\fg\pau{z\uppm}$ is a monomorphism of $\N$-fold algebras
with translation operator and
the image $Y(R)$ is local and translation covariant for $T$.
In particular, $(\fg,Y(R))$ is a local mode algebra.

\item[$\itii$]
The subspaces $\fg_{\geq}$ and $\fg_<$ are subalgebras of $\fg$,
we have $\fg=\fg_<\oplus\fg_{\geq}$, 
$\fg_{\geq}=\bfg_{\geq}$, $\fg_< =\bfg_<$, 
and the map
$\ka: R\to\fg_<, a\mapsto a_{-1}$, is a vector space isomorphism. 

\item[$\itiii$]
For an algebra $\fg'$, there exists a natural bijection
\begin{align}
\notag
\Alg(\fg(R),\fg')
\; &\lra\;
\set{\phi\in\SAlg{\N}_T(R,\fg'\pau{z\uppm})\mid \text{$\phi R$ is local}},
\\
\notag
\phi
\; &\longmapsto\;
\phi\circ Y;
\end{align}
the inverse is given by $\phi\mapsto (a_n\mapsto (\phi a)_n)$.
\end{enumerate}
}

\bigskip

Roughly speaking, $\itiii$ says that 
the functor $R\mapsto\fg(R)$ is left adjoint to the functor
$\fg\mapsto\fg\pau{z\uppm}$ and that the unit of adjunction
is $Y: R\to\fg(R)\pau{z\uppm}$.
Part $\itiii$ shows that the algebra $\fg(R)$ represents
a certain functor $\Alg\to\Set$ and is thus uniquely
determined up to a unique isomorphism.

\bigskip

\begin{pf}
$\iti$, $\itii$\:
For any $a,b\in R$, the distributions $a(z)$ and $b(z)$
satisfy the commutator formula with coefficients
$((a_{(n)}b)(z))_n$ because we have
\begin{align}
\notag
[a_n,b_m]
\; =\;
&\io((a\otimes t^n) {}_{(0)}(b\otimes t^m))
\\
\notag
=\;
&\io\sum_{i\in\N}(a_{(i)}b)\otimes (\del^{(i)}(t^n)t^m)
\\
\label{E:conformal Lie bracket}
=\;
&\sum_{i\in\N}\binom{n}{i}(a_{(i)}b)_{n+m-i}.
\end{align}
Proposition \ref{SS:locality} implies that
$a(z)$ and $b(z)$ are mutually local of order at most $o_R(a,b)$ and 
that $a(z)_{(n)}b(z)=(a_{(n)}b)(z)$.
Because
$\tT(a\otimes t^n)=(Ta)\otimes t^n+n\, a\otimes t^{n-1}$
we have $(Ta)_n=-n\, a_{n-1}$.
By definition, we have $Ta_n=(Ta)_n$.
Thus $Ta(z)=(Ta)(z)=\del_z a(z)$.

Equation \eqref{E:conformal Lie bracket} implies that
$\fg_{\geq}$ and $\fg_<$ are subalgebras of $\fg$.
Because $\tT(\tR)$ is spanned by the vectors 
$Ta\otimes t^n +n\, a\otimes t^{n-1}$
we have $\fg=\fg_<\oplus\fg_{\geq}$.
The identities $\fg_{\geq}=\bfg_{\geq}$ and $\fg_< =\bfg_<$
are trivial since $\overline{Y(R)}=Y(R)$. 

For a non-negative integer $n$, we have
$a_{-1-n}=(T^{(n)}a)_{-1}$.
This shows that $\ka: R\to\fg_<$ is surjective.
The map $p:\tR\to R, a\otimes t^n\mapsto T^{(-1-n)}(a)$,
induces a map $p:\fg\to R$ because
$p(\tT(a\otimes t^n))=p(T(a)\otimes t^n+n\, a\otimes t^{n-1})=
T^{(-1-n)}(Ta)+n\, T^{(-n)}(a)=0$.
Since $p(a_{-1})=a$ we see that $\ka$ is injective and 
hence so is $Y:R\to\fg\pau{z\uppm}$.

$\itiii$\:
It is clear that if $\phi:\fg(R)\to\fg'$ is an algebra morphism
then $\phi\circ Y$ is an $\N$-fold algebra morphism 
that is compatible with $T$ and such that $(\phi\circ Y)(R)$ is local.
Conversely, if $\phi: R\to\fg'\pau{z\uppm}$ is an 
$\N$-fold algebra morphism that is compatible with $T$ and 
such that $\phi(R)$ is local then
$\phi:a_n\mapsto (\phi a)_n$ is a well-defined linear map because
$\phi: (Ta)_n+n a_{n-1}\mapsto (\del_z \phi a)_n+n(\phi a)_{n-1}=0$.
Lemma \ref{SS:loclie} shows that 
$\phi:\fg(R)\to\fg'$ is an algebra morphism.
${}_{}$
\end{pf}

\bigskip

Let $R$ be bounded $\N$-fold algebra with translation operator $T$.
The local mode algebra $\fg(R)=(\fg(R),Y(R))$ is called
the local mode algebra of $R$.
\index{local mode algebra!of a conformal algebra}

If $\hk$ is a central element of $R$ 
then $\hk:=\hk_{-1}$ is a central element of $\fg(R)$
because $-n\, \hk_{n-1}=(T\hk)_n=0$ implies that 
$\hk_n=0\in\fg(R)$ for any $n\ne -1$ and 
if $\hk$ is non-zero then $\hk\notin\bfg_{\geq}(R)$ 
since $\fg(R)=\fg_<(R)\oplus\fg_{\geq}(R)$.
We have $\fg(R\oplus\K\hk)=\fg(R)\oplus\K\hk$
because $\hk(z)=\hk$.

\subsection{The Local Lie Algebra of a Conformal Algebra}
\label{SS:Lie of ConfA}

We prove that the functor $R\mapsto\fg(R)$ induces
an equivalence between the category of conformal algebras 
and a localization of the category of local Lie algebras.

\bigskip

{\bf Proposition.}\: {\it
The functor $\locLie\to\ConfA, \fg\mapsto R(\fg)$,
has the functor $R\mapsto\fg(R)$ as left adjoint.
The map $Y:R\to R(\fg(R))$ is an isomorphism 
between endofunctors of $\ConfA$. 
}

\bigskip

\begin{pf}
This follows from Proposition \ref{SS:conf affinization},
Remark \ref{SS:mod alg of n fold alg}, and 
Proposition \ref{SS:loc mod alg of n fold alg}.
\end{pf}

\bigskip

Let $\fg$ be a local Lie algebra and $R:=R(\fg)$.
The map $\fg(R)\to\fg$ given by 
$a(z)_n\mapsto a_n$ for $a(z)\in\bF\subset\fg\pau{z\uppm}$ and $n\in\Z$
is well-defined because $\bF$ is translation covariant according to
Proposition \ref{SS:derivative}\,$\itiii$.
Equation \eqref{E:conformal Lie bracket} shows that 
this natural map $\fg(R)\to\fg$ is a morphism of Lie algebras.
It is clear that this map is an epimorphism of local Lie algebras.
Applying the functor $R$ to the morphism $\fg(R)\to\fg$
we obtain an isomorphism $R(\fg(R))\to R(\fg)$ that is the
inverse of $Y$.

\bigskip

{\bf Definition.}\: 
A local Lie algebra $\fg$ is called {\bf regular} 
\index{regular local Lie algebra}
if the natural epimorphism $\fg(R(\fg))\to\fg$ is an isomorphism.

\bigskip

Let $\fg$ be a local Lie algebra and $R:=R(\fg)$.
An ideal $I$ of $\fg$ is called {\bf irregular} \index{irregular ideal}
if for any non-zero $a(z)\in\bF$ there exists an integer $n$ such that
$a_n\notin I$. 
If $\phi:\fg\to\fg'$ is a morphism of local Lie algebras
such that $R\phi$ is a monomorphism then the kernel $I$ of $\phi$
is irregular because if for $a(z)\in\bF$ we have 
$a_n\in I$ for any integer $n$ then $\phi a(z)=0$ and hence
$a(z)=0$.
Conversely,
if $I$ is an irregular ideal of $\fg$ 
and $\phi:\fg\to\fg/I$ is the quotient map then
$R\phi: R(\fg)\to R(\fg/I)$ is an isomorphism.

Thus the functors $R:\locLie\to\ConfA$
and $\fg:\ConfA\to\locLie$  
induce mutually inverse equivalences between
$\ConfA$ and $S^{-1}\locLie$
where $S^{-1}\locLie$
denotes the localization of $\locLie$ with respect to
morphisms $\phi:\fg\to\fg'$ 
whose kernel is irregular and for which $\phi:\bF\to\bF'$ is surjective.

\subsection{Free Conformal Algebras}
\label{SS:free conf}

We prove that for any set $S$ and any function $o:S\times S\to\N$
there exists the {\it free conformal algebra} $R(S,o)$ generated by $S$
such that $o_{R(S,o)}\leq o$.

\bigskip

A {\bf pole order bound} \index{pole order bound}
on a set $S$ is a function $o:S\times S\to\Z$. 
A morphism $\phi: (S,o)\to (T,o)$ of sets with pole order bounds 
is an even map $\phi:S\to T$ of sets such that 
$o(s,t)\geq o(\phi s,\phi t)$ for any $s, t\in S$.
A pole order bound $o$ on a set $S$ is called 
non-negative \index{pole order bound!non-negative}
if $o(S\times S)\subset\N$.

Let $S$ be a set with a non-negative pole order bound $o$.
Let $\fg(S,o)$ be the Lie algebra that is generated by $S\times\Z$
with relations
$$
\sum_{i\in\N}\:
\binom{o(s,t)}{i}\; [(s,n+i),(t,m+o(s,t)-i)]
\; =\; 0
$$
for any $s, t\in S$ and $n, m\in\Z$.
By construction, the $\fg(S,o)$-valued distributions
$s(z):=\sum_{n\in\Z}(s,n)z^{-n-1}$ for $s\in S$ are mutually local.
Let $R(S,o)$ be the conformal subalgebra of $\fg(S,o)\pau{z\uppm}$
that is generated by the distributions $s(z)$ where $s\in S$.

\bigskip

{\bf Proposition.}\: {\it
The forgetful functor $R\mapsto (R,o_R)$ from $\ConfA$ to the
category of sets with non-negative pole order bounds has 
the functor $(S,o)\mapsto R(S,o)$ as its left adjoint.
Moreover, we have $\fg(S,o)=\fg(R(S,o))$.
}

\bigskip

\begin{pf}
Let $S$ be a set with a non-negative pole order bound $o$
and $\fg:=\fg(S,o)$.
The map $s\mapsto s(z)$ defines a morphism 
$(S,o)\to (R(S,o),o_{R(S,o)})$ of sets with pole order bounds
because $s(z)$ and $t(z)$ are mutually local of order at most $o(s,t)$
for any $s, t\in S$. 

Let $R$ be a conformal algebra and 
$\phi:(S,o)\to (R,o_R)$ be a morphism of sets with pole order bounds.
For any $a, b\in R$, 
the $\fg(R)$-valued distributions $a(z)$ and $b(z)$
are mutually local of order $o_R(a,b)$
according to Proposition \ref{SS:loc mod alg of n fold alg}\,$\iti$.
Thus the elements $\phi(s)_n\in\fg(R)$ for $(s,n)\in S\times\Z$
satisfy the relations defining $\fg$.
Hence there exists a unique Lie algebra morphism 
$\fg\to\fg(R)$ such that $(s,n)\mapsto\phi(s)_n$.
This morphism $\fg\to\fg(R)$ induces a morphism 
$\fg\pau{z\uppm}\to\fg(R)\pau{z\uppm}$ of $\N$-fold algebras with
a translation operator such that $s(z)\mapsto (\phi s)(z)$ for any $s\in S$.
Thus we obtain a conformal algebra morphism 
$R(S,o)\to Y(R)$ whose composition with the inverse of 
$Y:R\to Y(R)$ satisfies $s(z)\mapsto\phi s$.

Note that if a conformal algebra $R$ is generated by a subset $T$
then the Lie algebra $\fg(R)$ is generated by the elements 
$t_n$ where $t\in T$ and $n\in\Z$ because
$(a_{(n)}b)_m=(a(z)_{(n)}b(z))_m=\sum_{i\in\N}(-1)^{n+i}\binom{n}{i}
[a_i,b_{m+n-i}]$ and $(Ta)_m=-m\, a_{m-1}$ 
for $a, b\in R$ and $n\in\N, m\in\Z$.

We have shown that there exists a Lie algebra morphism 
$\fg\to\fg(R(S,o))$ such that $(s,n)\mapsto (s(z))_n$.
Conversely, Proposition \ref{SS:loc mod alg of n fold alg}\,$\itiii$
implies that the inclusion $R(S,o)\subset\fg\pau{z\uppm}$ induces a 
Lie algebra morphism $\fg(R(S,o))\to\fg$ such that
$(s(z))_n\mapsto (s,n)$.
Because the Lie algebra $\fg$ is generated by the elements $(s,n)$ and 
the Lie algebra $\fg(R(S,o))$ is generated by the elements $(s(z))_n$
these two morphisms are inverse to each other.
\end{pf}

\bigskip

The conformal algebra $R(S,o)$ is called the
{\bf free conformal algebra} \index{free conformal algebra}
generated by $S$ with pole order bound $o$.

\subsection{Vertex Modules and Lie Algebra Modules}
\label{SS:vertex reps}

We define the notion of a {\it vertex module}
and prove that the categories of vertex modules over a conformal algebra $R$ 
and of modules over the local Lie algebra of $R$ are equivalent.

\bigskip

{\bf Definition.}\:
Let $R$ be an $\N$-fold algebra with a translation operator.
A vector space $M$ together with a map $Y: R\to \cF_z(M)$
is called a {\bf vertex module} \index{vertex module} over $R$ if
$Y$ is a morphism of 
$\N$-fold algebras with a translation operator and $Y(R)$ is local.

\bigskip

A vertex module over a centered conformal algebra $R$ of level $k$
\index{level!of a vertex module}
is by definition a vertex module $M$ over the conformal algebra $R$ 
such that $Y(\hk,z)=k 1(z)$.

\bigskip

{\bf Proposition.}\: {\it 
$\iti$\:
If $M$ is a module over a local Lie algebra $\fg$
then there exists a unique vertex module structure on $M$
over the conformal algebra of $\fg$ such that
$a_{(n)}b=a_n b$ for any $a(z)\in F, b\in M$, and $n\in\Z$.

$\itii$\:
Let $\fg$ be a local Lie algebra 
such that the elements $a_n$ for $a(z)\in F$ and $n\in\Z$
form a basis of $\fg$.
There exists a unique one-to-one correspondence between 
modules $M$ over $\fg$ and 
vertex modules $M$ over the conformal algebra of $\fg$ 
such that $a_n b=a_{(n)}b$ for any $a(z)\in F, b\in M$, and $n\in\Z$.

$\itiii$\:
For a conformal algebra $R$,
there exists a unique one-to-one corres\-pondence between
vertex modules $M$ over $R$ and 
modules $M$ over the local Lie algebra of $R$
such that $a_{(n)}b=a_n b$ for any $a\in R, b\in M$, and $n\in\Z$.
}

\bigskip

\begin{pf}
$\iti$\:
Let $R$ be the conformal algebra of $\fg$.
The representation $\rho:\fg\to\End(M)$ induces a
morphism of $\N$-fold algebras with translation operator
$\fg\pau{z\uppm}\to \End(M)\pau{z\uppm}$ that
restricts to a morphism $Y:R\to\End(M)\pau{z\uppm}$.
Since $R$ is a local subset of $\fg\pau{z\uppm}$
the image $Y(R)$ is local, too.
Because $\rho F$ is contained in $\cF_z(M)$, so is $Y(R)$.

$\itii$\:
Let $Y: R\to\cF_z(M)$ be a vertex module over $R:=R(\fg)$ and
$\rho:\fg\to\End(M)$ be the vector space morphism defined by
$a_n\mapsto a_{(n)}$ for any $a(z)\in F$ and $n\in\Z$.
Because $\bF=R$ and $Y:R\to \cF_z(M)$ is a morphism of 
$\N$-fold algebras with a translation operator
we get $\rho=Y:R\to\cF_z(M)$.
Together with locality of $Y(R)$ this implies by Lemma \ref{SS:loclie} 
that $\rho:\fg\to\End(M)$ is a representation.
Thus the claim follows from $\iti$.

$\itiii$\:
This follows from 
Proposition \ref{SS:Lie of ConfA}\,$\itiii$ and the fact that
for a morphism $Y:R\to\End(M)\pau{z\uppm}$ such that $Y(R)$ is local
we have $Y(R)\subset \cF_z(M)$ if and only if the corresponding 
representation $\fg\to\End(M)$ 
defines a module over the local Lie algebra $\fg$.
\end{pf}

\bigskip

{\bf Corollary.}\: {\it
There exists a one-to-one correspondence between modules $M$ over 
a vertex algebra $V$ and 
those modules $M$ over the centered local Lie algebra of $V$ 
of level $1$ that satisfy 
$(a_{(-1)}b)(z)c=\normord{a(z)b(z)}c$ for any $a, b\in V$ and $c\in M$.
}

\bigskip

\begin{pf}
A module $M$ over the vertex algebra $V$ gives rise to  
a vertex mo\-dule over the centered conformal algebra $V$ and hence 
to a module over $\fg(V)$.
Because $Y:V\to\cF_z(M)$ is a $\Z$-fold algebra morphism we have
$(a_{(-1)}b)(z)c=\normord{a(z)b(z)}c$.

Conversely, if $M$ is a module over $\fg(V)$ satisfying this additional
property then $M$ is a vertex module over the centered conformal algebra $V$.
Lemma \ref{SS:commu vertex} implies that $Y:V\to\cF_z(M)$ is a 
morphism of $\Z$-fold algebras.
\end{pf}

\section{Conformal Algebras and Vertex Algebras}
\label{S:conf and verta}

{\bf Summary.}\:
In section \ref{SS:verma va}
we construct a canonical vertex algebra structure
on certain Verma modules $V(\fg)$.
In section \ref{SS:univ va of g}
we prove that the vertex algebra $V(\fg)$ is characterized by a
universal property.
In section \ref{SS:vertex envelope}
we prove that the forgetful functor $\VertA\to\ConfA$
has a left adjoint $R\mapsto V(R)$.

In section \ref{SS:free va}
we prove the existence of free vertex algebras.
In section \ref{SS:linear opes}
we explain that the vertex algebra that is defined by linear OPEs $O$
is isomorphic to some $V(R(O))$.
In section \ref{SS:linear opes conf alg}
we prove that under certain conditions
$R(O)$ is isomorphic to $\K[T]\otimes\K^{(S)}$.
In section \ref{SS:free va module}
we prove the existence of free vertex algebra modules.

\subsection{Verma Modules and Vertex Algebras}
\label{SS:verma va}

We prove that on certain {\it Verma modules} over a local
Lie algebra there exists a unique vertex algebra structure.

\bigskip

Let $\fg$ be a Lie algebra, $\fg_+$ be a subalgebra of $\fg$,
and $\la:\fg_+\to\K$ be a one-dimensional representation of $\fg_+$;
in other words, $\la$ is a linear form whose kernel is an ideal.
The induced $\fg$-module 
$$
V^{\la}(\fg)
\; :=\;
U(\fg)\otimes_{U(\fg_+)}\K
\; =\;
U(\fg)/U(\fg)\set{a-\la(a)\mid a\in\fg_+}
$$
is called the {\bf Verma module} \index{Verma module}
over $\fg$ of level $\la$.
The image in $V^{\la}(\fg)$ of the identity $1$ of $U(\fg)$ is 
called the {\bf highest weight vector} \index{highest weight vector}
of $V^{\la}(\fg)$ and is denoted by $1$.

The Poincar\'e-Birkhoff-Witt Theorem implies that if
$(a_i)_{i\in I}$ is a basis of $\fg$ and $I_+$ is a subset of $I$
such that $(a_i)_{i\in I_+}$ is a basis of $\fg_+$ then 
the monomials $a_{i(1)}^{n_1}\dots a_{i(r)}^{n_r}$ 
form a basis of $V^{\la}(\fg)$ where 
$r$ and $n_i$ are non-negative integers and 
$i:\set{1, \dots, r}\to I\setminus I_+$ 
is an increasing function with respect to some chosen total order on $I$.
In particular, if $\hk$ is an element of $\fg_+$ such that
$\la(\hk)$ is non-zero and 
$\fg'$ is a vector space complement of $\fg_+$ in $\fg$
then the map $\fg'\oplus\K\hk\to V^{\la}(\fg), a\mapsto a1$,
is a monomorphism.

\bigskip

{\bf Proposition.}\: {\it
Let $\fg$ be a local Lie algebra and $\fg_+$ be a subalgebra of $\fg$
that contains $\fg_{\geq}$ and that is invariant with respect to $T$.
If $\la:\fg_+\to\K$ is a one-dimensional representation of $\fg_+$
such that $\la(T\fg_+)=0$ then 
there exists a unique vertex algebra structure 
on $V^{\la}(\fg)$ such that 
the highest weight vector is a right identity and
$Y(a_{-1}1,z)=\rho a(z)$ for any $a(z)\in F$.
}

\bigskip

\begin{pf}
Define $S:=\rho F$.
The derivation $T$ of $\fg$ induces a derivation $T$ of $U(\fg)$
such that $T1=0$.
For $a\in\fg_+$, we have $T(a-\la(a))=Ta\in\set{b-\la(b)\mid b\in\fg_+}$
because $\la(T\fg_+)=0$.
Thus $T$ induces an operator $T$ of $V^{\la}(\fg)$
such that $1$ is invariant for $T$. The set $S$ is creative 
because $\fg_{\geq}=T\fg_{\geq}\subset T\fg_+$
and $\la(T\fg_+)=0$ implies that $a_n 1=0$ for any $a(z)\in F$ and $n\in\N$.
The set $S$ is translation covariant because 
$F$ is translation covariant and because 
$T$ being a derivation implies that
$[T,\rho a]b=\rho(Ta)b$ for any $a\in\fg$ and $b\in V^{\la}(\fg)$.
It is clear that $S$ is generating and local.
Thus the claim follows from the existence theorem for OPE-algebras.
\end{pf}

\subsection{Universal Vertex Algebra over a Local Lie Algebra}
\label{SS:univ va of g}

We define the notion of a vertex algebra over a local
Lie algebra $\fg$ and prove that the Verma module $V(\fg)$
is the {\it universal} vertex algebra over $\fg$.

\bigskip

{\bf Definition.}\:
For a local Lie algebra $\fg$, 
a vertex algebra $V$ together with a $\fg$-module structure
is called a vertex algebra {\bf over} $\fg$ 
\index{vertex algebra!over a local Lie algebra}
if $\rho F$ is contained in the space of fields of $V$.

\bigskip

Let $\fg$ be a local Lie algebra.
A morphism $V\to W$ of vertex algebras over $\fg$
is a morphism $V\to W$ of vertex algebras that is also
a morphism of $\fg$-modules.
A vertex algebra $V$ over $\fg$ is called {\bf universal}
if for any vertex algebra $W$ over $\fg$ there
exists a unique morphism $V\to W$ of vertex algebras over $\fg$.

A vertex algebra together with a $\fg$-module structure
is a vertex algebra over $\fg$ if and only if
$Y(a_{-1}1,z)=\rho a(z)$ for any $a(z)\in F$.

Let $\fg$ be a centered local Lie algebra 
and $\la:\bfg_{\geq}\oplus\K\hk\to\K$
be a one-dimensional representation such that $\la(\bfg_{\geq})=0$;
thus $\la$ is determined by $k:=\la(\hk)$.
We denote by $V^k(\fg)$
the Verma module over $\fg$ of level $\la$.
Proposition \ref{SS:univ va of g} shows that there exists a unique
vertex algebra structure on $V^k(\fg)$ such that
the highest weight vector is a right identity and
$V^k(\fg)$ is a vertex algebra over $\fg$ of level $k$.

If $\fg$ is a local Lie algebra then 
$V(\fg):=V^0(\fg\oplus\K\hk)$ is a vertex algebra over $\fg$.

\bigskip

{\bf Proposition.}\: {\it
$\iti$\:
If $\fg$ is a centered local Lie algebra of level $k$
then $V^k(\fg)$ is the universal vertex algebra over $\fg$ of level $k$.

$\itii$\:
If $\fg$ is a local Lie algebra 
then $V(\fg)$ is the universal vertex algebra over $\fg$.
}

\bigskip

\begin{pf}
$\iti$\:
Let $W$ be a vertex algebra over $\fg$ of level $k$.
We have $\hk 1=k1$ and $\bfg_{\geq}1=0$
because $\rho\bF\subset\cF_Y$.
The universal property of the Verma module $V^k(\fg)$ implies that
there exists a unique morphism of $\fg$-modules $\psi:V^k(\fg)\to W$
such that $1\mapsto 1$. 
The canonical filtration of $U(\fg)$, that is induced from
the gradation of the tensor algebra of $\fg$, induces an
exhaustive increasing filtration $(V_l)_{l\in\N}$ of $V^k(\fg)$.
By induction on $l$, we will prove that
$\psi(a_{(n)}b)=\psi(a)_{(n)}\psi(b)$ for any
$a\in V_l, b\in V^k(\fg)$, and $n\in\Z$;
this shows that $\psi$ is a morphism of vertex algebras
and hence $V^k(\fg)$ is universal.

Remark \ref{SS:hol jacobi}\,$\itv$ implies that
for any integers $r, t, n$, and $m$ there exist
integers $N_i$ such that 
for any states $a, b, c$ of a vertex algebra 
such that $o'(a,b)\geq r$ and $o'(a,c)\geq t$ we have
$(a_{(n)}b)_{(m)}c=\sum_{i\in\N} N_i\, a_{(r+t-i)}b_{(n+m-r-t+i)}c$.
Thus if $a(z)\in F, b\in V_l, c\in V^k(\fg)$, and $n, m\in\Z$ 
then there exist integers $N$ and $N_i$ such that 
by induction and using that $\psi$ is a $\fg$-module morphism we get
\begin{align}
\notag
\psi((a_n b)_{(m)}c)
\; &=\;
\psi(((a_{-1}1)_{(n)}b)_{(m)}c)
\\
\notag
\; &=\;
\sum_{i\in\N} N_i\, \psi((a_{-1}1)_{(N-i)}b_{(n+m-N+i)}c)
\\
\notag
\; &=\;
\sum_{i\in\N} N_i\, \psi(a_{-1}1)_{(N-i)}\psi(b)_{(n+m-N+i)}\psi(c)
\\
\notag
\; &=\;
(\psi(a_{-1}1)_{(n)}\psi(b))_{(m)}\psi(c)
\\
\notag
\; &=\;
\psi(a_n b)_{(m)}\psi(c)
\end{align}
because $\psi(a_{-1}1)_{(n)}\psi(b)=(a_{-1}1)_{(n)}\psi(b)=
a_n\psi(b)=\psi(a_n b)$.

$\itii$\:
This follows from $\iti$. 
\end{pf}

\subsection{Enveloping Vertex Algebras}
\label{SS:vertex envelope}

We prove that the forgetful functor from vertex algebras to
conformal algebras has a left adjoint $R\mapsto V(R)$.
The vertex algebra $V(R)$ is called the {\it enveloping vertex algebra}
of $R$.

\bigskip

For a centered conformal algebra $R$ of level $k$,
the vertex algebra $V^k(\fg(R))$ is called the
{\bf enveloping vertex algebra} of $R$ and is denoted by $V^k(R)$.

For a conformal algebra $R$,
the vertex algebra $V(R):=V^0(R\oplus\K\hk)$ is called the
{\bf enveloping vertex algebra} of $R$.

\bigskip

{\bf Remark.}\: {\it
$\iti$\:
For a centered conformal algebra $R$ of level $k$,
the map $\io:R\to f_k V^k(R), a\mapsto a_{-1}1$, is a monomorphism
of centered conformal algebras. 

$\itii$\:
For a conformal algebra $R$,
the map $\io:R\to V(R), a\mapsto a_{-1}1$, is a monomorphism
of conformal algebras. 
}

\bigskip

\begin{pf}
$\iti$\:
Let $\fg:=\fg(R)$. The map $\io$ is a morphism
of conformal algebras because it is the composition of 
$Y:R\to\fg\pau{z\uppm}$, the morphism 
$\fg\pau{z\uppm}\to\cF_Y$ induced by the representation
$\fg\to\End(V^k(R))$, and the field-state correspondence $\cF_Y\to V^k(R)$. 
From $\fg=\fg_<\oplus\fg_{\geq}$ follows that the map
$\fg_<\to V^k(\fg), a\mapsto a1$, is injective.
By Proposition \ref{SS:loc mod alg of n fold alg}\,$\itii$ 
the map $R\to\fg_<, a\mapsto a_{-1}$, is a vector space
isomorphism.
Thus $\io:R\to V^k(R)$ is injective.
Finally, we have $\hk_{-1}1=\hk 1=k1$.

$\itii$\:
This follows from $\iti$.
\end{pf}

\bigskip

{\bf Proposition.}\: {\it
$\iti$\:
For a centered conformal algebra $R$ of level $k$, there exists a unique 
one-to-one correspondence between vertex modules $M$ over $R$ of level $k$
and modules $M$ over $V^k(R)$ such that $\io a_{(n)}b=a_{(n)}b$
for any $a\in R, b\in M$, and $n\in\Z$.

$\itii$\:
For a conformal algebra $R$, there exists a unique 
one-to-one correspondence between vertex modules $M$ over $R$ 
and modules $M$ over $V(R)$ such that $\io a_{(n)}b=a_{(n)}b$
for any $a\in R, b\in M$, and $n\in\Z$.
}

\bigskip

\begin{pf}
$\iti$\:
Let $Y:R\to\cF_z(M)$ be a vertex module over $R$.
By Remark \ref{SS:modules} the image $Y(R)$ generates a vertex subalgebra 
$V:=\sqbrack{Y(R)\cup\set{1(z)}}$ of $\cF_z(M)$.
The vector space $V$ together with the morphism 
$Y:R\to V\to\cF_z(V)$ is a vertex module over $R$
and hence a $\fg(R)$-module. 
By definition, $V$ is a vertex algebra over $\fg(R)$.
Proposition \ref{SS:univ va of g} implies 
that there exists a unique morphism of vertex algebras
$\psi:V(R)\to V$.
The composition of $\psi:V(R)\to V$ and the inclusion of $V$ into
$\cF_z(M)$ defines on $M$ a module structure over $V(R)$
such that $\io a_{(n)}b=a_{(n)}b$
for any $a\in R, b\in M$, and $n\in\Z$.

Conversely, if $M$ is a module over $V(R)$ then
the composition of $\io:R\to V(R)$ and 
$Y:V(R)\to\cF_z(M)$ defines on $M$ the structure of a vertex
module such that $\io a_{(n)}b=a_{(n)}b$.

$\itii$\:
This follows from $\iti$.
\end{pf}

\bigskip

{\bf Corollary.}\: {\it
$\iti$\:
For any $k\in\K$,
the functor $\ConfA_k\to\VertA, R\mapsto V^k(R)$, is 
the left adjoint of the forgetful functor $f_k$
from the category of non-zero vertex algebras to $\ConfA_k$.

$\itii$\:
The functor $\ConfA\to\VertA, R\mapsto V(R)$, is 
the left adjoint of the forgetful functor $\VertA\to\ConfA$.
}

\bigskip

\begin{pf}
$\iti$\:
Let $V$ be a vertex algebra.
The fact that $Y:V\to\cF_Y$ is an isomorphism shows that
to give a morphism $R\to f_k V$ of central conformal algebras
is equivalent to giving
a vertex module structure $Y:R\to\cF_z(V)$ on $V$ of level $k$
such that $Y(R)\subset\cF_Y$.
From the Proposition and the fact that $\cF_Y$ is a vertex algebra
follows that
to give a vertex module structure $Y:R\to\cF_z(V)$ on $V$ of level $k$
such that $Y(R)\subset\cF_Y$ is equivalent to giving
a module structure $Y:V^k(R)\to\cF_z(V)$ on $V$ such that 
$Y(V^k(R))\subset\cF_Y$, in other words, it is equivalent to giving
a morphism $V^k(R)\to V$ of vertex algebras.

$\itii$\:
This follows from $\iti$.
\end{pf}

\subsection{Free Vertex Algebras}
\label{SS:free va}

We prove that for any set $S$ and any function $o:S\times S\to\Z$
there exists the {\it free vertex algebra} $V(S,o)$ generated by $S$
such that $(o':V\times V\to\Z\cup\set{-\infty})\leq o$.

\bigskip

{\bf Proposition.}\: {\it
The forgetful functor $V\mapsto (V,o')$ from $\VertA$
to the category of sets with pole order bounds has a left
adjoint $(S,o)\mapsto V(S,o)$.
}

\bigskip

\begin{pf}
Let $S$ be a set with a pole order bound $o$.
Let $V(S,o)$ be the quotient of $V(R(S,o_+))$
by the ideal generated by $s_{(n)}t$
where $s, t\in S$ and $n\geq o(s,t)$.

Let $V$ be a vertex algebra.
To give a map $\phi:(S,o)\to (V,o')$ of sets with pole order bounds
is equivalent by Proposition \ref{SS:free conf} to giving a morphism 
$\phi:R(S,o_+)\to V$ of conformal algebras
such that $\phi s_{(n)}\phi t=0$
for any $s, t\in S$ and $n\geq o(s,t)$
which by Corollary \ref{SS:vertex envelope} is equivalent to giving a morphism
$\phi:V(S,o)\to V$ of vertex algebras.
\end{pf}

\bigskip

The vertex algebra $V(S,o)$ is called the
{\bf free vertex algebra} \index{free vertex algebra}
generated by $S$ with pole order bound $o$.

Let $S$ be a set with a pole order bound $o$.
We denote by $V(S)$ the free $\Z$-fold algebra
with operator generated by $S\sqcup\set{1}$.
We denote by $R(S)$ the free $\N$-fold algebra
with operator generated by $S$.
There exist canonical morphisms $\io: V(S)\to V(S,o)$ and 
$\io: R(S)\to R(S,o)$ and an inclusion $R(S)\subset V(S)$.

If $\cR$ is a subset of $V(S)$ then the quotient
of $V(S,o)$ by the ideal generated by $\io\cR$ is called
the vertex algebra {\bf generated} by $S$ with pole order bound $o$
and {\bf relations} $\cR$ and is denoted by $V(S,\cR)$.

Assume that $o$ is non-negative. 
If $\cR$ is a subset of $R(S)$ then the quotient
of $R(S,o)$ by the ideal generated by $\io\cR$ is called
the conformal algebra generated by $S$ with pole order bound $o$
and relations $\cR$ and is denoted by $R(S,\cR)$.
The central conformal algebra 
generated by $S$ with pole order bound $o$
and relations $\cR\subset R(S)\oplus\K\hk$
is defined in the same way.

\bigskip

{\bf Remark.}\: {\it
The vertex algebra $V(S,\cR)$ generated by a set $S$ with 
non-negative pole order bound $o$ and relations 
$\cR\subset R(S)$ is canonically isomorphic to
$V(R(S,\cR))$.
}

\bigskip

\begin{pf}
This is proven in the same way as the Proposition.
\end{pf}

\subsection{Vertex Algebras defined by Linear OPEs}
\label{SS:linear opes}

We define the notion of a family of OPEs indexed by a set $S$
and prove that a vertex algebra that is 
{\it defined by linear OPEs} is the enveloping vertex algebra of a 
conformal algebra that is a quotient of the $\N$-fold algebra
$\K[T]\otimes\K^{(S)}\oplus\K\hk$.

\bigskip

For a set $S$,
a map $O:S\times S\times\N\to V(S)$ is called a 
{\bf family of OPEs} \index{family of OPEs}
indexed by $S$ if there exists a pole order bound $o$ on $S$
such that $O(s,t,n)$ is zero for any $s, t\in S$ and 
$n\geq o(s,t)$; in this case we say that $o$ bounds $O$.

Let $O$ be a family of OPEs indexed by a set $S$ and 
$o$ be a non-negative pole order bound that bounds $O$.
The vertex algebra that is generated by $S$
with pole order bound $o$ and relations $s_{(n)}t=O(s,t,n)$, 
where $s, t\in S$ and $n\in\N$,
is called the vertex algebra {\bf defined by the OPEs} $O$
and is denoted by $V(S,O)$. 
Its universal property shows that $V(S,O)$ 
does not depend on the choice of $o$.

The family $O$ of OPEs is called a family of {\bf linear OPEs}
if the image of $O$ is contained in $\K[T]S\oplus\K 1$.
We identify $\K[T]S\oplus\K 1$ with a subset of $R(S)\oplus\K\hk$
by mapping $1$ to $\hk$.

Assume that $O$ is a family of linear OPEs.
The centered conformal algebra that is generated by $S$
with pole order bound $o$ and relations $s_{(n)}t=O(s,t,n)$, 
where $s, t\in S$ and $n\in\N$,
is called the conformal algebra {\bf defined by the linear OPEs} $O$
and is denoted by $R(S,O)$. 
By Remark \ref{SS:free va} we have $V(S,O)=V^1(R(S,O))$.

Let $R'(S,O)$ be the $\N$-fold algebra with 
a translation operator and a central element 
that is generated by $S$ with relations 
$s_{(n)}t=O(s,t,n)$ for any $s, t\in S$ and $n\in\N$.
There exists a unique morphism
$R'(S,O)\to R(S,O)$ of $\N$-fold algebras 
with a translation operator and a central element,
that is compatible with $S$.
This morphism is surjective because $R(S,O)$ is generated by $S$
as an $\N$-fold algebra with a translation operator and a central element.

\bigskip

{\bf Remark.}\: {\it
Let $O$ be a family of linear OPEs indexed by a set $S$.
Define an operator $T$ on $\K[T]\otimes\K^{(S)}\oplus\K\hk$ by 
$T(T^m\otimes s):=T^{m+1}\otimes s$ and $T\hk:=0$.
The unique vector space morphism 
$\ka:\K[T]\otimes\K^{(S)}\oplus\K\hk\to R'(S,O)$
that is compatible with $T, \hk$, and $S$ is an isomorphism.
Moreover, $R'(S,O)$ is bounded.
}

\bigskip

\begin{pf}
Let $s, t\in S$ and $n, m, m'\in\N$.
We define by induction on $m'$ 
an $\N$-fold algebra structure with central element $\hk$
on $\K[T]\otimes\K^{(S)}\oplus\K\hk$ by
\begin{equation}
\label{E:TTT}
T^m s_{(n)}T^{m'}t
=
\begin{cases}
\;
(-1)^m m! \binom{n}{m}\, O(s,t,n-m)
\quad \text{if $m'=0$ and} \\
\;
T(T^m s_{(n)}T^{m'-1}t) 
-
T^{m+1} s_{(n)}T^{m'-1}t
\quad \text{if $m'>0$}
\end{cases}
\end{equation}
where $T^n s:=T^n\otimes s$.
Using induction on $m'$ we show that $T$ is a translation operator
for $T^m s$ and $T^{m'} t$.
For $m'=0$, $T$ is a translation generator by definition and 
$T$ is a derivation because of the definition of 
$T^m s_{(n)}T^{m'}t$ in the case that $m'=1$.
For $m'>0$, it follows directly from the definition of 
$T^m s_{(n)}T^{m'}t$ and from the induction hypothesis. 
Thus the universal property of $R'(S,O)$ implies that
there exists a unique morphism 
$\ka': R'(S,O)\to \K[T]\otimes\K^{(S)}\oplus\K\hk$
of $\N$-fold algebras with a translation operator and a central element
that is compatible with $S$.

Conversely, 
from \eqref{E:TTT} follows by induction on $m'$ that $\ka$ is a morphism
of $\N$-fold algebras.
Because $\K[T]\otimes\K^{(S)}\oplus\K\hk$ and $R'(S,O)$ 
are $\N$-fold algebras with a translation operator and a central element, 
that are generated by $S$, 
the morphisms $\ka$ and $\ka'$ are mutually inverse.

Equation \eqref{E:TTT} shows that $R'(S,O)$ is bounded.
\end{pf}

\bigskip

{\bf Corollary.}\: {\it
The superaffinization of a Lie algebra $\fg$ with an invariant
symmetric pairing is a regular local Lie algebra.
The same is true for the affinization and the Clifford affinization.
}

\bigskip

\begin{pf}
By \eqref{E:lie opes} the OPEs of the distributions in 
$F\subset\hfg^s\pau{z\uppm}$ are linear. 
Thus there exists a natural epimorphism 
$p:R:=\K[T]\otimes (\fg\oplus\fg\theta)\oplus\K\hk\to\bF$ by the Remark.
This morphism is injective because the $(-1)$-st modes of 
$p(T^{(n)}\otimes a)$ and $p(T^{(n)}\otimes a\theta)$ are
$a_{-1-n}$ and $\ba_{-1/2-n}$, resp.
It follows directly that the morphisms
$\fg(R)\to\hfg^s, (p(T^n\otimes a))_m\mapsto T^n(a_m)$, and
$\hfg^s\to\fg(R), a_m\mapsto (p(1\otimes a))_m$, are well-defined and
inverse to each other, where $a\in\fg\oplus\fg\theta$ and $n, m\in\Z$.
\end{pf}

\subsection{Conformal Algebras defined by Linear OPEs}
\label{SS:linear opes conf alg}

We prove that if the elements of $S$ satisfy inside $R'(S,O)$
the Jacobi identity and skew-symmetry 
then $R'(S,O)=R(S,O)$.

\bigskip

{\bf Lemma.}\: {\it
Let $R$ be an $\N$-fold algebra with a translation operator $T$
and $a, b, c$ be elements of $R$.

\begin{enumerate}
\item[$\iti$]
If $a, b, c$ satisfy the Jacobi identity 
then $Ta, b, c$ and $a, Tb, c$ and $a, b, Tc$ satisfy
the Jacobi identity as well.

\item[$\itii$]
If $a$ and $b$ satisfy skew-symmetry then
$Ta, b$ and $a, Tb$ satisfy skew-symmetry as well.
\end{enumerate}
}

\bigskip

\begin{pf}
\iti\:
For $a\in R$, we define
$a(z):=\sum_{n\in\N}a_{(n)}z^{-n-1}\in\End(R)\pau{z\uppm}$.
The elements $a, b, c$ satisfy the associativity formula 
if and only if $(a_{(r)}b)(z)c=(a(z)_{(r)}b(z))_+c$ for any 
non-negative integer $r$.
Because $T$ is a translation ope\-rator and 
$\del_z(a(z)_+)=(\del_z a(z))_+$
the associativity formula for $a, b, c$ implies that
\begin{align}
\notag
(Ta_{(r)}b)(z)c
\; =\;
-r\, (a_{(r-1)}b)(z)c
\; &=\;
-r\, (a(z)_{(r-1)}b(z))_+c
\\
\notag
&=\;
(\del_z a(z)_{(r)}b(z))_+c
\\
\notag
&=\;
((Ta)(z)_{(r)}b(z))_+c
\end{align}
and
\begin{align}
\notag
(Ta_{(r)}b+a_{(r)}Tb)(z)c
\; &=\;
(T(a_{(r)}b))(z)c
\\
\notag
&=\;
\del_z (a_{(r)}b)(z)c
\\
\notag
&=\;
\del_z (a(z)_{(r)}b(z))_+c
\\
\notag
&=\;
(\del_z a(z)_{(r)}b(z))_+c \: +\: (a(z)_{(r)}\del_z b(z))_+c
\\
\notag
&=\;
((Ta)(z)_{(r)}b(z))_+c \: +\: (a(z)_{(r)}(Tb)(z))_+c.
\end{align}
Thus the associativity formula and hence by 
Remark \ref{SS:conformal algebras} the Jacobi identity hold for $Ta, b, c$
and $a, Tb, c$. 
Because $T$ is a derivation by applying $T$ to 
the Jacobi identity for $a, b, c$ we see that
the Jacobi identity is also satisfied for $a, b$, and $Tc$.

\itii\:
Because $T$ is a translation operator skew-symmetry for $a$ and $b$ 
implies that 
\begin{align}
\notag
b_{(n)}Ta 
\; &=\;
T(b_{(n)}a) \: +\: n\, b_{(n-1)}a
\\
\notag
&= \;
\paraab \sum_{i\in\N} \:
\big(
(-1)^{n+1+i}\, (i+1)\, T^{(i+1)}(a_{(n+i)}b)
\\
\notag
&\qquad\qquad\qquad\qquad\qquad
\: +\:
(-1)^{n+i}\, n\, T^{(i)}(a_{(n-1+i)}b)
\big)
\\
\notag
&=\;
\paraab\sum_{i\in\N}\:
(-1)^{n+1+i}\, T^{(i)}(Ta_{(n+i)}b)
\end{align}
since $T^{(i)}(Ta_{(n+i)}b)=-(n+i)\, T^{(i)}(a_{(n+i-1)}b)$.
Because $T$ is a derivation 
skew-symmetry for $a, b$ and $Ta, b$ implies skew-symmetry for $a, Tb$.
\end{pf}

\bigskip

{\bf Proposition.}\: {\it
Let $O$ be a family of linear OPEs indexed by a set $S$.
The canonical epimorphism $R'(S,O)\to R(S,O)$ is an isomorphism
if and only if 
the elements of $S\subset R'(S,O)$ satisfy 
the Jacobi identity and skew-symmetry.
}

\bigskip

\begin{pf}
The Lemma shows that if the elements of $S$ satisfy 
the Jacobi identity and skew-symmetry then $R'(S,O)$ is a 
centered conformal algebra.
Thus there exists a unique morphism of centered conformal algebras
$\ka: R(S,O)\to R'(S,O)$ that is compatible with $S$.
Because $R(S,O)$ and $R'(S,O)$ are both
centered conformal algebras that are generated by $S$
the canonical morphism $R'(S,O)\to R(S,O)$ is the inverse of $\ka$.
The converse is obvious.
\end{pf}

\subsection{Free Modules over Vertex Algebras}
\label{SS:free va module}

We prove that for any vertex algebra $V$ and any function $o:V\to\Z$
there exists the {\it free module} $M(V,o)$ over $V$.

\bigskip

{\bf Proposition.}\: {\it
Let $V$ be a vertex algebra and $o:V\to\Z$ be a map.
There exists a unique module $M(V,o)$ over $V$ together with
an element $1_o$ such that $a_{(n)} 1_o=0$ for any $a\in V$ and $n\geq o(a)$
and such that $M(V,o)$ and $1_o$ are universal with this property.
}

\bigskip

\begin{pf}
Let $\fg$ be the local Lie algebra of $V$ and 
$M$ be the quotient of $U(\fg)$ by the left ideal generated by
the elements $a_n$ where $a\in V$ and $n\geq o(a)$.
Thus for any $a\in V$ the distribution $a(z)$ is bounded on the
image $1_o\in M$ of $1\in U(\fg)$.
Remark \ref{SS:loclie} implies that $M$ is a module over the 
local Lie algebra $\fg$.
Let $M(V,o)$ be the quotient of $M$ by the submodule that is 
generated by the elements $1_n c-\de_{n,-1}c$ and
$(a_{(-1)}b)_n c-
\sum_{i\in\N}(a_{-1-i}b_{n+i}+\paraab\, b_{n-1-i}a_i)c$
where $a, b\in V, c\in M$, and $n\in\Z$.
Corollary \ref{SS:vertex reps} shows that $M(V,o)$ is a module over $V$. 
Let $1_o\in M(V,o)$ be the image of $1_o\in M$.
The module $M(V,o)$ together with the element $1_o$ clearly satisfies 
the universal property as stated.
\end{pf}

\bigskip

{\bf Lemma.}\: {\it
Let $a(z)$ and $b(z)$ be holomorphic fields on a vector space $V$,
$c\in V$, and $N, N_a, N_b, n\in\Z$.
If $a(z)$ and $b(z)$ are mutually local of order $N$,
$a(z)c\in z^{-N_a}V\pau{z}$, and $b(z)c\in z^{-N_b}V\pau{z}$
then $a(z)_{(n)}b(z)c\in z^{-M}V\pau{z}$ where 
$M:=N+N_a+N_b-n-1$.
}

\bigskip

\begin{pf}
Because the distributions $(z-w)^N a(z)b(w)c$ and 
$\paraab (z-w)^N b(w)a(z)c$ are equal they are contained in
$z^{-N_a}w^{-N_b}V\pau{z,w}$.
Thus we get
\begin{align}
\notag
a(w)_{(n)}b(w)
\; &=\;
\del_z^{(N-n-1)}((z-w)^N a(z)b(w)c)|_{z=w}
\\
\notag
&\in\;
\del_z^{(N-n-1)}(z^{-N_a}w^{-N_b}V\pau{z,w})|_{z=w}
\\
\notag
&\subset\;
z^{-N_a-N+n+1-N_b}V\pau{z}.
\end{align}
\end{pf}

\bigskip

{\bf Corollary.}\: {\it
Let $V$ be a vertex algebra that is generated by a subset $S$
and $o:S\to\Z$ be a map.
There exists a unique module $M$ over $V$ with an element $1_o$
such that $s_{(n)}1_o=0$ for any $s\in S$ and $n\geq o(s)$ and
such that $M$ and $1_o$ are universal with this property.
}

\bigskip

\begin{pf}
Define a map $o_V:V\times V\to\Z$ where
$o_V(a,b)$ is the pole order of $a(z)b(w)$ if $a(z)b(w)$ is non-zero
and $o_V(a,b)$ is zero otherwise.
Let $\cS$ be the set of sequences $\tis:=(s_1,n_1, \dots, s_r, n_r, s_{r+1})$
where $s_i\in S, n_i\in\Z$, and $r\geq -1$. 
Define a map $\io: \cS\to V$ by 
$\tis\mapsto s_{1(n_1)}(\dots (s_{r(n_r)}s_{r+1})\dots )$.
By induction, define a map $o:\cS\to\Z$ by 
$o(\emptyset):=0, o((s)):=o(s)$, and 
$$
o(s,n,s_1,n_1, \dots, s_{r+1})
\; :=\;
o(s)+o(\tis)+o_V(\io(s),\io(\tis))-n-1.
$$
Remark \ref{SS:hol jacobi}\,$\itv$ 
shows that the vector space $V$ is spanned by $\io(\cS)$. 
Choose a subset $B$ of $\cS$ such that $\io(B)$ is a basis of $V$
and $\io|_B$ is injective.
Define a map $o:V\to\Z$ by
$o(\sum_i \la_i \io(\tis_i)):=\max_i(o(\tis_i))$ where 
$\la_i\in\K$ and $\tis_i\in B$.
From the Lemma follows that $M:=M(V,o)$ satisfies the
universal property as asserted.
\end{pf}

\appendix
\chapter{Superalgebra}
\label{SS:superalgebra}

We give some basic definitions of superalgebra.
We introduce the notion of a superset because
we will consider super objects generated by supersets. 
We define a morphism of super vector spaces, 
in contrast e.g.~to \cite{deligne.morgan.susy},
as a morphism of the underlying vector spaces.
This has the advantage that left multiplication 
with any element of a superalgebra is 
an endomorphism of the underlying super vector space.
Because we thus do not restrict to even morphisms
it follows that our categories are actually 
supercategories, i.e.~the morphism sets are actually 
supersets.

\bigskip

A {\bf superset} is a set $S$ together with a decomposition
$S=S\even\sqcup S\odd$ indexed by $\Z/2=\set{\bar{0},\bar{1}}$.
Thus a superset is given by a set $S$ together with a map
$S\to\Z/2, a\mapsto\ta$.
We call $\ta$ the {\bf parity} of $a$
and we call $a$ {\bf even}, resp., {\bf odd},
if $\ta=\bar{0}$, resp., $\ta=\bar{1}$.
A morphism of supersets $S\to T$
is a morphism of sets $S\to T$ such that
$S_i\to T_i$ for $i\in\Z/2$.
We denote by $\Set$ the category of supersets.
\index{superset}
\index{parity}
\index{even}
\index{odd}

A {\bf super vector space} is a vector space $V$
together with a decomposition $V=V\even\oplus V\odd$.
A morphism $V\to W$ of super vector spaces
is a morphism $V\to W$ of vector spaces. 
By definition, a morphism $V\to W$ of super vector spaces 
has parity $p\in\Z/2$ if $V_i\to W_{i+p}$.
Thus the vector space $\Vect(V,W)$
of super vector space morphisms $V\to W$ 
has a natural structure of a super vector space.
\index{super vector space}

The category $\Vect\even$ of super vector spaces and even morphisms
together with the tensor product $V\otimes W$ defined by
\begin{equation}
\label{E:tensor super vs}
(V\otimes W)_i
\; :=\;
\bigoplus_{j\in\Z/2}\: V_j\otimes W_{i-j}
\end{equation}
is a closed symmetric monoidal category whose inner Homs 
are given by $\Vect(V,W)$;
see e.g.~\cite{maclane.working}, chapter VII.7,
for the definition of the notion of a closed symmetric monoidal category.
A {\bf linear supercategory} \index{linear supercategory}
is a category enriched in $\Vect\even$;
see e.g.~\cite{maclane.working}, chapter VII.7,
for the definition of the notion of an enriched category.
We denote by $\Vect$ the linear supercategory of super vector spaces.
The tensor product \eqref{E:tensor super vs}
makes $\Vect$ into a {\it monoidal supercategory}.
\index{monoidal supercategory}

A super vector space $A$ together with an even morphism
$A\otimes A\to A$ is called a {\bf superalgebra}.
We denote by $\Alg$ \index{Alg@$\Alg$}
the category of superalgebras.
For a super vector space $V$,
$\End_{\Vect}(V)$ is an associative 
unital superalgebra.
We call a superalgebra $A$ {\bf commutative} 
if $A$ is associative and unital and the diagram
$$
\xymatrix{
A\otimes A \ar[r] \ar[d]_{\tau} & A \ar[d]^{id} \\
A\otimes A \ar[r]               & A
}
$$
commutes where $\tau$ maps $a\otimes b$ to
$(-1)^{\ta\tb}\, b\otimes a$.
We call a superalgebra $A$ {\bf skew-symmetric}
if the above diagram commutes with $id$ replaced by $-id$.
\index{superalgebra}
\index{commutative superalgebra}
\index{superalgebra!commutative}
\index{skew-symmetric superalgebra}
\index{superalgebra!skew-symmetric}

For an element $a$ of a unital superalgebra $A$ and $n\in\K$, 
we have 
$\binom{a}{n}=(-1)^n\binom{-a-1+n}{n}$
and
$\binom{a+1}{n}=\binom{a}{n-1}+\binom{a}{n}$.
For non-negative integers $n$ and $m$, we have
$\sum_{i=n}^m\binom{i}{n}=\binom{m+1}{n+1}$.
If $a$ and $b$ are commuting even elements
of $A$
then for any $n\in\K$ we have
$\binom{a+b}{n}=\sum_{i\in\N}\binom{a}{i}\binom{b}{n-i}$
and
$(a+b)^{(n)}=\sum_{i\in\N}a^{(i)}b^{(n-i)}$.

A superalgebra $(\fg,[\; ,\; ])$ is called a
left {\bf super Leibniz algebra} 
if the left {\bf Leibniz identity} 
$$
[[a,b],c]
\; =\;
[a,[b,c]]
\: -\:
(-1)^{\ta\tb}\,
[b,[a,c]]
$$
is satisfied for any elements $a, b, c$ of $\fg$.
Note that $(\fg,a\otimes b\mapsto -(-1)^{\ta\tb}[b,a])$ is then
a right super Leibniz algebra.
We call a left super Leibniz algebra just a
super Leibniz algebra.
A superalgebra $\fg$ is a super Leibniz algebra if and only if 
$[a,\;\, ]$ is a derivation of $\fg$
for any element $a$ of $\fg$.
See \cite{loday.pirashvili.enveloping.leibniz.algebras.cohomology} 
for a discussion of Leibniz algebras. 
\index{Leibniz algebra}
\index{super Leibniz algebra}
\index{Leibniz identity}

A skew-symmetric super Leibniz algebra is called a
{\bf super Lie algebra}.
There is a forgetful functor from associative superalgebras
to super Lie algebras given by
$A\mapsto (A,[\; ,\; ])$ where
$[a,b]:=ab-(-1)^{\ta\tb} ba$.
This functor has a left adjoint $U$.
The superalgebra
$U(\fg)$ is called the {\bf enveloping algebra} of $\fg$.
\index{super Lie algebra}
\index{enveloping algebra}

\chapter{Distributions and Fields}
\label{C:formal distri}

\section{Distributions and Fields}
\label{S:formal distri}

We define {\it distributions} as continuous vector-valued
functions on a commutative algebra which is endowed with a topology.
In the next section and in the main body of the text
the commutative algebra will always be an algebra
of vertex or Laurent polynomials.
Still we formulate the basic notions and results about distributions
for arbitrary commutative algebras 
because in this way it becomes more apparent 
where the special properties of algebras of 
vertex polynomials enter the theory.
Moreover,
in this more general setting
our presentation parallels the 
discussion of distributions in analysis.

Another motivation for the more general approach
is that one may define {\it fields}
as distributions which take values in
an endomorphism algebra. 
This corresponds to the definition of 
a quantum field in axiomatic quantum field theory
as an operator-valued distribution,
see \cite{streater.wightman.pct.spin.all.that,
kazhdan.qft.princeton}.
It was noted by Chambert-Loir
\cite{chambertloir.distributions.champs}
that in the case of algebras of Laurent polynomials
this general definition of a field 
agrees with the usual one.
Finally,
Chambert-Loir's beautiful proof of the fact 
that local distributions 
can be expressed in terms of derivatives 
of the delta distribution,
see Proposition \ref{SS:locality},
is based on the notion of a kernel distribution.
It is useful to discuss this notion
in general terms.

\medskip

{\bf Summary.}\: 
In section \ref{SS:test distri}
we define the notions of a distribution 
and of a field.
In section \ref{SS:kernel delta}
we introduce kernel distributions and
prove some properties of the delta distribution.
In section \ref{SS:modes distr}
we define the notions of mode, support, and 
restriction of a distribution.

\medskip

{\bf Conventions.}\:
We denote by $V$ and $W$ vector spaces.
We endow $V$ with the discrete topology.

\subsection{Distributions, Test Functions, and Fields}
\label{SS:test distri}

Let $C$ be a commutative algebra
endowed with a topology, 
such that
multiplication with any element is continuous,
and let $U$ be a topological vector space.
A $U$-valued {\bf distribution} 
\index{distribution!on a commutative algebra} 
on $C$ is a
morphism $C\to U$ of topological vector spaces.
The elements of $C$ are called {\bf test functions}.
\index{test function}
The vector space $\cD(C,U)$
of $U$-valued distributions on $C$
is a module over $C$.
A morphism $U\to U'$ of topological vector spaces
induces
a morphism $\cD(C,U)\to\cD(C,U')$.

The {\bf weak topology} 
\index{weak topology}
on $\Vect(V,W)$ 
is defined by the family of 
neighborhoods of zero
$\set{a:V\to W\mid a(V')=0}$ 
where $V'$ runs over the finite-dimensional subspaces of $V$.
A {\bf field}
\index{field}
on $V$ over $C$ is an $\End(V)$-valued distribution on $C$
where $\End(V)$ is endowed with the weak topology.
We note that 
$\cD(C,\End(V))=\Vect(V,\cD(C,V))$.

Because $\Vect(C,U)$ has much better 
functorial properties than $\cD(C,U)$
we consider in this section only
the case where the topology of $C$ 
is discrete so that $\cD(C,U)=\Vect(C,U)$.
When the topology of $C$ is arbitrary 
we view $\cD(C,U)$ as a subspace of $\Vect(C,U)$.

\medskip

Let $C$ and $C'$ be commutative algebras.
The morphism $C'\to C\otimes C', a\mapsto 1\otimes a$,
induces a morphism of $C'$-modules
$\int_C:\cD(C\otimes C',V)\to\cD(C',V)$
which is called the {\bf residue} 
\index{residue}
morphism over $C$.
If $a$ is a distribution on $C$ and $b$ is a test function then
\begin{equation}
\label{E:distr int}
a(b)
\; =\;
(ab)(1)
\; =\;
\int_C ab.
\end{equation}

To give a linear form $\int:C\to\K$ 
is equivalent to giving
a morphism of $C$-modules $\io: C\to\cD(C,\K)$.
We have $\io(a)(b)=\int ab$ and $\int=\io(1)=\int_C\circ\io$.

If $\del$ is a derivation of $C$ 
then there exists a unique derivation $\del$ 
of the $C$-module $\cD(C,V)$
such that $\int_C\del\cD(C,V)=0$.
This follows from 
the {\bf integration-by-parts} formula
\index{integration-by-parts formula}
$$
\int_C \del(a)b
\; =\;
-\zeta^{\ta\tilde{\del}}\int_C a\del(b)
$$
which is satisfied for any distribution $a$ and
any test function $b$.
The morphism $C\to\cD(C,\K)$ induced by $\int$
is compatible with
$\del$ if and only if $\int\del C=0$
because both statements are equivalent to 
$\int\del(a)b=-\zeta^{\ta\tilde{\del}}\int a\del(b)$. 

In general, distributions do not form an algebra.
However, distributions can be multiplied if one
allows for a change of the algebra of test functions.
Indeed,
there exists a canonical morphism
$$
\cD(C,V)\otimes\cD(C',W)
\;\lra\;
\cD(C\otimes C',V\otimes W)
$$
which is called {\bf juxtaposition}
\index{juxtaposition} 
and is written $a\otimes b\mapsto ab$.
If $W$ is a module over an algebra $A$ then
juxtaposition induces a canonical morphism 
$\cD(C,A)\otimes W\to \cD(C,W)$.
In particular, there exists a canonical morphism
$\cD(C,C')\otimes\cD(C',V)\to\cD(C,\cD(C',V))=\cD(C\otimes C',V)$.

\subsection{Kernels and the Delta Distribution}
\label{SS:kernel delta}

Let $C$ and $C'$ be commutative algebras.
Consider the canonical isomorphism
$K:\Vect(C,\cD(C',V))\tra\cD(C\otimes C',V)$
defined by $K(\al)(a\otimes b):=\al(a)(b)$.
The distribution $K(\al)$ is called
the {\bf kernel} 
\index{kernel!distribution}
of $\al$.
From \eqref{E:distr int} follows that if $a\in C$ then
\begin{equation}
\label{E:kernel int}
\al(a)
\; =\;
\int_C K(\al)a.
\end{equation}
If $C'=\K$ then $K(\al)=\al:C\to V$.
The canonical isomorphism between
$\Vect(C,\cD(C',V))$ and 
$\Mod_{C'}(C\otimes C',\cD(C',V))$
will be used implicitly at times.

If 
$\varep:C\to C', a\mapsto \varep(a)\equiv a(\varep)$, 
is a morphism of algebras
then the kernel of $\varep$
is called the {\bf delta distribution} 
\index{delta distribution!on a commutative algebra}
on $C$ at $\varep$
and is denoted by $\de_{\varep}=\varep\in\cD(C,C')$. 
From \eqref{E:kernel int} follows that
if $a\in C$ then
$$
a(\varep)
\; =\;
\int_C \de_{\varep}a.
$$

If $a$ and $b$ are test functions in $C$
then
$(a\varep)(b)=\zeta^{\ta\tilde{\varep}}\varep(ab)
=\zeta^{\ta\tilde{\varep}}\varep(a)\varep(b)$
and thus 
$a\de_{\varep}=\zeta^{\ta\tilde{\varep}}\varep(a)\de_{\varep}$.

Let $\del$ and $\del'$ be even derivations of $C$ and $C'$
such that $\del'\varep(a)=\varep(\del a)$ for any $a\in C$.
If $c\in\cD(C',V)$ and $n\in\K$ 
then the product
$c(\del')^{(n)}\de_{\varep}\in\cD(C\otimes C',V)$
is the kernel of the morphism
$c(\del')^{(n)}\varep:C\to\cD(C',V)$.
For any $a\in C$ we have
$\del'\varep(a)=\varep(\del a)=-\del\varep(a)$
and thus $\del'\de_{\varep}=-\del\de_{\varep}$.
Let $z\in C$ such that $\del z=1$ and $w:=\varep(z)$.
For any $n\in\K$ and $a\in C$ we have
$(z-w)((\del')^{(n)}\varep)(a)
=
\varep(\del^{(n)}((z-w)a))
=
\varep(\del^{(n-1)}a)
$
and thus 
$(z-w)((\del')^{(n)}\de_{\varep})=(\del')^{(n-1)}\de_{\varep}$.

We call $\de_C:=\de_{\id_C}$ {\it the}
delta distribution on $C$.
Let $\int$ be a linear form on $C$
such that the induced morphism $C\to\cD(C,\K)$
is injective.
Then we may identify $\de_C$ with the 
$\K$-valued distribution on $C\otimes C$
defined by 
$a\otimes b\mapsto\int ab$.
From $\int ab=\paraab\int ba$ follows that
$\ga\de_C=\de_C$ where
$\ga$ is the automorphism of $\cD(C\otimes C,V)$
induced by the automorphism of $C\otimes C$
defined by $a\otimes b\mapsto\paraab b\otimes a$.

\subsection{Modes of a Distribution}
\label{SS:modes distr}

Let $C$ be a commutative algebra and
$(u_s)_{s\in S}$ be a basis of $C$.
If $a$ is a distribution on $C$ then
the vector $a_s:=a(u_s)$ 
is called the {\bf s-th mode} \index{mode} of $a$.
The {\bf support} 
\index{support!of a distribution}
of $a$ is the set of elements $s$ of $S$  
such that $a_s$ is non-zero.

Let $U$ be a topological vector space.
A family $(a_s)_{s\in S}$ in $U$ is called
{\bf summable} 
\index{summable!family of vectors}
with {\bf sum} 
\index{sum!of a summable family}
$a$ if $a$ is the unique limit of the
net $(\sum_{s\in T}a_s)_T$
where $T$ runs over the finite subsets of $S$.
The sum of a summable family $(a_s)_{s\in S}$ 
is denoted by $\sum_{s\in S}a_s$.
See \cite{bourbaki.top1}, Chapter III, \S 5,
for a discussion of summable families.

Let $\int$ be a linear form on $C$.
The canonical morphism 
$V\otimes C\to V\otimes\cD(C,\K)\to\cD(C,V)$
induced by $\int$
is written $a\otimes b\mapsto ab$. 
We endow the module $\cD(C,V)$ 
with the weak topology.
Assume that a basis $(v_s)_{s\in S}$ of $C$ exists
such that $\int v_s u_t=\de_{s,t}$.
Let $a\in\cD(C,V)$.
The family $(a_s v_s)_{s\in S}$ in $\cD(C,V)$
is summable with sum $a$.
In particular,
$V\otimes C$ is a dense submodule of $\cD(C,V)$.
The {\bf restriction} 
\index{restriction!of a distribution}
of $a$ to a subset $T$ of $S$ is the distribution 
$a|_T:=\sum_{s\in T}a_s v_s$.
The kernel of a vector space morphism 
$\al:C\to V$ is given by
\begin{equation}
\label{E:kernel modes}
K(\al)
\; =\;
\sum_{s\in S}\:
\al(u_s)\, v_s
\;\in\;
\cD(C,V).
\end{equation}

\section{Distributions and Fields on Vertex Polynomials}
\label{S:distri vertex}

{\bf Summary.}\: 
In section \ref{SS:distr variables}
we consider distributions on algebras of vertex
and Laurent polynomials.
In section \ref{SS:vertex series}   
we discuss vertex series and fields over vertex polynomials.
In section \ref{SS:power series expansions}
we discuss power series expansions of rational functions.

\medskip

{\bf Conventions.}\:
We denote by $V$ and $W$ vector spaces.
We endow $V$ with the discrete topology.

\subsection{Distributions on Vertex and Laurent Polynomials}
\label{SS:distr variables}

For the theory of OPE-algebras
the most important algebras of test functions
are the group rings $\K[\K^S]$ 
where $S$ is a finite even set.
The elements of $\K[\K^S]$ are
called {\bf vertex polynomials}.
\index{vertex polynomial}
We explain below the topology on $\K[\K^S]$.
If $z_1,\dots, z_r$ are the elements of $S$ 
then we denote $\K[\K^S]$
also by $\K[z_1^{\K},\dots,z_r^{\K}]$.
The elements of $S$ are called 
\index{variables}
{\bf variables}.

If we choose a total order $z_1>\dotsc >z_r$ on $S$
then distributions on $\K[\K^S]$ are denoted by
$a(z_1,\dots,z_r)$.
Conversely, if we denote a distribution by $a(z_1,\dots,z_r)$
then it is understood that $a(z_1,\dots,z_r)$ is a
distribution on $\K[\K^S]$ and that we have chosen 
the total order $z_1>\dotsc >z_r$ on $S$.
A total order on $S$ is used 
when we write a vector of $\K^S$ 
as an $r$-tuple $(n_1,\dots,n_r)$. 
It is also used implicitly in the definition of 
OPE-finiteness in section \ref{SS:ope finite}.

In the following we only discuss the case of one variable
if the extension to several variables 
is straightforward.

\medskip

In this subsection we consider 
distributions defined with respect to the 
discrete topology on $\K[z^{\K}]$.
In the next subsection
we endow $\K[z^{\K}]$ with a non-trivial topology.

We denote by $V\set{z}$
the module of $V$-valued distributions on $\K[z^{\K}]$
where $\K[z^{\K}]$ is endowed with the
discrete topology.
This notation was introduced in
\cite{frenkel.lepowsky.meurman.book}, Chapter 2.1.
The residue morphism $\int_{\K[z^{\K}]}:V\set{z}\to V$
is written $a(z)\mapsto\int a(z)dz\equiv \res_z(a(z))$.
Because $\K[z^{\K},w^{\K}]=\K[z^{\K}]\otimes\K[w^{\K}]$
we have $V\set{z,w}=V\set{z}\set{w}$.

We denote by $n\mapsto z^n$
the canonical morphism from 
the group $\K$ to the group ring $\K[z^{\K}]$.
Thus $(z^n)_{n\in\K}$ is a basis of $\K[z^{\K}]$
and $z^n z^m=z^{n+m}$.  

There exists a unique derivation $\del_z$ on 
$\K[z^{\K}]$ such that $\del_z z^n=n z^{n-1}$ for any $n\in\K$.
There exists a unique linear form $\int$ on $\K[z^{\K}]$
such that $\int\del_z\K[z^{\K}]=0$ and $\int z^{-1}=1$.
Because $\int z^n z^{-m-1}=\de_{n,m}$
we can identify 
$V[z^{\K}]:=V\otimes\K[z^{\K}]$ 
with a dense submodule of $V\set{z}$.

An automorphism $\al$ of the algebra $\K[z^{\K}]$
induces an automorphism of $V[z^{\K}]$.
If this automorphism extends 
to a bicontinuous automorphism $\al$ of $V\set{z}$ 
then $\al$ is written $a(z)\mapsto a(\al z)$.
For example,
the automorphism $z^n\mapsto z^{-n}$ of $\K[z^{\K}]$
induces an automorphism $a(z)\mapsto a(z^{-1})$
of $V\set{z}$.

Juxtaposition defines a
morphism
\index{juxtaposition} 
$$
V\set{z}\otimes W\set{w}
\;\lra\;
(V\otimes W)\set{z,w},
\quad
a(z)\otimes b(w)\mapsto a(z)b(w).
$$
In particular,
if $M$ is a module over an algebra $A$ 
then we obtain a morphism 
$A\set{z}\otimes M\set{w}\to M\set{z,w}$.

For historical reasons the modes $a_n$ of a distribution
$a(z)$ are defined with respect to the basis $(z^n)_{n\in\K}$
of $\K[z^{\K}]$ rather than with respect to the basis
$(z^{-n-1})_{n\in\K}$.
Therefore we have $\res_z(a(z))=a_0$
and $a(z)=\sum_{n\in\K}a_n z^{-n-1}$.

\medskip

We denote by $V\pau{z\uppm}$
the module of 
$V$-valued distributions on $\K[z\uppm]$
where $\K[z\uppm]$
is endowed with the discrete topology.
The inclusion  
$\K[z\uppm]\subset\K[z^{\K}]$
induces an epimorphism
$V\set{z}\to V\pau{z\uppm}$.
Its restriction to the subspace of distributions $a(z)$
such that $a(z)|_{\Z}=a(z)$ 
is an isomorphism onto $V\pau{z\uppm}$.
We thus identify $V\pau{z\uppm}$ with a subspace
of $V\set{z}$.
Distributions in $V\pau{z\uppm}$ are called {\bf integral}.
\index{distribution!integral}
\index{integral distribution}

\subsection{Vertex Series and Fields}
\label{SS:vertex series}

The {\bf linear topology} \index{linear topology} 
on $\K[z^{\K}]$ is defined by the neighborhoods of zero
$$
\sum_{n\in S}\:
z^n\K[z]
\; +\;
\sum_{n\in\K\setminus (S+\Z)}\:
z^n\K[z\uppm]
$$
where $S$ runs over the finite subsets of $\K$.
If $U$ and $U'$ are topological vector spaces
and $(U_s)$ and $(U'_{s'})$ are neighborhood bases
of zero of $U$ and $U'$ 
then we define a topology on $U\otimes U'$
by taking 
$(U_s\otimes U')+(U\otimes U'_{s'})$
as a neighborhood basis of zero.
The linear topology on $\K[z^{\K},w^{\K}]$
is defined to be the topology of the
tensor product 
$\K[z^{\K},w^{\K}]=\K[z^{\K}]\otimes\K[w^{\K}]$.
Note that $\K[z^{\K}]$ and 
$\K[z\uppm,w\uppm]$ are {\it not} topological algebras.

We denote by $V\sqbrack{z}$
the module of $V$-valued distributions on $\K[z^{\K}]$
where $\K[z^{\K}]$ is endowed with the linear topology.
Elements of $V\sqbrack{z}$
are called {\bf vertex series}.
\index{vertex series}
A distribution $a(z)$ is a vertex series 
if and only if the support of $a(z)$
is contained in finitely many cosets of $\Z_<\subset\K$.
The module $V\sqbrack{z,w}$ is contained in 
$V[(z/w)^{\K}]\set{z}$ and there exists
a morphism
$V\sqbrack{z,w}\to V\sqbrack{z}, 
a(z,w)\mapsto a(z,z)$.

We give another description of $V\sqbrack{z}$.
We endow $\K[z]$ with the linear topo\-lo\-gy.
Elements of the completion $V\pau{z}$ of 
$V[z]:=V\otimes\K[z]$ are called  
{\bf power series}.
\index{power series}
The canonical monomorphism
$V[z]\to V[z^{\K}]$ extends to a
continuous monomorphism
$V\pau{z}\to V\set{z}$
which induces an isomorphism
$\K[z^{\K}]\otimes_{\K[z]}V\pau{z}\tra V\sqbrack{z}$.

The morphism of multiplication
$V[z]\otimes_{\K[z]}W[z]\to (V\otimes W)[z]$
induces a morphism 
$V\pau{z}\otimes_{\K[z]}W\pau{z}\to (V\otimes W)\pau{z}$
which in turn induces a morphism
$$
V\sqbrack{z}\otimes_{\K[z^{\K}]}W\sqbrack{z}
\;\to\;
(V\otimes W)\sqbrack{z}, 
\qquad
a(z)\otimes b(z)
\;\mapsto\;
a(z)b(z).
$$
This morphism coincides with the morphism
$a(z)\otimes b(z)\mapsto a(z)b(w)|_{w=z}$.

\bigskip

Let $S$ and $\bS$ be finite even sets and 
$S\to\bS, z\mapsto\bz$, be a bijection.
We define $\vz:=(z,\bz)$ for any $z\in S$.
Despite the resulting ambiguity,
we write
$V\set{\vz}:=V\set{z,\bz}$
and
$a(\vz):=a(z,\bz)$.

We denote by $\cF_r(V)$
the module of fields 
on $V$ over 
$\K[\vz_1^{\K},\dots,\vz_r^{\K}]$
where $\K[\vz_1^{\K},\dots,\vz_r^{\K}]$ 
is endowed with the linear topology.
We have $\cF_r(V)=\Vect(V,V\sqbrack{\vz_1,\dots,\vz_r})$.
We denote by $\QEnd(V)$
the module of fields on $V$ over $\K[z\uppm]$
where $\K[z\uppm]$ is endowed with the linear topology.
\index{field}

\subsection{Power Series Expansions of Rational Functions}
\label{SS:power series expansions}

The morphism of fields $T_{z_1,\dots,z_r}$,
which we defined in section \ref{SS:prelim},
can be made more explicit as follows.
We first show that the morphism $T_z:\K(z)\to\K\lau{z}$
is given in terms of Taylor series.

\bigskip

{\bf Remark.}\: {\it
If $f(z)$ is a rational function and $n$ is an integer such that
$z^n f(z)$ is regular at $z=0$ then
\begin{equation}
\label{E:t is taylor}
T_z(f(z))
\; =\;
z^{-n}\sum_{m\in\N}\:
\del_z^{(m)}(z^n f(z))|_{z=0} \; z^m.
\end{equation}
}

% no bigskip

\begin{pf}
Denote the right-hand side of \eqref{E:t is taylor} by $T(f,n)$.
The product formula shows that 
$T(f,n+m)=T(f,n)$ for any non-negative integer $m$.
Thus $f(z)\mapsto T(f,n)$ defines a well-defined map
$T:\K(z)\to\K\lau{z}$ that is independent of the choice of $n$.
The product formula implies that $T(f,n)T(g,m)=T(fg,n+m)$.
Moreover, we have $T(z)=z$.
Thus $T=T_z$ because of the universal property of $T_z$.
\end{pf}

\bigskip

The morphism $T_{z_i}:\K(z_i)\to\K\lau{z_i}$ induces a morphism
$$
T_{z_i}:\:
\K\cB(z_1)\dots (z_i)\dots\cB(z_r)
\; \to\;
\K\cB(z_1)\dots\lau{z_i}\dots\cB(z_r)
$$
where for any $j\ne i$
the symbol $\cB(z_j)$ stands for either $(z_j)$ or $\lau{z_j}$.
Due to its universal property, the morphism $T_{z_1,\dots,z_r}$
is equal to the composition of the morphisms
$T_{z_1}, \dots, T_{z_r}$ in any order.
For example, for $r=2$ there exists a commutative diagram 
$$
\xymatrix{
\K(z,w)=\K(z)(w)
 \; \ar[rr]^{T_z}
   \ar[d]_{T_w}   
   \ar[rrd]^{T_{z,w}}  & &
\;\K\lau{z}(w)
   \ar[d]^{T_w}
\\
\K(z)\lau{w}
\;\ar[rr]^{T_z} & &
\;\K\lau{z}\lau{w}.
}
$$

The morphism $T_{z_1,\dots,z_r}$ commutes with $\del_{z_i}$
because the morphisms
$\del_{z_i}\circ T_{z_1,\dots,z_r}$ and
$T_{z_1,\dots,z_r}\circ \del_{z_i}$  
are both derivations that agree for $z_1, \dots, z_r$.

\chapter{Bibliographical Notes}
\label{C:biblio}

\section{Chapter \ref{C:intro}}
\label{S:biblio intro}

${}_{}$
\indent
{\bf Section \ref{S:intro overview}.}\:
The notion of an OPE-algebra was introduced by 
Kapustin and Orlov \cite{kapustin.orlov.vertex.algebras}, 
Definition 3.3. They use the term ``vertex algebra" for it
and they call vertex algebras ``chiral algebras".
The notion of a vertex algebra was introduced by
Borcherds \cite{borcherds.voa}, section 4.
Li \cite{li.localsystems},
Proposition 2.2.4 and Proposition 2.2.6,
proved that there are three equivalent definitions of a vertex algebra.

\medskip

{\bf Section \ref{S:intro op alg}.}\:
The idea of an operator algebra of quantum fields
was put forward by Polyakov \cite{polyakov.bootstrap} and 
Kadanoff \cite{kadanoff.ope,kadanoff.ceva.operator.algebra.2d.ising}.
The normal ordered product of quantum fields was 
introduced by Wick \cite{wick.theorem.normal.ordering}.
The use of the operator product expansion was pointed out
by Wilson \cite{wilson.ope1} and Kadanoff \cite{kadanoff.ope}.

\medskip

{\bf Section \ref{S:intro rev dual local}.}\:
Goddard \cite{goddard.meromorphic} calls 
equations \eqref{E:intro duali} and \eqref{E:intro locality}
``duality" and ``locality", resp.
Huang showed in \cite{huang.book} that vertex operator algebras
can be formulated geometrically.
The equivalence of the Jacobi identity
with locality and with duality and skew-symmetry
is due to Li \cite{li.localsystems},
Proposition 2.2.4 and Proposition 2.2.6.
Bakalov and Kac \cite{bakalov.kac.field.algebras}, Theorem 7.4,
showed that 
vertex algebras with a Virasoro vector can be defined in terms
of duality alone.
Both Li and Bakalov and Kac assume 
the existence of an identity and a translation operator.

The notion of a generalized vertex algebra 
was introduced by Dong and Lepowsky 
\cite{dong.lepowsky.book}.
Essentially the same concept was also defined by 
Feingold, Frenkel, and Ries
\cite{feingold.frenkel.ries.spinor.vertex.triality.e(1)8}
and by Mossberg \cite{mossberg.jacobi}.
Huang formulated the notion of an intertwining algebra 
in \cite{huang.g=0cft}.
The notion of a $G_n$-vertex algebra 
was introduced by Li \cite{li.higher.dim.analogues.vertex}.

\medskip

{\bf Section \ref{S:intro dual local}.}\:
The vertex operators of string theory were introduced by
Fubini and Veneziano
\cite{fubini.veneziano.vertex.operators}.

\medskip

{\bf Section \ref{S:intro examp}.}\:
For the work of Huang on duality and locality for intertwiners,
see \cite{huang.differ.equations.intertwining.operators} and 
references given there. 

The original works about 
cosets, orbifolds, simple current extensions, 
Gepner's $U(1)$-projection, and BRST cohomology
are
\cite{goddard.kent.olive.unitary.coset},
\cite{dixon.harvey.vafa.witten.orbifolds1},
\cite{schellekens.yankielowicz.extended.chiral.modular.inv.simple.curr},
\cite{gepner.models},
and
\cite{kato.ogawa.brsstring,felder.brst}, respectively.

\section{Chapter \ref{C:fields}}
\label{S:biblio fields}

${}_{}$
\indent
{\bf Section \ref{S:formal distributions}.}\:
The module $V\set{z}$ is discussed in
\cite{frenkel.lepowsky.meurman.book}, sections 2.1--2.2 and 8.1--8.3.
The term ``distribution" designating the elements of $V\set{z}$
goes back to Kac \cite{kac.beginners.first.ed}, equation (2.1.1),
who calls the elements of $V\pau{z\uppm}$ ``formal distributions".

The use of the operator product expansion 
in quantum field theory was pointed out by
Wilson \cite{wilson.ope1} and Kadanoff \cite{kadanoff.ope}.
The concept of a normal ordered product of quantum fields
is due to Wick \cite{wick.theorem.normal.ordering}.

Lemma \ref{SS:n-fold module holom distr} 
and Proposition \ref{SS:ope}
are due to Lian and Zuckerman 
\cite{lian.zuckerman.quantum}, Proposition 2.3.
The $n$-th products of the 
$\N$-fold module of holomorphic distributions
were originally defined for holomorphic fields 
by Li \cite{li.localsystems}, Lemma 3.1.4,
and by Lian and Zuckerman 
\cite{lian.zuckerman.quantum}, Definition 2.1.
The equivalence of the ``OPE" for $z>w$ \eqref{E:commut ope}
with the commutator formula \eqref{E:commut formula0}
was pointed out by 
Roitman \cite{roitman.free.conformal.vertex.algebras}, equation (1.3).
The normal ordered product $\normord{a(z)b(w)}$
was introduced by Lian and Zuckerman 
\cite{lian.zuckerman.quantum}, equation (2.2).

The holomorphic Jacobi identity was introduced 
by Frenkel, Lepowsky, and Meurman
\cite{frenkel.lepowsky.meurman.book}, section 8.10.
The proof we give of Proposition \ref{SS:Jacobi formal distributions}
is due to Matsuo and Nagatomo
\cite{matsuo.nagatomo.locality}, Corollary 3.2.2.
Instead of $\fg\set{\vz}$ they consider the 
$\Z$-fold algebra of holomorphic fields 
in which case the holomorphic Jacobi identity
is valid for any indices $r\in\N, s\in\Z$, and $t\in\N$.

The associativity formula already appears in
Borcherds' definition of a vertex algebra,
\cite{borcherds.voa}, section 4.
Moreover, Borcherds \cite{borcherds.voa}, section 8,
explains that a vertex algebra satisfies the commutator formula.
Proposition \ref{SS:pf hol jacobi}\,$\iti$
is due to Li \cite{li.localsystems}, equation (2.2.9).
Proposition \ref{SS:pf hol jacobi}\,$\itiii$
and the proof of Lemma \ref{SS:pf hol jacobi}
are due to Matsuo and Nagatomo \cite{matsuo.nagatomo.locality}, 
Proposition 4.3.3 and equation (4.3.1).

In \cite{kac.beginners.first.ed}, 
Theorem 4.8, equation (4.6.7), and Theorem 4.6,
the Jacobi identity is called Borcherds identity and 
the commutator formula is called Borcherds commutator formula and
when it is written in terms of OPEs it is called Borcherds OPE formula.
Kac \cite{kac.beginners}, Proposition 2.3\,$\itd$,
proves that the $\N$-fold algebra $\fg\pau{z\uppm}$
satisfies the commutator formula. 
Remark \ref{SS:conformal algebras} shows that his result implies
Proposition \ref{SS:Jacobi formal distributions}.

\medskip

{\bf Section \ref{S:vir symm}.}\:
Proposition \ref{SS:derivative}\,$\iti$, $\itiii$ and
Proposition \ref{SS:grad dilat cov}
are due to Li \cite{li.localsystems}, Lemma 3.1.7 and Lemma 3.1.8.
He considers the $\Z$-fold algebra of holomorphic fields.

The discovery of the Virasoro algebra is usually attributed 
to Gelfand and Fuchs \cite{gelfand.fuks.cohomologies.lie.vector.fields.circle} 
who computed the cohomology of the Lie algebra of smooth vector fields 
on the circle. The importance of this algebra was pointed out by
Virasoro \cite{virasoro.virasoro} who showed that an infinite number
of operators $L_n$ act on the Fock space of the bosonic string
that implement residual gauge conditions and 
single out the physical state space in the Gupta-Bleuler
quantization of the string. 
The operators $L_n$ generate the Virasoro algebra but
the central term was overlooked by Virasoro. 
That was later found by Fubini and Veneziano
\cite{fubini.veneziano.virasoro}.

\medskip

{\bf Section \ref{S:locality distr}.}\:
The one-variable version $\de(x)=\sum_{n\in\Z}x^n$ of the
delta distribution $\de(z,w)$ was introduced in
\cite{frenkel.lepowsky.meurman.book}, equation (2.1.22).
We have $\de(z,w)=z^{-1}\de(w/z)$.
The analogue of Proposition \ref{SS:delta distr}
with $\de(z,w)$ replaced by $\de(x)$
is obtained in \cite{frenkel.lepowsky.meurman.book},
Proposition 8.1.2.
Taylor's formula for holomorphic fields $a(z,w)$
and Propositions \ref{SS:locality},
\ref{SS:locality and ope}, and \ref{SS:skew-symmetry distr}
are due to Kac \cite{kac.beginners},
Proposition 3.1, Corollary 2.2, Theorem 2.3, and 
Proposition 2.3\,$\itc$.
The proof we give of Proposition \ref{SS:locality}
is due to Chambert-Loir 
\cite{chambertloir.distributions.champs},
Proposition 2.2.

The identity of holomorphic skew-symmetry was discovered by
Borcherds \cite{borcherds.voa}, section 4.
It is one of the axioms in his definition of a vertex algebra.
The original version of Dong's lemma for holomorphic locality
deals with the $n$-th products of the $\Z$-fold algebra of 
holomorphic fields and is due to Li \cite{li.localsystems},
Proposition 3.2.7, 
who acknowledges Dong for providing a proof of it.
The effective bound on the order is due to
Matsuo and Nagatomo 
\cite{matsuo.nagatomo.locality}, Proposition 2.1.5.

\medskip

{\bf Section \ref{S:module q fields}.}\:
The notion of a holomorphic field is due to Li \cite{li.localsystems},
Definition 3.1.1. 
Kac \cite{kac.beginners.first.ed}, equation (1.3.1), 
gave this notion the name ``field".
This term already appears in \cite{frenkel.kac.radul.wang.w},
section 3, where it is used for homogeneous $\End(V)$-valued
distributions in the case that $V$ is endowed with a gradation
that is bounded from below.
The $n$-th products of the $\Z$-fold algebra $\cF_z(V)$ 
of holomorphic fields were defined 
by Li \cite{li.localsystems}, Lemma 3.1.4,
and by Lian and Zuckerman \cite{lian.zuckerman.quantum}, Definition 2.1.
The normal ordered product $\normord{a(z)b(z)}$
was also introduced in \cite{frenkel.kac.radul.wang.w},
equation (3.2).
Proposition \ref{SS:holom loc implies itself} is due to
Li \cite{li.localsystems}, Proposition 3.2.9.

Proposition \ref{SS:creat}\,$\itii$ 
is an extension of a result of Li \cite{li.invariant.forms.voas},
Proposition 3.3\,$\ita$.
Proposition \ref{SS:creat}\,$\itiii$
is due to Matsuo and Nagatomo
\cite{matsuo.nagatomo.locality}, Lemma 5.1.2.
Proposition \ref{SS:id field}
is due to Li \cite{li.localsystems}, Lemma 3.1.6.

Independently from Lian-Zuckerman and Li, 
Meurman and Primc \cite{meurman.primc.annihilating.fields.sl2.combinat.ids},
Propositions 2.3 and 2.5, 
defined the $n$-th products of holomorphic fields, 
proved that these products preserve translation covariance,
dilatation covariance, and creativity, 
and proved Dong's lemma for holomorphic locality.

\medskip

{\bf Section \ref{S:algebra q fields}.}\:
The definition of the canonical non-integral powers
$(\mu\vz+\nu\vw)^{\vh}$ is due to Kapustin and Orlov
\cite{kapustin.orlov.vertex.algebras}, section 3.2.
Proposition \ref{SS:ope finite} is due to Kapustin and Orlov
\cite{kapustin.orlov.vertex.algebras}, Lemma B.1.
The idea of the proof of Remark \ref{SS:reduced ope}
is due to D.~Orlov.
Proposition \ref{SS:ope holom distr}
is due to Kapustin and Orlov
\cite{kapustin.orlov.vertex.algebras}, Proposition B.5.

\medskip

{\bf Section \ref{S:locality}.}\:
The definition of locality
and Proposition \ref{SS:ope commut}
are due to Kapustin and Orlov 
\cite{kapustin.orlov.vertex.algebras}, 
Definition 3.3, Proposition B.4, and end of appendix B.
The idea of the proof of Lemma \ref{SS:ope commut}\,$\itiii$
is due to D.~Orlov.
Goddard's uniqueness theorem 
was first derived by Goddard in the framework of
meromorphic conformal field theory
\cite{goddard.meromorphic}, Theorem 1.
Because local holomorphic fields are additively local
Dong's lemma for additive locality is an extension 
to the $\K^2$-fold algebra of fields
of a result of Kapustin and Orlov
\cite{kapustin.orlov.vertex.algebras},
Lemma B.7,
that only concerns the $\Z$-fold module 
$\cF(V)$ over $\cF_z(V)$.
Their result in turn generalizes 
the original Dong's lemma 
which only deals with the $\Z$-fold algebra $\cF_z(V)$.

\section{Chapter \ref{C:opea}}
\label{S:biblio opea}

${}_{}$
\indent
{\bf Section \ref{S:ope algebras}.}\:
Kapustin and Orlov 
\cite{kapustin.orlov.vertex.algebras}, Definition 3.3,
define a vertex algebra as a local bounded 
$\vbbK$-fold algebra together with a translation gene\-rator $\vT$ and 
an invariant strong identity such that $T$ and $\bT$ commute.
They prove \cite{kapustin.orlov.vertex.algebras}, Corollary B.3,
that $\vT$ is a translation operator.
Proposition \ref{SS:ko vertex} shows that their notion 
of a vertex algebra is equivalent to the notion of an OPE-algebra. 
Proposition \ref{SS:chiral alg}
is due to Kapustin and Orlov 
\cite{kapustin.orlov.vertex.algebras}, 
Proposition B.4 and Corollary B.6.
The analogue for vertex algebras of 
Proposition \ref{SS:space fields}\,$\iti$
is due to Matsuo and Nagatomo
\cite{matsuo.nagatomo.locality}, 
Theorem 6.3.3.
The analogue for vertex algebras of 
Proposition \ref{SS:space fields}\,$\itii$
was stated without proof by Lian and Zuckerman
\cite{lian.zuckerman.classicalq}, Theorem 5.6, 
under the additional assumption that $\cF$ contains the identity field.
The analogue was proven by Matsuo and Nagatomo 
\cite{matsuo.nagatomo.locality}, Theorem 5.3.2.
The existence theorem for OPE-algebras is a generalization of the 
existence theorem for vertex algebras which is due to 
Frenkel, Kac, Radul, and Wang \cite{frenkel.kac.radul.wang.w}, Proposition 3.1,
and Meurman and Primc 
\cite{meurman.primc.annihilating.fields.sl2.combinat.ids}, Theorem 2.6.
The proof we give is an extension of the proof of the existence theorem
for vertex algebras that was given by Matsuo and Nagatomo
\cite{matsuo.nagatomo.locality}, Corollary 6.4.1.

\medskip

{\bf Section \ref{S:addi opea}.}\:
Lemma \ref{SS:ope alg} and Proposition \ref{SS:ope alg} 
are extensions of two results of Li about $G_n$-vertex algebras
\cite{li.higher.dim.analogues.vertex}, 
Proposition 3.10 and Theorem 3.7.
The notion of a $G_n$-vertex algebra is a generalization of the
notion of a vertex algebra from $\Z$-fold algebras to
$\Z^n$-fold algebras; in other words,
integral distributions in one variable are replaced
by integral distributions in $n$ variables.
Li \cite{li.higher.dim.analogues.vertex}, Definition 3.1,
defines a $G_n$-vertex algebra as a bounded
$\Z^n$-fold algebra with an identity and translation
generator satisfying a multi-variable holomorphic Jacobi identity.
This Jacobi identity has $2^n+1$ terms.
In Theorem 3.7 Li proves that a $G_n$-vertex algebra 
can equivalently be defined in terms of boundedness,
existence of an identity, translation covariance,
and multi-variable holomorphic locality.
In Proposition 3.10 Li proves that a $G_n$-vertex algebra 
can equivalently be defined in terms of boundedness,
existence of an identity, translation covariance,
multi-variable holomorphic duality, and
multi-variable holomorphic skew-symmetry.

The Jacobi identity as defined in section \ref{SS:jacobi modes}
combines the Jacobi identity for generalized vertex algebras
\cite{dong.lepowsky.book,feingold.frenkel.ries.spinor.vertex.triality.e(1)8,
mossberg.jacobi}
and the Jacobi identity for $G_2$-vertex algebras.
Proposition \ref{SS:dual local pretex} 
is an extension of a result of Li 
\cite{li.higher.dim.analogues.vertex},
Proposition 3.5, that states that the
Jacobi identity in the definition of a $G_n$-vertex algebra
can be replaced by 
multi-variable holomorphic duality plus
multi-variable holomorphic locality.

\medskip

{\bf Section \ref{S:vertex algebras}.}\:
Proposition \ref{SS:vertex algebras}\,$\itiii$
with ``left identity" replaced by ``strong identity"
is exactly Borcherds' original definition of a vertex algebra 
\cite{borcherds.voa}, section 4.
Li \cite{li.localsystems},
Proposition 2.2.4 and Proposition 2.2.6,
proves that for a bounded $\Z$-fold algebra $V$
with an identity and a translation operator
the following statements are equivalent:
$\iti$\: $V$ satisfies the holomorphic Jacobi identity;
$\itii$\: $V$ is holomorphically local;
$\itiii$\: $V$ satisfies holomorphic duality and 
holomorphic skew-symmetry.
He also states \cite{li.localsystems}, Corollary 2.2.7,
that Borcherds' definition of a vertex algebra is equivalent to 
the definition in terms of properties $\iti$--$\itiii$.
In \cite{frenkel.huang.lepowsky.axiom}, Remark 2.2.4,
it is shown that a vertex operator algebra can be defined in terms of 
injectivity of $Y$ and existence of a left identity.
Proposition \ref{SS:vertex algebras}\,$\itvii$ is the analogue
of this result for vertex algebras.

Li \cite{li.localsystems} defines a module over a vertex algebra $V$
as a vector space $M$ together with a linear map $Y:V\to\cF_z(M)$
and a translation operator such that $Y(1,z)=1(z)$ and the 
Jacobi identity is satisfied. He notes, Proposition 2.3.3,
that in this definition the Jacobi identity can be replaced by duality.
Proposition \ref{SS:modules} is a variation of his result.
The first two claims of Remark \ref{SS:modules} 
are due to Li \cite{li.localsystems}, Theorem 3.2.10 and
Corollary 3.2.11.

In \cite{malikov.schechtman.derham}, equation (1.3.4),
it is claimed that the associativity formula implies the Jacobi
identity. The remark made after Proposition \ref{SS:commu vertex}
shows that this is wrong.

\section{Chapter \ref{C:str va}}
\label{S:biblio strva}

${}_{}$
\indent
{\bf Section \ref{S:conformal algebras}.}\:
The notion of a conformal algebra was introduced independently 
by Kac \cite{kac.beginners.first.ed}, Definition 2.7b,
and Primc \cite{primc.lie}, equation (3.3).
The claim made in \cite{kac.beginners.first.ed}, equation (2.7.2),
and also in \cite{kac.beginners}, after Definition 2.7,
that in the definition of a conformal algebra it suffices
to require that there exists a translation endomorphism, is wrong.
A survey about conformal algebras is given in 
\cite{kac.conf.algebras}.
Under natural additional assumptions
conformal algebras are amenable to classifications,
see \cite{kac.superconformal.algebras.quadrics,yamamoto.conformal,
kac.superconformal.algebras.quadrics.err,
fattori.kac.classifi.finite.simple.lie.conf.superalgebras}.
Remark \ref{SS:conformal algebras} is due to 
Primc \cite{primc.lie}, Lemma 6.1 and Lemma 6.2.
Proposition \ref{SS:conformal algebras} is due to 
Kac \cite{kac.beginners.first.ed}, beginning of section 2.7.

The notion of a local Lie algebra was 
introduced by Kac \cite{kac.beginners.first.ed},
Definition 2.7a and equation (4.7.1).
He calls them ``regular formal distribution Lie algebras".
The claim made in \cite{kac.beginners.first.ed} after equation (4.7.1)
and also in \cite{kac.beginners}, equation (4.7.1), 
that in the definition of a local Lie algebra $\fg$ it suffices 
to require that there exists an operator of $\fg$
for which $F$ is translation covariant, is wrong.
A notion very similar to the notion of a local Lie algebra
was defined by 
Dong, Li, Mason \cite{dong.li.mason.vertex.lie.poisson.algebras},
Definition 3.1.(1).

Section \ref{SS:affiniz} about superaffinization 
Kac \cite{kac.beginners.first.ed}, section 2.5.
Proposition \ref{SS:conf affinization} is due to 
Kac \cite{kac.beginners}, equation (2.7.4) and Remark 2.7d.

Remark \ref{SS:mod alg of n fold alg}
is due to Borcherds \cite{borcherds.voa}, section 4,
in the case of vertex algebras and to Primc \cite{primc.lie}, Lemma 3.1, 
and Kac \cite{kac.beginners.first.ed}, Remark 2.7a, in the general case.
Primc \cite{primc.lie}, Theorem 4.1, 
proves that if $R$ is a conformal algebra then $\fg(R)$ is a Lie
algebra with a derivation. 

Proposition \ref{SS:loc mod alg of n fold alg}\,$\iti$
is due to Kac \cite{kac.beginners}, after Remark 2.7d.
Proposition \ref{SS:loc mod alg of n fold alg}\,$\itii$
is due to Primc \cite{primc.lie}, Proposition 4.4 and Theorem 4.6,
and Kac \cite{kac.beginners}, before Lemma 2.7.
Dong, Li, Mason \cite{dong.li.mason.vertex.lie.poisson.algebras},
Lemma 5.3, 
obtained this result in the special case of vertex algebras.
Proposition \ref{SS:loc mod alg of n fold alg}\,$\itiii$
is due to Roitman \cite{roitman.free.conformal.vertex.algebras},
Proposition 1.3.(e).
Moreover, Roitman \cite{roitman.free.conformal.vertex.algebras},
Proposition 1.3, 
observes that Proposition \ref{SS:loc mod alg of n fold alg}
holds in the generality of $\N$-fold algebras.

The definition of an irregular ideal and the equivalence
between certain categories of local Lie algebras and 
conformal algebras is due to Kac \cite{kac.beginners}, Theorem 2.7.
Proposition \ref{SS:free conf} is due to
Roitman \cite{roitman.free.conformal.vertex.algebras}, Proposition 3.1.

\medskip

{\bf Section \ref{S:conf and verta}.}\:
Proposition \ref{SS:verma va} is due to 
Primc \cite{primc.lie}, Theorem 5.3, in the case $\fg=\fg(R)$
and to Kac \cite{kac.beginners.first.ed}, Theorem 4.7, and 
Dong, Li, Mason \cite{dong.li.mason.vertex.lie.poisson.algebras}, Theorem 4.8,
in the general case.
Remark \ref{SS:vertex envelope} and Corollary \ref{SS:vertex envelope}
are due to Primc \cite{primc.lie}, Proposition 5.4 and Theorem 5.5.
The proof we give of Corollary \ref{SS:vertex envelope},
using the universal property of the Verma module,
is different from Primc's and was found by Roitman
\cite{roitman.michael.phd}, Theorem 2.4.
Primc \cite{primc.lie}, Theorem 5.8, and
Dong, Li, Mason \cite{dong.li.mason.vertex.lie.poisson.algebras}, Theorem 4.8,
prove that any $\fg$-module is a $V(\fg)$-module.
This result follows from Proposition \ref{SS:vertex reps}\,$\iti$
and Proposition \ref{SS:vertex envelope}\,$\itii$.

Borcherds \cite{borcherds.voa}, section 4,
mentions that there exists a free vertex algebra associated to
any pole order bound and that it can be constructed as a subalgebra of 
the lattice vertex algebra of a lattice that is associated to the pole
order bound. 
Roitman \cite{roitman.free.conformal.vertex.algebras}, end of section 3.1, 
constructs the free vertex algebra corresponding to a non-negative
pole order bound as the enveloping vertex algebra of a free 
conformal algebra. 

Primc \cite{primc.lie}, Lemma 7.1, proves that there exists a 
unique $\N$-fold algebra structure with translation operator
on $\K[T]\otimes\K^{(S)}$ such that $s_{(n)}t=O(s,t,n)$.
This is the main part of Remark \ref{SS:linear opes}.
Moreover, Primc \cite{primc.lie}, Theorem 8.4, proves that if $O$ 
is a linear OPE that contains a conformal vector such that 
the elements of $S$ are quasi-primary, 
$\hk$ is the only element of weight zero, 
the weights are all non-negative, and 
there are only finitely many elements of a given weight then
$\K[T]\otimes\K^{(S)}\oplus\K\hk$ is a conformal algebra if and only if
the elements of $S$ satisfy the Jacobi identity and skew-symmetry.
Primc's proof makes essential use of the existence of a conformal
vector. Proposition \ref{SS:linear opes} is a generalization
and an extension of the result of Primc.
Kac \cite{kac.beginners.first.ed}, Theorem 2.7,
sketches a proof of the statement that
$\K[T]\otimes\K^{(S)}$ has a unique structure of a conformal
algebra if the elements in $S$ satisfy boundedness,
skew-symmetry, and the commutator formula.
This theorem is omitted in the second edition \cite{kac.beginners}.

\backmatter%%%%%%%%%%%%%%%%%%%%%%%%%%%%%%%%%%%%%%%%%%%%%%%%%%%%%%%

\bibliographystyle{alpha}
\bibliography{lit/maths,lit/physics,lit/topology,lit/collections}

\include{referenc}
\printindex

%%%%%%%%%%%%%%%%%%%%%%%%%%%%%%%%%%%%%%%%%%%%%%%%%%%%%%%%%%%%%%%%%%%%%%

\end{document}